\documentclass[11pt,twoside]{amsart}
\usepackage{amsmath}
\usepackage{amssymb}
\usepackage{amsfonts}
\usepackage{amscd}
\usepackage{color}
\usepackage[small,nohug,heads=vee]{diagrams}
\diagramstyle[labelstyle=\scriptstyle]

\newtheorem{thm}{Theorem}[section]

\newtheorem{prop}[thm]{Proposition}

\newtheorem{lemma}[thm]{Lemma}

\newtheorem{cor}[thm]{Corollary}

\newtheorem{claim}[thm]{Claim}

\theoremstyle{definition}
\newtheorem{defn}[thm]{Definition}

\newtheorem{example}[thm]{Example}

\newtheorem{notation}[thm]{Notation}

\newtheorem{convention}[thm]{Convention}
\newtheorem{rem}[thm]{Remark}

\newcommand\op{\operatorname}

\newcommand\abs{\op{abs}}
\newcommand\ample{\op{ample}}

\newcommand{\an}{{\op{an}}}

\newcommand\Aut{\op{Aut}}

\newcommand\can{\op{can}}

\newcommand\chara{\op{char}}

\newcommand\closed{\op{cl}}

\newcommand\codim{\op{codim}}

\newcommand\Cone{\op{Cone}}

\newcommand\Conv{\op{Conv}}
\newcommand\CUB{\op{CUB}}
\newcommand\Cut{\op{Cut}}
\newcommand\cut{\op{cut}}
\newcommand\cutlog{\op{cutlog}}

\newcommand\Del{\op{Del}}

\newcommand\DD{\op{DD}}

\newcommand\depth{\op{depth}}

\newcommand\Fan{\op{Fan}}

\newcommand\FC{\op{FC}}

\newcommand\FP{\op{FP}}

\newcommand\formal{\op{for}}

\newcommand\Gal{\op{Gal}}

\newcommand\GL{\op{GL}}

\newcommand\height{\op{ht}}
\newcommand\Hom{\op{Hom}}

\newcommand\id{\op{id}}

\newcommand\inv{\op{inv}}

\newcommand\init{\op{in}}

\newcommand\Isom{\op{Isom}}
\newcommand\Jac{\op{Jac}}

\newcommand{\LGRS}{{\op{LGRS}}}

\newcommand\lft{\op{lft}}

\newcommand\Max{\op{Max}}

\newcommand\nor{\op{nor}}

\newcommand\Orb{\op{Orb}}
\newcommand\orb{\op{orb}}

\newcommand\Pic{\op{Pic}}

\newcommand\Proj{\op{Proj}}

\newcommand\quot{\op{quot}}
\newcommand\red{\op{red}}
\newcommand\Reg{\op{Reg}}

\newcommand\res{\op{res}}
\newcommand\rank{\op{rank}}
\newcommand{\rig}{{\op{rig}}}
\newcommand{\Rig}{{\op{Rig}}}
\newcommand\Sat{\op{Sat}}
\newcommand\Sch{\op{Sch}}
\newcommand\sep{\op{sep}}
\newcommand\Semi{\op{Semi}}
\newcommand\Sing{\op{Sing}}
\newcommand\Sk{\op{Sk}}

\newcommand\specialization{\op{sp}}
\newcommand\Spec{\op{Spec}}
\newcommand\Spf{\op{Spf}\,}
\newcommand\Sprig{\op{Sp}\,}

\newcommand\spl{\op{spl}}

\newcommand\Sup{\op{Sup}}

\newcommand\unr{\op{unr}}

\newcommand\Vor{\op{Vor}}
\newcommand\vor{\op{vor}}
\newcommand\wt{\op{wt}}

\newcommand{\tf}{{\tilde f}}

\newcommand{\tp}{{\tilde p}}

\newcommand{\tw}{{\tilde w}}
\newcommand{\tx}{{\tilde x}}

\newcommand{\trho}{{\tilde\rho}}

\newcommand{\tpsi}{{\tilde\psi}}

\newcommand{\tmu}{{\tilde\mu}}

\newcommand{\tGamma}{{\tilde\Gamma}}
\newcommand{\tLambda}{\tilde\Lambda}
\newcommand{\tA}{{\tilde A}}

\newcommand{\tC}{{\tilde C}}
\newcommand{\tD}{{\tilde D}}

\newcommand{\tG}{{\tilde G}}

\newcommand{\tL}{{\tilde L}}

\newcommand{\tP}{{\tilde P}}
\newcommand{\tQ}{{\tilde Q}}
\newcommand{\tS}{{\tilde S}}
\newcommand{\tT}{{\tilde T}}

\newcommand{\tV}{{\tilde V}}
\newcommand{\tW}{{\tilde W}}
\newcommand{\tX}{{\tilde X}}

\newcommand{\tcL}{{\tilde\cL}}

\newcommand{\tcV}{{\tilde\cV}}

\newcommand{\ep}{\epsilon}

\newcommand{\bara}{\bar{a}}
\newcommand{\barb}{\bar{b}}

\newcommand{\barp}{{\bar{p}}}

\newcommand{\barG}{\bar G}

\newcommand{\barZ}{\bar Z}

\newcommand{\bG}{{\mathbf G}}

\newcommand{\bL}{{\mathbf L}}

\newcommand{\bN}{{\mathbf N}}
\newcommand{\bP}{{\mathbf P}}
\newcommand{\bQ}{{\mathbf Q}}
\newcommand{\bR}{{\mathbf R}}

\newcommand{\bZ}{{\mathbf Z}}
\newcommand{\cA}{{\mathcal A}}
\newcommand{\cB}{{\mathcal B}}

\newcommand{\cD}{{\mathcal D}}

\newcommand{\cF}{{\mathcal F}}
\newcommand{\cG}{{\mathcal G}}

\newcommand{\cL}{{\mathcal L}}
\newcommand{\cM}{{\mathcal M}}
\newcommand{\cN}{{\mathcal N}}

\newcommand{\cO}{{\mathcal O}}

\newcommand{\cV}{{\mathcal V}}
\newcommand{\cW}{{\mathcal W}}
\newcommand{\cX}{{\mathcal X}}

\newcommand{\fraka}{{\mathfrak a}}
\newcommand{\fm}{{\mathfrak m}}

\newcommand{\fb}{{\mathfrak b}}

\ifx\undefined\pf
\newcommand{\pf}{\proof}
\fi
\ifx\undefined\cal
\newcommand{\cal}{\mathcal}
\fi
\ifx\undefined\Bbb
\newcommand{\Bbb}{\mathbb}
\fi
\ifx\undefined\bold
\newcommand{\bold}{\mathbf}
\fi
\ifx\undefined\frak
\newcommand{\frak}{\mathfrak}
\fi
\begin{document}

\title[Relative compactification]
{Relative compactification of semiabelian N\'eron models, I}
\author[Kentaro Mitsui  and Iku Nakamura]{K.~ Mitsui and I.~Nakamura}
\address{Department of Mathematics, 
Graduate School of Science, Kobe University, Hyogo 657-8501, Japan} 
\email{mitsui@math.kobe-u.ac.jp}
\address{Department of Mathematics, 
Hokkaido University, Sapporo, 060, Japan}
\email{nakamura@math.sci.hokudai.ac.jp}
\thanks{The first author is partially supported by JSPS KAKENHI (No. 17K14167, 17H06127, 17H02832, 21K03179). The second author was partially 
supported by the Grants-in-aid 
 for Scientific Research (C) (No. 17K05188, 22K03261), JSPS. 
This work was partially supported also by the Research Institute 
for Mathematical Sciences, an International Joint Usage/Research Center located in Kyoto University.}
\begin{abstract}Let $R$ be a complete discrete valuation ring, $k(\eta)$ its fraction field, 
$S:=\Spec R$, 
$(G_{\eta},\cL_{\eta})$ a polarized   
abelian variety  over $k(\eta)$ 
with $\cL_{\eta}$ ample cubical    
 and $\cG$ the N\'eron model of $G_{\eta}$ over $S$. 
Suppose that $\cG$ is totally degenerate semiabelian over $S$. 
Then there exists a (unique) relative compactification 
$(P,\cN)$ of $\cG$ such that 
($\alpha$) $P$ is Cohen-Macaulay with $\codim_P(P\setminus\cG)=2$ and  
($\beta$) $\cN$ is ample invertible 
with $\cN_{|\cG}$ cubical and $\cN_{\eta}=\cL^{\otimes n}_{\eta}$ 
for some positive integer $n$. 
\end{abstract}
\maketitle
\setcounter{tocdepth}{1}
\tableofcontents

\section{Introduction}
\label{sec:introduction}
Let $R$ be a complete discrete valuation ring (abbr. CDVR), 
$k(\eta)$ (resp. $k(0)$)
the fraction (resp. residue) field of $R$, $S:=\Spec R$ 
and $\eta$ (resp. $0$) the generic (resp. closed) point  of $S$. 
Let $(G_{\eta},\cL_{\eta})$ be a polarized 
abelian variety over $k(\eta)$ with $\cL_{\eta}$ 
symmetric ample cubical, 
$\cG$ the N\'eron model 
of $G_{\eta}$ 
and $G:=\cG^0$ the identity component of $\cG$. 
A triple $(P,i,\cN)$ is called a 
{\it relative compactification of $\cG$} 
(extending $(G_{\eta},\cL_{\eta})$) if 
\begin{enumerate}
\item[(rc1)]$P$ is an irreducible proper flat $S$-scheme;  
\item[(rc2)]
$i:\cG\hookrightarrow P$ is an open immersion with  
$P_{\eta}=i(\cG_{\eta})=i(G_{\eta})$;
\item[(rc3)]$\cN$ is an ample invertible $\cO_P$-module 
with $i^*\cN_{\eta}\simeq\cL^{\otimes l}_{\eta}$ for some $l\geq 1$. 
\end{enumerate}

\cite[3.5]{K98} constructs relative compactifications of a given semiabelian  
N\'eron model over $R$. 
His compactifications $P$ 
are {\it regular} with $\codim_P(P\setminus i(\cG))_{\red}=1$ in general, 
hence $P_0:=P\times_S0$ can have extra components other than 
the closures of components of $\cG_0$. 

In contrast with it, 
we construct a {\it unique but possibly singular} 
relative compactification of a given semiabelian N\'eron model. 
Our work consists of two parts, Part I (this article) 
and Part II (\cite{Nakamura24}). Part I (resp. Part II) 
treats the totally degenerate case (resp. 
the partially degenerate case). This article proves  
the following 
\begin{thm}\label{thm:main thm}\ 
If $\cG$ is totally degenerate semiabelian  over $S$,  
then there exists a relative compactification 
$(P,i,\cN)$ of $\cG$ extending $(G_{\eta},\cL_{\eta})$
 such that 
\begin{enumerate}
\item[(a)] $P$ is Cohen-Macaulay;
\item[(b)] $i(\cG)=P\setminus\Sing(P_0)$ with   
$\codim_{P}\Sing(P_0)\geq 2$ where $\Sing(P_0)$ denotes 
 the singular locus of $P_0$;
\item[(c)] $i^*\cN$ is cubical;
\item[(d)] 
$\cG$ acts on $P$ so that $i$ is $\cG$-equivariant.
\end{enumerate}
\end{thm}

We call $(P,i,\cN)$ in Theorem~\ref{thm:main thm}
a {\it cubical compactification of $\cG$.} 
Theorem~\ref{thm:main thm} is proved 
in \S\S~\ref{sec:relative compactifications}\,-\ref{sec:nonsplit case}. 
A similar theorem in the partially degenerate case  
and the uniqueness of $P$ satisfying the conditions (a)-(c) 
are proved in \cite{Nakamura24}.   
\cite[\S~7]{Mumford72} gives an example of 
$(P,i,\cN)$ with $i^*\cN$ non-cubical. 
It would be interesting to construct $(P,i,\cN)$
as geometric moduli like \cite{OS79}. 
See Remarks~\ref{rem:no cubic str on cLflat on cG}-\ref{rem:Jac_phi}.

This article (Part I) is organized as follows. 
In \S~\ref{sec:preliminaries}, 
we briefly review comparison theorems of 
local rings and their completions from \cite[IV$_2$]{EGA}, 
N\'eron models from \cite{BLR90} 
and related notions and results.   
From \S\S~\ref{sec:FC data}\,-\ref{sec:nonsplit case}, 
we consider a semiabelian scheme $\cG$ 
over a CDVR 
such that $\cG_0$ has a trivial abelian part. 
In \S\S~\ref{sec:FC data}-\ref{sec:Mumford fam}, 
we recall  
Faltings-Chai's degeneration data (FC data)  
and (simplified) 
Mumford families 
of degenerating abelian varieties 
from \cite{FC90} and \cite{AN99}. \par 
In \S~\ref{sec:graded algebras}, 
starting from a Mumford family 
$(P,\cL)$ associated with an FC datum $\xi$, 
we introduce an eFC datum $\xi^e$
extending $\xi$, two graded algebras 
$R^{\quot}(\xi^e)$ and $R^{\sharp}(\xi^e_{2N})$, and secondly 
we define a N\'eron-FC kit $\xi^{\natural}$\  
(abbr. a N\'eFC kit) of $\xi^e_{2N}$ and   
an infinite series of graded algebras 
$R_{l}(\xi^{\natural})$\ $(l=1,2,\cdots)$  
naturally associated with $\xi^{\natural}$. 
In \S~\ref{sec:Voronoi polytopes}, 
we study some combinatorics related to $R_{l}(\xi^{\natural})$. 
In \S~\ref{sec:rel complete models}, we construct the 
projective models $(P_{l},\cL_{l})$ 
associated with $\xi^{\natural}$, and will see 
that $(P_{l},\cL_{l})$ 
is independent of $l$ (Theorem~\ref{thm:scheme normalization isom})    
if the Voronoi polytope $\Sigma_{l}(0)$ is integral    
(Definition~\ref{defn:integral Sigma-Voronoi cell}). 

In \S~\ref{sec:str of Pl}, we describe the structure 
of ``\,the semi-universal covering $\tP_{l}$ 
of $P_{l}$\,'' in terms of $\Sigma_{l}(0)$. 
In \S~\ref{sec:relative compactifications}, we prove 
that $P_{l}$ is a relative compactification 
of the pullback $G_{l}$ of the N\'eron model of an abelian variety 
$G_{\eta}$ if $\Sigma_{l}(0)$ is integral
(Theorem~\ref{thm:summary of Pl Gl}). We prove that $(P_{l},\cL_{l})$ 
satisfies the conditions (a)-(c) of Theorem~\ref{thm:main thm}.
In \S~\ref{sec:action of Gl}, we prove that 
$G_{l}$ acts on $P_{l}$ extending the group law of $G_{l}$. 
The proof uses the theory of rigid analytic spaces, 
of which we give a brief survey  in \S~\ref{subsec:rigid analytic spaces} 
based on \cite{Bosch14}. 
\S\S~\ref{sec:proof of main thm}\,-\ref{sec:nonsplit case} are devoted 
to proving Theorem~\ref{thm:main thm} 
after a brief review on 
Galois descent based on \cite{Conrad} and \cite{W79}. 
In \S~\ref{sec:examples}, we give four examples for comparison. 

See \S\S~\ref{subsec:notation convention}/\ref{subsec:notation for torus emb}/\ref{subsec:theorems notation}
 and Notation~\ref{notation:Notation Rinit Sinit keta=Kmin(xi)} 
for general notation.

\medskip
{\small
{\it Acknowledgement.} We are very grateful to Professor T.~Suzuki for his advices during the preparation of this article.}
\section{Preliminaries}
\label{sec:preliminaries}

\subsection{The conditions (R$_k$) and (S$_k$)}
\label{subsec:Rk Sk}

Let $A$ be a noetherian ring. We recall the conditions 
\cite[(17.I)]{Matsumura70}: 
\begin{align*}
(\op{R}_k)&\qquad A_{p}\ \text{is regular 
for any prime ideal $p$ of $A$ with $\height(p)\leq k$};\\
(\op{S}_k) &\qquad \depth(A_{p})\geq \inf(k,\op{ht}(p))
\quad \text{for any prime ideal $p$ of $A$}.
\end{align*}

\begin{rem}\label{rem:normal/regular etc}
(\cite[17.I)]{Matsumura70})
Let $A$ be any noetherian ring. 
Then 
\begin{enumerate}
\item $A$ is reduced iff $A$ is (S$_1$) and (R$_0$); 
\item
$A$ is normal iff $A$ is (S$_2$) and (R$_1$); 
\item
$A$ is regular iff $A$ is (R$_k$) for all $k$; 
\item
$A$ is Cohen-Macaulay iff $A$ is (S$_k$) for all $k$. 
\end{enumerate}
\end{rem}

\begin{defn}Let $A$ be a noetherian ring. The  ring 
$A$ is called {\it catenary} if for any prime ideals 
$p$ and $q$ with $q\subset p$, we have $\height(p)=\height(q)+\height(p/q)$.  
The ring $A$ is called {\it universally catenary} if any polynomial $A$-algebra $A[x_1,\cdots,x_n]$ is catenary. See \cite[IV$_2$, 5.6.1/5.6.2]{EGA}. \par
The ring $A$ is called
{\it excellent}  \cite[(34.A)]{Matsumura70} if 
\begin{enumerate}
\item[(i)] $A$ is universally catenary;
\item[(ii)]
the formal fiber of the canonical morphism $\Spec A^{\wedge}_p\to \Spec A_p$ is  geometrically regular for any prime ideal $p$ of $A$ where $A^{\wedge}_p$ is the $pA_p$-adic completion of $A_p$;
\item[(iii)] for any finitely generated algebra $B$ over $A$, the set $\Reg(\Spec B)$ of regular points of $\Spec B$ is open in $\Spec B$.
\end{enumerate} 
\end{defn}

\begin{prop}
\label{prop:basics of excellent}{\rm (\cite[IV$_2$, 7.8.3]{EGA})}
The following are true:
\begin{enumerate}
\item any noetherian complete local ring is excellent;
\item if $A$ is an excellent ring, then so is  
the fraction ring $S^{-1}A$ for any multiplicatively closed 
subset $S$ of $A$;
\item if $A$ is an excellent ring, so is 
any $A$-algebra of finite type;
\item 
let $A$ be an excellent local ring, $I$ an ideal of $A$, and 
$A^{\wedge}$ the $I$-adic completion of $A$; then
$A$ is (R$_k$) (resp.  (S$_k$)) iff $A^{\wedge}$ is (R$_k$) (resp.  (S$_k$)).
\end{enumerate}
\end{prop}

\begin{cor}
\label{cor:if Chat=Bhat Rk Sk}
Let $R$ be a noetherian complete local ring, $I$ an ideal of $R$ 
and $A$ an $R$-algebra of finite type. Let $p$ be 
a prime ideal of $A$ containing $I$, $B:=A_p$, $C$ 
an excellent local $R$-algebra and  
$B^{\wedge}$ (resp. $C^{\wedge}$) 
the $I$-adic completion of $B$ (resp. $C$). 
Then 
\begin{enumerate}
\item $B$ is excellent;
\item 
$B$ is (R$_k$) (resp. (S$_k$)) iff $B^{\wedge}$ is (R$_k$) (resp. (S$_k$));
\item  if $C^{\wedge}\simeq B^{\wedge}$,  then 
$B$ is (R$_k$) (resp. (S$_k$)) iff $C$ is (R$_k$) (resp. (S$_k$)).
\end{enumerate}
\end{cor}
\begin{proof}By Proposition~\ref{prop:basics of excellent}~(1)-(3), 
the rings $R$, $A$ and $B$ are excellent. 
The rest follows from 
Proposition~\ref{prop:basics of excellent}~(4).
\end{proof}

\subsection{N\'eron models}
\label{sec:Neron models}
Here we quote a few standard facts about N\'eron models from 
\cite{BLR90}. 
Let $R$ be a Dedekind domain, $S:=\Spec R$,  
$\eta$ the generic point of $S$ and $k(\eta)$ the fraction field of $R$.
\begin{defn}
\label{defn:Neron mapping prop} {\rm (\cite[1.2/1]{BLR90})}
Let $X_{\eta}$ be a smooth separated $k(\eta)$-scheme of finite type. 
An $S$-scheme $X$ is defined to be a {\it N\'eron model of $X_{\eta}$} if 
$X$ is a smooth separated $S$-scheme of finite type,   
which has the following property, called N\'eron mapping property:

{\sl for each smooth $S$-scheme $Y$ and each $k(\eta)$-morphism $u_{\eta}:Y_{\eta}\to X_{\eta}$, there is a unique $S$-morphism $u:Y\to X$ 
extending $u_{\eta}$.} 
\end{defn}

\begin{prop}\label{prop:neron model of gr scheme}
{\rm (\cite[1.2/6]{BLR90})}
Let $X$ be a N\'eron model 
of a group $k(\eta)$-scheme $X_{\eta}$.  
Then the 
group $k(\eta)$-scheme structure of $X_{\eta}$ extends 
 uniquely to a group $S$-scheme structure  
of $X$. 
\end{prop}

\begin{cor}\label{cor:morphism of neron model of gr scheme}
Let $X_{\eta}$ be a group $k(\eta)$-scheme with N\'eron model $X$, $Y$ a smooth group $S$-scheme and $f\colon Y\to X$ an $S$-morphism such that $f\times_S\eta$ is a morphism of group $k(\eta)$-schemes. 
Then $f$ is a morphism of group $S$-schemes.
\end{cor}
\begin{proof}
Let $m_Y:Y\times_SY\to Y$ (resp. $m_X:X\times_SX\to X$) 
be multiplication of $Y$ (resp. $X$) and 
$f_{\eta}:=f\times_S\eta$. Note that $m_X$ exists
by Proposition~\ref{prop:neron model of gr scheme}.
By assumption, 
we have $m_{X,\eta}\circ(f_{\eta}\times_Sf_{\eta})
=f_{\eta}\circ m_{Y,\eta}$.   
Hence we have $m_{X}\circ(f\times_Sf)=f\circ m_{Y}$, which proves Corollary.
\end{proof}

\begin{thm}{\rm (\cite[1.2/8, 1.3/2]{BLR90})}
\label{thm:existence of neron model of AV}The following are true:
\begin{enumerate}
\item any abelian variety over $k(\eta)$ admits a N\'eron model;
\item any abelian scheme $X$ over $S$ is a N\'eron model of its generic fiber. 
\end{enumerate}
\end{thm}

\begin{prop}
\label{prop:semi abelian into neron}{\rm (\cite[7.4/3]{BLR90})}
Let $A_{\eta}$ be an abelian variety over $k(\eta)$ 
with N\'eron model $A$, and 
let $G$ be a semiabelian 
 $S$-scheme with connected fibers such that $G_{\eta}\simeq A_{\eta}$.  
Then the canonical morphism $G\to A$ (induced from N\'eron mapping property) is an open immersion; it is an isomorphism between their identity components.
\end{prop}

\subsection{Notation and convention}
\label{subsec:notation convention}
\begin{notation}\label{not:notation convention}
In what follows, we use the following notation 
unless otherwise mentioned. 
Let $R$ be a (C)DVR, 
$I$ the maximal ideal of $R$, $S:=\Spec R$, 
$\eta$ the generic point of $S$, $0$ the closed point of $S$, 
$k(\eta)$ the fraction field of $R$, 
$k(\eta)^{\times}:=k(\eta)\setminus\{0\}$, 
$\Omega:=\overline{k(\eta)}$ an algebraic  
closure of  $k(\eta)$, $R_{\Omega}$ the integral closure of $R$ in $\Omega$, 
$R^{\times}$ (resp. $R_{\Omega}^{\times}$) the group of units in $R$ 
(resp. $R_{\Omega}$), 
and $G_{\abs}:=\Aut(\Omega/k(\eta))$. 
Let $k(0):=R/I$ and let 
$\overline{k(0)}$ be the residue field of $R_{\Omega}$ and 
$k(0)^{\sep}$ the separable closure of $k(0)$ in $\overline{k(0)}$.
Let $s$ be a uniformizer of $R$ and 
$v_s:k(\eta)^{\times}\to\bZ$ 
the valuation of $k(\eta)$ with $v_s(s)=1$. 

Let $K$ be an extension field of $k(\eta)$. Then 
we denote the integral closure of $R$ in $K$ by $R_K$ 
and set $S_K:=\Spec R_K$.  

Let $Z$ be any (formal) $S$-scheme, $\cL$ any invertible sheaf on $Z$ 
and $N$ any $R$-module. For any $S$-scheme $T$, we denote $Z\times_ST$\  
 (resp. $Z\times_S0$,  $Z\times_S{\eta}$, $\cL\times_S0$, $\cL\times_S{\eta}$) 
by $Z_T$\ (resp. $Z_0$, $Z_{\eta}$, $\cL_0$, $\cL_{\eta}$). 
We also denote their $I$-adic completions 
by $Z^{\wedge}$, $\cL^{\wedge}$ and $N^{\wedge}$ 
(or occasionally $N^{I\op{-adic}}$) respectively. 
For an $S$-morphism  $f:Z\to Z'$ of $S$-schemes, 
let $f_0:=f\times_S0$ and $f_{\eta}:=f\times_S\eta$. 
We denote the $I$-adic completion of $f$ 
by $f^{\wedge}:Z^{\wedge}\to (Z')^{\wedge}$. \par
Let $\bN$ be the set of all positive integers.
\end{notation}

\begin{convention}\label{conv:convention}
\ \ Throughout this article, we 
assume that $(G_{\eta},\cL_{\eta})$ is a polarized abelian variety 
with $\cL_{\eta}$ ample and rigidified in the sense that 
a $k(\eta)$-isomorphism ($=$ a rigidification) 
$\sigma_{\eta}:e_{\eta}^*\cL_{\eta}\simeq \bG_{m,\eta}$ is given  
for the unit section  $e_{\eta}:\Spec k(\eta)\to G_{\eta}$. 
Similarly we mean 
by a {\em  semi-abelian $S$-scheme} a (polarized) semiabelian 
$S$-scheme $\pi:(G,\cL)\to S$ with $\pi$-ample  
invertible sheaf $\cL$ 
on $G$ such that $(G_{\eta},\cL_{\eta})$ is a polarized  
abelian variety over $k(\eta)$ and $\cL$ is rigidified in the sense that 
an $S$-isomorphism $\sigma:e^*\cL\simeq\bG_{m,S}$ is given 
for the unit section $e:S\to G$. 
We say that $\cL$ is ample instead of $\pi$-ample 
if no confusion is possible. We denote 
a $\bG_{m}$-torsor $\cL^{\times}:=\cL\setminus \{\text{zero section}\}$ 
associated with $\cL$ by the same $\cL$ if no confusion is possible. 
\end{convention}

\subsection{The component group of a N\'eron model}
\label{subsec;ample cubical sheaf over cG}
Let $(G_{\eta},\cL_{\eta})$ be an abelian variety of dimension $g$ 
over $k(\eta)$, $\cG$ the {\it N\'eron model of $G_{\eta}$} 
and $e$ the unit of $\cG$.  
It is clear that $\cG_{\eta}=G_{\eta}$.

Since $\cG$ is 
by definition \cite[1.2/1]{BLR90} 
smooth separated and of finite type over $S$,  
there exists,  by \cite[IV$_3$, (15.6.5), p.~238]{EGA}, an open subgroup 
$S$-scheme $\cG^0$ of $\cG$ such that every geometric fiber of 
$\cG^0$ is connected. 
Let $G:=\cG^0$. We call it the {\it connected N\'eron model 
of} $G_{\eta}$. 

Let $\Phi:=\cG_0/G_0$ and 
$N:=|\Phi|=|\Phi(\overline{k(0)})|$.
 Then $\Phi$ is a finite \'etale $k(0)$-scheme, 
so that $\Phi(\overline{k(0)})=\Phi(k(0)^{\sep})$ by \cite[6.2 Theorem]{W79}.
We call $\Phi$ the {\it component group of $\cG$.}

\subsection{Cubical structures}
\label{sec;cubical str}

Let   
$X$ be a commutative group $S$-scheme, 
$X^n:=X\times_SX\times_S\cdots \times_SX$ ($n$-times), 
$e_{X^n}$ the unit of $X^n$  
 and $\cL$ an invertible sheaf on $X$, which is viewed as a  
$\bG_{m}$-torsor over $X$. In what follows, we switch freely 
between invertible sheaves and $\bG_{m}$-torsors over $X$. 
 Let $m:X^2\to X$ be 
multiplication of $X$ and $e:=e_{X^1}$. 
For any subset $I$ of $\{1,2,3\}$, we denote by $m_I:X^3\to X$ 
the morphism  sending $(x_1,x_2,x_3)\mapsto \sum_{i\in I}x_i$. 
We also denote by $p_i:X^2\to X$ the $i$-th projection $(i=1,2)$.
We define 
\begin{align*}
\Theta(\cL)&=
\bigotimes_{\emptyset\neq I\subset \{1,2,3\}}m_I^*(\cL)^{(-1)^{|I|+1}},\\
\Lambda(\cL)&=m^*(\cL)\otimes p_1^*(\cL)^{-1}
\otimes p_2^*(\cL)^{-1}
\end{align*}  
where $\Theta(\cL)=\cD_3(L)$ and $\Lambda(\cL)=\cD_2(\cL)$ 
in \cite[I, 2.1]{MB85}. By \cite[I, 1.2.4.1]{MB85}, if 
$\Theta(\cL)\simeq\bG_{m,X^3}$, 
then $\Lambda(\cL)$ admits a structure of a biextension.
\begin{defn}\label{defn:cubical str}
A {\it cubical structure of $\cL$} is an isomorphism 
$$\tau:=\tau_{l}:\Theta(\cL)\simeq\bG_{m,X^3}
$$such that $\Lambda(\cL)$ is a {\it symmetric 
biextension} of $X\times_SX$ by $\bG_{m,S}$.  
See \cite[2.2]{Breen83}/\cite[I, 2.5.4]{MB85}.  
There is a canonical isomorphism 
$e^*_{X^3}\Theta(\cL)\simeq e^*\cL$ by \cite[2.1.3]{Breen83}/
\cite[I, 2.4.1]{MB85}. 
Hence a cubical structure of $\cL$ induces 
a rigidification of $\cL$ (along $e$). \par
Let $\cM$ and $\cN$
 be cubical invertible sheaves on $X$. 
A {\it cubical homomorphism from $\cM$ to $\cN$} is defined to be 
a homomorphism $f:\cM\to \cN$ of $\cO_X$-modules such that 
$\tau_{\cM}=\tau_{\cN}\circ\Theta(f)$. 
We denote by $\Hom_{\CUB}(\cM,\cN)$ the set of all cubical homomorphisms 
 from $\cM$ to $\cN$.
See \cite[I, 2.4.5]{MB85}. 
\end{defn}
\begin{lemma}\label{lemma:cub str}
Let $R$ be a discrete valuation ring and $k(\eta)$ its fraction field. 
Let $(G_{\eta},\cL_{\eta})$ be an abelian variety 
with $\cL_{\eta}$ rigidified, $\cG$ 
the N\'eron model of $G_{\eta}$, $G:=\cG^0$, $\cG^{\sharp}$ an open subgroup 
$S$-scheme of $\cG$ with $G\subset\cG^{\sharp}$, 
$N^{\sharp}:=|\cG^{\sharp}_0/G_0|$ 
and $n\in\bN$ divisible by $N^{\sharp}$.  
Then  
\begin{enumerate}
\item\label{item:cubic str on G} 
any rigidified invertible sheaf on $G$ has a unique 
cubical structure;
\item\label{item:extension to G} $\cL_{\eta}$ is uniquely extended 
to an ample cubical invertible sheaf on $G$;
\item\label{item:unique cubic extension to cG/Geta} 
$\Hom_{\CUB}(\cM,\cN)\simeq\Hom_{\CUB}(\cM_{\eta},\cN_{\eta})$ for any cubical invertible sheaves $\cM$ and $\cN$  on $\cG^{\sharp}$;
\item\label{item:cubic extension to cG/Geta} 
there exists a unique ample cubical invertible sheaf $\cM$ 
on $\cG^{\sharp}$ such that 
$\cM_{\eta}=\cL^{\otimes 2n}_{\eta}$ if $n$ is even, 
while $\cM_{\eta}=\cL^{\otimes n}_{\eta}$ if $n$ is odd.  
\end{enumerate}
\end{lemma}
\begin{proof}See \cite[I, 2.6/II, 1.1~(i)]{MB85} for 
(\ref{item:cubic str on G})/(\ref{item:unique cubic extension to cG/Geta}).
Since $\Gamma(G_{\eta},\cL_{\eta})\neq 0$, there exists 
a nontrivial section $\sigma\in\Gamma(G_{\eta},\cL_{\eta})$. 
The closure of the zero locus $\sigma=0$ defines a divisor of $G$, which 
defines an invertible sheaf $\cL$ on $G$ extending $\cL_{\eta}$.  
Since $G$ is a smooth group $S$-scheme with connected fibers, 
$\cL$ admits a unique cubical structure by \cite[I, 2.6]{MB85}. 
By \cite[VI, 2.1]{MB85}, $\cL$ is ample on $G$. 
This proves (\ref{item:extension to G}).   By \cite[II, 1.2.1]{MB85}, 
there exists a cubical invertible sheaf $\cM$ 
on $\cG^{\sharp}$ such that 
$\cM_{\eta}=\cL_{\eta}^{\otimes 2n}$ if $n$ is even, 
while $\cM_{\eta}=\cL_{\eta}^{\otimes n}$ if $n$ is odd.   
By \cite[VI, 2.1]{MB85}, $\cM$ is ample on $\cG^{\sharp}$, while 
$\cM$ is unique by (\ref{item:unique cubic extension to cG/Geta}).
\end{proof}

\section{FC data}
\label{sec:FC data}

\subsection{Grothendieck's stable reduction}
\label{subsec:Grothendieck stable red}
Suppose that we are given a polarized abelian variety  
$(G_{\eta},\cL_{\eta})$ over $k(\eta)$ with $\cL_{\eta}$ 
ample cubical. 
Then by Grothendieck's stable reduction theorem
\cite[IX, Th.~3.6]{SGA7},   
$G_{\eta}$  can be extended 
to a semiabelian $S$-scheme $G$ of $G_{\eta}$ 
as the connected N\'eron model 
after replacing $k(\eta)$ by its finite 
extension $K'$ and $R$ by its integral closure $R'$ 
in $K'$ if necessary. 
\footnotemark\footnotetext{In what follows, we say this simply 
``by taking a finite extension of $k(\eta)$''. 
The integral closure of $R$ in any finite extension 
of $k(\eta)$ is again a CDVR. }  
Since $G_0$ is irreducible,  
$\cL_{\eta}$ extends uniquely to $G$ as an ample 
cubical invertible sheaf $\cL$ by 
Lemma~\ref{lemma:cub str}~(\ref{item:extension to G}).  
The closed fiber $G_0$ is a semiabelian variety over $k(0)$, namely
an extension of an abelian variety $A_0$ of dimension $g'$ 
by a torus $T_0$ of dimension $g''$ where $g'+g''=g$. 
See \cite[pp.~37-45]{FC90} for more details. 

\begin{rem}\label{rem:plan}
In \S\S~\ref{sec:FC data}-\ref{sec:examples}  
we limit ourselves to the totally degenerate case, 
that is, the case where $A_0=0$ and $g=g''$.  
By \cite[\S~7.3]{W79} and \cite[X, \S~1]{SGA3}, 
we may assume that 
$T_0$ is a split $k(0)$-torus:
\begin{equation}
\label{eq:T0 split torus}
T_0\simeq \bG_{m,k(0)}^{g''}
\end{equation}
by taking a finite unramified extension 
of $k(\eta)$ \footnote{We assume that its residue field extension is separable and any uniformizer of $R$ is a uniformizer of $R'$, {\it i.e.},  
the natural morphism $\varpi:S'\to S$ is \'etale where 
$S':=\Spec R'$.}\label{footnote:etale}
 if necessary. 
\end{rem}

\subsection{Faltings-Chai's degeneration data}
\label{subsec:Faltings-Chai deg data}

Since $G_0\simeq\bG_{m,k(0)}^g$ by (\ref{eq:T0 split torus}), 
$G^{\wedge}\simeq 
(\bG^g_{m,S})^{\wedge}$ by  
\cite[III, 2.8, p.~118]{SGA3}.
Let $\cL$ be the unique extension of $\cL_{\eta}$ to $G$. 
Since any line bundle 
on $(\bG^g_{m,S})^{\wedge}$ is trivial, so is 
$\cL^{\wedge}$. 
\par
 Let $X$ be a free $\bZ$-module of rank $g$ and $(w^x;x\in X)$ the coordinates 
of $(\bG^g_{m,S})^{\wedge}$ such that 
$w^{x+y}=w^x\cdot w^y$ $(\forall x,y\in X)$. It follows that 
any $\theta\in\Gamma(G^{\wedge},\cL^{\wedge})
=\Gamma(G^{\wedge},\cO_{G^{\wedge}})$ 
is an $I$-adically convergent formal power series in $w^x$. 
Since $\Gamma(G,\cL)$ is an $R$-free module of finite rank 
by Raynaud (see \cite[VI, 1.4.2]{MB85}), we have natural homomorphisms:
{\small\begin{equation}\label{eq:Gamma(Geta,Leta) as subspace}
\Gamma(G_{\eta},\cL_{\eta})\simeq\Gamma(G,\cL)\otimes_Rk(\eta)
\hookrightarrow\Gamma(G^{\wedge},\cL^{\wedge})\otimes_Rk(\eta)
\hookrightarrow\prod_{x\in X}k(\eta)\cdot w^x.
\end{equation}}

For $\theta\in\Gamma(G_{\eta}, \cL_{\eta})$,
$\theta\in\Gamma(G,\cL)$ iff 
$\theta\in \Gamma(G^{\wedge},\cL^{\wedge})$. Hence 
$\Gamma(G,\cL)$ is an $R$-submodule of 
$\Gamma(G_{\eta},\cL_{\eta})$ consisting 
of formal power series $\theta=\sum_{x\in X}\sigma_x(\theta)w^x$
with $\sigma_x(\theta)\in R$ $(\forall x\in X)$.

\begin{thm}
\label{thm:FC datum}
Suppose that $A_0=0$ and $G_0$ is 
a split $k(0)$-torus $T_0$ in \S~\ref{subsec:Grothendieck stable red}. 
Let $X$ be  
a free $\bZ$-module of rank $g$ such that 
$T_0\simeq\Spec k(0)[X]$. 
Then there exist a submodule $Y$ of $X$ of finite index, 
a function $a:Y\to k(\eta)^{\times}$ and 
a bilinear function $b:Y\times X\to k(\eta)^{\times}$ 
\footnote{To be exact, $b(y+y',x)=b(y,x)b(y',x)$ and 
$b(y,x+x')=b(y,x)b(y,x')$.} 
such that  
\begin{enumerate}
\item[(i)]\label{item:a0} 
$a(0)=1$;
\item[(ii)]\label{a,b} $b(y,z)=b(z,y)=a(y+z)a(y)^{-1}a(z)^{-1}$ 
$\  (\forall y,z\in Y)$;
\item[(iii)]\label{item:positive definite} $b(y,y)\in I \ (\forall y\in Y\setminus\{0\})$, 
and for 
every $n\geq 0$, $a(y)\in I^n$ for all but finitely many $y\in Y$; 
\item[(iv)]\label{item:k(eta) module} 
$\Gamma(G_{\eta},\cL_{\eta})$ is 
identified with the $k(\eta)$-vector subspace of formal Fourier series 
$\theta=\sum_{x\in X}\sigma_x(\theta)w^x$
such that 
\begin{gather*}
\sigma_{x+y}(\theta)
=a(y)b(y,x)\sigma_x(\theta),\ 
\sigma_x(\theta)\in k(\eta)\ (\forall x\in X, \forall y\in Y).
\end{gather*}
\end{enumerate}
\end{thm}

 Note that $Y$, $a$ and $b$ 
are uniquely determined by the conditions (i)-(iv).  
See \cite{R70} and \cite[pp.~37-42]{FC90}. 
See also \cite[pp.\ 232 - 233]{Nakamura16}.

\begin{rem}Since $\Gamma(G_{\eta},\cL_{\eta})$ is of rank $|X/Y|$ by 
Theorem~\ref{thm:FC datum}~(iv), 
$(G_{\eta},\cL_{\eta})$ is principally polarized iff $X=Y$. 
\end{rem}

\begin{defn}\label{defn:deg data of G L}
For a polarized semiabelian 
$S$-scheme $(G,\cL)$ with $G_0$ 
a $k(0)$-split torus, 
we define $\FC(G,\cL)$ by 
$$\FC(G,\cL)=(X,Y,a,b,A,B)$$ 
where $A(y):=v_s a(y)$ and $B(y,x):=v_s b(y,x)$\ 
$(\forall x\in X, \forall y\in Y)$.
We call $\FC(G,\cL)$ the {\it Faltings-Chai's degeneration datum of $(G,\cL)$}.
\end{defn}

\begin{lemma}
\label{lemma:FC(G_L) L symmetric}
If $\cL$ is symmetric in Theorem~\ref{thm:FC datum}, {\it i.e.}, 
$[-\id_G]^*\cL\simeq\cL$, then  
$a(y)=a(-y)$, $b(y,y)=a(y)^2$ and $A(y)=A(-y)=B(y,y)/2$ $(\forall y\in Y)$. 
In particular, $B$ is even on $Y\times Y$ in the sense that $B(y,y)$ is 
even  $(\forall y\in Y)$.
\end{lemma}
\begin{proof} Let $\theta\in\Gamma(G_{\eta},\cL_{\eta})$ and $\theta':=[-\id_G]^*\theta$. 
The Fourier 
expansion of $\theta$ is uniquely given 
by $\theta=\sum_{x\in X}\sigma_x(\theta)w^x$ 
such that $\sigma_{x+y}(\theta)=a(y)b(y,x)\sigma_x(\theta)$ $(\forall x\in X, \forall y\in Y)$. Meanwhile $\theta'=\sum_{x\in X}\sigma_x(\theta)w^{-x}\in\Gamma(G_{\eta},\cL_{\eta})$ because $[-\id_G]^*\cL\simeq\cL$ by assumption, so that $\sigma_x(\theta)=\sigma_{-x}(\theta')$ $(\forall x\in X)$. It follows 
\begin{equation}
\label{eq:sigma_x+y = ab sigma_x}
\sigma_{x+y}(\theta)=\sigma_{-x-y}(\theta')=a(-y)b(y,x)\sigma_{-x}(\theta')=a(-y)b(y,x)\sigma_{x}(\theta) 
\end{equation}for any $x\in X$ and $y\in Y$.  
By Eq.~(\ref{eq:sigma_x+y = ab sigma_x}), 
we obtain $a(-y)=a(y)$ and $A(y)=A(-y)$ $(\forall y\in Y)$. 
By Theorem~\ref{thm:FC datum}~(i)-(ii), 
$b(y,y)=a(y)a(-y)a(0)^{-1}=a(y)^2$. 
Hence $B(y,y)/2=A(y)=A(-y)$. This completes the proof. 
\end{proof}

\begin{rem}
\label{rem:symmetric L}
$\cL\otimes [-\id_{G}]^*\cL$ is symmetric and rigidified.
\end{rem}

\begin{defn}
\label{defn:FC datum}
A sextuple  
\begin{equation}\label{eq:FC datum}
\xi=(X, Y, a, b, A, B)\ \ ({\rm abbr.}\ \ (X,Y,a,b))
\end{equation}
is called a ({\it totally degenerate})  
{\it FC datum over $R$} (or {\it over $S$})
if the following conditions are satisfied:
\begin{enumerate}
\item[(i)] $X$ is a free $\bZ$-module; 
$Y$ is a $\bZ$-submodule of $X$ of finite index;
\item[(ii)] $b:Y\times X\to k(\eta)^{\times}$ 
is a bilinear function; 
\item[(iii)] $a:Y\to k(\eta)^{\times}$ is a function such that 
$$a(0)=1,\ b(y,z)=a(y+z)a(y)^{-1}a(z)^{-1}\quad (\forall y,z\in Y);$$
\item[(iv)] $A(y)=v_s a(y)$ and $B(y,x)=v_s b(y,x)$\ 
$(\forall x\in X, \forall y\in Y)$;
\item[(v)] $B$ is even symmetric and positive definite on $Y\times Y$.
\end{enumerate}

Note that $\FC(G,\cL)$ is an FC datum. 
If $a(y)=a(-y)$ $(\forall y\in Y)$, then 
we call $\xi$ {\it symmetric}. 
 If $\xi=\FC(G,\cL)$ with $\cL$ symmetric, 
then Lemma~\ref{lemma:FC(G_L) L symmetric} is true, so that 
$\xi$ is {\it symmetric}. 
If $X=Y$, then we call $\xi$ {\it principal} and denote it by $\xi=(X,a,b,A,B)$  or $\xi=(X,a,b)$.   
For a given FC datum $\xi$ in Eq.~(\ref{eq:FC datum}), 
we denote $a=a_{\xi}$, $b=b_{\xi}$, 
$A=A_{\xi}$ and $B=B_{\xi}$ if necessary. 
\end{defn}

\begin{defn}\label{defn:bara barb}
We define   
$$\bara(y)=a(y)s^{-A(y)}\in R^{\times},\ 
\barb(y,x)=b(y,x)s^{-B(y,x)}\in R^{\times}.$$ 
Then we have a $k(0)^{\times}$-valued bilinear form  
by Definition~\ref{defn:FC datum}~(iii):    
\begin{gather*}\label{cond:bilinear a0(x)}
\barb(y,z):=\bara(y+z)\bara(y)^{-1}\bara(z)^{-1}\ (\forall y, z\in Y).
\end{gather*}
\end{defn}

\begin{defn}
\label{defn:pullback FC datum}
Let $k(\eta')$ be a finite extension of $k(\eta)$, $R'$ the integral closure of $R$ in $k(\eta')$, $S':=\Spec R'$, $\varpi:S'\to S$ the natural morphism    
and {$t$ a uniformizer of $R'$.} 
For a given FC datum $\xi$ over $R$, 
we define the {\it pullback $\varpi^*\xi$ of $\xi$ by $\varpi$}  
(by abuse of notation $\varpi^*$) as follows:
\begin{gather*}
\varpi^*\xi=(X, Y, a', b', A', B'),\\
a'(y):=\varpi^*a(y),\ b'(y,x):=\varpi^*b(y,x),\\
 A'(y):=v_ta'(y),\ B'(y,x):=v_tb'(y,x)\ (\forall  x\in X, \forall y\in Y)
\end{gather*}where $\varpi^*:k(\eta)\to k(\eta')$ is the inclusion.  
Then $\varpi^*\xi$ 
is an FC datum over $R'$ with $A'(x)=v_t(s)A(x)$ and $B'(y,x)=v_t(s)B(y,x)$. 
\end{defn}

\begin{defn}
\label{defn:restriction of FC datum}
Let $\xi=(X,Y,a,b,A,B)$ be an FC datum and $X'$ (resp. $Y'$) 
a submodule of $X$ (resp. $Y$) of finite index such that $Y'\subset X'$. We define the {\it restriction $\xi'$ of $\xi$ to $(X',Y')$} to be an FC datum 
$(X',Y',a',b',A',B')$ over $R$ such that 
\begin{gather*}
a'(y):=a(y),\ b'(y,x):=b(y,x),\\
 A'(y):=A(y),\ B'(y,x):=B(y,x)\ (\forall x\in X', \forall y\in Y').
\end{gather*}
 
We denote $\xi'$ by $\res_{X',Y'}\xi$ or simply by $(X',Y',a,b,A,B)$. 
\end{defn}

\subsection{Fourier series (1)}
\label{subsec:Fouries 1}
Let $(G,\cL)$ be a polarized 
semiabelian $S$-scheme with $G_0$ a split $k(0)$-torus. 
The polarization morphism $\lambda(\cL_{\eta}^{\otimes m})=m\lambda(\cL_{\eta}):G_{\eta}\to G^t_{\eta}$ of $\cL_{\eta}^{\otimes m}$ extends to $G$  by Proposition~\ref{prop:semi abelian into neron}, which we denote by $\lambda(\cL^{\otimes m})=m\lambda(\cL):G\to G^t$. Since 
\begin{gather*}
G_0=\Spec k(0)[w^x;x\in X],\ 
G_0^t=\Spec k(0)[w^y;y\in Y],
\end{gather*} a monomorphism 
$\phi_{\cL^{\otimes m}}:Y\to X$ is induced from 
$m\lambda(\cL)$. 
Taking the above remark into account, we obtain, 
by Theorem~\ref{thm:FC datum},
\begin{cor}\label{cor:am bm}
There exist a function 
$a_m:mY\to k(\eta)^{\times}$ and 
a bilinear function $b_m:mY\times X\to k(\eta)^{\times}$ 
such that  
\begin{enumerate}
\item[(i)]\label{item:am(0)} 
$a_m(0)=1$;
\item[(ii)]\label{item:am,bm} 
$b_m(my,mz)=a_m(my+mz)a_m(my)^{-1}a_m(mz)^{-1}$ $(\forall y,z\in Y)$;
\item[(iii)]\label{item:positive definite Bm} 
$b_m(my,my)\in I \ (\forall y\ne 0)$, 
and for 
every $n\geq 0$, $a_m(my)\in I^n$ for all but finitely many $y\in Y$; 
\item[(iv)]\label{item:k(eta) module Gamma_Lm} 
$\Gamma(G_{\eta},\cL_{\eta}^{\otimes m})$ is 
identified with the $k(\eta)$-vector subspace of formal Fourier series 
$\theta=\sum_{x\in X}\sigma^{(m)}_x(\theta)w^x$
such that 
{\small \begin{gather*}\quad\quad
\sigma^{(m)}_{x+my}(\theta)
=a_m(my)b_m(my,x)\sigma^{(m)}_x(\theta),\  
\sigma^{(m)}_x(\theta)\in k(\eta)\ 
(\forall x\in X, \forall  y\in Y).
\end{gather*}}
\end{enumerate}
\end{cor}

\begin{lemma}\label{lemma:Gamma(Geta,Leta^m m=2^l)}
Let $m:=2^{l}$ $(l\geq 0)$. Then 
with the notation in Theorem~\ref{thm:FC datum},  
$$a_m(my)=a(y)^m,\ b_m(my,x)=b(y,x)\ \ (\forall x\in X, \forall y\in Y). 
$$
\end{lemma}
\begin{proof}
By \cite[p.~42, line 27]{FC90}, we see
$$a_2(2y)b_2(2y,x)=a(y)^2b(y,x)\ (\forall x\in X, \forall y\in Y). 
$$

This is the case $l=1$. We obtain Lemma by the induction on $l$. 
\end{proof}

\begin{thm}
\label{thm:G isom G' iff FC same}
Let $(G,\cL)$ (resp. $(G',\cL')$) be 
a polarized semiabelian $S$-scheme with $G_0$ 
(resp. $G'_0$) a split $k(0)$-torus. 
If $\FC(G,\cL)=\FC(G',\cL')$, then 
$(G,\cL)\simeq (G',\cL')$ as polarized group $S$-schemes.
\end{thm}
\begin{proof}Let $\FC(G,\cL)=\FC(G',\cL')=(X,Y,a,b,A,B)$. 
Since $G_0$ and $G'_0$ are split $k(0)$-tori, 
\begin{equation}
\label{eq:isom Gformal G'formal TXformal}
G^{\wedge}\simeq \Spf (R[u^x; x\in X])^{\wedge},
\ (G')^{\wedge}\simeq\Spf (R[v^x; x\in X])^{\wedge}
\end{equation} 
for some lattice $X$ by \cite[IX, \S~3]{SGA3}. 
Let $e$ (resp. $e'$) be the unit section of $G$ (resp. $G'$) 
over $S$, each being defined by $u^x=1$\ (resp. $v^x=1$) $(\forall x\in X)$.

Let $m:=2^{l}$ $(l\geq 0)$. 
 By our assumption we can choose an isomorphism 
$\Phi:G^{\wedge}\simeq (G')^{\wedge}$ such that 
$\Phi^*(v^x)=u^x$ $(\forall x\in X)$. 
Let $V_m:=\Gamma(G_{\eta},\cL_{\eta}^{\otimes m})$ and 
$V'_m:=\Gamma(G'_{\eta},(\cL'_{\eta})^{\otimes m})$.  
By Theorem~\ref{thm:FC datum}~(iv) and 
Lemma~\ref{lemma:Gamma(Geta,Leta^m m=2^l)}, 
{\small 
\begin{equation}
\label{eq:Geta Letam explicit}
\begin{aligned}
V'_m
&=\left\{
\theta'=\sum_{x\in X}\sigma^{(m)}_x(\theta')v^x;
\begin{matrix}\sigma^{(m)}_{x+my}(\theta')
=a(y)^mb(y,x)\sigma^{(m)}_x(\theta') 
\\  
\sigma^{(m)}_x(\theta')\in k(\eta)\, (\forall x\in X, \forall y\in Y)
\end{matrix}\right\},\\
V_m
&=\left\{
\theta=\sum_{x\in X}\sigma^{(m)}_x(\theta)u^x;
\begin{matrix}\sigma^{(m)}_{x+my}(\theta)=a(y)^mb(y,x)
\sigma^{(m)}_x(\theta)
\\ 
\sigma^{(m)}_x(\theta)\in k(\eta)\ (\forall x\in X, \forall y\in Y)
\end{matrix}\right\}.
\end{aligned}
\end{equation}}

By Eq.~(\ref{eq:Geta Letam explicit}),  $\Phi$ induces an isomorphism
$\Psi_m:V'_m\to V_m$ defined by 
$\Psi_m(\theta')=\sum_{x\in X}\sigma^{(m)}_x(\theta')u^x$ 
for $\theta'=\sum_{x\in X}\sigma^{(m)}_x(\theta')v^x\in V'_m$.

As is well known, 
$\cL_{\eta}^{\otimes m}$ and $(\cL'_{\eta})^{\otimes m}$ 
are very ample for $m\geq 4$. Let $J_m$ be a set of  
representatives of $X/mY$.  
Then the function field of $G_{\eta}$ (resp. $G'_{\eta}$) 
is the field generated over $k(\eta)$ by 
the quotients $\theta_{\alpha}/\theta_{\beta}$ (resp. $\theta'_{\alpha}/\theta'_{\beta}$), where $\alpha$ and $\beta$ 
range over $J_m$, and 
\begin{gather*}
\theta_{\alpha}:=\sum_{y\in Y}a(y)^mb(y,\alpha)u^{\alpha+my},\ 
\theta'_{\alpha}:=\sum_{y\in Y}a(y)^mb(y,\alpha)v^{\alpha+my}.
\end{gather*}

Hence the function field of $G'_{\eta}$ is isomorphic to that of $G_{\eta}$, that is,  $G_{\eta}$ is birationally equivalent to 
$G'_{\eta}$. Since $\phi(e^{\wedge})=(e')^{\wedge}$,  
$\phi$ induces a birational map 
$f:G_{\eta}\to G'_{\eta}$ such that $f(e_{\eta})=e'_{\eta}$. 
Suppose that $f$ is not a morphism. From the final blowing-up in the process of elimination of indeterminacy of $f$, we have a rational curve of 
the blown-up $G_{\eta}$ which is sent to a rational curve of $G'_{\eta}$ \cite[p.~600]{Mori79}. 
This is a contradiction because the abelian variety $G'_{\eta}$ 
contains no rational curves. Hence $f$ is an isomorphism. 

Let $\cG$ (resp. $\cG'$) be the N\'eron model of $G_{\eta}$ (resp. $G'_{\eta}$) whose identity component is $G$ (resp. $G'$). 
Since $f(e_{\eta})=e'_{\eta}$, by N\'eron mapping property, 
$f$ is extended to an $S$-isomorphism $h:\cG\to\cG'$ with $h(e)=e'$. Hence 
the restriction of $h$ to $G$ is an isomorphism between $G$ and $G'$, 
which we denote by the same $h$. 
Hence there is an isomorphism 
$h^{\wedge}:G^{\wedge}\simeq (G')^{\wedge}$ 
with $h^{\wedge}(e^{\wedge})=(e')^{\wedge}$ which induces $\psi_m$. 
This implies that 
$h^*:V'_m\to V_m$ is an isomorphism $(m=2^{l}, \forall l\geq 0)$.  
 Since $\Gamma(G'_{\eta},\cL'_{\eta})\neq\{0\}$, 
$h$ induces an isomorphism 
$(G_{\eta},\cL_{\eta})\simeq (G'_{\eta},\cL'_{\eta})$. 
Since $G_0$ (resp. $G_0'$) is irreducible, 
the extension of $\cL_{\eta}$ (resp. $\cL'$)  to $G$ (resp. $G'$) is unique, so that $f:(G,\cL)\to (G',\cL')$ is an isomorphism as polarized group $S$-schemes.  \end{proof}

\begin{rem}Theorem~\ref{thm:G isom G' iff FC same} is true 
also when $G$ has nontrivial abelian part. 
See \cite[3.7]{Nakamura24}.
\end{rem}

\section{The simplified Mumford families}
\label{sec:Mumford fam}
\subsection{Notation}
\label{subsec:notation for torus emb}
Let $X$ be a free $\bZ$-module of rank $g$, 
$X^{\vee}:=\Hom_{\bZ}(X,\bZ)$, $X_{\bR}:=X\otimes_{\bZ}\bR$, 
$X_{\bR}^{\vee}:=X^{\vee}\otimes_{\bZ}\bR$ and 
$X\times X^{\vee}\ni (x,u)\mapsto u(x)\in \bZ$ the natural pairing.  
Let $(m_i;i\in [1,g])$ be a basis of $X$ and $(f_i;i\in [1,g])$ a dual basis 
with $f_i(m_j)=\delta_{i,j}$.  
A split $\bZ$-torus is by definition
{\small$$T_X:=\Spec \bZ[w^x;x\in X]=\Spec \bZ[w_k^{\pm 1}; k\in[1,g]]$$}
\hspace*{-0.15cm}for some free $\bZ$-module $X$,   
where $w_k$ denotes the monomial $w^{m_k}$ 
corresponding to $m_k$, and we denote $\prod_{k=1}^gw_k^{x_k}$ by $w^x$ 
for $x=\sum_{k=1}^gx_km_k\in X$. 

For a rational convex cone $\sigma$ 
of $X_{\bR}^{\vee}$, we define a $T_X$-embedding by
\begin{equation*}
\label{eq:TMsigma} V(\sigma):=T_X(\sigma):=\Spec \bZ[\sigma^{\vee}\cap X].
\end{equation*} 
We also define a torus embedding for a fan $\cF$ in $X^{\vee}$ 
following \cite[Chap.~1]{Oda85} and \cite[Chap.~I]{TE73}, which we denote by $V(\cF)$.
\par

Let $m_0:=1$ be a basis of $\bZ$, $\tilde{X}:=\bZ m_0\oplus X$ and 
$f_0\in \tilde{X}^{\vee}$ with $f_0(m_0)=1$ 
and $f_0(x)=0$ $(\forall x\in X)$. 
Let  $\langle\phantom{\ \cdot},\phantom{\,\cdot}\rangle:
\tilde{X}\times \tilde{X}^{\vee}\to\bZ$ be the natural pairing
\begin{equation}
\label{eq:pairing of Msharp} 
\langle x_0m_0\oplus x, u_0f_0\oplus u\rangle=u_0x_0+u(x)
\quad (x\in X, u\in X^{\vee})
\end{equation} 
 and $\tilde{X}_{\bR}^{\vee}:=\tilde{X}^{\vee}\otimes_{\bZ}\bR$.  
Let $w_0:=w^{m_0}$  be  
the monomial corresponding to $m_0$ and $w^x=\prod_{k=0}^gw_k^{x_k}$ for 
$x=\sum_{k=0}^gx_km_k\in\tX$. 
We define  
$$T_{\tilde{X}}=\Spec \bZ[w^x;x\in \tX]=\Spec \bZ[w_k^{\pm 1}; k\in[0,g]].$$

\subsection{Fans and torus embeddings}
\label{subsec:fan and torus emb}
Let $R$ be a CDVR, $s$ 
its uniformizer and $S:=\Spec R$. 
Let $\iota:S\to \Spec\bZ[w_0]$ be the morphism with $\iota^*w_0=s$. 
In what follows, we always identify $w_0=s$ by $\iota^*$. 
Let 
$\sigma$ be any rational convex cone of $\tX^{\vee}_{\bR}$.  
Let $T_{\tilde X}(\sigma):=\Spec\bZ[\sigma^{\vee}\cap\tX]$ and  
\begin{equation*}
 U(\sigma):=\iota^*T_{\tilde X}(\sigma)=\Spec R[\sigma^{\vee}\cap \tilde{X}].
\end{equation*}

We call $U(\sigma)$ an {\it $S$-torus embedding 
associated with $\sigma$.}  For a fan $\cF$ in $\tX^{\vee}$, we denote 
the pullback by $\iota$ of  
$T_{\tilde{X}}(\cF)$ by $U(\cF)$, which we call  an {\it $S$-torus embedding 
associated with $\cF$.}

\begin{defn}
\label{defn:fan over S}
We call a fan $\cF$ in $\tX^{\vee}$ a {\it fan over $S$} ({\it or} $R$) 
if any $\sigma\in\cF$ satisfies the following conditions: 
\begin{enumerate}
\item[(i)] $\sigma^{\vee}\ni m_0$\ \ ($\Leftrightarrow \sigma\subset\bR_{\geq 0}f_0+X^{\vee}_{\bR}$);
\item[(ii)] $\sigma\cap X^{\vee}=(0)$\ \ ($\Leftrightarrow \sigma\cap X^{\vee}_{\bR}=(0)$);  
\item[(iii)] $\bZ m_0+\sigma^{\vee}\cap\tX=\tX$\ \ ($\Leftrightarrow \bR m_0+\sigma^{\vee}=\tX_{\bR}$);
\item[(iv)] $U(\sigma)_{\eta}\simeq T_{X,k(\eta)}:=\Spec k(\eta)[X]$.
\end{enumerate}
 
\end{defn}

\begin{lemma}Under the condition (i) 
of Definition~\ref{defn:fan over S}, 
the three conditions (ii),(iii) 
and (iv) are equivalent to each other. 
\end{lemma}
\begin{proof}
Note that {\small$U(\sigma)_{\eta}
=\Spec R[\sigma^{\vee}\cap \tX][1/s]=\Spec R[\bZ m_0+\sigma^{\vee}\cap \tX]
$.}  Hence  (iii) $\Leftrightarrow$ {\small$U(\sigma)_{\eta}=\Spec R[\tX]\simeq T_{X,k(\eta)}$} $\Leftrightarrow$ (iv).  \par
Next, we assume (ii). If $\sigma=(0)$, 
(iii) is clear. So we assume $\sigma\neq (0)$. 
Let $z\in\tX$ and $\sigma$ the convex closure of $(u_i; i\in I)$ 
for some $u_i\in \tX^{\vee}$. 
By (i) and (ii), $u_i=a_if_0+v_i$ for some $a_i>0$ and $v_i\in X^{\vee}$. 
Let $a:=\max(|u_i(z)|;i\in I)\ (\geq 0)$. Since 
$u_i(am_0+z)=aa_i+u_i(z)\geq a+u_i(z)\geq 0$, we have 
$am_0+z\in\sigma^{\vee}$. Hence $z\in \bZ m_0+\sigma^{\vee}\cap\tX$. 
This implies (iii). 

Finally, we prove that (iii) implies (ii). 
Let $u\in\sigma\cap X^{\vee}$. Then $u(m_0)=0$ and $u(x)\geq 0$\ 
$(\forall x\in \sigma^{\vee})$. 
By (iii), $\tX=\bZ m_0+\sigma^{\vee}\cap\tX$, so that 
$u(x)\geq 0$\ $(\forall x\in \tX)$. Hence $u(x)=0$\ $(\forall x\in X)$, that is, $u=0$. This proves (ii). 
\end{proof}

\begin{defn}
\label{defn]Sat}
For a subset $K$ of $\tX$, 
$\Semi(K)$ is defined to be 
the semisubgroup of $\tX$ generated by $K$. 
For a semi(sub)group $\Sigma$ of $\tX$, the {\it saturation  
$\Sat(\Sigma)$ of $\Sigma$} is defined by\ 
$\Sat(\Sigma)=\{x\in\tX; rx\in\Sigma\ (\exists\, r\in\bZ_{>0})\}$.  
For two subsets $K$ and $K'$ of $\tX_{\bR}$, 
we define $K+K'=\{\kappa+\kappa'; \kappa\in K, \kappa'\in K'\}$ and 
$-K=\{-\kappa; \kappa\in K\}$.  For a subset $T$ 
of $\tX_{\bR}$, we define $\Conv(T)$ to be the {\it convex closure of\,\, $T$} 
{\it in  $\tX_{\bR}$,} and $\Cone(T)$ to be the {\it cone over} $\Conv(T):=
\{rx; r\geq 0, x\in\Conv(T)\}$. 
Let $\Delta$ be a bounded closed polytope 
in $X_{\bR}$ with $0$ origin. 
If $0$ is a vertex of $\Delta$, then we define 
$\Cone(\Delta,0)=\Cone(\Delta)$. 
If $\alpha\in\Delta$ is a vertex of $\Delta$, then 
we define $\Cone(\Delta,\alpha)=\alpha+\Cone(\Delta-\alpha,0)$, which we call 
the cone of $\Delta$ at $\alpha$.

For a subset $T$ or a polytope 
$\Delta$ of $\tX^{\vee}_{\bR}$, we can define 
the convex closure of $T$, the cone over $T$ and 
the cone of $\Delta$ at some vertex $\alpha$ of $\Delta$ similarly.
\end{defn}

\begin{defn}\label{defn:weight}
Let $k(\eta)[X]$ be the group ring of $X$ over $k(\eta)$ and 
$\Sigma$ a semi(sub)group of $\tX$.
We define $R[\Sigma]$ to be the $R$-subalgebra of $k(\eta)[X]$ 
generated by $\Sigma$ where $w_0$ and $s$ are identified.
We define $\wt(s^aw^x)=am_0+x\in\tilde X$, 
which we call the {\it weight of $s^aw^x$.} 
\end{defn}

\begin{defn}
\label{defn:toric ring basic}
Let $A$ be an $R$-subalgebra of $k(\eta)[X]$ generated 
by some monomials $s^aw^x$. Then we define 
\begin{align*}
\Semi(A)&=\text{the semigroup in $\tX$ generated by $am_0+x$}\ 
\text{with $s^aw^x\in A$},\\
\Cone(A)&
=\text{the cone over $\Conv(\Semi(A))$}. 
\end{align*} 

 Then $A=R[\Semi(A)]$ and $\Sat(\Semi(A))=\Cone(A)\cap\tX$. 
Hence the normalization $\tA$ 
\footnote{It is the integral closure 
of $A$ in  $k(\eta)[X]$.} of $A$ is given by 
$\tA=R[\Cone(A)\cap\tX]$.     
In what follows, we call a cone over $\bR_{\geq 0}$ a cone
if no confusion is possible. 

\end{defn}
\begin{defn}\label{defn:Cut cF}
Let $\cF$ be any fan in $\tX$ over $S$. 
For $\sigma\in\cF$, we define 
\begin{align*}
\Cut(\sigma)&=-f_0+\sigma\cap(f_0+X^{\vee}_{\bR}).
\end{align*}  

Let $\Cut(\cF):=\{\Cut(\tau);\tau\in\cF\}$. Since $\sigma\in\cF$ 
contains no nontrivial linear subspace, $\Cut(\sigma)$ is a bounded convex polytope in $X_{\bR}^{\vee}$. We denote by $\Sk(\Cut(\sigma))$ 
(resp. $\Sk^q(\Cut(\sigma))$) the set of faces (resp. the set of $q$-dimensional faces)  of $\Cut(\sigma)$. 
\end{defn}

\begin{lemma}
\label{lemma:cF and CutcF}
Let $\cF$ be a fan in $\tX$ over $S$ and 
$\cF^{\times}:=\cF\setminus (0)$. 
Then the following are true:
\begin{enumerate}
\item\label{item:Cut bijective} the map 
$\Cut:\cF^{\times}\to\Cut(\cF)$ is bijective;
\item\label{item:dim Cut} if $\sigma\in\cF$, then 
$\dim\sigma=\dim\Cut(\sigma)+1$ where we define 
$\dim(\emptyset)=-1$;
\item\label{item:sigma=ConeCut} if $\sigma\in\cF$, then 
$\sigma=\Cone(f_0+\Cut(\sigma))$;
\item\label{item:sigma_cap_tau} if $\sigma,\tau\in\cF^{\times}$, then 
$\Cut(\sigma\cap\tau)=\Cut(\sigma)\cap\Cut(\tau)$.
\end{enumerate}
\end{lemma}
\begin{proof}Easy. 
\end{proof}

\begin{rem}\label{rem:def of faces}
Let $\sigma$ be a rational convex cone in $\tX_{\bR}$. 
Any face of $\sigma$ of codimension $\geq 2$
is the non-empty intersection of 
some of one-codimensional faces of $\sigma$. 
The same  is true as well for convex polytopes. 
\end{rem}

\subsection{The simplified Mumford families}
\label{subsec:a Mumford family} 
Let $R^{\sharp}$ be the polynomial algebra over $k(\eta)[X]$ 
in an indeterminate $\vartheta$:
\begin{gather*}R^{\sharp}:=k(\eta)[X][\vartheta]
=k(\eta)[w^x;x\in X][\vartheta] 
\end{gather*}
where we consider $R^{\sharp}$ as a graded algebra by defining 
\begin{gather*}
\deg\vartheta=1,\ \deg(cw^x)=0\ \ (c\in k(\eta), x\in X).
\end{gather*}

\begin{defn}
\label{defn:graded alg Aflat}
Let $\xi=(X, a, b, A, B)$ be a {\it principal} FC datum over $R$.  
We define  
a graded $R$-subalgebra $R^{\flat}(\xi)$ of $R^{\sharp}$ by
\begin{equation}
\label{eq:FC datum and Rflat}
R^{\flat}(\xi)=R[\xi_x\vartheta; x\in X],\ \xi_x=a(x)w^x.
\end{equation}

We also define an action  
$S_y^*$ $(y\in X$) 
on $R^{\sharp}$ by 
\begin{equation}\label{eq:Sy on Rsharp}
S_y^*(cw^x\vartheta)=ca(y)b(y,x)w^{x+y}\vartheta\ \ (c\in R, x\in X), 
\end{equation}which induces an action  
 on $R^{\flat}(\xi)$:
\begin{gather*}
\label{eq:action Sy}
S_y^*(ca(x)w^x\vartheta)=ca(x+y)w^{x+y}\vartheta.
\end{gather*}
\end{defn}
\begin{defn}
\label{defn:Xwedge Lwedge}
Let $\xi=(X,a,b,A,B)$ be a principal FC datum over $R$, 
$Y$ {\it a submodule of $X$ of finite index}, 
$(\tQ,\tcL):=(\Proj R^{\flat}(\xi),\cO_{\tQ}(1))$. 
Let $\tP$ be the normalization of $\tQ$,  and we denote 
by the same letter $\tcL$ 
the pullback to $\tP$ of $\tcL$ on $\tQ$ if no confusion is possible.
Let 
\begin{equation*}
U_{\alpha}:=\Spec R[\xi_{\alpha+\beta}/\xi_{\alpha}; \beta\in X]\quad (\alpha\in X),
\end{equation*}and $W_{\alpha}$ the normalization of $U_{\alpha}$. 
 Then $\tP$ (resp. $\tQ$) admits an affine open covering $(W_{\alpha};\alpha\in X)$ (resp. $(U_{\alpha};\alpha\in X)$). 
  
The ring homomorphism 
$S_y^*$ of $R^{\flat}(\xi)$ induces an $S$-automorphism 
$S_y$ of $(\tQ,\tcL)$, $(\tP,\tcL)$ 
and  $(\tP^{\wedge},\tcL^{\wedge})$. 
Let $\Lambda$ be a set of all representatives of $X\op{mod} Y$, 
 $\cD$ the union of all $W_{\alpha}$ $(\alpha\in\Lambda)$.
The action of $Y$ on $\tP_0$ (via $(S_y;y\in Y)$) is proper in the sense that 
there is a finite subset $\Gamma$ of $Y$ such that 
\begin{equation}\label{eq:action of Y is proper}
S_y(\overline{\cD_0})\cap S_z(\overline{\cD_0})=\emptyset\quad\text{if $y-z\not\in \Gamma$ and $y,z\in Y$}
\end{equation}
where  $\overline{\cD_0}$ is the closure of 
$\cD_0$ in $\tP_0$. Hence 
there is a formal quotient 
$(\tP^{\wedge},\tcL^{\wedge})/Y$ 
of $(\tP^{\wedge},\tcL^{\wedge})$ by $Y$, and similarly there is a formal quotient $(\tQ^{\wedge},\tcL^{\wedge})/Y$ 
of $(\tQ^{\wedge},\tcL^{\wedge})$ by $Y$.
\end{defn}

\begin{defn}
\label{defn:fan(xi)}
For $\alpha\in X$, we define 
\begin{align*}\tau_{\alpha}&=\left\{\begin{matrix}
x_0f_0+x\in {\tilde X}_{\bR}^{\vee};\ 
x_0\geq 0, x\in X_{\bR}^{\vee}\\ 
 (A(\beta)+B(\alpha,\beta))x_0+x(\beta)\geq 0\ (\forall \beta\in X)
\end{matrix}
\right\}.
\end{align*} 

Let
$\Fan(\xi):=\{\tau_{\alpha}\ \text{and their faces};\ \alpha\in X\}$. 
\end{defn} 

\begin{lemma}
\label{lemma:Fanxi over S}{\rm (\cite[3.8~(v)]{AN99})}
The $S$-scheme $\tP$ is isomorphic to the $S$-torus embedding $U(\Fan(\xi))$ associated with 
the fan $\Fan(\xi)$.
\end{lemma}
\begin{proof}Definition~\ref{defn:fan over S}~(i) 
follows from Definition~\ref{defn:fan(xi)}.
Let $u\in\tau_{\alpha}\cap X^{\vee}$. Since $\tau_{\alpha}^{\vee}$ is generated over $\bR_{\geq 0}$ by $z(\beta):=(A(\beta)+B(\alpha,\beta))m_0+\beta$\  
$(\forall \beta\in X)$, we have $u(\beta)=u(z(\beta))\geq 0$\ 
$(\forall \beta\in X)$.  Hence $u=0$, so that 
$\tau_{\alpha}\cap X^{\vee}=(0)$. 
Any $\sigma\in\Fan(\xi)$ is a face 
of $\tau_{\alpha}$ for some $\alpha\in X$. Hence 
$\sigma\cap X^{\vee}\subset \tau_{\alpha}\cap X^{\vee}=(0)$. 
This proves Definition~\ref{defn:fan over S}~(ii), so that 
$\Fan(\xi)$ is a fan in $\tX$ over $S$. Since $W_{\alpha}\simeq U(\tau_{\alpha})$, $\tP$ is isomorphic to $U(\Fan(\xi))$.   
\end{proof}

\begin{lemma}\label{lemma:toric is excellent}
Let $\Sigma$ be a finitely generated semigroup in  
$\tilde X$ 
with $m_0\in \Sigma$, $A(\Sat(\Sigma)):=R[\Sat(\Sigma)]/(w^{m_0}-s)R[\Sat(\Sigma)]$, $p$ any prime ideal of $A(\Sat(\Sigma))$ containing $I$, $B:=A(\Sat(\Sigma))_p$ and $C$ an excellent local $R$-algebra. 
Then 
\begin{enumerate}
\item\label{item:B Bwedge isom} 
$B$ and $B^{\wedge}$ are excellent normal and Cohen-Macaulay;
\item\label{item:C normal CM}  if $C^{\wedge}\simeq B^{\wedge}$,
then $C$ is normal and Cohen-Macaulay. 
\end{enumerate}
\end{lemma}
\begin{proof}By Gordan's lemma \cite[p.\ 7]{TE73}, 
$\Sat(\Sigma)$ is finitely generated. Since $\Sat(\Sigma)$ is saturated, 
$A(\Sat(\Sigma))$ is normal and Cohen-Macaulay 
by \cite[Theorem~14, p.\ 52]{TE73}, hence 
so is $B$. By Corollary~\ref{cor:if Chat=Bhat Rk Sk}~(1), $B$ is excellent. 
By Proposition~\ref{prop:basics of excellent}~(1), 
$B^{\wedge}$ is excellent, while 
it is normal and Cohen-Macaulay by Corollary~\ref{cor:if Chat=Bhat Rk Sk}~(2).
This proves (1). 
By (1) and Corollary~\ref{cor:if Chat=Bhat Rk Sk}~(3), 
$C$ is normal and Cohen-Macaulay. This proves (2).
\end{proof}

\begin{cor}
\label{cor:hatB normal CM}
Let  
$A_{\alpha}:=R[\xi_{\alpha+x}/\xi_{\alpha};x\in X]$ and 
$B_{\alpha}$ the integral closure of $A_{\alpha}$. 
Let $p$ (resp. $q$) be any prime ideal of $A_{\alpha}$ 
(resp. $B_{\alpha}$) 
containing $I$, $A^{\alpha}:=(A_{\alpha})_p$,   
$B^{\alpha}:=(B_{\alpha})_q$ and $C$ any  excellent local $R$-algebra. 
Then 
\begin{enumerate}
\item 
$(A^{\alpha})^{\wedge}$ is reduced, while 
$(B^{\alpha})^{\wedge}$ is normal and Cohen-Macaulay;
\item if $C^{\wedge}\simeq (A^{\alpha})^{\wedge}$ (resp. $C^{\wedge}\simeq (B^{\alpha})^{\wedge}$), then $C$ is reduced (resp. normal and Cohen-Macaulay). 
\end{enumerate}
\end{cor}
\begin{proof}Since $A^{\alpha}$ is an $R$-algebra of finite type, $A^{\alpha}$ is excellent by Proposition~\ref{prop:basics of excellent}~(3). Since  $A^{\alpha}$ is reduced, so is $(A^{\alpha})^{\wedge}$ by Corollary~\ref{cor:if Chat=Bhat Rk Sk}~(2). By \cite[Lemma~1, p.\ 5]{TE73}, $B_{\alpha}\simeq A(\Sat(\Sigma))$ for a semigroup $\Sigma$ in $\tX$. 
Hence $B^{\alpha}$ is excellent normal and Cohen-Macaulay 
by Lemma~\ref{lemma:toric is excellent}. 
Hence so is $(B^{\alpha})^{\wedge}$ 
by Corollary~\ref{cor:if Chat=Bhat Rk Sk}~(2). This proves (1). (2) is clear from  Corollary~\ref{cor:if Chat=Bhat Rk Sk}~(3).
\end{proof}

Let $\xi=(X,a,b,A,B)$ be a principal FC datum over $R$, 
$Y$ a submodule of $X$ of finite index
and $\tQ:=\Proj R^{\flat}(\xi)$. By \cite[\S~3]{AN99} 
and \cite[III${}_1$, 5.4.5]{EGA}, 
there exist {\em  flat projective $S$-schemes} 
$(P,\cL)$ and $(Q,\cL)$   
such that 
\begin{gather*}
(P^{\wedge},\cL^{\wedge})\simeq (\tP^{\wedge},\cO_{\tP^{\wedge}}(1))/Y,\ (Q^{\wedge},\cL^{\wedge})\simeq (\tQ^{\wedge},\cO_{\tQ^{\wedge}}(1))/Y.
\end{gather*} 
 Since $U_{\alpha}$ and $W_{\alpha}$ 
are flat over $S$, so are $Q$ and $P$ 
by Corollary~\ref{cor:hatB normal CM}.

\begin{defn}
\label{defn:Mumford family PQxi}Let $\xi=(X,a,b,A,B)$ 
be a principal FC datum over $R$, $Y$ a submodule of $X$ 
of finite index and  
$\xi':=\res_{X,Y}\xi=(X,Y,a,b,A,B)$.  
We call the above flat $S$-scheme $(P,\cL)$ 
the {\it Mumford family 
associated with $\xi'$,} which we denote
by $(P(\xi'),\cL(\xi'))$ if necessary. Similarly, we denote 
$(Q,\cL)$ by $(Q(\xi'),\cL(\xi'))$ if necessary. 
\end{defn}

\begin{prop}
\label{summary:Mumford family}
Let $R$ be a CDVR, $S$, 
$\eta$, $0$ and $k(\eta)$ the same as before. 
Let $\xi=(X,a,b,A,B)$ be a principal FC datum over $R$, 
$Y$ a submodule of $X$ of finite index, 
$\xi':=\res_{X,Y}\xi$, $(P,\cL):=(P(\xi'),\cL(\xi'))$ and  
$(Q,\cL):=(Q(\xi'),\cL(\xi'))$.
Then 
\begin{enumerate}
\item\label{item:normalization} $Q$ is reduced, 
and $P$ is the normalization of $Q$; 
\item\label{item:P0 reduced} $P(f^*\xi')_{0'}$ is reduced 
by taking a finite Galois cover $f:S'\to S$ 
ramifying only over $0$   
if necessary where $0'$ is the closed point of $S'$;  
\item\label{item:if P0 reduced} if $P_0$ is reduced, then 
there is an open $S$-subscheme $G$ 
of $P$ which is also a semiabelian $S$-scheme such that 
\begin{enumerate}
\item\label{item:4 semiabelian}   
$G_{\eta}=P_{\eta},\ G_0\simeq T_{X,k(0)}=\Spec k(0)[w^x; x\in X]$;
\item\label{item:4 Neron} $G$ is 
the connected N\'eron model 
\footnote{See Proposition~\ref{prop:semi abelian into neron}.} 
of $P_{\eta}$, but the $S$-immersion of $G$ into $P$ is not unique: 
for any irreducible component $W$ 
of $P_0\setminus\Sing(P_0)$, we can choose $G$ such that $G_0=W$. 
\end{enumerate}
\end{enumerate}
\end{prop}
\begin{proof}
\footnote{This proposition is true when $\xi$ is replaced by 
any split object \cite[\S~3.2]{Nakamura24} 
in the partially degenerate case, as is easily shown.}
 See \cite[3.10-3.12, 3.17, 3.18]{AN99} for 
(\ref{item:P0 reduced}).  
 (\ref{item:if P0 reduced}) follows from 
\cite[4.7]{Mumford72} 
and \cite[III, \S~4]{FC90} by \cite{AN99}. 
We shall prove (\ref{item:normalization}).
By Corollary~\ref{cor:hatB normal CM},  
$Q$ is reduced, and $P$ 
is normal. By \cite[3.10-3.12]{AN99} and (\ref{item:P0 reduced}), 
there exists a finite Galois covering 
$S'$ of $S$ ramifying only at $0'$
such that the closed fiber $P'_{0'}$ of $(P',\cL')$ is reduced where 
 $h:(P',\cL')\to (Q',\cL')$ is the normalization of 
$(Q',\cL'):=(Q,\cL)_{S'}$. Let $H:=\Gal(S'/S)$. Then 
$P'/H$ is the normalization of $Q'/H$. Since $Q'/H\simeq Q$, $P'/H$ is 
the normalization of $Q$.

Let  
$g:(P,\cL)\to (Q,\cL)$ be the natural morphism.  
  Since $P$ is normal and $P'/H$ is 
the normalization of $Q$, there exists an $S$-morphism
$$\pi:P\to P'/H$$
which factors $g$. 
Let $(\tQ',\tcL'):=(\tQ,\tcL)_{S'}$, and 
let $(\tP',\tcL')$ be the normalization of 
$(\tQ',\tcL')$.  Then the natural morphism $(\tQ',\tcL')\to (\tQ,\tcL)$ induces a morphism $\tmu:(\tP',\tcL')\to (\tP,\tcL)$ of their normalizations, a morphism $\tilde{\mu}^{\wedge}:((\tP')^{\wedge},(\tcL')^{\wedge})
\to (\tP^{\wedge},\tcL^{\wedge})$ of their $I$-adic completions,  
a morphism of their formal quotients by $Y$ and an $H$-invariant morphism 
$\mu:(P',\cL')\to (P,\cL)$ of their algebraizations. 
This $\mu$ induces a morphism 
$$\lambda:P'/H\to P$$
such that $\pi\lambda=\id_{P'/H}$ because 
$\pi^{\wedge}\lambda^{\wedge}=\id_{(P'/H)^{\wedge}}$ is the normalization of the composite $(=\id_{\tQ'/H}):\tQ'/H\overset{\simeq}{\to} \tQ\overset{\simeq}{\to} \tQ'/H$. Similarly, $\lambda\pi=\id_P$. It follows 
that $P\simeq P'/H$, which is the normalization of $Q$. 
This proves (\ref{item:normalization}). 
\end{proof}

\begin{rem}Proposition~\ref{summary:Mumford family}~(\ref{item:4 semiabelian}) tells that $P_{\eta}$ is a polarized abelian variety over $k(\eta)$ if $P_0$ 
is reduced, while by (\ref{item:P0 reduced}), $P_{\eta}$ is a polarized abelian torsor over $k(\eta)$ even if $P_0$ is not reduced. 
\end{rem}

\begin{defn}If $P_0$ 
is reduced, by Proposition~\ref{summary:Mumford family}, 
$P$ contains a semiabelian $S$-scheme $G$, 
which we denote by $G(\xi')$ if necessary. We also denote the restriction of $\cL(\xi')$ to $G(\xi')$ by the same letter. By Proposition~\ref{prop:semi abelian into neron}, $G$ is isomorphic 
to the connected N\'eron model of $G_{\eta}=P_{\eta}$. 
\end{defn}

\begin{rem}
\label{rem:modify by nilpotency}
Let $\xi=(X,a,b,A,B)$ be a principal FC datum, $Y$ any submodule of $X$ of finite index and $\xi':=\res_{X,Y}\xi$. 
If $B=E_8$, then $P(\xi')_0$ is nonreduced everywhere 
by \cite[3.12]{AN99}, so that $P(\xi')$ 
contains no semiabelian $S$-scheme $G$. 
See also \cite[11.11]{NSugawara06}. 
\end{rem}

\subsection{Fourier series (2)}
\label{subsec:Fouries 2}
Let $\xi=(X,a,b,A,B)$ be a principal FC datum, 
$Y$ a submodule of $X$ of finite index, 
$\xi':=\res_{X,Y}\xi$  
and $(P,\cL):=(P(\xi'),\cL(\xi'))$. 
The $S$-scheme $\tP$ is covered by the affine open sets  
$(W_{\alpha};\alpha\in X)$ in Definition~\ref{defn:Xwedge Lwedge} where   
$\Gamma(\cO_{W_{\alpha}})
=R[\tau_{\alpha}^{\vee}\cap\tX]=R[\zeta_{x,\alpha};x\in X]$ and  
 $\zeta_{x,\alpha}=s^{A(x)+B(x,\alpha)}w^x$ 
by \cite[3.8~(ii)]{AN99}. Then   
\begin{equation}\label{eq:Sy zetaxc}
S_y^*(\zeta_{x,\alpha})=\barb(y,x)\zeta_{x,\alpha+y}\ (\forall y\in Y)
\end{equation}by \cite[3.16]{AN99}.
We define 
$\eta_m(x)$ for $m\in\bN$ and $x\in X$ as follows. 
Associated with $B$, we have Delaunay cells of $X_{\bR}$, each being a bounded convex polytope. The union of all the Delaunay cells is $X_{\bR}$. See \cite[p.~662]{Nakamura99}. Hence there exists a Delaunay cell $\sigma$ 
such that $\frac{x}{m}\in\sigma\cap\frac{X}{m}$. By 
choosing any vertex $\alpha\in\sigma\cap X$, 
we define $\eta_m(x)=\zeta_{x-m\alpha,\alpha}(\xi_{\alpha}^*)^m$  
where $\xi^*_{\alpha}:=s^{A(\alpha)}w^{\alpha}$ (\cite[p.~672]{Nakamura99}). \footnote{This $\xi_{\alpha}^*$ is the same as $\xi_{\alpha}$ in \cite[3.4, p.~670]{Nakamura99}.}   
This is well-defined.  
By Eq.~(\ref{eq:Sy zetaxc}), we see
\begin{equation}\label{eq;etam}
S_y^*(\eta_m(x)\vartheta^m)
=\bara(y)^m\barb(y,x)\eta_m(x+my)\vartheta^m\ \ (\forall y\in Y).
\end{equation}  See also \cite[17th line, p.~673]{Nakamura99}.

\begin{lemma}\label{lemma:H0(Geta,Leta^m)}
Let $\xi=(X,a,b,A,B)$ be a principal FC datum, 
$Y$ a submodule of $X$ of finite index, $\xi':=\res_{X,Y}\xi$ and 
{\small $(G,\cL):=(G(\xi'),\cL(\xi'))$}. Then  
$\Gamma(G_{\eta},\cL_{\eta}^{\otimes m})$ is identified with 
the following  $k(\eta)$-vector space for any $m\geq 1$:
\begin{align*}
&\left\{\theta=\sum_{x\in X}\sigma_x(\theta)w^x; 
\begin{matrix}
\sigma_{x+my}(\theta)=a(y)^mb(y,x)\sigma_x(\theta)\\
\sigma_x(\theta)\in k(\eta)\ (\forall x\in X, \forall y\in Y)
\end{matrix}\right\}.
\end{align*}
\end{lemma}
\begin{proof}\
Let $(P,\cL):=(P(\xi'),\cL(\xi'))$. 
Since $H^q(P_0,\cL_0^{\otimes m})=0$ $(\forall q,m\geq 1)$ by 
\cite[3.9]{Nakamura99}, by \cite[Cor.~2, p.48]{Mumford12}
$\Gamma(P,\cL^{\otimes m})$ is $R$-free and  
\begin{equation*}\label{eq:Gamma(P,Lm)otimes k0, otimes keta}
\begin{aligned}
\Gamma(P_d,\cL_d^{\otimes m})
&=\Gamma(P,\cL^{\otimes m})\otimes_Rk(d)\ \  (d\in\{0,\eta\}).
\end{aligned}
\end{equation*} 

Let $j:\Gamma(P,\cL^{\otimes m})
\hookrightarrow \Gamma(P^{\wedge},(\cL^{\wedge})^{\otimes m})$ 
be the natural monomorphism. Since 
$\Gamma(P^{\wedge},(\cL^{\wedge})^{\otimes m})$ is $R$-free, we have 
$\Gamma(P^{\wedge},(\cL^{\wedge})^{\otimes m})\otimes_R k(0)
=\Gamma(P_0,\cL_0^{\otimes m})$, so that $j$ is 
an isomorphism by Nakayama's lemma. 
Moreover $\Gamma(P_0,\cL_0^{\otimes m})$ 
is given explicitly 
in terms of $\eta_m(x)$ $(x\in X)$ by 
\cite[3.9]{Nakamura99}. 
Hence $\Gamma(P^{\wedge},(\cL^{\wedge})^{\otimes m})$ 
is an $R$-module consisting of formal series of $\eta_m(x)$:
\begin{align*}
\Gamma(P^{\wedge},(\cL^{\wedge})^{\otimes m})
&=\left\{\sum_{x\in X}c(x)\eta_m(x);
\begin{matrix}c(x+my)=\bara(y)^m\barb(y,x)c(x)\\
c(x)\in R\ (\forall x\in X,\forall y\in Y)
\end{matrix}
\right\}.  
\end{align*}

Let $\theta:=\sum_{x\in X}c(x)\eta_m(x)$. 
Since $\sigma_x(\theta)w^x=c(x)\eta_m(x)$, by Eq.~(\ref{eq;etam}), 
\begin{align*}
&c(x+my)=\bara(y)^m\barb(y,x)c(x)
\\
&\Leftrightarrow 
S_y^*(c(x)\eta_m(x)\vartheta^m)=c(x+my)\eta_m(x+my)\vartheta^m\\
&\Leftrightarrow 
S_y^*(\sigma_x(\theta)w^x\vartheta^m)
=\sigma_{x+my}(\theta)w^{x+my}\vartheta^m\\
&\Leftrightarrow \sigma_{x+my}(\theta)=a(y)^mb(y,x)\sigma_x(\theta).
\end{align*} 

Since $j$ is an isomorphism, we see
\begin{align*}\Gamma(G_{\eta},\cL_{\eta}^{\otimes m})
&=\Gamma(P_{\eta},\cL_{\eta}^{\otimes m})
=\Gamma(P,\cL^{\otimes m})\otimes_R k(\eta)\\
&=\Gamma(P^{\wedge},(\cL^{\wedge})^{\otimes m})\otimes_R k(\eta).
\end{align*} 

It follows that $\Gamma(G_{\eta},\cL_{\eta}^{\otimes m})$ is 
identified  with 
\begin{align*}
&\left\{\theta=\sum_{x\in X}\sigma_x(\theta)w^x;
\begin{matrix}\sigma_{x+my}(\theta)=a(y)^mb(y,x)\sigma_x(\theta)\\
\sigma_x(\theta)\in k(\eta)\ (\forall x\in X,\forall y\in Y)
\end{matrix}
\right\}.\end{align*}

This completes the proof. 
\end{proof}

\begin{lemma}\label{lemma:H0(G,Lm) over keta as Fourier}
Let $(G,\cL)$ be a semiabelian $S$-scheme 
such that $G_0\simeq T_{X,k(0)}$, and  $\xi:=\FC(G,\cL)=(X,Y,a,b,A,B)$. 
Then $\Gamma(G_{\eta},\cL_{\eta}^{\otimes m})$ is identified with 
the following $k(\eta)$-vector space for any $m\geq 1$:
\begin{align*}
&\left\{
\theta=\sum_{x\in X}\sigma_x(\theta)w^x; 
\begin{matrix}
\sigma_{x+my}(\theta)=a(y)^mb(y,x)\sigma_x(\theta)\\
\sigma_x(\theta)\in k(\eta)\ (\forall x\in X, \forall y\in Y)
\end{matrix}\right\}.
\end{align*}
\end{lemma}
\begin{proof}Let $R_i$ be a CDVR\ $(i=1,2)$, 
$S_i:=\Spec R_i$, $0_i$ (resp. $\eta_i$) 
the closed (resp. generic) point of $S_i$ and $k(\eta_2)$ 
the fraction field of $R_2$. 

There exists a finite Galois covering 
$\varpi_1:S_1\to S$ ramifying only at $0$ and 
 a {\it principal} FC datum $\xi_1$ 
over $S_1$ such that 
$\varpi_1^*\xi=\res_{X,Y}\xi_1$. See \cite[3.3]{Nakamura99}. 
Then by Proposition~\ref{summary:Mumford family}~(\ref{item:P0 reduced}), 
there exists a finite  Galois covering 
$\varpi_2:S_2\to S_1$ ramifying only at $0_1$ 
such that $S_2$ is a Galois covering of $S$ via 
$\varpi:=\varpi_1\circ\varpi_2$ and 
$P(\xi_2)_{0_2}:=P(\xi_2)\times_{S_2}0_2$ 
is reduced where $\xi_2:=\res_{X,Y}\varpi_2^*\xi_1$.  
By Proposition~\ref{summary:Mumford family}~(\ref{item:if P0 reduced}), 
there exists an open $S_2$-subscheme $(G',\cL')$
of $(P',\cL')$ satisfying the conditions 
(\ref{item:4 semiabelian}) and (\ref{item:4 Neron})  
where $(P',\cL'):=(P(\xi_2),\cL(\xi_2))$. \par

 Let $\Gamma:=\Gal(k(\eta_2)/k(\eta))$. 
 Since $\xi_2=\varpi^*\xi$, 
we have 
$(G,\cL)\times_{S} S_2\simeq (G',\cL')$ 
by Theorem~\ref{thm:G isom G' iff FC same}. Hence  
$\Gamma(G_{\eta},\cL_{\eta}^{\otimes m})\otimes_{k(\eta)}k(\eta_2)
=\Gamma(G'_{\eta_2},(\cL'_{\eta_2})^{\otimes m})$.  It follows that  
$\Gamma(G_{\eta},\cL_{\eta}^{\otimes m})=\Gamma(G'_{\eta_2},(\cL'_{\eta_2})^{\otimes m})^{\Gamma}$.  
By Lemma~\ref{lemma:H0(Geta,Leta^m)}, 
{\small\begin{align*}
\Gamma(G'_{\eta_2},(\cL'_{\eta_2})^{\otimes m})^{\Gamma}&=\left\{
\theta=\sum_{x\in X}\sigma_x(\theta)w^x; 
\begin{matrix}
\sigma_{x+my}(\theta)=a(y)^mb(y,x)\sigma_x(\theta)\\
\sigma_x(\theta)\in k(\eta)\ (\forall x\in X, \forall y\in Y)
\end{matrix}\right\}.
\end{align*}}

 This completes the proof. 
\end{proof}

\begin{cor}\label{cor:FC(GLm)}
Let $m\geq 1$. Then 
with the notation in Corollary~\ref{cor:am bm},  
$$a_m(my)=a(y)^m,\ b_m(my,x)=b(y,x)\ \ (\forall y\in Y, \forall x\in X). 
$$
\end{cor}
\begin{proof}
Clear from Lemma~\ref{lemma:H0(G,Lm) over keta as Fourier}. 
See Lemma~\ref{lemma:Gamma(Geta,Leta^m m=2^l)}. 
\end{proof}

\begin{rem}Corollary~\ref{cor:FC(GLm)}  
proves \cite[3.10]{Nakamura99}.
\end{rem}

\begin{defn}\label{defn:xim}
Let $\xi=(X,Y,a,b,A,B)$ be an FC datum and  
$\xi_m:=(X_m,Y_m,a_m,b_m,A_m,B_m)$ (abbr.  
$(X,mY,a_m,b_m)$) where   
\begin{gather*}
X_m:=X,\ Y_m:=mY,\ a_m(my):=a(y)^m,\, b_m(my,x):=b(y,x)\\
A_m(my):=mA(y),\ B_m(my,x):=B(y,x)\ \  
 (x\in X, y\in Y). 
\end{gather*} 

Then $\xi_m$ is an FC datum. If 
$\xi=\FC(G,\cL)$, then $\xi_m=\FC(G,\cL^{\otimes m})$ 
by Corollary~\ref{cor:FC(GLm)}. 
\end{defn}

\section{The graded algebras}
\label{sec:graded algebras}
In what follows, throughout this article, 
the FC data we treat are assumed to be {\em symmetric.} 
In this section, starting with a (symmetric) FC datum $\xi$, 
we define/construct an eFC datum $\xi^e$ 
extending $\xi$ and a N\'eFC kit $\xi^{\natural}$) 
of $\xi^e_{2N}$). Then we define three graded algebras associated 
with $\xi^e$ or $\xi^{\natural}$:
$$R^{\quot}(\xi^e),\ R^{\sharp}(\xi^e_{2N}),\ R_{l}(\xi^{\natural})\ 
(l\in\bN)$$
in \S~\ref{subsec:graded alg Rnatural(xie)}, \S~\ref{subsec:graded alg Rsharp(xie)} and \S~\ref{subsec:graded alg R(xinatural)} respectively.
Our relative compactification of N\'eron model 
is constructed in 
\S~\ref{sec:relative compactifications} by using $\Proj R_{l}(\xi^{\natural})$.

\subsection{The eFC data}
\label{subsec:eFC data}

\begin{defn}\label{defn:normal extension}Let $F$ and $L$ be fields, 
An algebraic (field) extension $L/F$ is called {\it normal} 
if every irreducible polynomial in $F[X]$ that has a root in $L$ 
completely factors into linear factors over $L$. 
Let $\overline{F}$ be an algebraic closure of $F$ which contains $L$. 
Then (an algebraic extension) $L/F$ is normal iff any $F$-automorphism of 
$\overline{F}$ maps $L$ into $L$. 
Note that if $\chara F=0$ and $L/F$ is a finite normal extension, 
then $L/F$ is Galois. 
\end{defn}

Let $\xi=(X,Y,a,b,A,B)$ be an FC datum over a CDVR $R$. 
\begin{defn}\label{defn:beta}
There is a canonical isomorphism $j:\Hom_{\bZ}(Y\times X,\bZ)\simeq\Hom_{\bZ}(Y, X^{\vee})$. We define a monomorphism $\beta:=\beta_{\xi}:
=j(B)\in\Hom_{\bZ}(Y,X^{\vee})$. To be explicit,     
$\beta(y)(x):=B(y,x)$\ \ $(y\in Y, x\in X)$.
\end{defn}
\begin{defn}\label{defn:Nxi} 
Let $N:=N_{\xi}:=|X^{\vee}/\beta(Y)|=|\det(B(y_i,x_j))|$, 
which we denote by $|B|$ or $|B_{\xi}|$ 
where $(x_i;i\in [1,g])$ (resp. $(y_i;i\in [1,g])$) 
is a $\bZ$-basis of $X$ 
(resp. $Y$). 
\end{defn}
\begin{defn}
\label{defn:phiB}Let $N:=N_{\xi}$. 
We define 
$$\mu:=\mu_{\xi}:=\beta^{-1}N\id_{X^{\vee}}:X^{\vee}\to Y.$$

This is well-defined. Indeed, since 
$|X^{\vee}/\beta(Y)|=N$, we have $NX^{\vee}\subset\beta(Y)$.   
For any $x\in X^{\vee}$, there exists $y\in Y$ such that $Nx=\beta(y)$. Since $\beta$ is injective, this $y\ (=\mu(x))$ is unique, and hence 
$\mu$ is well-defined.
\end{defn}

\begin{lemma}
\label{lemma:be 2nd}Let $\xi=(X,Y,a,b,A,B)$ be an FC datum 
over a CDVR $R$ with $t$ uniformizer. By choosing a finite 
normal extension $K'$ of $k(\eta)$ if necessary, there exists 
a bilinear form $b^e:X^{\vee}\times X\to (K')^{\times}$ such that 
\begin{equation}\label{eq:be_betayx=b_yx}
\begin{aligned}
b^e&(\beta(y), x)=b(y, x),\ v_tb^e(u,x)=u(x)\ 
(\forall x\in X, \forall y\in Y, \forall u\in X^{\vee}).
\end{aligned}
\end{equation}
\end{lemma}
\begin{proof} 
Let $\barb(y, x):=t^{-B(y,x)}b(y,x)$. 
Suppose that $\barb^e:X^{\vee}\times X\to R_{\Omega}^{\times}$ 
is a bilinear form satisfying the following equations:
\begin{equation}\label{eq:barbe_betayx=barb_yx}
\barb^e(\beta(y), x)=\barb(y, x)\ (\forall x\in X, \forall y\in Y).\end{equation}
We set $b^e(u,x):=t^{u(x)}\barb^e(u,x)$. Then  
$b^e(\beta(y),x)=t^{B(y,x)}\barb^e(\beta(y), x)
=t^{B(y,x)}\barb(y,x)=b(y,x)$.  
Thus it suffices to construct $\barb^e$ 
satisfying  (\ref{eq:barbe_betayx=barb_yx}). 

 Let $e_i:=1$ $(i\in[1,r])$ and let 
$e_{r+1}|\cdots|e_g$ be elementary divisors of $X^{\vee}/\beta(Y)$. 
Then there exist a $\bZ$-basis $(u_i;i\in[1,g])$ of $X^{\vee}$ and    
a $\bZ$-basis $(y_j;j\in[1,g])$  of $Y$ 
such that $e_iu_i=\beta(y_i)$ $(\forall i\in[1,g])$. 
Let $(x_j;j\in[1,g])$ be a $\bZ$-basis of $X$.
To construct $b^e$, it suffices to solve 
the equations in $\Omega$:
\begin{equation}\label{eq:Galois eq for be}
\barb^e(u_i,x_j)^{e_i}=\barb(y_i,x_j)\in R^{\times}\ (i\in[1,g])
\end{equation} 
for every fixed $j\in [1,g]$. 
Any solution of 
(\ref{eq:Galois eq for be}) 
gives a solution $\barb^e$ of 
(\ref{eq:barbe_betayx=barb_yx}). Then we extend 
$\barb^e$ to $X^{\vee}\times X$ additively. 

Let $K$ be the subfield of $\Omega$ generated by all the solutions 
$\barb^e(u_i,x_j)$ of (\ref{eq:Galois eq for be}). 
 Then $K$ is a finite normal extension of $k(\eta)$. This completes the proof. 
\end{proof}
\begin{defn}\label{defn:min Galois ext Kxi for bexi}
Let $\xi$ be a symmetric FC datum over $R$. 
Let $K_{\min}(\xi)$ be the minimal normal extension of $k(\eta)$ such that 
$b^e(u,x)\in K_{\min}(\xi)$ $(\forall u\in X^{\vee},\forall x\in X)$. 
Let $R_{\min}(\xi)$ be the integral closure of $R$ in $K_{\min}(\xi)$, 
$S_{\min}(\xi):=\Spec R_{\min}(\xi)$ and $\eta_{\min}(\xi)$ 
the generic point of $S_{\min}(\xi)$. 
Note that 
$$K_{\min}(\xi_m)=K_{\min}(\xi),\ R_{\min}(\xi_m)=R_{\min}(\xi),\  
 \eta_{\min}(\xi_m)=\eta_{\min}(\xi)$$ 
for any $m\geq 1$ by Definition~\ref{defn:xim} 
and Eq.~(\ref{eq:be_betayx=b_yx}). 
\end{defn}

\begin{notation}
\label{notation:Notation Rinit Sinit keta=Kmin(xi)}
For later convenience, 
we introduce the following notation.
Let $R_{\init}$ be a CDVR with $t:=s_{\init}$ uniformizer, 
$S_{\init}:=\Spec R_{\init}$, 
$\eta_{\init}$ (resp. $0_{\init}$) 
the generic point (resp. the closed point) of $S_{\init}$ 
 and $k(\eta_{\init})$ (resp. $k(0_{\init})$)
the fraction field (resp. the residue field) of $R_{\init}$.

In what follows, throughout this article,  
let $\xi=(X,Y,a,b,A,B)$ be 
an FC datum over $S_{\init}$ with $A=v_t(a)$ and $B=v_t(b)$.   
Let $\beta:=\beta_{\xi}:=j(B)\in\Hom_{\bZ}(Y,X^{\vee})$, 
$N:=|B|$ and $\mu:=\beta^{-1}N\id_{X^{\vee}}
\in\Hom_{\bZ}(X^{\vee},Y)$.  

Let $k(\eta):=K_{\min}(\xi)$, $R:=R_{\min}(\xi)$, 
$S:=S_{\min}(\xi)=\Spec R$, $\eta:=\eta_{\min}(\xi)$, 
$s$ a uniformizer of $R$ 
and $e(\xi)$ the ramification index of 
$k(\eta)/k(\eta_{\init})$. Hence we have $v_t=v_s/e(\xi)$ 
on $k(\eta_{\init})$, by which  
we extend $v_t$ to $k(\eta)$. 
We keep this notation unless otherwise mentioned.
\end{notation}

\begin{defn} 
\label{defn:extended FC datum}Let $\xi=(X,Y,a,b)$ 
be a symmetric FC datum over $R_{\init}$
Then we define a sextuple 
\begin{equation}
\label{eq:eFC datum}
\xi^e=(X,Y,a,b^e,A,B)\ \ ({\rm abbr.}\ \ (X,Y,a,b^e))
\end{equation} 
to be a ({\it totally degenerate}) {\it eFC datum 
over} $R$ ({\it extending} $\xi$) 
if \begin{enumerate}
\item[(a)] $b^e:X^{\vee}\times X\to k(\eta)^{\times}$ is 
a bilinear function;
\item[(b)] $\xi(\xi^e):=(X,Y,a,b,A,B)$ is an FC datum 
(equal to $\xi$) over $R_{\init}$ with $b(y,x):=b^e(\beta(y),x)$
\ $(y\in Y,  x\in X)$.  
\end{enumerate}
\end{defn}

\begin{defn}\label{defn:xi_m_e}For a symmetric FC datum $\xi=(X,Y,a,b)$ 
over $R_{\init}$, 
we have $\xi_m=(X,mY,a_m,b_m)$ for $m\geq 1$ by Definition~\ref{defn:xim}. 
We define 
$$\xi^e_m=(X, mY, a_m, b^e_m),\ \ b^e_m=b^e,\ \beta_m=\beta_{\xi_m}.$$ 
Since $\beta_m$ satisfies 
$\beta_m(my)=\beta(y)$ 
by Corollary~\ref{cor:FC(GLm)}, we have    
$$b^e_m(\beta_m(my),x)=b^e_m(\beta(y),x)=b^e(\beta(y),x)=b(y,x)=b_m(my,x).$$
Hence $\xi^e_m$ is an eFC datum over $R=R_{\min}(\xi)$ 
extending $\xi_m$. By Lemma~\ref{lemma:be 2nd},  
any eFC datum over $R$ extending $\xi$ is 
{\it conjugate under} $\Aut(k(\eta)/k(\eta_{\init}))$. 
Hence {\it so} is any eFC datum extending $\xi_m$ over $R$.
\end{defn}

\begin{lemma}\label{lemma:Be and Bephi}
Let $\xi^e$ be an eFC datum over $R$ extending an FC datum $\xi$ 
over $R_{\init}$, $u,v\in X^{\vee}$, $x\in X$  
and $F(u,v):=u(\mu_{\xi}(v))$. Then 
\begin{enumerate}
\item\label{item:Be} $B^e(u,x):=v_t b^e(u,x)=u(x)$;
\item\label{item:F pos def} 
$F$ is symmetric bilinear and positive definite on $X^{\vee}\times X^{\vee}$. 
\end{enumerate}
\end{lemma}
\begin{proof}Let $\mu:=\mu_{\xi}$. 
 By Lemma~\ref{lemma:be 2nd}, we have  
$v_tb^e(u,x)=u(x)$. 
This proves (1). Next we shall prove (2). 
It suffices to prove the claim on $\beta(Y)\times\beta(Y)$.
Let $u=\beta(x)$ and $v=\beta(y)$ for $x,y\in Y$. Since $B$ is symmetric,  
$F(u,v)=\beta(x)(\mu\beta(y))=NB(x,y)=\beta(y)(\mu\beta(x))
=F(v,u)$, 
while  
$F(u,u)=\beta(x)(\mu\beta(x))=NB(x,x)\geq 0$ 
by the positivity of $B$.  
If $u\neq 0$, then $x\neq 0$, so that $NB(x,x)>0$. 
This proves (2). 
\end{proof}

\subsection{The graded algebra $R^{\quot}(\xi^e)$}
\label{subsec:graded alg Rnatural(xie)}
We use the notation in \S~\ref{subsec:eFC data}.
 
\begin{defn}
\label{defn:deltay} 
We define a commutative $k(\eta)[X]$-algebra
\begin{equation}
\label{defn:R}
R(\xi^e)=k(\eta)[X][\vartheta_u; u\in X^{\vee}]
\end{equation}where $\vartheta_u$ is an indeterminate 
$(\forall u\in X^{\vee})$. 
This is a graded algebra with $\deg(\vartheta_u)=1$ 
and $\deg(cw^{\alpha})=0$ $(c\in k(\eta), u\in X^{\vee},\alpha\in X)$. 
The algebra $R^{\sharp}:=k(\eta)[X][\vartheta]$ 
is identified with a $k(\eta)[X]$-subalgebra of 
$R(\xi^e)$ by setting $\vartheta_0=\vartheta$. 
Next, we define $k(\eta)$-endomorphisms $\delta_u^*$ and $S_y^*$ 
of $R(\xi^e)$ by
\begin{equation}\label{eq:defn of delta_y}
\begin{aligned}
\delta_u^*(cw^{\alpha})&=cb^e(u,\alpha)w^{\alpha},\ 
\delta_u^*(\vartheta_v)=\vartheta_{u+v}, \\ 
S_y^*(cw^{\alpha})&=cb(y,\alpha)w^{\alpha},\ 
S_y^*(\vartheta_v)=\vartheta_{v+\beta(y)} 
\end{aligned}
\end{equation}where $c\in k(\eta)$, $u,v\in X^{\vee}$, 
$\alpha\in X$ and $y\in Y$. Hence 
\begin{equation}\label{eq:delta_u Sy}
\delta^*_u\delta^*_v=\delta^*_{u+v},\  \delta^*_{\beta(y)}=S_y^* 
\quad(\forall u, v\in X^{\vee}, \forall y\in Y).
\end{equation} 
\end{defn}

\begin{defn}\label{defn:Rnatural_xie}
We define a quotient algebra of $R(\xi^e)$ by 
\begin{equation*}
\label{defn:Rnatural}
R^{\quot}(\xi^e)=R(\xi^e)/(\vartheta_{v+\beta(y)}-a(y)b^e(v,y)w^y
\vartheta_v; v\in X^{\vee}, y\in Y).
\end{equation*}
This is a graded $k(\eta)[X]$-algebra generated by $\vartheta_v$ 
$(v\in X^{\vee}/\beta(Y))$ whose grading is inherited from that of 
$R(\xi^e)$. Moreover $\delta_u^*$ and $S_y^*$ induce 
$R$-endomorphisms of $R^{\quot}(\xi^e)$ because 
\begin{align*}
\delta_u^*(\vartheta_{v+\beta(y)})&=\vartheta_{u+v+\beta(y)}
=a(y)b^e(u+v,y)w^y\vartheta_{u+v}\\
&=a(y)b^e(v,y)b^e(u,y)w^y\vartheta_{u+v}
=\delta_u^*(a(y)b^e(v,y)
w^y\vartheta_v).
\end{align*}

Since both $\delta_u^*$ and $S_y^*$ are well-defined 
on $R^{\quot}(\xi^e)$,
the equalities (\ref{eq:defn of delta_y}) and (\ref{eq:delta_u Sy}) 
are true on $R^{\quot}(\xi^e)$.
\end{defn}

\subsection{The graded algebra  $R^{\sharp}(\xi^e_{2N})$}
\label{subsec:graded alg Rsharp(xie)}
\begin{defn}
\label{defn:theta=prod delta_ over Delta}Let 
$\Psi=\{v_i\in X^{\vee}; i\in[1,N]\}$ be a set of representatives 
in $X^{\vee}$ 
of $X^{\vee}/\beta(Y)$ such that $0\in\Psi$. Let 
$\theta_{\pm}:=\prod_{v\in (\pm\Psi)}\vartheta_v$ 
and $\theta:=\theta_{+}\theta_{-}.$  
Recall an FC datum $\xi^e_{2N}$ in 
Definition~\ref{defn:extended FC datum} with $m=2N$.  
We define a graded $k(\eta)[X]$-subalgebra $R^{\sharp}(\xi^e_{2N})$ 
of $R^{\quot}(\xi^e)$ by 
\begin{gather*}\label{defn:Rsharp}
R^{\sharp}(\xi^e_{2N})=k(\eta)[X][\theta],\ 
\deg(\theta)=1,\ \deg(c)=0\ (c\in k(\eta)[X]).
\end{gather*}
\end{defn}

\begin{lemma}
\label{lemma:deltaz*theta}
 $R^{\sharp}(\xi^e_{2N})$ is a subalgebra of $R^{\quot}(\xi^e)$ 
stable under $(\delta_u^*;u\in X^{\vee})$. Let 
$\epsilon(u):=b^e(u,\mu(u))$.  
Then $\delta_u^*\theta=\epsilon(u)w^{2\mu(u)}\theta.$
\end{lemma}
\begin{proof}
Let $\Lambda:=[1,N]$ and $u\in X^{\vee}$. Then 
for every $i\in\Lambda$, there exists a unique pair $j(i)\in\Lambda$ and 
$y_i\in Y$ such that $u+v_i=v_{j(i)}+\beta(y_i)$. Then $\sum_{i\in\Lambda}y_i
=\beta^{-1}Nu=\mu(u)$ because 
\begin{align*}
Nu+\sum_{i\in\Lambda}v_i&=\sum_{i\in\Lambda}(u+v_i)=\sum_{i\in\Lambda}(v_{j(i)}+\beta(y_i))\\
&=\sum_{i\in\Lambda}v_{j(i)}+\beta\sum_{i\in\Lambda}y_i=\sum_{i\in\Lambda}v_i+\beta\sum_{i\in\Lambda}y_i.
\end{align*}

Now we define $\epsilon(u)=\epsilon_{+}(u)\epsilon_{-}(u)$ where 
\begin{equation}\label{eq:computation of epsilon_plus_minus}
\epsilon_{+}(u):=\prod_{i\in\Lambda}a(y_i)b^e(v_{j(i)},y_i),\ 
\epsilon_{-}(u):=\prod_{i\in\Lambda}a(y_i)b^e(-v_i,y_i).
\end{equation}

Hence we see 
\begin{align*}
\delta_u^*\vartheta_{v_i}&=\vartheta_{u+v_i}=\vartheta_{v_{j(i)}+\beta(y_i)}
=a(y_i)b^e(v_{j(i)},y_i)w^{y_i}\vartheta_{v_{j(i)}},\\
\delta_u^*\theta_{+}&=\prod_{i\in\Lambda}\delta_u^*\vartheta_{v_i}=\prod_{i\in\Lambda}a(y_i)b^e(v_{j(i)},y_i)w^{y_i}
\vartheta_{v_{j(i)}}
=\epsilon_{+}(u)w^{\mu(u)}\theta_{+}.
\end{align*}

Since $u+(-v_{j(i)})=(-v_i)+\beta(y_i)$, similarly we have 
\begin{align*}
\delta_u^*\theta_{-}&=\prod_{i\in\Lambda}\delta_u^*\vartheta_{-v_{j(i)}}=
\prod_{i\in\Lambda}\vartheta_{-v_i+\beta(y_i)}
=\epsilon_{-}(u)w^{\mu(u)}\theta_{-}.
\end{align*}

By Eq.~(\ref{eq:computation of epsilon_plus_minus}) and 
Lemma~\ref{lemma:FC(G_L) L symmetric} 
\begin{align*}
\epsilon(u)&
=\prod_{i\in\Lambda}a(y_i)^2b^e(v_{j(i)}-v_i,y_i)\\
&=\prod_{i\in\Lambda}a(y_i)^2b^e(u-\beta(y_i),y_i)
=b^e(u,\mu(u)).
\end{align*}
\end{proof}

\begin{lemma}
\label{lemma:basic formulae}
The following are true in $R^{\sharp}(\xi^e_{2N})$:
\begin{enumerate}
\item\label{item:Sy walpha theta} $S_y^*(w^{\alpha}\theta)=a(y)^{2N}b(y,\alpha)
w^{\alpha+2Ny}\theta$; 
\item\label{item:deltau walpha theta} $\delta_u^*(w^{\alpha}\theta)=
\epsilon(u)b^e(u,\alpha)w^{\alpha+2\mu(u)}\theta$;
\item\label{item:epsilon beta= a_2N} $\epsilon(\beta(y))=a(y)^{2N}$ and 
$b^e(\beta(y),\alpha)=b(y,\alpha)$;
\item\label{item:be} $b^e(u,\mu(v))=b^e(v,\mu(u))$;
\item\label{item:epsilon u+v} $\epsilon(u+v)
=\epsilon(u)\epsilon(v)b^e(u,\mu(v))^2$;
\item\label{item:epsilon u+betay} $\epsilon(u+\beta(y))
=\epsilon(u)a(y)^{2N}b^e(u,y)^{2N}$
\end{enumerate} where $\alpha\in X$, $y\in Y$ and $u,v\in X^{\vee}$.
\end{lemma}
\begin{proof}Let $m_{\Psi}:=\sum_{v\in\Psi}v$.
By Lemma~\ref{lemma:deltaz*theta}, 
{\small\begin{gather*}
S_y^*(\vartheta_v)=\vartheta_{v+\beta(y)}
=a(y)b^e(v,y)w^y\vartheta_v,\\
S_y^*(\theta_{\pm})=\left(\prod_{v\in\pm\Psi}a(y)b^e(v,y)w^y\right)
\theta_{\pm}
=a(y)^Nb^e(\pm m_{\Psi},y)w^{Ny}\theta_{\pm},\\
S_y^*(w^{\alpha}\theta)=a(y)^{2N}b(y,\alpha)
w^{\alpha+2Ny}\theta,\ 
\delta_u^*(w^{\alpha}\theta)=
\epsilon(u)b^e(u,\alpha)w^{\alpha+2\mu(u)}\theta.
\end{gather*}}
\hskip -0.15cm These prove 
(\ref{item:Sy walpha theta})-(\ref{item:epsilon beta= a_2N}). 
By Eq.~(\ref{eq:delta_u Sy}) and 
Lemma~\ref{lemma:deltaz*theta}, 
$\epsilon(u+v)\epsilon(u)^{-1}\epsilon(v)^{-1}
=b^e(u,\mu(v))^2=b^e(v,\mu(u))^2=b^e(u,\mu(v))b^e(v,\mu(u))$. 
This proves (\ref{item:be}) and (\ref{item:epsilon u+v}). 
See also Lemma~\ref{lemma:be 2nd}.
(\ref{item:epsilon u+betay}) follows from (\ref{item:epsilon beta= a_2N}) 
and (\ref{item:epsilon u+v}). 
\end{proof}

\begin{cor}
\label{cor:E(x)}
Let $E(u):=v_t\epsilon(u)$, 
$x\in X$, $y\in Y$ and $u,v\in X^{\vee}$. Then 
\begin{enumerate}
\item\label{item:E=umuu} $E(u)=u(\mu(u))=F(u,u)$, $E(\beta(y))=2NA(y)$;
\item\label{item:E sum} $E(u+v)=E(u)+E(v)+2u(\mu(v))$; 
\item\label{item N=Bmu} $Nu(x)=B(\mu(u),x)$.
\end{enumerate}
\end{cor}
\begin{proof}
Lemma follows from 
Lemmas~\ref{lemma:Be and Bephi}/\ref{lemma:deltaz*theta}/\ref{lemma:basic formulae} 
and $\beta\mu=N\id_{X^{\vee}}$.  
\end{proof}

\begin{rem}The algebra $R^{\sharp}(\xi^e_{2N})$ 
has actions $\delta^*_u$ and $S^*_y$ 
$(u\in X^{\vee},y\in Y)$, each being described  
by Lemma~\ref{lemma:basic formulae} 
in terms of the eFC datum 
$\xi^e_{2N}$ rather than $\xi^e$. 
See Definition~\ref{defn:2nd neron FC datum}~(iii). 
\end{rem}

\subsection{The graded algebra $R(\xi^{\natural})$}
\label{subsec:graded alg R(xinatural)}
\begin{defn}
\label{defn:2nd neron FC datum} With the notation 
in \S~\ref{subsec:eFC data}, 
we define a sextuple  
\begin{equation}\label{eq:NeFC kit xinatural}
\xi^{\natural}=(X,Y,\epsilon,b^e,E,\Sigma)
\end{equation}to be a {\it N\'eron-FC kit} (abbr. a {\it N\'eFC kit}) 
{\it over $R$} of $\xi^e_{2N}$ 
if 
\begin{enumerate}
\item[(i)]  $\epsilon:X^{\vee}\to k(\eta)^{\times}$ and 
$b^e:X^{\vee}\times X\to k(\eta)^{\times}$ are functions such that 
\begin{gather*}
 \epsilon(u)=b^e(u,\mu(u)),\ 
E(u)=v_t\epsilon(u)=u(\mu(u))\ \   
(\forall u\in X^{\vee});
\end{gather*}
\item[(ii)] $\Sigma:=\{\alpha\in X; E(u)+u(\alpha)\geq 0\ 
(\forall u\in X^{\vee})\}$;
\item[(iii)] $\xi^e_{2N}:=\xi^e_{2N}(\xi^{\natural}):
=(X,2NY,a_{2N},b^e)$ 
is an eFC datum over $R$ extending $\xi_{2N}$   
where $a_{2N}(2Ny)=a(y)^{2N}=\epsilon(\beta(y))$ $(\forall y\in Y)$.  
\end{enumerate} 

By Lemma~\ref{lemma:be 2nd}, {\it there does exist} 
a N\'eFC kit $\xi^{\natural}$ 
over $R$ of $\xi^e_{2N}$. 
By Lemma~\ref{lemma:Be and Bephi}~(2), 
$\Sigma$ is a finite subset of $X$. Meanwhile (iii) implies that 
Lemma~\ref{lemma:basic formulae}~(3)-(6) 
and Corollary~\ref{cor:E(x)} are true for $\xi^{\natural}$.  
\end{defn}

\begin{defn}
\label{defn:R(xinatural)}Let 
$\xi^{\natural}$ be a N\'eFC kit given by Eq.~(\ref{eq:NeFC kit xinatural})   
and $\theta:=\theta_{+}\theta_{-}$.    
We define an $R$-subalgebra $R(\xi^{\natural})$ 
of $R^{\sharp}(\xi^e_{2N})=k(\eta)[X][\theta]$ by 
\begin{gather*}
R(\xi^{\natural})
=R[\xi_{\alpha,v}\theta; \alpha\in\Sigma, v\in X^{\vee}],\ 
\xi_{\alpha,v}
=\epsilon(v)b^e(v,\alpha)w^{\alpha+2\mu(v)}.
\end{gather*}

By Lemma~\ref{lemma:basic formulae} 
we have 
\begin{gather*}
\delta_u^*(\xi_{\alpha,v}\theta)=\xi_{\alpha,v+u}\theta,\  
S_y^*(\xi_{\alpha,v}\theta)=\xi_{\alpha,v+\beta(y)}\theta\ (\forall u,v\in X^{\vee}, \forall y\in Y).
\end{gather*} 
Hence $R(\xi^{\natural})$ is stable under  
$\delta_u^*$ and $S_y^*$,  
so that $\delta_u^*$ (resp. $S_y^*$)  induces 
an $S$-automorphism $\delta_u$ (resp. $S_y$) 
of $\Proj R(\xi^{\natural})$ and $\cO_{\Proj R(\xi^{\natural})}(1)$. 
\end{defn}

\begin{defn}(\cite[p.~209]{Mumford72})
Let $\cB$ be any $R$-algebra and $a\in T_{X,R}(\cB)$. 
Then we define an $R$-automorphism $T_a^*$ of 
$\cB\otimes_RR(\xi^{\natural})$ by 
$$T_a^*(w^x\theta)=w^x(a)\otimes w^x\theta\ (\forall x\in X)  
$$where $w^x(a):=a^*w^x$.  
Since $T_a^*\circ T_b^*=T_{a+b}^*$ $(\forall a,b\in T_{X,R}(\cB))$ and $T_a^*$ is functorial in $a$ and $\cB$,  $T_{X,R}$ acts on $\Proj R(\xi^{\natural})$ and $\cO_{\Proj R(\xi^{\natural})}(1)$. 
\end{defn}

\subsection{The graded algebra $R_{l}(\xi^{\natural})$}
\label{subsec:Rl(xinatural)}
We use the notation in \S~\ref{subsec:graded alg R(xinatural)}.
\begin{defn}\label{defn:xinatural_l} 
We define the {\it $l$-th N\'eFC kit $\xi^{\natural}_{l}$ 
induced from $\xi^{\natural}$} by 
\begin{gather*}
\xi^{\natural}_{l}=(X,Y,\epsilon_{l},b^e,E_{l},
\Sigma_{l})\ (l\geq 1),
\ \xi^{\natural}_1=\xi^{\natural},\\
\Sigma_{l}(0)=\{\alpha\in X_{\bR}; E_{l}(u)+u(\alpha)\geq 0\ 
(\forall u\in X^{\vee})\},\\ 
\epsilon_{l}=\epsilon^{l}=b^e(u,l\mu(u)),\ 
E_{l}=l E=l v_t(\epsilon),\ \Sigma_{l}=\Sigma_{l}(0)\cap X.
\end{gather*}
\end{defn}

\begin{defn}
\label{defn:R_l(xinatural)}
We define a graded $R$-subalgebra $R_{l}(\xi^{\natural})$ 
of $k(\eta)[X][\theta_{l}]$ by 
\begin{gather*}
R_{l}(\xi^{\natural})
=R[\xi_{l,\alpha,v}\theta_{l}; \alpha\in\Sigma_{l}, v\in X^{\vee}]
\end{gather*}where $\theta_{l}:=\theta^{l}$ 
and 
$\xi_{l,\alpha,v}:=\epsilon_{l}(v)b^e(v,\alpha)w^{\alpha+2l\mu(v)}
$ $(\alpha\in\Sigma_{l}, v\in X^{\vee})$. Then 
\begin{gather*}
\delta_u^*(\xi_{l,\alpha,v}\theta_{l})
=\xi_{l,\alpha,v+u}\theta_{l}\   
(\forall u\in X^{\vee}).
\end{gather*}
Thus $R_{l}(\xi^{\natural})$ is stable 
under $(\delta_u^*;u\in X^{\vee})$ 
and $(S_y^*;y\in Y)$ where $S_y^*=\delta^*_{\beta(y)}$. 

The $R$-automorphism $\delta_u^*$ (resp. $S_y^*$) 
of $R_{l}(\xi^{\natural})$ induces 
an $S$-automorphism $\delta_u$ (resp. $S_y$) 
of $\Proj R_{l}(\xi^{\natural})$ and $\cO_{\Proj R_{l}(\xi^{\natural})}(1)$. 
Similarly $T_{X,R}$ acts 
on $\Proj R_{l}(\xi^{\natural})$ and $\cO_{\Proj R_{l}(\xi^{\natural})}(1)$. 
Definition~\ref{defn:R_l(xinatural)} is derived from 
Definition~\ref{defn:R(xinatural)}.
\end{defn}

 The following is obvious from Lemma~\ref{lemma:basic formulae}. 
\begin{lemma}
\label{lemma:lth basic formulae}
The following are true in $R_{l}(\xi^{\natural})$:
\begin{enumerate}
\item\label{item:lth Sy walpha theta} 
$S_y^*(w^{\alpha}\theta_{l})=a(y)^{2Nl}b(y,\alpha)
w^{\alpha+2Nl y}\theta_{l}$; 
\item\label{item:lth deltau walpha theta} $\delta_u^*(w^{\alpha}\theta_{l})=
\epsilon_{l}(u)b^e(u,\alpha)w^{\alpha+2l\mu(u)}\theta_{l}$;
\item\label{item:lth epsilon beta= a_2N} 
$\epsilon_{l}(\beta(y))=a(y)^{2Nl}$ and 
$b^e(\beta(y),\alpha)=b(y,\alpha)$;
\item\label{item:lth be} $b^e(u,\mu(v))=b^e(v,\mu(u))$;
\item\label{item:lth epsilon u+v} $\epsilon_{l}(u+v)
=\epsilon_{l}(u)\epsilon_{l}(v)b^e(u,\mu(v))^{2l}$;
\item\label{item:lth epsilon u+betay} 
$\epsilon_{l}(u+\beta(y))=\epsilon_{l}(u)a(y)^{2Nl}b^e(u,y)^{2Nl}$ 
\end{enumerate}where $\alpha\in X$, $y\in Y$ and $u,v\in X^{\vee}$.
\end{lemma} 
\begin{cor}
\label{cor:El(x)}
Let $y\in Y$, $u,v\in X^{\vee}$ and 
$E_{l}:=v_t\epsilon_{l}$ $(l\geq 1)$. 
Then 
\begin{enumerate}
\item\label{item:El=lumuu} 
$E_{l}(u)=l u(\mu(u))$, $E_{l}(\beta(y))=2Nl A(y)$;
\item\label{item:El sum} $E_{l}(u+v)=E_{l}(u)+E_{l}(v)+2lu(\mu(v))$. 
\end{enumerate}
\end{cor}

The following is clear.
\begin{lemma}\label{lemma:Rlxinatural=Rxinaturall}
$R_{l}(\xi^{\natural})\simeq R(\xi^{\natural}_{l})$
 for any $l\geq 1$.
\end{lemma}

\section{Voronoi polytopes}
\label{sec:Voronoi polytopes}
\subsection{The Voronoi polytopes of $\xi^{\natural}_{l}$}
\label{subsec:Voronoi cells}
Let $\xi=(X,Y,a,b,A,B)$ and $\xi^{\natural}_{l}
=(X,Y,\epsilon_{l},b^e,E_{l},\Sigma_{l})$ be 
the same as in \S~\ref{subsec:eFC data}  
and \S~\ref{subsec:Rl(xinatural)} respectively. 
\begin{defn}
\label{defn:Voronoi polytopes}We define a distance 
$\|z\|=\sqrt{B(z,z)}$ on $X_{\bR}$. For any $c\in X^{\vee}$, 
we define a Voronoi polytope
$\Sigma_{l}(c)$ by 
$$\Sigma_{l}(c)=\{x\in X_{\bR}; \|x-2l\mu(u)\|\geq 
\|x-2l\mu(c)\|\ (\forall u\in X^{\vee})\}.
$$ 
\end{defn}

The set $\Vor_{l,B}$  of 
all $\Sigma_{l}(c)$ $(c\in X^{\vee})$ and their closed faces 
 forms a polyhedral decomposition of $X_{\bR}$ 
into bounded convex polytopes. We denote $\Vor_{l,B}$ by $\Vor_{l}$  
if no confusion is possible. 
Every $\Sigma_{l}(c)$ is the set consisting of all $x\in X_{\bR}$ 
such that $2l\mu(c)$ is the nearest to $x$ 
among $2l\mu(X^{\vee})$.  
Each $\Sigma_{l}(c)$ is a bounded convex polytope such that 
$$X_{\bR}=\bigcup_{c\in X^{\vee}}\Sigma_{l}(c),\  
\Sigma_{l}(c)=\Sigma_{l}(0)+2l\mu(c).
$$

 Every pair $\Sigma_{l}(c)$ and $\Sigma_{l}(d)$\ $(c\neq d)$ 
have no interior points in common. 
For a polytope $W$ of $X_{\bR}$, we denote by $W^0$ 
the relative interior of $W$. If necessary, we denote $\Sigma_{l}(c)$ by 
$\Sigma_{l}(\xi^{\natural})(c)$.

\begin{lemma} 
\label{lemma:Voronoi_cap_X=Sigmal} 
$\Sigma_{l}=\Sigma_{l}(0)\cap X$, $0\in\Sigma_{l}(0)^0$ and  
\begin{align*}
\Sigma_{l}(c)&=\{x\in X_{\bR};  E_{l}(u)+u(x-2l\mu(c))\geq 0
\ (\forall u\in X^{\vee})\}.
\end{align*}
\end{lemma}
\begin{proof} 
Note that $\Sigma_{l}(0)=\{x\in X_{\bR};  E_{l}(u)+u(x)\geq 0
\ (\forall u\in X^{\vee})\}$. 
We prove $0\in \Sigma_{l}(0)^0$. Since $E_{l}$ is positive definite by Lemma~\ref{lemma:Be and Bephi}~(2), the minimal eigenvalue $a$ of $E_{l}$ is positive.  Let $(f_i;i\in [1,g])$ be a basis of $X^{\vee}$  
and $u=\sum_{i=1}^gu_if_i\in X^{\vee}$.   
If $u\neq 0$ and $|f_i(x)|<a/2$ $(\forall i)$, then  
$$|u(x)|=\biggl|\sum_{i=1}^gu_if_i(x)\biggr|\leq 
(a/2)\sum_{i=1}^gu_i^2< E_{l}(u),$$  
whence $0\in \{x\in X_{\bR}; |f_j(x)|<a/2\ (\forall j)\}\subset \Sigma_{l}(0)$. It follows $0\in \Sigma_{l}(0)^0$. The rest is clear.
\end{proof}

\begin{cor}
\label{cor:Sigma_l}$\Sigma_{l}$  
is a finite subset of $X$ such that $\Sigma_{l}+2l\mu(X^{\vee})=X$.
\end{cor}
\begin{proof}
Since $2l\mu(X^{\vee})\subset X$, we obtain 
by Lemma~\ref{lemma:Voronoi_cap_X=Sigmal},  
\begin{align*}X&=X\cap \bigcup_{c\in X^{\vee}}(2l\mu(c)+\Sigma_{l}(0))
=\bigcup_{c\in X^{\vee}}(2l\mu(c)+\Sigma_{l}(0)\cap X)\\
&=\bigcup_{c\in X^{\vee}}(2l\mu(c)+\Sigma_{l})
=2l\mu(X^{\vee})+\Sigma_{l}.
\end{align*}
\end{proof}

\begin{defn}\label{defn:integral Sigma-Voronoi cell}
A convex polytope $\Delta$ in $X_{\bR}$  is 
said to be {\it integral} if 
\begin{enumerate}
\item[(i)]$\Delta=\Conv(\Delta\cap X)$  
and $\Delta=-\Delta\ni 0$;
\item[(ii)] $\Delta\cap X$ contains a basis of $X$.
\end{enumerate}

If $\Delta$ is integral, then $\Delta\cap X$ 
is a star of \cite[2.2]{Mumford72}. 
\end{defn}

\begin{lemma}
\label{lemma:constructing integral Sigma_n}
There exists $l_0\in\bN$ such that   
$\Sigma_{l_0l}(0)$ is integral $(\forall l\in\bN)$.
\end{lemma}
\begin{proof}
Since $\Sk^0(\Sigma_1(0))\subset X\otimes_{\bZ}\bQ$, 
there exists $l_0\in\bN$ such that 
$\Sk^0(l_0\Sigma_1(0))$ is contained in $X$, and 
$l_0\Sigma_1(0)\cap X$ contains a basis of $X$.  Since $l_0\Sigma_1(0)^0$ contains $0$, $\Sigma_{l_0}(0)=l_0\Sigma_1(0)$ satisfies (i)-(ii) of Definition~\ref{defn:integral Sigma-Voronoi cell}, and so does $\Sigma_{l_0l}(0)$ for any $l\geq 1$. Hence $\Sigma_{l_0l}(0)$ is integral. 
\end{proof}

\begin{lemma}
\label{lemma:comparison of val}Let $l\in\bN$,
$\alpha\in X$, $a, b\in X^{\vee}$ and $\gamma\in\Sigma_{l}$ 
such that $\alpha+2l\mu(a)=\gamma+2l\mu(b)$. 
Then $\alpha\in\Sigma_{l}$ iff 
there exists a unit $c\in R^{\times}$ such that 
$\xi_{l,\alpha,a}=c\xi_{l,\gamma,b}$.
\end{lemma} 
\begin{proof}Since $\alpha+2l\mu(a)=\gamma+2l\mu(b)$ iff 
$\alpha=\gamma+2l\mu(b-a)$, the proof is reduced to the case $a=0$ because $\xi_{l,\alpha,0}=c\xi_{l,\gamma,b-a}$ $(\exists c\in R^{\times})$ iff $\xi_{l,\alpha,a}\theta_{l}=\delta_{a}^*(\xi_{l,\alpha,0}\theta_{l})=\delta_{a}^*(c\xi_{l,\gamma,b-a}\theta_{l})=c\xi_{l,\gamma,b}\theta_{l}$  $(\exists c\in R^{\times})$. \par
In what follows, we assume $a=0$. Let $z\in X^{\vee}$.  We define 
\begin{equation}\label{eq:defn Dl(x,y)}
D_{l}(\gamma, 2l\mu(z))
=v_t (\xi_{l,\gamma,z}/w^{\gamma+2l\mu(z)})=E_{l}(z)+z(\gamma).
\end{equation} 
By Lemma~\ref{lemma:Voronoi_cap_X=Sigmal}, $D_{l}(\gamma, 2l\mu(z))\geq 0$ for any $z\in X^{\vee}$. To prove Lemma, 
it suffices to prove that $D_{l}(\gamma, 2l\mu(b))=0$ 
iff $\gamma+2l\mu(b)\in\Sigma_{l}$. \par
First, we prove the if part. Assume $\alpha\in\Sigma_{l}$. 
Since $\alpha=\gamma+2l\mu(b)\in\Sigma_{l}(b)$, we have
$\alpha\in\Sigma_{l}(0)\cap\Sigma_{l}(b)$. 
Then $E_{l}(b)=b(\alpha)$ because  
\begin{align*}
0&=\|\alpha-2l\mu(b)\|^2-\|\alpha\|^2=-4l B(\mu(b),\alpha)+4l^2B(\mu(b),\mu(b))\\
&=-4Nl b(\alpha)+4Nl^2 b(\mu(b))=4Nl (-b(\alpha)+E_{l}(b)).
\end{align*}
Since $b(\gamma)=b(\alpha)-2E_{l}(b)=-b(\alpha)$, we have  
$D_{l}(\gamma, 2l\mu(b))=E_{l}(b)+b(\gamma)=0$.

Conversely, we shall 
prove that $\alpha\in\Sigma_{l}$ if $D_{l}(\gamma, 2l\mu(b))=0$. 
Since $\alpha=\gamma+2l\mu(b)$, we have 
 $0=D_{l}(\gamma,2l\mu(b))=E_{l}(b)+b(\gamma).$ 
Then 
\begin{align*}
E_{l}(z)+z(\alpha)
&=E_{l}(b)+b(\gamma)+E_{l}(z)+z(\gamma)+2l z(\mu(b))\\
&=E_{l}(b+z)+(b+z)(\gamma)\geq 0\quad (\forall z\in X^{\vee})
\end{align*}by Corollary~\ref{cor:El(x)}~(\ref{item:El sum}) 
and $\gamma\in\Sigma_{l}$. 
It follows that $\alpha\in\Sigma_{l}$. 
\end{proof}

\begin{defn}
\label{defn:Dl(x)}
Via (\ref{eq:defn Dl(x,y)}), 
we define a function $D_{l} : X\to \bZ$ by 
\begin{equation*}
D_{l}(x):=D_{l}(\gamma,2l\mu(z))
=E_{l}(z)+z(\gamma)
\end{equation*}where $x=\gamma+2l\mu(z)$ for 
$\gamma\in\Sigma_{l}$ and $z\in X^{\vee}$.  This is well-defined by Lemma~\ref{lemma:comparison of val}. Hence 
$\xi_{l,\gamma,z}$ is an $R^{\times}$-multiple 
of $t^{D_{l}(x)}w^x$ and $s^{e(\xi)D_{l}(x)}w^x$. 
\end{defn}

\begin{cor}
\label{cor:Dlalpha=0 iff alpha in Sigma}
 $D_{l}(\alpha)\geq 0$ $(\forall\alpha\in X)$, while 
$D_{l}(\alpha)=0$ iff $\alpha\in\Sigma_{l}$.
\end{cor}
\begin{proof}Clear from the proof of Lemma~\ref{lemma:comparison of val}.
\end{proof}

\begin{lemma}\label{lemma:formula of Dl(x)}
Let $C_{l}(y):=\frac{1}{4Nl}B(y,y)$ $(y\in X)$, $\gamma\in\Sigma_{l}$, $z\in X^{\vee}$ and $x:=\gamma+2l\mu(z)$. Then 
$D_{l}(x)=C_{l}(x)-C_{l}(\gamma)\leq C_{l}(x)$, with equality holding iff 
$\gamma=0$ and $x\in 2l\mu(X^{\vee})$. 
\end{lemma}
\begin{proof}It is easy to see that 
$E_{l}(z)+z(\gamma)=C_{l}(x)-C_{l}(\gamma)$. 
Since $C_{l}(\gamma)\geq 0$, we have $D_{l}(x)\leq C_{l}(x)$, 
with equality holding iff $\gamma=0$ and  
$x\in 2l\mu(X^{\vee})$. 
\end{proof}

\subsection{Examples}
\begin{example}
\label{example:E8}
Let $B$ be an even
\footnote{$B$ is defined to be {\it even} 
if $B(x,x)$ is even for any $x\in X$.} 
unimodular positive 
bilinear form on a lattice $X$ of rank $g\geq 8$. 
Let $\xi^e=(X,a,b^e,A,B)$ be a principal eFC datum. 
By Definition~\ref{defn:theta=prod delta_ over Delta}, 
$X^{\vee}/\beta(Y)=0$, $N=|B|=1$, $e(\xi)=1$ 
and $\mu\beta=\id_Y$.  
The Voronoi polytope of $B$ is defined to be 
\begin{align*}
V(0)&=\{x\in X_{\bR}; \|x-a\|\geq 
\|x\|\ (\forall a\in X)\}\\
&=\{x\in X_{\bR};B(a,a)+2B(a,x)\geq 0\ (\forall a\in X)\}.
\end{align*}

Meanwhile, we obtain
\begin{align*}
\Sigma_{l}(0)&
=\{x\in X_{\bR}; E_{l}(u)+u(x)\geq 0\ (\forall u\in X^{\vee})\}\\
&=\{x\in X_{\bR}; B(y,y)+2B(y,x/2l)\geq 0\ (\forall y\in Y)\}
\end{align*}because $\beta:Y\to X^{\vee}$ is an isomorphism. Hence 
$\Sigma_{l}(0)=2l V(0)$. \par
  
 Let $B:=E_8$,  $(\alpha_i;i\in [1,8])$ 
a positive root system of $E_8$, $X:=\bigoplus_{i=1}^8\bZ \alpha_i$ and 
$\alpha_0=\sum_{i=1}^8\lambda_i\alpha_i$ 
its maximal root. To be more explicit, 
\begin{equation*}
\alpha_0
=2\alpha_1+3\alpha_2+4\alpha_3+6\alpha_4
+5\alpha_5+4\alpha_6+3\alpha_7+2\alpha_8.
\end{equation*}
Let $(\omega_i;i\in[1,8])$ be the weights 
(as the elements of $X_{\bQ}$) 
with $B(\omega_i,\alpha_k)=\delta_{i,k}$, and 
$F:=\Conv(\omega_i/\lambda_i;i\in[1,8])$. 
Then $V(0)$ is   
  the union of $wF$ $(w\in W(E_8))$ where $W(E_8)$ is 
the Weyl group of $E_8$. 
See \cite[9.3, p.~263]{NSugawara06}. 
Since $|B|=1$, $(\omega_i;i\in[1,8])$ is a basis of $X$. 
Hence $\Sigma_{l}(0)=2l V(0)$ is integral    
iff $l$ is divisible by $30$. 
\end{example}

\begin{example}
\label{example:two variables}
Let $\xi=(X,a,b,A,B)$ be a principal FC datum 
with rank of $X$ equal to 2.  
As is well-known, any integral symmetric positive $2\times 2$ matrix 
is, up to $\GL(2,\bZ)$-equivalence, one of the following:
$$\begin{pmatrix}p+r&-r\\
-r&q+r
\end{pmatrix}\ (p,q>0,r\geq 0).$$

Therefore there exists a basis $(m_i; i=1,2)$ of $X$
such that the matrix $(B(m_i,m_j))$ is as above. Let 
$(f_i\in X^\vee; i=1,2)$ with $f_i(m_j)=\delta_{i,j}$,  
$y=y_1m_1+y_2m_2$ and $u=u_1f_1+u_2f_2$.  It is easy to see 
{\small\begin{gather*}
\beta(y)=((p+r)y_1-ry_2)f_1+(-ry_1+(q+r)y_2)f_2,\\
\mu(u)=((q+r)u_1+ru_2)m_1+(ru_1+(p+r)u_2)m_2,\\
\Sigma_1(0)=\left\{x_1m_1+x_2m_2\in X_{\bR}; 
\begin{matrix}|x_1|\leq q+r,\ |x_2|\leq p+r, 
 |x_1-x_2|\leq 
p+q\end{matrix}
\right\}.
\end{gather*}}  
Thus $\Sigma_1(0)$ is the convex closure of 
{\small\begin{gather*}
\pm ((q+r)m_1+(p+r)m_2),\ \pm ((q+r)m_1+(-p+r)m_2),\\
\pm ((-q+r)m_1+(p+r)m_2), 
\end{gather*}}which contains $m_1$ and $m_2$. It follows that 
$\Sigma_1(0)$ is integral.    
Hence $\Sigma_{l}(0)$ is integral $(\forall l\geq 1)$ because 
$\Sigma_{l}(0)=l\Sigma_{1}(0)$. 
\end{example}

\section{Relatively complete models}
\label{sec:rel complete models}
We use 
Notation~\ref{notation:Notation Rinit Sinit keta=Kmin(xi)}. 
Let $\xi=(X,Y,a,b,A,B)$ and $\xi^{\natural}_{l}
=(X,Y,\epsilon_{l},b^e,E_{l},\Sigma_{l})$ be the same as in 
\S~\ref{subsec:eFC data} and \S~\ref{subsec:Rl(xinatural)} 
respectively. 
Let $t$ (resp. $s$) be a uniformizer of $R_{\init}$ (resp. $R:=R_{\min}(\xi)$) 
and $e(\xi):=v_s(t)$. 
Let $\Sigma_{l}(0)$ be the Voronoi polytope for $c=0$ 
in \S~\ref{subsec:Voronoi cells}.  
In \S\S~\ref{sec:rel complete models}-\ref{sec:action of Gl}  
we assume 
\begin{equation}
\label{assump:Sigmal(0) integral}
\text{$\Sigma_{l}(0)$ is integral.}
\end{equation}

\subsection{The $S$-scheme $\tP_{l}$}
\label{subsec:tP_l}
Let $(\tQ_{l},\tcL_{l}):=(\Proj R_{l}(\xi^{\natural}),\cO_{\tQ_{l}}(1))$ and $\pi_{l}:\tP_{l}\to\tQ_{l}$ 
the normalization of $\tQ_{l}$. We denote the pullback of $\tcL_{l}$ to $\tP_{l}$ by the same $\tcL_{l}$ if no confusion is possible. 
For any $\alpha\in\Sigma_{l}$ and $u\in X^{\ \vee}$, let 
\begin{equation*}\label{eq:B_lalphau}
\begin{aligned}
A_{l, \alpha, u}&:=R[\xi_{l,\beta,v}/\xi_{l,\alpha,u}; 
\beta\in\Sigma_{l}, v\in X^{\vee}],\\
B_{l, \alpha, u}&:=\text{the normalization of $A_{l,\alpha,u}$}.
\end{aligned}
\end{equation*}  Hence $\Cone(A_{l, \alpha, u})=\Cone(B_{l, \alpha, u})$ by Definition~\ref{defn:toric ring basic}. 

\begin{lemma}
\label{lemma:finite generation} $A_{l, \alpha, u}$ is finitely generated over $R$ $(\forall \alpha\in\Sigma_{l}, \forall u\in X^{\vee})$. 
\end{lemma}
\begin{proof}See \cite[2.4, p.~290]{Mumford72}. 
Since $A_{l,\alpha,u}\simeq A_{l,\alpha,0}$, it suffices to prove Lemma 
when $u=0$. Since $v_s(t)=e(\xi)$ and $D_{l}(\alpha)=0$, 
by Corollary~\ref{cor:Sigma_l}, 
\begin{align*}
A_{l,\alpha,0}&=R[t^{D_{l}(x)}w^{x-\alpha};x\in\Sigma_{l}+2l\mu(X^{\vee})]=R[t^{D_{l}(x)}w^{x-\alpha};x\in X].
\end{align*}

We define $C_{l}(x,y)=\frac{1}{4Nl}B(x,y)$ $(x, y\in X)$. 
Then $C_{l}(x)=C_{l}(x,x)$ in Lemma~\ref{lemma:formula of Dl(x)}. 
Let $M:=\max\{C_{l}(\gamma);\gamma\in \Sigma_{3l}\}$. 
We define a subset $\Delta$ of $X$ by 
$$\Delta=\{x\in X; C_{l}(x)\leq C_{l}(x-\lambda)+2M\ 
(\forall\lambda\in \Sigma_{2l}-\alpha)\}.
$$
Since $\Sigma_{2l}-\alpha$ contains the star $\Sigma_{l}$,  
$\Delta$ is finite by the inclusion relation:
$$\Delta\subset \{x\in X; |2C_{l}(x,\lambda)|\leq 3M\ (\forall\lambda\in\Sigma_{l})\}.$$  

Let $\sigma_{l,x}:=t^{D_{l}(x)}w^{x-\alpha}$ $(x\in X)$. 
Let $A^*_{l,\alpha,0}$ be the $R$-subalgebra of 
$A_{l,\alpha,0}$ generated by 
$(\sigma_{l,v};v\in \Delta\cup \Sigma_{2l})$. 
We prove $\sigma_{l,x}\in A^*_{l,\alpha,0}$ $(\forall x\in X)$ 
by the induction on $D_{l}(x)$.  
If $D_{l}(x)=0$, then $x\in\Sigma_{l}$ 
by Corollary~\ref{cor:Dlalpha=0 iff alpha in Sigma}, 
so that  $\sigma_{l,x}\in A^*_{l,\alpha,0}$. 
Let $n\in\bN$, and assume that 
$\sigma_{l,y}\in A^*_{l,\alpha,0}$ if $D_{l}(y)<n$. 
We take any $x$ with $D_{l}(x)=n$. 
We write $x=\gamma+2l\mu(z)$ for some 
$\gamma\in\Sigma_{l}$ and $z\in X^{\vee}$. Then 
$D_{l}(x)=C_{l}(x)-C_{l}(\gamma)$ by Lemma~\ref{lemma:formula of Dl(x)}. 
Suppose that $x\not\in\Delta$. 
Then there exists $\lambda\in\Sigma_{2l}-\alpha$ such that 
$C_{l}(x)>  C_{l}(x-\lambda)+2M$. 
Let $y:=x-\lambda$ and $\mu:=\lambda+\alpha$. Then  
$C_{l}(x)> C_{l}(y)+2M$. Since $\gamma\in\Sigma_{l}$,  we have 
$C_{l}(\gamma)\leq M$. 
Since $M\geq C_{l}(\mu)\geq D_{l}(\mu)$, 
\begin{gather*}
D_{l}(x)=C_{l}(x)-C_{l}(\gamma)
>C_{l}(y)+M\geq D_{l}(y)+D_{l}(\mu),
\end{gather*}whence $n=D_{l}(x)>D_{l}(y)$ and 
$\sigma_{l,x}=t^c\sigma_{l,y}\sigma_{l,\mu}$ 
for some $c\geq 1$.  Hence $\sigma_{l,y}\in A^*_{l,\alpha,0}$ 
by the induction hypothesis, so that 
$\sigma_{l,x}\in A^*_{l,\alpha,0}$ by $\mu\in\Sigma_{2l}$. 
It follows that $A_{l,\alpha,0}=A^*_{l,\alpha,0}$, 
which is finitely generated over $R$. \end{proof}

By the proof of Lemma~\ref{lemma:finite generation}, 
$\Delta\cup\Sigma_{3l}$ is a bounded  subset of $X$. Hence we can choose $M^*>0$ such that $M^*\geq 16M$ and 
$$\Delta\cup\Sigma_{3l}\subset \{x\in X; C_{l}(x)\leq M^*\}.$$ 

\begin{lemma}\label{lemma:Iadic conv}Let $M^*$ be as above and $x\in X$. If $C_{l}(x)\geq 2^{n}M^*$ for some integer $n\geq 0$, 
then $\sigma_{l,x}\in I^nA_{l,\alpha,0}$. 
\end{lemma}
\begin{proof}We prove this by the induction on $n$. If $n=0$, then 
$\sigma_{l,x}\in A_{l,\alpha,0}$ by Lemma~\ref{lemma:finite generation}. 
Assume that $\sigma_{l,y}\in I^{n-1}A_{l,\alpha,0}$ for $y\in X$ 
if $C_{l}(y)\geq 2^{n-1}M^*$. 
Suppose that $C_{l}(x)\geq 2^{n}M^*$ for $x\in X$. 
Hence $x\not\in\Delta$.
By the proof of Lemma~\ref{lemma:finite generation}, 
there exists $\lambda\in \Sigma_{2l}-\alpha$ such that 
$\sigma_{l,x}=t^{c}\sigma_{l,y}\sigma_{l,\mu}$ for some $c\geq 1$  
where $y=x-\lambda$ and $\mu=\lambda+\alpha$. 
Since $B$ is positive definite, 
$C_{l}(x,\lambda)^2\leq C_{l}(x)C_{l}(\lambda)\leq MC_{l}(x)$. 
Since $C_{l}(x)\geq 2^{n}M^*\geq 16M$, 
{\small\begin{align*}
C_{l}(y)&\geq C_{l}(x)-2C_{l}(x,\lambda)\geq C_{l}(x)-2\sqrt{MC_{l}(x)}
\geq (1/2)C_{l}(x)\geq 2^{n-1}M^*.
\end{align*}} 

\vskip -0.4cm 
By the induction hypothesis,  
$\sigma_{l,y}\in I^{n-1}A_{l,\alpha,0}$, so that 
$\sigma_{l,x}\in I^{n}A_{l,\alpha,0}$.
\end{proof}

\begin{cor}\label{cor:Iadic convergence}
Let $m\geq 1$ and 
$\gamma\in\Sigma_{l m}$. We define 
$\vartheta^{(m)}_{l,\gamma}$ by
\begin{equation}\label{eq:defn of vartheta_m_l_gamma}
\vartheta^{(m)}_{l,\gamma}
=\sum_{y\in Y}a(y)^{2Nlm}b(y,\gamma)w^{\gamma+2Nlmy}. 
\end{equation}
Then 
$\vartheta^{(m)}_{l,\gamma}/\xi_{l,\alpha,u}^m
\in A_{l,\alpha,u}^{\wedge}$. 
\end{cor}
\begin{proof} Let $\gamma\in\Sigma_{l m}$. We choose 
and fix an $m$-tuple $(\gamma_i\in\Sigma_{l};i\in[1,m])$ 
with $\gamma=\sum_{i=1}^m\gamma_i$. 
Then $a(y)^{2Nlm}b(y,\gamma)w^{\gamma+2Nlmy}=\prod_{i=1}^m\xi_{l,\gamma_i,\beta(y)}$. Let   
$\sigma^{*}_{l,\gamma_i,\beta(y)}:
=\xi_{l,\gamma_i,\beta(y)}/\xi_{l,\alpha,u}
=\delta_u^*(\xi_{l,\gamma_i,\beta(y)-u}/\xi_{l,\alpha,0})$ and $\sigma_{l,x}:=t^{D_{l}(x)}w^{x-\alpha}$ $(x\in X)$. 
Then $\vartheta^{(m)}_{l,\gamma}/\xi_{l,\alpha,u}^m=\sum_{y\in Y}\prod_{i=1}^m\sigma^{*}_{l,\gamma_i,\beta(y)}$. 
Since $\xi_{l,\gamma_i,\beta(y)-u}/\xi_{l,\alpha,0}$ is 
an $R^{\times}$-multiple of $\sigma_{l,\gamma_i+2l\mu(\beta(y)-u)}$,  
we have $\sigma_{l,\gamma_i+2l\mu(\beta(y)-u)}\in I^nA_{l,\alpha,0}$ and 
hence $\sigma^{*}_{l,\gamma_i,\beta(y)}\in I^nA_{l,\alpha,u}$ by Lemma 7.2 
if $C_{l}(\gamma_i+2l\mu(\beta(y)-u))\geq 2^nM^*$. 
Hence    
{\small\begin{align*}\vartheta^{(m)}_{l,\gamma}/\xi_{l,\alpha,u}^m
&\equiv \sum_{\begin{matrix}
y\in Y\ \text{s.t. for $\exists \gamma_i$} \\
\text{$C_{l}(\gamma_i+2Nly-2l\mu(u))<2^nM^*$}
\end{matrix}}\prod_{i=1}^m\sigma^{*}_{l,\gamma_i,\beta(y)}\ 
\op{mod}\, I^{nm}A_{l,\alpha,0}^{\wedge} 
\end{align*}}for any $n\geq 1$, 
which is a finite sum.     
Hence $\vartheta^{(m)}_{l,\gamma}/\xi_{l,\alpha,u}^m
\in A^{\wedge}_{l,\alpha,u}$.
\end{proof}

Let $U_{l,\alpha, u}:=\Spec A_{l,\alpha, u}$ 
and $W_{l,\alpha, u}:=\Spec B_{l,\alpha, u}$. 
Then $\tQ_{l}$ (resp. $\tP_{l}$) 
is an $S$-scheme with an affine covering $(U_{l,\alpha,u};\alpha\in\Sigma_{l}, u\in X^{\vee})$\ (resp. $(W_{l,\alpha, u};\alpha\in\Sigma_{l}, u\in X^{\vee})$) such that $\pi_{l}(W_{l,\alpha, u})=U_{l,\alpha, u}$. 
See Definition~\ref{defn:R_l(xinatural)}.

If $\alpha+2l\mu(u)=\beta+2l\mu(v)$, then 
by Lemma~\ref{lemma:comparison of val}, 
$$A_{l,\alpha, u}=A_{l,\beta, v},\ B_{l,\alpha, u}=B_{l,\beta, v}.
$$

\begin{defn}
\label{defn:torus embedding of tildecX}
We define a fan $\Fan(\xi^{\natural}_{l})$ in $\tilde X^{\vee}$  by 
\begin{gather*}
\Fan(\xi^{\natural}_{l})
=\left\{\tau_{l,\alpha,u}\ \text{and their faces}; 
\alpha\in\Sigma_{l}, u\in X^{\vee}\right\},\\
\tau_{l,\alpha,u}={\left\{\begin{matrix}
v_0f_0+v\in \bR f_0\oplus X_{\bR}^{\vee};\ 
v_0\geq 0, v\in X_{\bR}^{\vee}, \forall y\in X\\ 
e(\xi)(D_{l}(y)-D_{l}(\alpha+2l\mu(u)))v_0+v(y-(\alpha+2l\mu(u)))\geq 0
\end{matrix}
\right\}.}
\end{gather*}

Since $A_{l,\alpha,u}$ is finitely generated over $R$ by Lemma~\ref{lemma:finite generation},  $\tau_{l,\alpha,u}$ has finitely many faces. Hence 
$\Fan(\xi^{\natural}_{l})$ is a fan in $\tilde X^{\vee}$ 
in the sense of \cite[p.~2]{Oda85}. 
\end{defn}

\begin{lemma}\label{lemma:Fanl is fan over S}
$\Fan(\xi^{\natural}_{l})$ is a fan in $\tX^{\vee}$ over $S$, and  
$\tP_{l}$ is isomorphic to the $S$-torus embedding 
$U(\Fan(\xi^{\natural}_{l}))$ with 
$U(\tau_{l,\alpha,u})=W_{l,\alpha, u}$\ $(\forall\alpha\in\Sigma_{l}, 
\forall u\in X^{\vee})$. 
\end{lemma}

\begin{proof}We shall check 
Definition~\ref{defn:fan over S}~(i) and (ii) for $\Fan(\xi^{\natural}_{l})$. 
The condition~\ref{defn:fan over S}~(i) 
follows from Definition~\ref{defn:torus embedding of tildecX}. 
Let $\tau:=\tau_{l,\alpha,u}$\ for $\alpha\in\Sigma_{l}$ 
and $u\in X^{\vee}$. We shall prove that $\tau\cap X^{\vee}=(0)$. 
Let $v\in\tau\cap X^{\vee}$. 
Since $v(y-(\alpha+2l\mu(u)))\geq 0$\ 
$(\forall y\in X)$ by Definition~\ref{defn:torus embedding of tildecX}, 
we have $v(y)=0$\ $(\forall y\in X)$, 
so that $v=0$ and 
$\tau\cap X^{\vee}=(0)$. 
Next, let $\sigma\in\Fan(\xi^{\natural}_{l})$. 
Since $\sigma$ is a face 
of $\tau_{l,\alpha,u}$ for some $\alpha\in X$ and $u\in X^{\vee}$,  
$\sigma\cap X^{\vee}\subset \tau_{l,\alpha,u}\cap X^{\vee}=(0)$. 
This proves the condition~\ref{defn:fan over S}~(ii). 
Hence $\Fan(\xi^{\natural}_{l})$ is a fan in $\tX^{\vee}$ over $S$. \par

By the remarks in Definitions~\ref{defn:toric ring basic} 
and \ref{defn:torus embedding of tildecX},
$$\Semi(B_{l,\alpha,u})\cap\tX
=\Sat(\Semi(A_{l,\alpha,u}))=\Cone(A_{l,\alpha,u})\cap\tX=\tau_{l,\alpha,u}^{\vee}\cap\tX.$$ 
Hence $W_{l,\alpha,u}
=U(\tau_{l,\alpha,u})$, 
so that 
$\tP_{l}$ is isomorphic to $U(\Fan(\xi^{\natural}_{l}))$. 
\end{proof}

\subsection{The completeness of $\tP_{l}$}
\label{subsec:completeness}
\begin{table}[ht]
 \begin{centering}
 \renewcommand{\arraystretch}{1.4}
 \renewcommand{\arraycolsep}{0.6em}
 $\begin{array}{|c|c|c|c|c|c|c|}
 \hline
\text{Here} & \theta_{l}&\mu_{l}&X&X^{\vee}
&u,v\in X^{\vee}&\delta_u,\delta_v\\
\text{\cite{Mumford72}}&\theta&\phi&X&Y&y,z\in Y&S_y,S_z\\
\hline
\text{Here}&\Sigma_{l}&\alpha\in\Sigma_{l}&w^{\alpha}&b^e(u,\alpha)
&b^e(u,\mu_{l}(v))&\epsilon_{l}(u)\\
\text{\cite{Mumford72}}&\Sigma&\alpha\in\Sigma&\cX^{\alpha}
&\cX^{\alpha}(y)&\cX^{\phi(z)}(y)&\cX^{\phi(y)}(y)\\
 \hline
 \end{array}$
 \caption{Correspondence of data}
 \label{table:comparison}
 \end{centering}
 \end{table}

 For the N\'eFC kit $\xi^{\natural}_{l}$, 
there is an obvious correspondence between 
Definition~\ref{defn:R_l(xinatural)} and 
\cite[2.3]{Mumford72} as in Table~\ref{table:comparison}.
Hence $R_{l}(\xi^{\natural})$ can be regarded as a particular case 
of $R_{\phi,\Sigma}$ in \cite[2.3]{Mumford72}. 
\begin{lemma}\label{lemma:completeness}Let $K'$ be any algebraic field extension of $k(\eta)$ and $R'$ the integral closure of $R$ in $K'$. Then 
$\tP_{l}(K')=\tP_{l}(R')$. 
\end{lemma}
\begin{proof}See \cite[III, 3.3]{FC90}. 
We may assume that $K'=k(\eta)$ and $R'=R$. 
Let $f\in\tQ_{l}(k(\eta))$ 
be a $k(\eta)$-morphism.  
Since $\tQ_{l,k(\eta)}\simeq T_{X,k(\eta)}$, $f$ induces a 
$k(\eta)$-homomorphism $f^*:k(\eta)[w^x;x\in X]\to k(\eta)$. 
Let $\nu:=v_t\circ f^*$. 
Since $\nu(t)=1$, $\nu(\xi_{l,\beta,u})
=E_{l}(u)+u(\beta)+\nu(w^{\beta+2l\mu(u)})$ 
by Definition~\ref{defn:R_l(xinatural)}. 
Since $F$ is positive definite 
by Lemma~\ref{lemma:Be and Bephi}~(\ref{item:F pos def}), 
$\Delta_{\beta}:=\{u\in X^{\vee}: \nu(\xi_{l,\beta,u})\leq 0\}$ 
is finite $(\forall \beta\in\Sigma_{l})$. 
Since $\nu(\xi_{l,0,0})=0$, we obtain 
$$
\min(\nu(\xi_{l,\beta,u});\beta\in\Sigma_{l}, u\in X^{\vee})=
\min(\nu(\xi_{l,\beta,u});\beta\in\Sigma_{l}, u\in \Delta_{\beta}). 
$$ Hence there exists $(\alpha_0, u_0)$ such that 
$\nu(\xi_{l,\alpha_0, u_0})$ is the minimum. Then 
$$\nu(\xi_{l,\beta,u}/\xi_{l,\alpha_0,u_0})\geq 0\ (\forall \beta\in\Sigma_{l}, \forall u\in X^{\vee}).$$ Therefore $f^*(A_{l,\alpha_0,u_0})\subset R$, so that $f\in U_{l,\alpha_0,u_0}(R)$. 
Hence $\tQ_{l}(k(\eta))=\tQ_{l}(R)$. 
Since $\tP_{l}$ is finite over $\tQ_{l}$, $\tP_{l}(k(\eta))=\tP_{l}(R)$.  
\end{proof} 

\subsection{The $S$-scheme $P_{l}$}
\label{subsec:theta on Pl}
The $S$-scheme $\tP_{l}$ (resp. $\tQ_{l}$) 
has an action $(S_y;y\in Y)$ of $Y$. 
 Let $(\tP_{l}^{\wedge},\tcL_{l}^{\wedge})/Y$ 
(resp. $(\tQ_{l}^{\wedge},\tcL_{l}^{\wedge})/Y$) be 
the formal quotient by $Y$ of $(\tP_{l}^{\wedge},\tcL_{l}^{\wedge})$ (resp. $(\tQ_{l}^{\wedge},\tcL_{l}^{\wedge})$). 
By \cite[III$_1$, 5.4.5]{EGA}, there exist  projective $S$-schemes $(P_{l},\cL_{l})$ and 
 $(Q_{l},\cL_{l})$ such that 
\begin{gather*}
(P_{l}^{\wedge},\cL_{l}^{\wedge})\simeq 
(\tP_{l}^{\wedge},\tcL_{l}^{\wedge})/Y,\  
(Q_{l}^{\wedge},\cL_{l}^{\wedge})\simeq 
(\tQ_{l}^{\wedge},\tcL_{l}^{\wedge})/Y.
\end{gather*} 

\begin{defn}
\label{defn:twisted Mumford family}
We call $(P_{l},\cL_{l})$ the {\it twisted Mumford family 
associated with $\xi^{\natural}_{l}$}. We denote $(P_{l},\cL_{l})$ 
(resp. $(Q_{l},\cL_{l})$) 
by $(P_{l}(\xi^{\natural}),\cL_{l}(\xi^{\natural}))$ 
(resp. $(Q_{l}(\xi^{\natural}),\cL_{l}(\xi^{\natural}))$) 
if necessary. Recall that $R_{l}(\xi^{\natural})\simeq 
R(\xi^{\natural}_{l})$ by Lemma~\ref{lemma:Rlxinatural=Rxinaturall}. 
\end{defn}

\begin{defn}\label{defn:Y invariance}
Let $k(\eta)[[X]]:=\prod_{x\in X}k(\eta)w^x$.
Since $(S_y;y\in Y)$ acts on $k(\eta)[X][\theta]$ by 
Definition~\ref{defn:theta=prod delta_ over Delta}, 
it acts on $k(\eta)[[X]]\theta_{l}^m$ 
where $\theta_{l}=\theta^{l}$. 
Let $\psi\in k(\eta)[[X]]\theta_{l}^m$. 
Then we say that $\psi$ is {\it $Y$-invariant} if $S_y^*(\psi)=\psi$ 
$(\forall y\in Y)$. 
Let $(k(\eta)[[X]]\theta_{l}^m)^{Y\op{-inv}}$ be the $k(\eta)$-vector subspace of $k(\eta)[[X]]\theta_{l}^m$    
consisting of $Y$-invariant elements.\end{defn}

\begin{rem}\label{rem:subspace Wm}
Let $m\in\bN$ and 
\begin{gather*}
W_{1,m}:=\Gamma(\tQ_{l}^{\wedge},(\tcL_{l}^{\otimes m})^{\wedge})\otimes_Rk(\eta),\ \  
W_{2,m}:=\Gamma(\tP_{l}^{\wedge},(\tcL_{l}^{\otimes m})^{\wedge})\otimes_Rk(\eta).
\end{gather*}

Since $\Sigma_{l}+2l\mu(X^{\vee})=X$ by Corollary~\ref{cor:Sigma_l}, 
we have $B_{l,\alpha,u}\otimes_Rk(\eta)
=A_{l,\alpha,u}\otimes_Rk(\eta)=k(\eta)[X]$, so that 
\begin{align}\label{eq:Ahat=Bhat}
B_{l,\alpha,u}^{\wedge}\otimes_Rk(\eta)
&=A_{l,\alpha,u}^{\wedge}\otimes_Rk(\eta)
\subset k(\eta)[[X]].
\end{align}
It is clear from (\ref{eq:Ahat=Bhat}) that $W_{1,m}=W_{2,m}$, 
which we denote by $W_m$.
\end{rem}
\begin{defn}\label{defn:i(f) Y inv}Let $f\in W_m$. 
Then $f$  is a collection 
$(f_{l,\alpha,u}\in A_{l,\alpha,u}^{\wedge}\otimes_Rk(\eta); 
\alpha\in\Sigma_{l}, u\in X^{\vee})$ such that 
$f_{l,\alpha,u}
=(\xi_{l,\beta,v}/\xi_{l,\alpha,u})^mf_{l,\beta,v}$\ 
$(\forall \alpha,\beta\in\Sigma_{l}, \forall u, v\in X^{\vee})$.
We define  
\begin{equation}\label{eq:defn of iqfq}
i(f)=f_{l,\alpha,u}\xi_{l,\alpha,u}^m\theta_{l}^m
\in k(\eta)[[X]]\theta_{l}^m,
\end{equation}
which is independent of the choice of $(\alpha,u)$.
 We say that $f\in W_{m}$ is {\it $Y$-invariant} if 
$i(f)$ is $Y$-invariant in the sense of Definition~\ref{defn:Y invariance}. 
Let $W_{m}^{Y\op{-inv}}$ be the $k(\eta)$-vector subspace 
of $W_m$ consisting of $Y$-invariant elements.
\end{defn}

\begin{lemma}
\label{lemma:Gamma(Qhat Lhatm) over keta}
Let $m\in\bN$ and $N:=|B_{\xi^{\natural}}|$. 
Then 
\begin{align*}
W_{m}^{Y\op{-inv}}
&\simeq (k(\eta)[[X]]\theta_{l}^m)^{Y\op{-inv}}
\simeq\Gamma(G_{\eta},\cL_{\eta}^{\otimes 2Nlm}).
\end{align*} 
\end{lemma}
\begin{proof}Let $U_m:=(k(\eta)[[X]]\theta^m_{l})^{Y\op{-inv}}$ and 
$V_{2Nlm}:=\Gamma(G_{\eta},\cL_{\eta}^{\otimes 2Nlm})$. 
By Remark~\ref{rem:subspace Wm}, 
$W_m^{Y\op{-inv}}\simeq i(W_m^{Y\op{-inv}})\subset U_m$. 
Let $\psi:=\sum_{x\in X}\sigma_x(\theta)w^x\theta_{l}^m\in U_m$ and 
$\theta:=j(\psi):=\psi/\theta_{l}^m$.    
By Lemma~\ref{lemma:lth basic formulae}~(\ref{item:lth Sy walpha theta}), 
$S_y^*(\psi)=\psi$ is 
equivalent to  
\begin{equation}
\label{eq:expansion of theta}
\sigma_{x+2Nlmy}(\theta)=a(y)^{2Nlm}b(y,x)
\sigma_x(\theta)\ (\forall x\in X). 
\end{equation} 

This proves $j(U_m)=V_{2Nlm}$ 
by Lemma~\ref{lemma:H0(G,Lm) over keta as Fourier}. 
Let $\theta\in V_{2Nlm}$. 
By Eq.~(\ref{eq:expansion of theta}), 
$\theta$ is written as 
$\theta=\sum_{\gamma\in 
\Sigma_{lm}}c_{\gamma}\vartheta^{(m)}_{l,\gamma}$ for 
some $c_{\gamma}\in k(\eta)$ where
$\vartheta^{(m)}_{l,\gamma}$ is given by 
Eq.~(\ref{eq:defn of vartheta_m_l_gamma}). 
Let $f_{l,\alpha,u}:=\theta/\xi^{m}_{l,\alpha,u}$  
and  $f:=(f_{l,\alpha,u}; \alpha\in\Sigma_{l},u\in X^{\vee})$. 
Since $f_{l,\alpha,u}\in A_{l,\alpha,u}^{\wedge}\otimes_Rk(\eta)$ 
by Corollary~\ref{cor:Iadic convergence}, 
we have $f\in W_m$, so that $\theta=(j\circ i)(f)\in (j\circ i)(W_m)$. 
Hence $V_{2Nlm}\subset (j\circ i)(W_m)\subset j(U_m)=V_{2Nlm}$.  
Since $i$ and $j$ are injective, 
this proves Lemma.
\end{proof}

\begin{lemma}\label{lemma:theta on Pl}
There exists $m_0\in\bN$ such that 
for any $m\geq m_0$ 
{\small\begin{equation*}
\label{eq:Gamma(Gleta Leta^m)}
\Gamma(P_{l,\eta},\cL_{l,\eta}^{\otimes m})=\left\{
\theta=\sum_{x\in X}\sigma_x(\theta)w^x;
\begin{matrix}
\sigma_{x+2Nlmy}(\theta)
=a(y)^{2Nlm}b(y,x)\sigma_x(\theta)\\
\sigma_x(\theta)\in k(\eta)\ (\forall x\in X, \forall y\in Y)
\end{matrix}
\right\}.
\end{equation*}}
\end{lemma}
\begin{proof}Since $\cL_{l}$ and $\cL_{l}^{\wedge}$ 
are ample, by Serre's vanishing theorem \cite[II, 5.17]{Hartshorne77}, 
there exists $m_0\in\bN$ such that for any $m\geq m_0$
\begin{gather*}
H^q(P_{l,0},\cL_{l,0}^{\otimes m})=H^q(P_{l},\cL_{l}^{\otimes m})
=H^q(P^{\wedge}_{l},(\cL^{\wedge}_{l})^{\otimes m})=0 
\ (\forall q\geq 1).
\end{gather*} 

By \cite[Cor.~3, p.\ 50]{Mumford12},  
\begin{equation}\label{eq:Pl0 Pl0wedge}
H^0(P_{l},\cL_{l}^{\otimes m})\otimes_Rk(0)\simeq 
H^0(P_{l,0},\cL_{l,0}^{\otimes m})\simeq  
H^0(P^{\wedge}_{l},(\cL^{\wedge}_{l})^{\otimes m})\otimes_Rk(0).
\end{equation}
Since $P_{l}$ is normal and projective, 
both $\Gamma(P_{l},\cL_{l}^{\otimes m})$ and 
$\Gamma(P_{l}^{\wedge},(\cL_{l}^{\wedge})^{\otimes m})$ 
are finite torsion-free and hence free $R$-modules.  
 Let $s\in \Gamma(P_{l},\cL_{l}^{\otimes m})$. 
Since the pullback of $s$ to $\tP_{l}^{\wedge}$ is $S_y^*$-invariant 
$(\forall y\in Y)$, we have a monomorphism $\iota^*$
\begin{gather*}
\Gamma(P_{l},\cL_{l}^{\otimes m})\overset{\iota^*}{\hookrightarrow} 
\Gamma(P_{l}^{\wedge},(\cL^{\wedge}_{l})^{\otimes m})
\simeq \Gamma(\tP_{l}^{\wedge},(\tcL^{\wedge}_{l})^{\otimes m})
^{\text{$Y$-inv}}.
\end{gather*}

By Eq.~(\ref{eq:Pl0 Pl0wedge}), we see 
$\Gamma(P_{l},\cL_{l}^{\otimes m})
\simeq \Gamma(P_{l}^{\wedge},(\cL_{l}^{\wedge})^{\otimes m})$ 
by applying Nakayama's lemma. 
Hence 
$\Gamma(P_{l,\eta},\cL_{l,\eta}^{\otimes m})\simeq
\Gamma(P_{l},\cL_{l}^{\otimes m})\otimes_R k(\eta)\simeq 
\Gamma(P_{l}^{\wedge},(\cL_{l}^{\wedge})^{\otimes m})\otimes_R k(\eta).$

Hence Lemma follows from Lemma~\ref{lemma:Gamma(Qhat Lhatm) over keta}.
\end{proof}

\subsection{The $R$-algebra $B_{l,\alpha,u}$}
\label{subsec:Blalphau}
Let $\Sk(\Sigma_{l}(0))$ (resp. $\Sk^q(\Sigma_{l}(0))$) 
be the set of all faces of $\rho$ (resp. all $q$-dimensional 
faces of $\rho$). 

\begin{lemma}
\label{lemma:Defn of A_lrho0}  
Let $\rho\in\Sk(\Sigma_{l}(0))$, 
$\rho^0$ the relative interior of $\rho$ and $u\in X^{\vee}$. Assume that $\alpha\in\rho^0\cap\Sigma_{l}$ and $\beta\in\rho\cap\Sigma_{l}$. Then 
$\xi_{l,\alpha,u}/\xi_{l,\beta,u}\in A_{l,\alpha,u}$ and 
$A_{l,\alpha,u}\supset A_{l,\beta,u}$. 
If moreover $\beta\in\rho^0\cap\Sigma_{l}$, then 
$A_{l,\alpha,u}=A_{l,\beta,u}$. 
\end{lemma}
\begin{proof} 
Since $\alpha\in\rho^0\cap\Sigma_{l}$, 
$\beta\in\rho\cap\Sigma_{l}$ and 
$\rho$ is convex, 
there exist some $\beta_i\in\rho\cap\Sigma_{l}$ 
and $s_i\in\bQ$ $(i\in[1,r])$ for some finite $r$ 
such that $\beta_1=\beta$ and  
$\alpha=\sum_{i=1}^rs_i\beta_i,\ s_i>0,\ \sum_{i=1}^rs_i=1.$
We write $s_i=n_i/n$ 
for $n,n_i\in\bN$. 
Then $\sum_{i=1}^rn_i(\beta_i-\alpha)=0$, so that  
$$\xi_{l,\alpha,u}/\xi_{l,\beta,u}=
(\xi_{l,\beta,u}/\xi_{l,\alpha,u})^{n_1-1}\prod_{i=2}^r
(\xi_{l,\beta_i,u}/\xi_{l,\alpha,u})^{n_i}\in A_{l,\alpha,u}.
$$

Therefore $A_{l,\beta,u}\subset A_{l,\alpha,u}$ because
$$\xi_{l,\gamma,v}/\xi_{l,\beta,u}=(\xi_{l,\gamma,v}/\xi_{l,\alpha,u})
(\xi_{l,\alpha,u}/\xi_{l,\beta,u})\in A_{l,\alpha,u}\ 
(\forall \gamma\in\Sigma_{l}, \forall v\in X^{\vee}).$$  

If moreover $\beta\in\rho^0\cap\Sigma_{l}$,  then 
$A_{l,\alpha,u}\subset A_{l,\beta,u}$, whence  
$A_{l,\alpha,u}=A_{l,\beta,u}$. 
\end{proof}

\begin{defn}\label{defn:Alrho0}
For $\rho\in\Sk(\Sigma_{l}(0))$ and $u\in X^{\vee}$, we define
\begin{align*}
A_{l,\rho,u}&=R[\xi_{l,\gamma,v}/\xi_{l,\beta,u}; \beta\in\rho\cap\Sigma_{l}, \gamma\in\Sigma_{l},v\in X^{\vee}],\\
B_{l,\rho,u}&=\text{the normalization of $A_{l,\rho,u}$.}
\end{align*}
\end{defn}

\begin{lemma}\label{lemma:Alrho0 as intersection} 
Let $\rho\in\Sk(\Sigma_{l}(0))$ and $u\in X^{\vee}$. 
Then 
\begin{align*}
A_{l,\rho,u}
&=R[\xi_{l,\gamma,v}/\xi_{l,\beta,u};\beta\in\Sk^0(\rho), 
\gamma\in\Sigma_{l},v\in X^{\vee}].
\end{align*}

If $\alpha\in\rho^0\cap\Sigma_{l}$, then $A_{l,\rho,u}=A_{l,\alpha,u}$ and 
$B_{l,\rho,u}=B_{l,\alpha,u}$.
\end{lemma}
\begin{proof} 
Let 
\begin{gather*}
A^*:=A_{l,\rho,u},\ 
B^*:=R[\xi_{l,\gamma,v}/\xi_{l,\beta,u};\beta\in\Sk^0(\rho), 
\gamma\in\Sigma_{l},v\in X^{\vee}].\end{gather*}

It is clear that $B^*\subset A^*$.  
It remains to prove $A^*\subset B^*$. 
Suppose $\beta\in\rho\cap\Sigma_{l}$ but $\beta\not\in\Sk^0(\rho)$. 
Let $\sigma$ be the minimal face of $\rho$ such that $\beta\in\sigma^0$. 
There exist $\beta_i\in\Sk^0(\sigma)\ (\subset\Sk^0(\rho))$ and 
$n_i\in\bN$\ $(i\in[1,r])$ for some finite $r$ such that  
$\sum_{i=1}^rn_i(\beta_i-\beta)=0.$ Then  
$$\xi_{l,\beta_1,u}/\xi_{l,\beta,u}=(\xi_{l,\beta,u}/\xi_{l,\beta_1,u})^{n_1-1}\prod_{i=2}^r(\xi_{l,\beta,u}/\xi_{l,\beta_i,u})^{n_i}\in B^*.
$$ 

Hence $\xi_{l,\gamma,v}/\xi_{l,\beta,u}=(\xi_{l,\gamma,v}/\xi_{l,\beta_1,u})(\xi_{l,\beta_1,u}/\xi_{l,\beta,u})\in B^*$ \ $(\forall\gamma\in\Sigma_{l},\forall v\in X^{\vee})$.  It follows that $A^*\subset B^*$. 
This proves the first assertion. \par 
Next, suppose 
$\alpha\in\rho^0\cap\Sigma_{l}$. 
Since $A_{l,\rho,u}$ is 
generated by all $A_{l,\beta,u}$ 
$(\beta\in\rho\cap\Sigma_{l}$, $A_{l,\rho,u}\subset A_{l,\alpha,u}$ 
by Lemma~\ref{lemma:Defn of A_lrho0}, while $A_{l,\rho,u}\supset A_{l,\alpha,u}$ is clear. Hence $A_{l,\rho,u}=A_{l,\alpha,u}$, so that 
$B_{l,\rho,u}=B_{l,\alpha,u}$. 
This completes the proof. 
\end{proof}

\begin{cor}
\label{cor:Wlalphau as covering}The $S$-scheme $\tQ_{l}$ 
(resp. $\tP_{l}$) admits an affine covering 
 $(U_{l,\beta,u};\beta\in\Sk^0(\Sigma_{l}(0)) , u\in X^{\vee})$ (resp. $(W_{l,\beta,u};\beta\in\Sk^0(\Sigma_{l}(0)) , u\in X^{\vee})$). 
\end{cor}
\begin{proof}If 
$\alpha\in\rho^0\cap\Sigma_{l}$ for $\rho\in\Vor_{l}$, 
then by Lemma~\ref{lemma:Alrho0 as intersection}, 
$$U_{l,\alpha,u}=\bigcap_{\beta\in\Sk^0(\rho)} U_{l,\beta,u},\ 
W_{l,\alpha,u}=\bigcap_{\beta\in\Sk^0(\rho)} W_{l,\beta,u}.
$$  

If there exists no $\rho\in\Sk^q(\Vor_{l})$ $(q\geq 1)$ 
such that $\alpha\in\rho^0\cap\Sigma_{l}$, then $\alpha\in\Sk^0(\Sigma_{l}(0))$. This completes the proof. 
\end{proof}

\begin{defn}We define an $R$-subalgebra $R_{l}^*$ 
of $R_{l}(\xi^{\natural})$ by 
\begin{gather*}
R_{l}^*=R_{l}^*(\xi^{\natural}):=R[\xi_{l,\alpha,u}\theta_{l};\alpha\in\Sk^0(\Sigma_{l}(0)), u\in X^{\vee}],\\
C_{l,\alpha,u}=
R[\xi_{l,\gamma,v}/\xi_{l,\alpha,u}; 
\gamma\in\Sk^0(\Sigma_{l}(0)),v\in X^{\vee}]\\ 
(\alpha\in\Sk^0(\Sigma_{l}(0)), u\in X^{\vee}). 
\end{gather*} 

Let $\tT_{l}:=\Proj R_{l}^*$. 
The $S$-scheme $\tT_{l}$ admits an affine covering 
$$T_{l,\alpha,u}:=\Spec C_{l,\alpha,u}\ 
 (\alpha\in\Sk^0(\Sigma_{l}(0)), u\in X^{\vee}).$$ 

Since $R_{l}^*$ is an $R$-subalgebra of $R_{l}(\xi^{\natural})$, the inclusion $R_{l}^*\hookrightarrow R_{l}(\xi^{\natural})$ induces an $S$-morphism $h_{l}:\tQ_{l}\to \tT_{l}$ such that $h_{l}(U_{l,\alpha,u})=T_{l,\alpha,u}$\ $(\forall\alpha\in\Sk^0(\Sigma_{l}(0)), \forall u\in X^{\vee})$. 
\end{defn}
\begin{defn}
\label{defn:Clrho0}For $\rho\in\Sk(\Sigma_{l}(0))$ 
and $u\in X^{\vee}$,  we define 
\begin{align*}
C_{l,\rho,u}&=
R[\xi_{l,\gamma,v}/\xi_{l,\beta,u}; \beta\in\rho\cap\Sigma_{l}, 
\gamma\in\Sk^0(\Sigma_{l}(0)),v\in X^{\vee}]. 
\end{align*}
\end{defn}

\begin{lemma}\label{lemma:C rho and beta}
 Let $\rho\in\Sk(\Sigma_{l}(0))$ and $u\in X^{\vee}$. Then 
\begin{align*}
C_{l,\rho,u}
&=R[\xi_{l,\gamma,v}/\xi_{l,\beta,u}; \beta\in\Sk^0(\rho), \gamma\in\Sk^0(\Sigma_{l}(0)),v\in X^{\vee}].
\end{align*}

If 
$\alpha\in\rho^0\cap\Sigma_{l}$, then $C_{l,\rho,u}=C_{l,\alpha,u}$.  
\end{lemma}
\begin{proof}Proof is similar to those of Lemmas~\ref{lemma:Defn of A_lrho0}/\ref{lemma:Alrho0 as intersection}. 
\end{proof}

\begin{lemma}\label{lemma:integral closure of Alrho0} 
Let $\rho\in\Sk(\Sigma_{l}(0))$ and $u\in X^{\vee}$. Then $B_{l,\rho,u}$ is 
the normalization of $C_{l,\rho,u}$.
\end{lemma} 
\begin{proof} 
Let 
$B^*:=R[\xi_{l,\gamma,v}/\xi_{l,\beta,u};\beta\in\Sk^0(\rho), 
\gamma\in\Sigma_{l},v\in X^{\vee}]$ and $C^*:=C_{l,\rho,u}$. 
Let $B^{\natural}$ (resp. $C^{\natural}$) be the integral closure of $B^*$ (resp. $C^*$) in $k(\eta)[X]$. 
By Lemma~\ref{lemma:C rho and beta}, $C^*\subset B^*$, 
hence $C^{\natural}\subset B^{\natural}$. 

Next,  
we shall prove $B^*\subset C^{\natural}$. 
Let $\gamma\in\Sigma_{l}$. Let 
$\sigma$ be the minimal face of 
$\Sigma_{l}(0)$ such that $\gamma\in\sigma^0$. Then there exist 
$\gamma_i\in\Sk^0(\sigma)$ and $s_i\in\bQ$ $(i\in [1,m])$ 
for some finite $m$ such that $\gamma=\sum_{i=1}^ms_i\gamma_i$, 
$s_i>0$ and $\sum_{i=1}^ms_i=1$. We write $s_i=n_i/n$ 
with $n,n_i\in\bN$, so that   
$n\gamma=\sum_{i=1}^mn_i\gamma_i$ and $n=\sum_{i=1}^mn_i$.  
For any $\beta\in\Sk^0(\rho)$, we have  
$$\prod_{i=1}^m(\xi_{l,\gamma_i,v}/\xi_{l,\beta,u})^{n_i}=(\xi_{l,\gamma,v}/\xi_{l,\beta,u})^n.$$  

Hence $\xi_{l,\gamma,v}/\xi_{l,\beta,u}$ is integral over $C^*$, 
so that $B^*\subset C^{\natural}$. It follows 
that $B^{\natural}=C^{\natural}$. 
Since $A_{l,\rho,u}=B^*$ by Lemma~\ref{lemma:Alrho0 as intersection}, 
$B_{l,\rho,u}=B^{\natural}=C^{\natural}$, 
which is the normalization of $C^*=C_{l,\rho,u}$. 
This completes the proof. 
\end{proof}

\begin{cor}\label{cor:cones Alrho0 etc}
Let $\rho\in\Sk(\Sigma_{l}(0))$, $\alpha\in\rho^0\cap\Sigma_{l}$ and 
$u\in X^{\vee}$. Then 
\begin{gather*}
\Cone(A_{l,\rho,u})=\Cone(B_{l,\rho,u})=\Cone(C_{l,\rho,u}),\\
\Cone(A_{l,\rho,u})=\Cone(A_{l,\alpha,u})=
\Cone(B_{l,\alpha,u})=\Cone(C_{l,\alpha,u}).
\end{gather*}
\end{cor}
\begin{proof}This follows from Definition~\ref{defn:toric ring basic} 
and Lemmas~\ref{lemma:Alrho0 as intersection}/\ref{lemma:C rho and beta}/\ref{lemma:integral closure of Alrho0}. 
\end{proof}

\begin{cor}
\label{cor:normalization of tTl}
The $S$-morphism $h_{l}\pi_{l}:\tP_{l}\to \tT_{l}$ is the normalization.
\end{cor}
\begin{proof} Both $A_{l,\alpha,u}$ and $C_{l,\alpha,u}$ are $R$-subalgebras of $k(\eta)[X]$ such that $C_{l,\alpha,u}\subset A_{l,\alpha,u}$, and the homomorphism $h_{l,\alpha,u}^*:C_{l,\alpha,u}\to A_{l,\alpha,u}$ induced from $h_{l}$ is the inclusion map. The normalization $\pi_{l}:\tP_{l}\to \tQ_{l}$ induces a homomorphism $\pi_{l,\alpha,u}^*:A_{l,\alpha,u}\to B_{l,\alpha,u}$, so that the composite $h_{l}\pi_{l}$ induces the inclusion map $(h_{l}\pi_{l})^*:C_{l,\alpha,u}\subset B_{l,\alpha,u}$. 
Hence Corollary follows from Lemma~\ref{lemma:integral closure of Alrho0}. 
\end{proof}

\subsection{The isomorphism $P_{l}\simeq P_{ll'}$}
\begin{lemma}
\label{lemma:Sigmall' = l'Sigmal}
For any $l'\in\bN$, 
\begin{enumerate}
\item $\Sigma_{ll'}(0)=l'\Sigma_{l}(0)$; 
\item $\Sk^0(\Sigma_{ll'}(0))=l'\Sk^0(\Sigma_{l}(0)):=\{l'\alpha; \alpha\in \Sk^0(\Sigma_{l}(0))\}$.
\end{enumerate} 
\end{lemma}
\begin{proof}Clear. 
\end{proof}

\begin{lemma}
\label{lemma:cone Alrho0 and All'rho0} For any $l'>0$,
$\Cone(B_{l,\alpha,u})=\Cone(B_{ll',l'\alpha,u})$. 
\end{lemma}
\begin{proof} 
By Corollary~\ref{cor:cones Alrho0 etc}, 
$\Cone(C_{ll',l'\alpha,u})=\Cone(A_{ll',l'\alpha,u})$ for any $l'$. Hence  
{\small\begin{align*}
&\Cone(B_{ll',l'\alpha,u})=\Cone(A_{ll',l'\alpha,u})
=\Cone(C_{ll',l'\alpha,u})\\
&=\begin{matrix}
\text{the cone over $\Conv(\wt(\xi_{ll',\gamma,v}/\xi_{ll',l'\alpha,u});\gamma\in\Sk^0(\Sigma_{ll'}(0)), v\in X^{\vee})$}
\end{matrix}\\
&\overset{*}{=}\begin{matrix}
\text{the cone over 
 $\Conv(\wt(\xi_{ll',l'\gamma,v}/\xi_{ll',l'\alpha,u});
\gamma\in\Sk^0(\Sigma_{l}(0)), v\in X^{\vee})$}
\end{matrix}\\
&\overset{**}{=}\begin{matrix}
\text{the cone over  
 $\Conv(\wt(\xi_{l,\gamma,v}/\xi_{l,\alpha,u});
\gamma\in\Sk^0(\Sigma_{l}(0)), v\in X^{\vee})$}
\end{matrix}\\
&=\Cone(C_{l,\alpha,u})=\Cone(A_{l,\alpha,u})=\Cone(B_{l,\alpha,u}), 
\end{align*}}where the equality with * 
(resp. **) follows from Lemma~\ref{lemma:Sigmall' = l'Sigmal}~(2) (resp. from the identity $\xi_{ll',l'\gamma,v}/\xi_{ll',l'\alpha,u}=
(\xi_{l,\gamma,v}/\xi_{l,\alpha,u})^{l'}$). 
\end{proof}

\begin{thm}\label{thm:scheme normalization isom}
The following are true:
\begin{enumerate}
\item there exists an $S$-isomorphism
$\tpsi_{l,ll'}:(\tP_{l},\tcL_{l}^{\otimes l'})
\simeq (\tP_{ll'},\tcL_{ll'})$
such that $\tpsi_{l,ll'}(W_{l,\alpha,u})=W_{ll',l'\alpha,u}$ 
$(\forall \alpha\in\Sk^0(\Vor_{l}),\forall u\in X^{\vee})$;
\item $\tpsi_{l,ll'}$ induces an $S$-isomorphism 
$\psi_{l,ll'}:(P_{l},\cL_{l}^{\otimes l'})
\simeq (P_{ll'},\cL_{ll'}).$
\end{enumerate} 
\end{thm}

We denote $(P_{l},\cL_{l})$ by 
$(P_{l}(\xi^{\natural}),\cL_{l}(\xi^{\natural}))$  
if necessary. 
\begin{proof}First we define a homomorphism 
$\phi^*_{l,ll'}:R_{ll'}^*\to R_{l}^*$ by
$$\phi^*_{l,ll'}(\xi_{ll',l'\gamma,v}\theta_{ll'})
=(\xi_{l,\gamma,v}\theta_{l})^{l'}\ (\forall\gamma\in\Sk^0(\Sigma_{l}(0)), \forall v\in X^{\vee}).
$$ 

Note that both $R_{ll'}^*$ and $R_{l}^*$ are $R$-subalgebras of $R(\xi^{\natural})$ where $\xi_{ll',l'\gamma,v}\theta_{ll'}=\xi_{l,\gamma,v}^{l'}\theta^{ll'}=(\xi_{l,\gamma,v}\theta_{l})^{l'}$. This shows that 
$\phi^*_{l,ll'}$ is just the restriction of $\id_{R(\xi^{\natural})}$ to $R_{ll'}^*$. Hence $\phi^*_{l,ll'}$ is well-defined, 
which induces an $S$-morphism 
$$\phi_{l,ll'}:\tT_{l}\to \tT_{ll'}.$$

By Corollary~\ref{cor:normalization of tTl}, 
there exists an $S$-morphism 
$\tpsi_{l,ll'}:\tP_{l}\to \tP_{ll'}$ 
such that the following diagram is commutative:
\begin{diagram}
\tP_{l}&\rTo^{\tpsi_{l,ll'}}&\tP_{ll'}\\
\dTo^{h_{l}\pi_{l}}&&\dTo_{h_{ll'}\pi_{ll'}}\\
\tT_{l}&\rTo_{\phi_{l,ll'}}&\tT_{ll'}.
\end{diagram}

So $\phi_{l,ll'}$ induces 
an $S$-isomorphism $\tpsi_{l,ll',\alpha,u}:W_{l,\alpha,u}\to W_{ll',l'\alpha,u}$ by Corollary~\ref{cor:cones Alrho0 etc} 
and Lemma~\ref{lemma:cone Alrho0 and All'rho0}. 
 Hence 
$\tpsi_{l,ll'}:\tP_{l}\to \tP_{ll'}$ is an $S$-isomorphism.  
It is easy to see that $\tcL_{l}^{\otimes l'}
=\tpsi_{l,ll'}^*\tcL_{ll'}$. This proves (1). 

Since $S_y$ $(y\in Y)$ commutes with $\tpsi_{l,ll'}$, 
the formal quotients by $Y$ of 
$\tP_{l}^{\wedge}$ and $\tP_{ll'}^{\wedge}$ are isomorphic. 
By \cite[III${}_1$, 5.4.1]{EGA}, $\tpsi_{l,ll'}$ induces an $S$-isomorphism 
$\psi_{l,ll'}:(P_{l},\cL_{l}^{\otimes l'})
\simeq (P_{ll'},\cL_{ll'})$.
This proves (2). \end{proof}

\section{The structure of $\tP_{l}$}
\label{sec:str of Pl}

We keep the same notation as in \S~\ref{sec:rel complete models}.

\subsection{The $S$-scheme $W_{l,\alpha,u}$}
\label{subsec:Wlaplha0}

\begin{defn}\label{defn:C_l_beta}
 For $\beta\in\Sigma_{l}$, we define
\begin{align*}
C_{l}^{\beta}&=\text{the semi(sub)group of $X$ generated by 
$(\gamma-\beta;\gamma\in\Sigma_{l})$},\\
\Cone(C^{\beta}_{l})&:=\text{the cone over $\Conv(C_{l}^{\beta})$}
=\sum_{\gamma\in\Sigma_{l}}\bR_{\geq 0}(\gamma-\beta).
\end{align*} 
\end{defn}
\begin{defn}\label{defn:v_beta} 
For $\alpha\in\Sigma_{l}$,  we define 
$$\Sigma_{l}^{\alpha}=(\alpha+2l\mu(X^{\vee}))\cap\Sigma_{l}.$$
For each $\beta\in\Sigma_{l}^{\alpha}$, there exists a unique $v_{\beta}\in X^{\vee}$ such that 
$\beta=\alpha-2l\mu(v_{\beta})$. 
\end{defn}

\begin{lemma}
\label{lemma:integral Sigma}
Let $\alpha\in\Sigma_{l}$. Then 
$\Sigma_{l}^{\alpha}=\{\alpha\}\Leftrightarrow|\Sigma_{l}^{\alpha}|=1\Leftrightarrow\alpha\in\Sigma_{l}(0)^0$. 
\end{lemma}
\begin{proof}Every $\Sigma_{l}(c)$ $(c\in X^{\vee})$ 
is a closed convex polytope such that 
$$\Sigma_{l}(c)^0\cap \Sigma_{l}(d)=\emptyset\ (c\neq d),\ 
X_{\bR}=\bigcup_{d\in X^{\vee}}\Sigma_{l}(d).$$ 
Hence the boundary $\partial(\Sigma_{l}(0))$ of 
$\Sigma_{l}(0)$ is the union of 
$\Sigma_{l}(0)\cap\Sigma_{l}(d)$ $(d\neq 0)$. 
Now we prove Lemma. Since $\alpha\in\Sigma_{l}^{\alpha}$, it is clear that 
$\Sigma_{l}^{\alpha}=\{\alpha\}$ iff $|\Sigma_{l}^{\alpha}|=1$. 
Suppose $\alpha\in\Sigma_{l}\cap\partial\Sigma_{l}(0)$. Hence $\alpha\in\Sigma_{l}\cap(2l\mu(d)+\Sigma_{l})$ for some $d\in X^{\vee}\setminus\{0\}$. Therefore $\alpha, \alpha-2l\mu(d)\in\Sigma_{l}^{\alpha}$, that is, 
$|\Sigma_{l}^{\alpha}|\geq 2$. 
Conversely, if $|\Sigma_{l}^{\alpha}|\geq 2$, 
then there is $\beta\in\Sigma_{l}^{\alpha}$ with $\beta\neq\alpha$. 
Let $\beta=\alpha-2l\mu(d)$. Then $\alpha\in\Sigma_{l}\cap(2l\mu(d)+\Sigma_{l})\subset\partial\Sigma_{l}(0)$. This completes the proof. 
\end{proof}

\begin{lemma}\label{lemma:union of Cgamma_l}
 Let $\alpha\in\Sigma_{l}$. Then 
\begin{enumerate}
\item\label{item:union is XR}$X_{\bR}=\bigcup_{\beta\in\Sigma_{l}^{\alpha}} 
\Cone(C^{\beta}_{l})$;
\item\label{item:cone Clbeta} 
$\Cone(C^{\beta}_{l})$ is a $g$-dimensional cone with vertex 
$0$ $(\forall \beta\in\Sigma_{l}^{\alpha})$; 
\item\label{item:Clbeta  cap Cbeta'} $\Cone(C^{\beta}_{l})^0\cap\Cone(C^{\gamma}_{l})=\emptyset \ (\forall\beta,\gamma\in\Sigma_{l}^{\alpha};
\beta\neq\gamma)$.
\end{enumerate}
\end{lemma}
\begin{proof}
Recall  
$X_{\bR}=\bigcup_{c\in X^{\vee}} \Sigma_{l}(c)$ and $\Sigma_{l}(c)^0\cap \Sigma_{l}(d)=\emptyset$\ $(c\neq d)$. Since $\Sigma_{l}(0)=\Conv(\Sigma_{l})$ 
by assumption, 
$\beta+\Cone(C^{\beta}_{l})$ is a $g$-dimensional cone with vertex $\beta$
 which is spanned by $\Sigma_{l}(0)$ at $\beta$. 
Hence $X_{\bR}$ ($=$ the tangent space of $X_{\bR}$ at the origin) 
is the union of 
$\Cone(C^{\beta}_{l})$ $(\beta\in\Sigma_{l}^{\alpha})$. This proves (\ref{item:union is XR}). The rest is clear.  
\end{proof}

\begin{lemma}\label{lemma:Alalpha0=Adagger_lalpha0}
Let $\alpha,\gamma\in\Sigma_{l}$, $\beta\in\Sigma_{l}^{\alpha}$  
 and $u\in X^{\vee}$. Then 
\begin{enumerate}
\item\label{item:wt alpha} 
$\wt(\xi_{l,\beta,v_{\beta}})=\wt(\xi_{l,\alpha,0})=\alpha$;
\item\label{item:lgamma vbeta u/lbeta vbeta u in A l alpha u}
 $\xi_{l,\gamma,v_{\beta}+u}
/\xi_{l,\beta,v_{\beta}+u}\in A_{l,\alpha,u}$;
\item\label{item:wt lgamma vbeta u - lbeta vbeta u} $\wt(\xi_{l,\gamma,v_{\beta}+u}/\xi_{l,\beta,v_{\beta}+u})
=e(\xi)(v_{\beta}+u)(\gamma-\beta)m_0+\gamma-\beta$. 
\end{enumerate}
\end{lemma}
\begin{proof}
Since $\alpha=\beta+2l\mu(v_{\beta})$, 
by Lemma~\ref{lemma:comparison of val}, 
$\xi_{l,\beta,v_{\beta}}$ is a multiple of $\xi_{l,\alpha,0}$ 
by a unit of $R$, so that $\wt(\xi_{l,\beta,v_{\beta}})
=\wt(\xi_{l,\alpha,0})=\alpha$.  
This proves (\ref{item:wt alpha}). It also follows that 
$\xi_{l,\gamma,v_{\beta}}/\xi_{l,\beta,v_{\beta}}\in A_{l,\alpha,0}$. Let $\zeta:=\xi_{l,\gamma,v_{\beta}+u}/\xi_{l,\beta,v_{\beta}+u}$. Then 
$\zeta=\delta_{u}^*(\xi_{l,\gamma,v_{\beta}}/\xi_{l,\beta,v_{\beta}})\in A_{l,\alpha,u}$, which proves (\ref{item:lgamma vbeta u/lbeta vbeta u in A l alpha u}). Next, (\ref{item:wt lgamma vbeta u - lbeta vbeta u}) follows from  
$\wt(\zeta)
=\wt(\delta_{v_{\beta}+u}^*(\xi_{l,\gamma,0}/\xi_{l,\beta,0}))$. 
\end{proof}

\begin{lemma}\label{lemma:vbeta(gamma-beta) geq 0}
Let $\alpha,\gamma\in\Sigma_{l}$ and $\beta\in\Sigma_{l}^{\alpha}$. Then 
$v_{\beta}(\gamma-\beta)\geq 0$, with equality holding 
iff $\gamma\in(-\alpha+\beta)+\Sigma_{l}$.
\end{lemma}
\begin{proof}By Corollary~\ref{cor:Dlalpha=0 iff alpha in Sigma}, $D_{l}(\alpha-2l\mu(v_{\beta}))=D_{l}(\beta)=0$, so that   
\begin{align*}0\leq D_{l}(\gamma-\beta+\alpha)
&=D_{l}(\gamma+2l\mu(v_{\beta}))
+D_{l}(\alpha-2l\mu(v_{\beta}))\\
&=E_{l}(v_{\beta})+v_{\beta}(\gamma)
+E_{l}(-v_{\beta})-v_{\beta}(\alpha)\\
&=v_{\beta}(\gamma-\alpha+2l\mu(v_{\beta}))
=v_{\beta}(\gamma-\beta).
\end{align*} 
By Corollary~\ref{cor:Dlalpha=0 iff alpha in Sigma}, 
{\small $v_{\beta}(\gamma-\beta)=0\Leftrightarrow D_{l}(\gamma-\beta+\alpha)=0
\Leftrightarrow \gamma-\beta+\alpha\in\Sigma_{l}$.} 
\end{proof}

\begin{cor}\label{cor:(beta,beta')=(alpha,beta)}Let $\alpha\in\Sigma_{l}$, 
$\beta,\beta'\in\Sigma_{l}^{\alpha}$ and $x\in\Cone(C_{l}^{\beta})$. Then 
\begin{enumerate}
\item\label{item:vbeta geq vbeta'}$v_{\beta}(x)\geq v_{\beta'}(x)$, 
 with equality holding iff 
$x\in\Cone(C_{l}^{\beta'})$; 
\item\label{item:vbeta geq valpah=0} 
$v_{\beta}(x)\geq 0$, with equality holding iff 
$x\in\Cone(C_{l}^{\alpha})$.
\end{enumerate}
\end{cor}
\begin{proof}
We apply Lemma~\ref{lemma:vbeta(gamma-beta) geq 0} 
by choosing $\beta'$ for $\alpha$. Let $v^*:=v_{\beta}-v_{\beta'}$. 
Then $\beta'=\beta+2l\mu(v^*)$, so that  
we obtain $v^*(\gamma-\beta)\geq 0\ (\forall \gamma\in\Sigma_{l})$, 
with equality holding iff $\gamma\in(-\beta'+\beta)+\Sigma_{l}$. 
Let $x\in\Cone(C_{l}^{\beta})$. We write 
$x=\sum_{i=1}^rn_i(\gamma_i-\beta)$ for some $n_i>0$ and $\gamma_i\in\Sigma_{l}$.  Then it follows that $v^*(x)\geq 0$, with equality holding iff  
$\gamma_i\in(-\beta'+\beta)+\Sigma_{l}$ for any $i$, namely,  
 $\gamma_i-\beta=\delta_i-\beta'$ 
for some $\delta_i\in\Sigma_{l}$, so that $x=\sum_{i=1}^rn_i(\delta_i-\beta')\in \Cone(C_{l}^{\beta'})$. This proves (\ref{item:vbeta geq vbeta'}). 
 If $\beta'=\alpha$, then $v_{\beta'}=0$, so that $v_{\beta}(x)\geq 0$, with equality holding iff $x\in\Cone(C_{l}^{\alpha})$. 
This proves (\ref{item:vbeta geq valpah=0}).
\end{proof}

\begin{cor}\label{cor:C_beta}Let $\alpha\in\Sigma_{l}$,  
$\beta\in\Sigma_{l}^{\alpha}$ and $x\in X_{\bR}$. Then 
$x\in\Cone(C_{l}^{\beta})$ iff $v_{\beta}(x)\geq v_{\beta'}(x)$  
$(\forall \beta'\in\Sigma_{l}^{\alpha})$.
\end{cor}
\begin{defn}\label{defn:defn of Adagger_lalpha0}
Let $\alpha\in\Sigma_{l}$, $\beta\in\Sigma_{l}^{\alpha}$ 
and $u\in X^{\vee}$. 
We define 
\begin{gather*}\Delta_{l,\alpha,u}^{\beta}
=\left\{\text{$\wt(\xi_{l,\lambda,v_{\beta}+u}/\xi_{l,\beta,v_{\beta}+u});  
 \lambda\in\Sigma_{l}$}\right\},\ 
\Delta_{l,\alpha,u}=\bigcup_{\beta\in\Sigma_{l}^{\alpha}}
\Delta_{l,\alpha,u}^{\beta},\\
\Cone(\Delta_{l,\alpha,u})=\text{the cone over 
$\Conv(\Delta_{l,\alpha,u})$}.
\end{gather*}
\end{defn}

\begin{prop}\label{prop:isom Clbeta and Semi(Alalpha0^beta)}
Let $\alpha\in\Sigma_{l}$, $\beta\in\Sigma_{l}^{\alpha}$ 
and $u\in X^{\vee}$. We define a map  
$$h^{\beta}_{\alpha,u}:C_{l}^{\beta}\to \Semi(\Delta_{l,\alpha,u}^{\beta})$$ 
by $h^{\beta}_{\alpha,u}(x)=e(\xi)(v_{\beta}+u)(x)m_0+x$ 
$(x\in C_{l}^{\beta})$. Then 
\begin{enumerate}
\item\label{item:h_beta_alpha isom} 
$h^{\beta}_{\alpha,u}$ is an isomorphism of semigroups;
\item\label{item:cone B_l_alpha_u} $\Cone(B_{l,\alpha,u})=\bR_{\geq 0}m_0+\bigcup_{\beta\in\Sigma_{l}^{\alpha}}\Cone(h^{\beta}_{\alpha,u}(C_{l}^{\beta}))$; 
to be more precise, 
if $\tx=x_0m_0+x$ for some $x_0\in\bR$ and $x\in X_{\bR}$,  
 then 
$\tx\in\Cone(B_{l,\alpha,u})$ iff there exists $\beta\in\Sigma_{l}^{\alpha}$ 
with $x\in\Cone(C_{l}^{\beta})$ 
and $x_0\geq e(\xi)(v_{\beta}+u)(x)$;   
\item\label{item:x0 geq u(x)} if $x_0m_0+x\in\Cone(B_{l,\alpha,u})$ for some 
$x_0\in\bR$ and $x\in X_{\bR}$, then $x_0\geq e(\xi)u(x)$, 
with equality holding iff 
$x_0m_0+x\in\Cone(h^{\alpha}_{\alpha,u}(C_{l}^{\alpha}))$.
\end{enumerate}
\end{prop}
\begin{proof}(\ref{item:h_beta_alpha isom}) follows from 
Lemma~\ref{lemma:Alalpha0=Adagger_lalpha0}
~(\ref{item:wt lgamma vbeta u - lbeta vbeta u}). 
 Next, we shall prove 
(\ref{item:cone B_l_alpha_u}).
Let $$C:=\bR_{\geq 0}m_0+\bigcup_{\beta\in\Sigma_{l}^{\alpha}}\Cone(h^{\beta}_{\alpha,u}(C_{l}^{\beta})).$$ 

Let $\tx',\tx''\in C$. 
Then there are $\beta',\beta''\in\Sigma_{l}^{\alpha}$ such that 
$\tx'=x'_0m_0+x'\in C$, $\tx''=x_0''m_0+x''\in C$, 
$x'\in\Cone(C_{l}^{\beta'})$, 
$x''\in\Cone(C_{l}^{\beta''})$, $x'_0/e(\xi)\geq (v_{\beta'}+u)(x')$ and 
$x''_0/e(\xi)\geq (v_{\beta''}+u)(x'')$. There exists 
$\beta\in\Sigma_{l}^{\alpha}$ by 
Lemma~\ref{lemma:union of Cgamma_l}~(\ref{item:union is XR}) such that $x'+x''\in\Cone(C_{l}^{\beta})$. By Corollary~\ref{cor:(beta,beta')=(alpha,beta)}~(\ref{item:vbeta geq vbeta'}), 
$(x'_0+x''_0)/e(\xi)\geq (v_{\beta'}+u)(x')+(v_{\beta''}+u)(x'')
\geq (v_{\beta}+u)(x'+x'')$. Hence 
$\tx'+\tx''\in \bR_{\geq 0}m_0+
\Cone(h^{\beta}_{\alpha,u}(C_{l}^{\beta}))\subset C$. 
It follows that $C$ is convex.

Since $C$ is convex and 
$C\subset\Cone(A_{l,\alpha,u})=\Cone(B_{l,\alpha,u})$ 
by Lemma~\ref{lemma:Alalpha0=Adagger_lalpha0}
~(\ref{item:lgamma vbeta u/lbeta vbeta u in A l alpha u})-(\ref{item:wt lgamma vbeta u - lbeta vbeta u}), it suffices to prove that
$C\supset\Semi(A_{l,\alpha,u})$.  
Let $\tx\in\Semi(A_{l,\alpha,u})$. Then $\tx$ is of the form 
$\tx:=x_0m_0+x:=\wt(\xi_{l,\gamma,v+u}/\xi_{l,\alpha,u})$, where 
$\gamma\in\Sigma_{l}$, $v\in X^{\vee}$, $x=\gamma-\alpha+2l\mu(v)$ and 
\begin{gather*}
x_0/e(\xi)=E_{l}(v+u)+(v+u)(\gamma)-E_{l}(u)-u(\alpha).
\end{gather*}
We shall prove $\tx\in C$. 
By Lemma~\ref{lemma:union of Cgamma_l}~(\ref{item:union is XR}), 
$x\in\Cone(C_{l}^{\beta})$ 
for some $\beta\in\Sigma_{l}^{\alpha}$. Then it suffices to prove 
$x_0/e(\xi)\geq (v_{\beta}+u)(x)$.  
Let $v':=v-v_{\beta}$. Then $x=\gamma-\beta+2l\mu(v')$. 
By Corollary~\ref{cor:E(x)}~(2) and  
Corollary~\ref{cor:Dlalpha=0 iff alpha in Sigma}, 
\begin{equation}\label{eq:inequality x0 geq (vbeta+u)(x)}
\begin{aligned}
x_0/e(\xi)-(v_{\beta}+u)(x)&=E_{l}(v)+v(\gamma)-v_{\beta}(x)\\
&=E_{l}(v')+v'(\gamma)+E_{l}(v_{\beta})+v_{\beta}(\beta)\\
&=D_{l}(\gamma+2l\mu(v'))+D_{l}(\alpha)\geq 0.
\end{aligned}
\end{equation}
This shows $\tx\in C$. Hence $C\supset\Semi(A_{l,\alpha,u})$, which proves 
(\ref{item:cone B_l_alpha_u}). 

Next we prove (\ref{item:x0 geq u(x)}). Let $\tx:=x_0m_0+x\in\Cone(B_{l,\alpha,u})$ for some $x\in X_{\bR}$. By Lemma~\ref{lemma:union of Cgamma_l}~(1), 
there exists $\beta\in\Sigma_{l}^{\alpha}$ such that 
 $x\in\Cone(C_{l}^{\beta})$. Then 
$x_0/e(\xi)\geq (v_{\beta}+u)(x)\geq u(x)$ 
by Eq.~(\ref{eq:inequality x0 geq (vbeta+u)(x)}) and 
Corollary~\ref{cor:(beta,beta')=(alpha,beta)}~(\ref{item:vbeta geq valpah=0}). 
If $x_0=e(\xi)u(x)$, then $v_{\beta}(x)=0$, 
so that $x\in\Cone(C_{l}^{\alpha})$ by 
Corollary~\ref{cor:(beta,beta')=(alpha,beta)}~(\ref{item:vbeta geq valpah=0}). 
It follows $\tx\in\Cone(h^{\alpha}_{\alpha,u}(C_{l}^{\alpha}))$. 
Conversely, 
if $\tx\in\Cone(h^{\alpha}_{\alpha,u}(C_{l}^{\alpha}))$, then $x_0=e(\xi)(v_{\alpha}+u)(x)=e(\xi)u(x)$. This proves (\ref{item:x0 geq u(x)}).   
\end{proof}

\begin{lemma}
\label{lemma:cone taulalpha0}Let $\alpha\in\Sigma_{l}$ and $u\in X^{\vee}$. 
 Then 
\begin{enumerate}
\item\label{item:cone tau_l_alpha_u}
$\tau_{l,\alpha,u}=\Cone(B_{l,\alpha,u})^{\vee}
=(\bR_{\geq 0}f_0+X^{\vee}_{\bR})\cap
\Cone(\Delta_{l,\alpha,u})^{\vee}$;
\item\label{item:case Sigma_0_0} 
$\tau_{l,\alpha,u}=\bR_{\geq 0}(f_0-e(\xi)u)$ 
if $\alpha\in\Sigma_{l}(0)^0$;
\item\label{item:case boundary Sigma_0}  
$\tau_{l,\alpha,u}=\Cone(\Delta_{l,\alpha,u})^{\vee}$ and $m_0\in\Cone(\Delta_{l,\alpha,u})$
if $\alpha\in\partial\Sigma_{l}(0)$.
\end{enumerate}
\end{lemma}
\begin{proof}
Since $\Cone(\Delta_{l,\alpha,u}^{\beta}) \subset\Cone(\Delta_{l,\alpha,u})\subset \Cone(B_{l,\alpha,u})$ $(\forall \beta\in\Sigma_{l}^{\alpha})$,  by  Proposition~\ref{prop:isom Clbeta and Semi(Alalpha0^beta)}~(\ref{item:cone B_l_alpha_u}), we see 
$\tau_{l,\alpha,u}^{\vee}=\Cone(B_{l,\alpha,u})=\bR_{\geq 0}m_0+\Cone(\Delta_{l,\alpha,u})$.   
This proves (\ref{item:cone tau_l_alpha_u}).   
By Lemma~\ref{lemma:integral Sigma}, if $\alpha\in\Sigma_{l}(0)^0$, then 
$\Sigma^{\alpha}_{l}=\{\alpha\}$, 
$v_{\alpha}=0$ and 
$\Cone(C_{l}^{\alpha})=X_{\bR}$. 
So by Proposition~\ref{prop:isom Clbeta and Semi(Alalpha0^beta)}
~(\ref{item:h_beta_alpha isom}),   
\begin{align*}
\Cone(\Delta_{l,\alpha,u})
&=\Cone(h^{\alpha}_{\alpha,u}(C_{l}^{\alpha}))=\left\{e(\xi)u(x)m_0+x;
x\in X_{\bR}\right\},\\
\Cone(\Delta_{l,\alpha,u})^{\vee}&=
{\left\{
\begin{matrix}
z_0f_0+z; z_0\in\bR, z\in X^{\vee}_{\bR} \\
(e(\xi)z_0u+z)(x)\geq 0\ (\forall x\in X_{\bR})
\end{matrix}
\right\}=\bR(f_0-e(\xi)u).}
\end{align*}  Hence $\tau_{l,\alpha,u}=\bR_{\geq 0}(f_0-e(\xi)u)$ 
by (\ref{item:cone tau_l_alpha_u}).  
This proves (\ref{item:case Sigma_0_0}). \par
Next, we shall prove (\ref{item:case boundary Sigma_0}). By assumption, 
$\alpha\in\partial\Sigma_{l}(0)$. 
By Lemma~\ref{lemma:union of Cgamma_l}, there exist a nonempty subset $J$ of $\Sigma_{l}^{\alpha}$ and $(x_{\beta}\in C_{l}^{\beta}\cap\Cone(C_{l}^{\beta})^0;\beta\in J)$ with $\sum_{\beta\in J} x_{\beta}=0$.  Hence 
$$\sum_{\beta\in J}h_{\alpha,u}^{\beta}(x_{\beta})
=e(\xi)\sum_{\beta\in J}(v_{\beta}+u)(x_{\beta})m_0+\sum_{\beta\in J}x_{\beta}
=e(\xi)\sum_{\beta\in J}v_{\beta}(x_{\beta})m_0.$$   
By Definition~\ref{defn:defn of Adagger_lalpha0},  
$\sum_{\beta\in J}h_{\alpha,u}^{\beta}(x_{\beta})\in
\Cone(\Delta_{l,\alpha,u})$. 
By Corollary~\ref{cor:(beta,beta')=(alpha,beta)}
~(\ref{item:vbeta geq valpah=0}), 
$\sum_{\beta\in J}v_{\beta}(x_{\beta})>0$, hence 
$m_0\in\Cone(\Delta_{l,\alpha,u})$,  so that (\ref{item:case boundary Sigma_0}) follows from (\ref{item:cone tau_l_alpha_u}). 
\end{proof}

\begin{thm}
\label{thm:W_l,alpha,0} 
Let $\alpha\in\Sigma_{l}$ and $u\in X^{\vee}$. 
Then 
\begin{enumerate}  
\item\label{item:Sigma(0)0} if $\alpha\in\Sigma_{l}(0)^0$, then 
$W_{l,\alpha,u}=W_{l,0,u}=\Spec R[s^{e(\xi)u(x)}w^x;x\in X]$;
\item\label{item:boundary_Sigma(0)0} if $\alpha\in\partial\Sigma_{l}(0)$, then every irreducible component of $(W_{l,\alpha,u})_0$ is 
the closure of $(W_{l,0,v_{\beta}+u})_0$ in $(W_{l,\alpha,u})_0$ 
for some $\beta\in\Sigma_{l}^{\alpha}$, and vice versa  
where $(W_{l,0,v_{\beta}+u})_0=\delta_{-v_{\beta}}((W_{l,0,u})_0)$ and $\beta=\alpha-2l\mu(v_{\beta})$;
\item\label{item:boundary_Sigma(0)0 ring str} if $\alpha\in\partial\Sigma_{l}(0)$, then 
\begin{align*}
B_{l,\alpha,u}\otimes_R k(0)&=
k(0)[s^{e(\xi)(v_{\beta}+u)(x)}w^{x}; 
\beta\in\Sigma_{l}^{\alpha}, 
x\in\Cone(C_{l}^{\beta})\cap X]
\end{align*}
with the fundamental relations in $B_{l,\alpha,0}\otimes_R k(0)$ 
given by
\begin{align*}
&s^{e(\xi)(v_{\beta'}+u)(x)}w^{x}\cdot s^{e(\xi)(v_{\beta''}+u)(y)}w^{y}\\
&=\begin{cases}s^{e(\xi)(v_{\beta}+u)(x+y)}w^{x+y}
&\text{if 
$x,y\in\Cone(C_{l}^{\beta})\cap X$\ 
$(\exists \beta\in\Sigma_{l}^{\alpha})$}\\
0&\text{otherwise}
\end{cases}
\end{align*}
where $\beta',\beta''\in\Sigma_{l}^{\alpha}$, $x\in\Cone(C_{l}^{\beta'})\cap X$ and $y\in\Cone(C_{l}^{\beta''})\cap X$. 
\end{enumerate} 
\end{thm}
\begin{proof} Let $\alpha\in\Sigma_{l}(0)^0$. Then 
$B_{l,\alpha,u}=R[\tau_{l,\alpha,u}^{\vee}\cap\tX]
=R[s^{e(\xi)u(x)}w^x;x\in X]$ by Lemma~\ref{lemma:cone taulalpha0}~(2), hence $W_{l,\alpha,u}=W_{l,0,u}$.  This proves (1). 
Next, we shall prove (\ref{item:boundary_Sigma(0)0 ring str}). Since $\alpha\in\partial\Sigma_{l}(0)$, $|\Sigma_{l}^{\alpha}|\geq 2$ by Lemma~\ref{lemma:integral Sigma}.
By Proposition~\ref{prop:isom Clbeta and Semi(Alalpha0^beta)}~(\ref{item:cone B_l_alpha_u}),
\begin{equation}
\label{eq:Bl_alpha_u explicit}
B_{l,\alpha,u}
=R[s^{e(\xi)(v_{\beta}+u)(x)}w^x; \beta\in\Sigma_{l}^{\alpha}, x\in \Cone(C_{l}^{\beta})\cap X].
\end{equation}

Let $\beta',\beta''\in\Sigma_{l}^{\alpha}$, $x\in\Cone(C_{l}^{\beta'})\cap X$ and $y\in\Cone(C_{l}^{\beta''})\cap X$.  Suppose $x+y\in\Cone(C_{l}^{\beta})\cap X$ for some $\beta\in\Sigma_{l}^{\alpha}$. By Corollary~\ref{cor:(beta,beta')=(alpha,beta)}~(\ref{item:vbeta geq vbeta'}), $v_{\beta'}(x)\geq v_{\beta}(x)$ and  
$v_{\beta''}(y)\geq v_{\beta}(y)$, so that 
$v_{\beta'}(x)+v_{\beta''}(y)\geq v_{\beta}(x+y)$. 
It follows that  
$$s^{e(\xi)(v_{\beta'}+u)(x)}w^{x}\cdot s^{e(\xi)(v_{\beta''}+u)(y)}w^{y}
=s^as^{e(\xi)(v_{\beta}+u)(x+y)}w^{x+y}
$$for some integer $a\geq 0$.  If $a=0$,  then 
$v_{\beta'}(x)=v_{\beta}(x)$ and  
$v_{\beta''}(y)=v_{\beta}(y)$, so that 
$x, y\in\Cone(C_{l}^{\beta})\cap X$ by Corollary~\ref{cor:(beta,beta')=(alpha,beta)}~(\ref{item:vbeta geq vbeta'}).  This proves (\ref{item:boundary_Sigma(0)0 ring str}). 

Finally we shall prove (\ref{item:boundary_Sigma(0)0}). 
By (1), Eq.~(\ref{eq:Bl_alpha_u explicit})  and Corollary~\ref{cor:(beta,beta')=(alpha,beta)}~(\ref{item:vbeta geq valpah=0}), \footnote{See \cite[lines 3-5, p.~24]{TE73}.}
{\small\begin{align*}
W_{l,\alpha,u}\cap W_{l,0,u}&=
\Spec R\left[
s^{e(\xi)(v_{\beta}+u)(x)}w^{x}, s^{e(\xi)u(y)}w^y; \begin{matrix}
\beta\in\Sigma_{l}^{\alpha}, y\in X\\
x\in\Cone(C_{l}^{\beta})\cap X
\end{matrix}
\right]\\
&=
\Spec R[s^{e(\xi)u(y)}w^y;  y\in X]=W_{l,0,u}\simeq T_{X,R}.
\end{align*}}

\vskip -0.5cm  It follows that  
$(W_{l,\alpha,u})_0\cap (W_{l,0,u})_0\simeq (W_{l,0,u})_0$, which is 
an irreducible open subscheme of $(W_{l,\alpha,u})_0$. 
Let $\beta\in\Sigma_{l}^{\alpha}\setminus\{\alpha\}$. Similarly since  
$(W_{l,0,u})_0$ is an irreducible open subscheme of $(W_{l,\beta,u})_0$, 
 $(W_{l,0,v_{\beta}+u})_0=\delta_{-v_{\beta}}((W_{l,0,u})_0)$ is an irreducible open subscheme of 
$\delta_{-v_{\beta}}((W_{l,\beta,u})_0)=(W_{l,\alpha,u})_0$. 
By (3), the nonsingular part of 
$(W_{l,\alpha,u})_0$ is the union of $((W_{l,0,v_{\beta}+u})_0;\beta\in\Sigma_{l}^{\alpha})$. 
This proves (\ref{item:boundary_Sigma(0)0}), which completes the proof.
\end{proof}

\subsection{$\Vor_{l}(\xi^{\natural})$}
\label{subsec:Vor_l_xinatural}

\begin{defn}\label{defn:Voronoi decomp}
Let  
$\Sigma_{l}(c)\ (c\in X^{\vee})$ 
be the Voronoi polytopes in \S~\ref{subsec:Voronoi cells}. 
We define the Voronoi decomposition $\Vor_{l}:=\Vor_{l,B}$ of $X_{\bR}$ by 
$$\Vor_{l,B}:=\left\{\Sigma_{l}(c)\ (c\in X^{\vee})\ \text{and their faces}\right\},  
$$which we denote by 
$\Vor_{l}(\xi^{\natural})$ if necessary. 
  $\Vor_{l}$ is a 
polyhedral decomposition of $X_{\bR}$ invariant under translations 
by $2l\mu(X^{\vee})$, hence 
invariant by $2l\mu\beta(Y)=2Nl Y$. 
 Let $\Vor_{l}/2Nl Y$ be the quotient 
of $\Vor_{l}$ by $2Nl Y$.  
\end{defn}

\begin{defn}
\label{defn:Skq}
Let $\Sk^q(\Vor_{l})$ be the set of 
all $q$-dimensional polytopes in 
$\Vor_{l}$. Then $\Vor_{l}=\coprod_{q=0}^g\Sk^q(\Vor_{l})$. 
Let $\rho,\rho'\in\Vor_{l}$. Then $\rho\cap\rho'\in\Vor_{l}$, while 
$\rho\subset\rho'$ iff $\rho$ is a face of $\rho'$. 
Let $\Sk(\rho)$ (resp. $\Sk^q(\rho)$) 
be the set of all faces of $\rho$ (resp. all $q$-dimensional 
faces of $\rho$). 
By the assumption (\ref{assump:Sigmal(0) integral}), 
$\Sk^0(\rho)\subset\rho\cap X$ 
and $\rho=\Conv(\Sk^0(\rho))$. We also define 
$$\overline{\Vor}_{l}=\Vor_{l}/2Nl Y,\ \Sk^q(\overline{\Vor}_{l}):=\Sk^q(\Vor_{l})/2Nl Y.$$ 
\end{defn}

\begin{rem}
\label{rem:vo_lll'}
If $\Sigma_{l}(0)$ is integral, then so is $\Sigma_{ll'}(0)$. 
So we consider $\Vor_{ll'}$ and $\Sk^q(\Vor_{ll'})$ 
as well for any $l'\geq 1$.
Note that 
\begin{align*}
\Sk^q(\Vor_{ll'})&=\{l'\sigma; \sigma\in \Sk^q(\Vor_{l})\}.
\end{align*}

Let $\vor_{l,ll'}:X_{\bR}\to X_{\bR}$ be multiplication by $l'$. 
Then  
\begin{gather*}\vor_{l,ll'}\Sk^q(\Vor_{l})=\Sk^q(\Vor_{ll'}),\\ 
\vor_{l,ll'}(\sigma+2l\mu(u))=\vor_{l,ll'}(\sigma)+2ll'\mu(u)\ (\forall u\in X^{\vee}).
\end{gather*}

 Hence the induced map is bijective:
$$\overline{\vor}_{l,ll'}:\overline{\Vor}_{l}\to \overline{\Vor}_{ll'}.
$$
\end{rem}

\begin{lemma}
\label{lemma:well defined of tau_lrho}
Let $\rho',\rho''\in\Sk^q(\Sigma_{l}(0))$ and $u',u''\in X^{\vee}$ such that 
$\rho'+2l\mu(u')=\rho''+2l\mu(u'')$. Then 
$$A_{l,\rho',u'}=A_{l,\rho'',u''},\ B_{l,\rho',u'}=B_{l,\rho'',u''}.
$$
\end{lemma}
 \begin{proof} It suffices to prove 
$A_{l,\rho',u'}=A_{l,\rho'',u''}$. 
Let $\beta'\in\rho'\cap\Sigma_{l}$ and $\beta'':=\beta'+2l\mu(u'-u'')$. 
Then $\beta''\in\rho''\cap X\subset\Sigma_{l}$, and 
by Lemma~\ref{lemma:comparison of val}, there exists a unit 
$c_{\beta',\beta''}\in R$ such that $\xi_{l,\beta',u'}=c_{\beta',\beta''}\xi_{l,\beta'',u''}$.  
By Definition~\ref{defn:Alrho0},   
\begin{align*}
A_{l,\rho',u'}&=R[\xi_{l,\gamma,v}/\xi_{l,\beta',u'}; \beta'\in\rho'\cap\Sigma_{l}, \gamma\in\Sigma_{l},v\in X^{\vee}]\\
&=R[\xi_{l,\gamma,v}/\xi_{l,\beta'',u''}; \beta''\in\rho''\cap\Sigma_{l}, \gamma\in\Sigma_{l},v\in X^{\vee}]=A_{l,\rho'',u''}.
\end{align*}

 This completes the proof. 
\end{proof}

\subsection{$\Fan(\xi^{\natural}_{l})$}
\label{subsec:Fan_l_xinatural}

\begin{defn}\label{defn:Bla,BlDelta}
If $a=\alpha+2l\mu(u)$ for $\alpha\in\Sigma_{l}$ and 
$u\in X^{\vee}$, then we define 
\begin{gather*}
\zeta_{l,a}=t^{D_{l}(a)}w^a,\ \tau_{l,a}=\tau_{l,\alpha,u},\ 
A_{l,a}=A_{l,\alpha,u},\ B_{l,a}=B_{l,\alpha,u},\\
 U_{l,a}=\Spec A_{l,a},\ W_{l,a}=\Spec B_{l,a}.
\end{gather*}
 
These are well-defined by Definition~\ref{defn:Dl(x)} and 
Lemma~\ref{lemma:well defined of tau_lrho}. 
For any $\Delta\in\Vor_{l}$, there exists $\rho\in\Sk(\Sigma_{l}(0))$ and $u\in X^{\vee}$ such that $\Delta=\rho+2l\mu(u)$. Then we define 
\begin{gather*}
\tau_{l,\Delta}=\bigcap_{\alpha\in\Sk^0(\rho)}\tau_{l,\alpha,u}=
\bigcap_{a\in\Sk^0(\Delta)}\tau_{l,a},\\ 
U_{l,\Delta}=\bigcap_{\alpha\in\Sk^0(\rho)}U_{l,\alpha,u}
=\bigcap_{a\in\Sk^0(\Delta)}U_{l,a},\\
W_{l,\Delta}=\bigcap_{\alpha\in\Sk^0(\rho)}W_{l,\alpha,u}
=\bigcap_{a\in\Sk^0(\Delta)}W_{l,a},\\
A_{l,\Delta}=\Gamma(U_{l,\Delta}, \cO_{U_{l,\Delta}})=A_{l,\rho,u},\ 
B_{l,\Delta}=\Gamma(W_{l,\Delta}, \cO_{W_{l,\Delta}})=B_{l,\rho,u}.
\end{gather*}

Note that 
$B_{l,\Delta}=R[\tau_{l,\Delta}^{\vee}\cap\tX]$, which is generated 
by $(B_{l,a}; a\in\Sk^0(\Delta))$, so that 
$\Cone(B_{l,\Delta})$ is generated by $(\Cone(B_{l,a}); a\in\Sk^0(\Delta))$.
\end{defn}

\begin{defn}
\label{defn:CDelta XDelta}Let $u\in X^{\vee}$ and  
$\Delta\in\Sk(\Sigma_{l}(u))$.  Then 
we define 
\begin{align*}
C_u(\Delta)&=\text{the cone over 
$\Conv(d-a; d\in\Sigma_{l}(u)\cap X, a\in\Sk^0(\Delta))$},\\
X(\Delta)&=\begin{matrix}
\text{the saturation in $X$ of the $\bZ$-submodule of $X$}\\    
\text{generated by $\Delta-\Delta$}
\end{matrix}
\end{align*}where  
$C_u(\emptyset)=X(\emptyset):=\emptyset$.  
Note that $C_u(\Delta)=C_0(\Delta-2l\mu(u))$. 
\end{defn}

\begin{cor}\label{cor:u_geq u(x) iff tx CuDelta}Let 
$u\in X^{\vee}$,  
$\Delta\in\Sk(\Sigma_{l}(u))$ and 
$\tx:=x_0m_0+x\in\Cone(B_{l,\Delta})$ for some $x_0\in\bR$ 
and $x\in X_{\bR}$.  
Then $x_0\geq e(\xi)u(x)$, 
with equality holding iff $\tx\in f(C_{u}(\Delta))$ where  
$f:=e(\xi)um_0+\id_{X_{\bR}}$.
\end{cor}
\begin{proof}  
Since $\Cone(B_{l,\Delta})=\Cone(B_{ll',l'\Delta})$, 
by taking $ll'$ for $l$ if necessary,  
we may assume that there exists $c\in\Delta^0\cap X$. 
Then we write $c=\alpha+2l\mu(u)$ 
for some $\alpha\in\Sigma_{l}$.   
Then $\Cone(B_{l,\Delta})=\Cone(B_{l,\alpha,u})$ 
by Lemma~\ref{lemma:Alrho0 as intersection}. 

Next, we prove  
$\Cone(C_{l}^{\alpha})=C_u(\Delta)$. 
Let $a\in\Sk^0(\Delta)$. 
Since $c\in\Delta^0$, 
there exist $n_i\in\bN$ and  
$a_i\in \Sk^0(\Delta)$\ $(i\in[1,n])$ for some $n\geq 2$ 
such that $a_1=a$ and 
$\sum_{i=1}^nn_i(c-a_i)=0$. Hence 
$a-c=(1/n_1)\sum_{i=2}^nn_i(c-a_i)\in C_u(\Delta)$ 
because $c-a_i\in C_u(\Delta)$. 
Since $b-c=(b-a)+(a-c)\in C_u(\Delta)$ 
$(\forall b\in\Sigma_{l}(u)\cap X)$, 
we have $\Cone(C_{l}^{\alpha})\subset C_u(\Delta)$.  
Meanwhile $c-a=(1/n_1)\sum_{i=2}^nn_i(a_i-c)\in\Cone(C_{l}^{\alpha})$.  
Hence $b-a=(b-c)+(c-a)\in\Cone(C_{l}^{\alpha})$  
$(\forall b\in\Sigma_{l}(u)\cap X)$, so that $C_u(\Delta)\subset\Cone(C_{l}^{\alpha})$. It follows $\Cone(C_{l}^{\alpha})=C_u(\Delta)$. 
 
Let $\tx=x_0m_0+x\in\Cone(B_{l,\Delta})=\Cone(B_{l,\alpha,u})$. 
By Proposition~\ref{prop:isom Clbeta and Semi(Alalpha0^beta)}~
(\ref{item:x0 geq u(x)}), $x_0\geq e(\xi)u(x)$, 
with equality holding iff 
$\tx=f(x)$ and $x\in\Cone(C_{l}^{\alpha})=C_u(\Delta)$. 
This completes the proof.
\end{proof}
\begin{example}\label{example:Bl_Sigma(u)}
Let $\Delta:=\Sigma_{l}(u)$ $(u\in X^{\vee})$. 
Then $B_{l,\Delta}=B_{l,0,u}$. By Lemma~\ref{lemma:cone taulalpha0}, 
$\Cone(B_{l,0,u})=\tau_{l,0,u}^{\vee}$ and  
\begin{gather*}
B_{l,\Delta}=R[s^{e(\xi)u(x)}w^x; x\in X],\ B_{l,\Delta}\otimes_R k(0)
=k(0)[s^{e(\xi)u(x)}w^x; x\in X],\\     
\Cone(B_{l,\Delta})=
\{x_0m_0+x\in\tX;x_0\geq e(\xi)u(x), x\in C_u(\Delta)\cap X\}
\end{gather*}by $C_u(\Delta)\cap X=X$, which is compatible with Corollary~\ref{cor:u_geq u(x) iff tx CuDelta}.
\end{example}

\begin{lemma}\label{lemma:Cone BlDelta_cap X=-Delta+lin(Delta)}
Let $u\in X^{\vee}$, $\Delta\in\Sk(\Sigma_{l}(u))$ and 
$f:=e(\xi)um_0+\id_{X_{\bR}}$. Then 
the following are true:
\begin{enumerate}
\item\label{item:a+XDelta cuts out Delta}
 $(a+X(\Delta)_{\bR})\cap\Sigma_{l}(u)=\Delta$ for any $a\in\Delta$;
\item\label{item:cone C(Delta)} 
$\Cone(B_{l,\Delta})\cap f(X_{\bR})=f(C_u(\Delta))
=\bigcup_{\alpha\in\Sk^0(\rho)} 
h^{\alpha}_{\alpha,u}(C_{l}^{\alpha})$ 
where $\rho:=\Delta-2l\mu(u)\in\Sk(\Sigma_{l}(0))$; 
\item\label{item:max linear X(Delta)} $f(X(\Delta)_{\bR})$ is 
the maximal $\bR$-linear subspace of 
$\Cone(B_{l,\Delta})$;
\item\label{item: chara of DeltacapX} $a\in\Delta\cap X$ iff  
$a\in\Sigma_{l}(u)\cap X$ and  
$f(\Sigma_{l}(u))\subset f(a)+\Cone(B_{l,\Delta})$.
\end{enumerate} 
\end{lemma}

\begin{proof}
Let $a\in\Delta$. Then $a+X(\Delta)_{\bR}$ is a linear subspace 
of $X_{\bR}$ containing $\Delta$, hence 
$(a+X(\Delta)_{\bR})\cap\Sigma_{l}(u)=\Delta$ 
because $\Sigma_{l}(u)$ is convex. This proves 
(\ref{item:a+XDelta cuts out Delta}).  
Next, we prove (\ref{item:cone C(Delta)}). 
If $\tx\in\Cone(B_{l,\Delta})\cap f(X_{\bR})$, then 
$\tx=e(\xi)u(x)m_0+x$ for some $x\in X_{\bR}$, whence 
$x\in C_u(\Delta)$ by Corollary~\ref{cor:u_geq u(x) iff tx CuDelta}.

Next we prove (\ref{item:max linear X(Delta)}). 
Since $\Cone(B_{l,\Delta})=\Cone(B_{ll',l'\Delta})$, 
by taking $ll'$ for $l$ if necessary, 
we may assume that there exists $a\in\Delta^0\cap X$. 
Then we write $a=\alpha+2l\mu(u)$ 
for some $\alpha\in\Sigma_{l}$.   
Let $L$ be a linear subspace of  
$\Cone(B_{l,\Delta})$ and $\tx\in L$. 
Since $\Cone(B_{l,\Delta})=\Cone(B_{l,\alpha,u})$, 
we have $\pm\tx\in \Cone(B_{l,\alpha,u})$. 
By Proposition~\ref{prop:isom Clbeta and Semi(Alalpha0^beta)}
~(\ref{item:cone B_l_alpha_u}), there exists 
$\beta\in\Sigma_{l}^{\alpha}$ and $x\in\Cone(C_{l}^{\beta})$ 
such that $\tx=x_0m_0+x$ and $x_0\geq e(\xi)(v_{\beta}+u)(x)$. 
There also exists $\beta'\in\Sigma_{l}^{\alpha}$ 
such that $-x\in\Cone(C_{l}^{\beta'})$ 
and $-x_0\geq e(\xi)(v_{\beta'}+u)(-x)$. 
 Hence $v_{\beta'}(x)\geq v_{\beta}(x)$, 
while $v_{\beta}(x)\geq v_{\beta'}(x)$ 
by Corollary~\ref{cor:(beta,beta')=(alpha,beta)}
~(\ref{item:vbeta geq vbeta'}). 
It follows that  $v_{\beta}(x)=v_{\beta'}(x)$,  
$x_0=e(\xi)(v_{\beta}+u)(x)$ and $\pm x\in\Cone(C_{l}^{\beta})\cap \Cone(C_{l}^{\beta'})$. By Corollary~\ref{cor:(beta,beta')=(alpha,beta)}~(\ref{item:vbeta geq valpah=0}), $v_{\beta'}(\pm x)\geq 0$, 
whence $v_{\beta}(x)=v_{\beta'}(x)=0$ 
and $\tx=e(\xi)u(x)m_0+x=f(x)$. It follows that 
$L\subset f(\Cone(C_{l}^{\beta}))
\subset\Cone(B_{l,\Delta})\cap f(X_{\bR})=f(C_u(\Delta))$. 
Since $X(\Delta)_{\bR}$ is 
the maximal linear subspace of $C_u(\Delta)$, 
$f(X(\Delta)_{\bR})$ is the maximal linear subspace of 
$f(C_u(\Delta))$. Hence $L\subset f(X(\Delta)_{\bR})$. 
This proves (\ref{item:max linear X(Delta)}).

Finally we shall prove (\ref{item: chara of DeltacapX}). 
Suppose $a\in\Delta\cap X$. 
Let $b\in\Sigma_{l}(u)\cap X$ and $\alpha:=a-2l\mu(u)$. 
Since $b-a\in C_{l}^{\alpha}$ and 
$f=h^{\alpha}_{\alpha,u}$ on $C_{l}^{\alpha}$, 
by Proposition~\ref{prop:isom Clbeta and Semi(Alalpha0^beta)}
~(\ref{item:cone B_l_alpha_u}),
$$f(b-a)=e(\xi)u(b-a)m_0+b-a\in\Cone(B_{l,a})\subset\Cone(B_{l,\Delta}).$$ 
It follows $f(\Sigma_{l}(u)\cap X)\subset f(a)+\Cone(B_{l,\Delta})$. 
This proves the only-if-part of (\ref{item: chara of DeltacapX}). 
Next, we shall prove the if-part of (\ref{item: chara of DeltacapX}). 
Suppose that  $a\in\Sigma_{l}(u)\cap X$ and  
$f(\Sigma_{l}(u))\subset f(a)+\Cone(B_{l,\Delta})$.  
Let $b\in\Sk^0(\Delta)$. 
Then $f(b-a)\in\Cone(B_{l,\Delta})$, while 
$f(a-b)\in\Cone(B_{l,\Delta})$ by the only-if-part of (\ref{item: chara of DeltacapX}). 
It follows that $\pm f(a-b)\in\Cone(B_{l,\Delta})$. 
Hence $\bR f(a-b)$ is a linear subspace contained in $\Cone(B_{l,\Delta})$. 
By (\ref{item:max linear X(Delta)}),
$f(a-b)\in f(X(\Delta)_{\bR})$, so that $a-b\in X(\Delta)_{\bR}$   
because $f$ is injective. 
It follows from (\ref{item:a+XDelta cuts out Delta}) 
that $a\in\Sigma_{l}(u)\cap(b+X(\Delta)_{\bR})=\Delta$. 
This  proves the if-part of (\ref{item: chara of DeltacapX}). 
\end{proof}

\begin{lemma}\label{lemma:tauDelta=BlDeltavee}Let $a,b\in X$ and 
$\Delta,\Delta'\in\Vor_{l}$. Then
\begin{enumerate}
\item\label{item:tau_l_Delta}$\tau_{l,a}=\Cone(B_{l,a})^{\vee}$ and  
$\tau_{l,\Delta}=\Cone(B_{l,\Delta})^{\vee}$;
\item\label{item:dim_X_Delta} $\dim\tau_{l,\Delta}=g+1-\dim_{\bR}\Delta$; 
\item\label{item:tau_l,a=tau_l_Delta} $\tau_{l,b}=\tau_{l,\Delta}$\ 
if $b\in\Delta^0\cap X$;
\item\label{item:tau_Delta tau_Delta'} $\tau_{l,\Delta}\supset\tau_{l,\Delta'}\Leftrightarrow 
\Delta\subset\Delta'$. 
\end{enumerate}
\end{lemma}
\begin{proof}The first half of (\ref{item:tau_l_Delta}) follows from 
$\tau_{l,\alpha,u}=\Cone(B_{l,\alpha,u})^{\vee}$ $(\forall\alpha\in\Sigma_{l},\forall u\in X^{\vee})$. 
Since $\Cone(B_{l,\Delta})$ is generated by 
$(\Cone(B_{l,a}); a\in\Sk^0(\Delta))$,   
\begin{align*}
\Cone(B_{l,\Delta})^{\vee}&=\bigcap_{a\in\Sk^0(\Delta)}\Cone(B_{l,a})^{\vee}=\bigcap_{a\in\Sk^0(\Delta)}\tau_{l,a}=\tau_{l,\Delta},
\end{align*}which proves (\ref{item:tau_l_Delta}). In general, for a cone $C$ in $\tX_{\bR}$, 
$$\dim C^{\vee}=g+1-\max_{
\begin{matrix}L\subset C\\
L:\text{$\bR$-lin. subsp.}
\end{matrix}}\dim L.$$
Hence {\small$\dim\tau_{l,\Delta}=\dim\Cone(B_{l,\Delta})^{\vee}=g+1-\dim_{\bR}f(X(\Delta)_{\bR})$} by Lemma~\ref{lemma:Cone BlDelta_cap X=-Delta+lin(Delta)}~(\ref{item:max linear X(Delta)}), while 
$\dim_{\bR}f(X(\Delta)_{\bR})=\dim_{\bR}\Delta$. 
This proves (\ref{item:dim_X_Delta}). \par

Next, we prove (\ref{item:tau_l,a=tau_l_Delta}). 
Let $b\in\Delta^0\cap X$. 
Then there exist a face $\rho$ of $\Sigma_{l}(0)$, 
$\alpha\in\Sigma_{l}$ and $u\in X^{\vee}$ such that 
$\Delta=\rho+2l\mu(u)$,  
$b=\alpha+2l\mu(u)$ and $\alpha\in\rho^0\cap X$. 
By Lemma~\ref{lemma:Alrho0 as intersection}, 
$B_{l,b}=B_{l,\alpha,u}=B_{l,\rho,u}=B_{l,\Delta}$. Hence 
$\tau_{l,b}=\Cone(B_{l,b})^{\vee}=\Cone(B_{l,\Delta})^{\vee}=\tau_{l,\Delta}$ by  (\ref{item:tau_l_Delta}).  
This proves (\ref{item:tau_l,a=tau_l_Delta}).

Finally we prove (\ref{item:tau_Delta tau_Delta'}). 
It is clear  
that $\tau_{l,\Delta'}\subset\tau_{l,\Delta}$ 
if $\Delta'\supset\Delta$. 
Now assume $\tau_{l,\Delta'}\subset\tau_{l,\Delta}$.  
We choose $u\in X^{\vee}$ such that  
$\Delta'\subset\Sigma_{l}(u)$. 
Then $\tau_{l,0,u}=\tau_{l,\Sigma_{l}(u)}\subset\tau_{l,\Delta'}\subset\tau_{l,\Delta}=\bigcap_{a\in\Sk^0(\Delta)}\tau_{l,a}$. 
Let $a\in X$. By Corollary~\ref{cor:Dlalpha=0 iff alpha in Sigma},  
\begin{equation*}
\begin{aligned}
\tau_{l,0,u}\subset\tau_{l,a}&\Leftrightarrow 
\bR_{\geq 0}f_0=\tau_{l,0,0}\subset\tau_{l,a-2l\mu(u)}\\
&\Leftrightarrow D_{l}(a-2l\mu(u))=0
\Leftrightarrow a\in\Sigma_{l}(u).
\end{aligned}
\end{equation*}
Hence $\Sk^0(\Delta)\subset\Sigma_{l}(u)$. 
Since 
$\Delta=\Conv(\Sk^0(\Delta))$, we have 
$\Delta\subset\Sigma_{l}(u)$. 

Let $a\in\Delta\cap X$. Since $\Delta\subset\Sigma_{l}(u)$, $f(\Sigma_{l}(u))\subset f(a)+\Cone(B_{l,\Delta})$ by Lemma~\ref{lemma:Cone BlDelta_cap X=-Delta+lin(Delta)}~(\ref{item: chara of DeltacapX}), and 
hence $f(\Sigma_{l}(u))\subset f(a)+\Cone(B_{l,\Delta'})$ because 
$\Cone(B_{l,\Delta'})\supset \Cone(B_{l,\Delta})$.
 By Lemma~\ref{lemma:Cone BlDelta_cap X=-Delta+lin(Delta)}~(\ref{item: chara of DeltacapX}), we obtain 
$a\in\Delta'\cap X$, so that 
$\Delta\cap X\subset\Delta'\cap X$. Hence $\Delta\subset\Delta'$.  This proves (\ref{item:tau_Delta tau_Delta'}). 
\end{proof}
\subsection{$\Cut(\Fan(\xi^{\natural}_{l}))$}

\begin{lemma}\label{lemma:Cut and Vor_l}Let $\tau\in\Fan(\xi^{\natural}_{l})$. Then 
\begin{enumerate}
\item\label{item:Cut tau_l_a_u} $\Cut(\tau_{l,\alpha,u})=\Cut(\tau_{l,\alpha,0})-e(\xi)u$\ $(\forall\alpha\in\Sigma_{l}, \forall u\in X^{\vee})$; 
\item\label{item:Cut tau_l_0_u} $\Cut(\tau_{l,0,u})=\{-e(\xi)u\}$\ $(\forall u\in X^{\vee})$;
\item\label{item:Xvee R} $X^{\vee}_{\bR}=\bigcup_{a\in \Sk^0(\Vor_{l})}\Cut(\tau_{l,a})$;
\item\label{item:dim g Cut(tau)}$\dim \Cut(\tau)=g$ iff 
$\Cut(\tau)=\Cut(\tau_{l,a})$ 
for some $a\in\Sk^0(\Vor_{l})$;
\item\label{item:dim g-1 Cut(tau)}$\dim \Cut(\tau)=g-1$ iff there is an adjacent pair $\Cut(\tau_{l,a})$ and $\Cut(\tau_{l,b})$ 
such that $a,b\in\Sk^0(\Vor_{l})$ and 
$\Cut(\tau)=\Cut(\tau_{l,a})\cap\Cut(\tau_{l,b})$. 
\end{enumerate}
\end{lemma}
\begin{proof}Since $\delta_u^*:B_{l,\alpha,0}\to B_{l,\alpha,u}$ sends 
$s^{x_0}w^x\mapsto s^{x_0}b^e(u,x)w^x$, 
$\delta_u^*$ induces an isomorphism $\tau_{l,\alpha,0}^{\vee}\simeq\tau_{l,\alpha,u}^{\vee}$ sending $x_0m_0+x\mapsto (x_0+e(\xi)u(x))m_0+x$ in view of 
$v_sb^e(u,x)=e(\xi)u(x)$. 
Hence $z_0f_0+(z+e(\xi)z_0u)\in\tau_{l,\alpha,0}$ 
iff $z_0f_0+z\in\tau_{l,\alpha,u}$, 
so that $z+e(\xi)u\in\Cut(\tau_{l,\alpha,0})$ 
iff $z\in\Cut(\tau_{l,\alpha,u})$. 
This proves (\ref{item:Cut tau_l_a_u}).  
By Lemma~\ref{lemma:cone taulalpha0}~(\ref{item:case Sigma_0_0}), 
$\tau_{l,0,u}=\bR_{\geq 0}(f_0-e(\xi)u)$. 
This proves (\ref{item:Cut tau_l_0_u}). 

Next we prove (\ref{item:Xvee R}).  Let 
$z\in X^{\vee}_{\bR}$, $b\in X$ and 
$H(x):=D_{l}(x)+z(x)$ $(x\in X)$. 
By 
Lemma~\ref{lemma:tauDelta=BlDeltavee}~(\ref{item:tau_l_Delta}), 
$z\in\Cut(\tau_{l,b})\Leftrightarrow f_0+z\in\tau_{l,b}
\Leftrightarrow f_0+z\in \Cone(A_{l,b})^{\vee}$. 
Since $A_{l,b}$ is generated by 
$s^{D_{l}(x)}w^x/s^{D_{l}(b)}w^b$ 
$(x\in X)$, we see $f_0+z\in \Cone(A_{l,b})^{\vee}
\Leftrightarrow H(x)\geq H(b)\ (\forall x\in X)$. 
The function $H$ has the minimum in $X$, 
as is shown in the same manner as in Lemma~\ref{lemma:completeness}. 
Hence there exists $a\in X$ such that 
$H(a)=\min_{x\in X}H(x)$. Then $z\in\Cut(\tau_{l,a})$. 
If $a\in\Sk^0(\Vor_{l})$, then (\ref{item:Xvee R}) is proved. Otherwise, there exists $\Delta\in\Vor_{l}$ such that $a\in\Delta^0\cap X$. Let $c\in\Sk^0(\Delta)$. Then $z\in\Cut(\tau_{l,a})=\Cut(\tau_{l,\Delta})\subset\Cut(\tau_{l,c})$ by Lemma~\ref{lemma:tauDelta=BlDeltavee}~(\ref{item:tau_l,a=tau_l_Delta})-(\ref{item:tau_Delta tau_Delta'}).   This proves (\ref{item:Xvee R}). \par 
 By Lemma~\ref{lemma:tauDelta=BlDeltavee}, 
$\dim\Cut(\tau_{l,a})=g$ 
$(\forall a\in\Sk^0(\Vor_{l}))$.  If 
$\Cut(\tau)=g$, then $\Cut(\tau)=\Cut(\tau_{l,a})$ 
for some $a\in\Sk^0(\Vor_{l})$ 
by (\ref{item:Xvee R}).
This proves (\ref{item:dim g Cut(tau)}). 
If $\Cut(\tau)=g-1$, then $\Cut(\tau)$ is a one-codimensional 
face of $\Cut(\tau_{l,a})$ for some $a\in\Sk^0(\Vor_{l})$. 
By (\ref{item:Xvee R}), 
there is $b\in\Sk^0(\Vor_{l})\setminus\{a\}$ such that 
$\Cut(\tau)\subset\Cut(\tau_{l,a})\cap\Cut(\tau_{l,b})$. 
  Since 
$\dim\Cut(\tau)=\dim\Cut(\tau_{l,a})\cap\Cut(\tau_{l,b})$, we see  
$\Cut(\tau)=\Cut(\tau_{l,a})\cap\Cut(\tau_{l,b})$ 
by Remark~\ref{rem:def of faces}.   
 This proves (\ref{item:dim g-1 Cut(tau)}).
\end{proof}

\begin{defn}\label{defn:Delta_tau}
For any $\tau\in\Fan(\xi^{\natural}_{l})\setminus \{(0)\}$, we define 
$$\Delta(\tau)=\bigcap_{(-v)\in\Sk^0(\Cut(\tau))}\Sigma_{l}(v).
$$
\end{defn}

\begin{prop}
\label{prop:Delta to tauDelta corresp} 
We define a map
$$\cut_{l}:\Vor_{l}(\xi^{\natural})\to\Cut(\Fan(\xi^{\natural}_{l}))$$
by $\cut_{l}(\Delta)=\Cut(\tau_{l,\Delta})$. 
Then
\begin{enumerate}
\item\label{item:tau=tau_Delta}
 $\Delta(\tau)\in\Vor_{l}(\xi^{\natural})$ and $\cut_{l}(\Delta(\tau))=\Cut(\tau)$\ $(\forall \tau\in\Fan(\xi^{\natural}_{l}))$; 
\item\label{item:bij} $\cut_{l}$ is a bijection 
between $\Vor_{l}(\xi^{\natural})$ and $\Cut(\Fan_{l}(\xi^{\natural}))$ 
such that 
\begin{enumerate}
\item\label{item:tau_Delta sub tau_Delta'} $\sigma\subset\tau\Leftrightarrow\Cut(\sigma)\subset\Cut(\tau)\Leftrightarrow\Delta(\sigma)\supset\Delta(\tau)$;
\item\label{item:dim tau_Delta} $\dim\Cut(\sigma)+\dim\Delta(\sigma)=g$. 
\end{enumerate}
\end{enumerate}
\end{prop}

\begin{proof}
Let $\tau\in\Fan(\xi^{\natural}_{l})\setminus\{(0)\}$. 
By Remark~\ref{rem:def of faces} and 
Lemma~\ref{lemma:Cut and Vor_l}~(\ref{item:dim g Cut(tau)})-(\ref{item:dim g-1 Cut(tau)}),   
$$\tau=
\bigcap_{
\begin{matrix}
\tau\subset\tau_{l,a}, a\in\Sk^0(\Vor_{l})\end{matrix}}\tau_{l,a}.
$$ 

Let $a\in\Sk^0(\Vor_{l})$. 
 By Lemmas~\ref{lemma:tauDelta=BlDeltavee}
~(\ref{item:tau_l,a=tau_l_Delta})-(\ref{item:tau_Delta tau_Delta'}) 
and \ref{lemma:Cut and Vor_l}~(\ref{item:Cut tau_l_0_u}),  
\begin{align*}
\tau\subset\tau_{l,a}
&\Leftrightarrow \Cut(\tau)\subset\Cut(\tau_{l,a})
\Leftrightarrow \Sk^0(\Cut(\tau))\subset\Sk^0(\Cut(\tau_{l,a}))\\
&\Leftrightarrow  -e(\xi)v\in\Cut(\tau_{l,a})\ \ 
(\forall (-e(\xi)v)\in\Sk^0(\Cut(\tau)))\\
&\Leftrightarrow \Cut(\tau_{l,0,v})\subset\Cut(\tau_{l,a})\ \ 
(\forall (-e(\xi)v)\in\Sk^0(\Cut(\tau)))\\
&\Leftrightarrow \tau_{l,\Sigma_{l}(v)}\ (=\tau_{l,0,v})\subset\tau_{l,a}
\ \ (\forall (-e(\xi)v)\in\Sk^0(\Cut(\tau)))\\ 
&\Leftrightarrow \Sigma_{l}(v)\ni a\ \ 
(\forall (-e(\xi)v)\in\Sk^0(\Cut(\tau)))\\
&\Leftrightarrow a\in \bigcap_{-e(\xi)v\in\Sk^0(\Cut(\tau))}
\Sigma_{l}(v)\quad (=\Delta(\tau)).
\end{align*} 

It follows that $\tau\subset\tau_{l,a}$ 
iff $a\in\Delta(\tau)$. 
For $a\in\Sk^0(\Vor_{l})$, $a\in\Delta(\tau)$ iff $a\in \Sk^0(\Delta(\tau))$. 
Hence (\ref{item:tau=tau_Delta}) follows from  
$$\tau=\bigcap_{\begin{matrix}
\tau\subset\tau_{l,a}, 
a\in\Sk^0(\Vor_{l})
\end{matrix}}\tau_{l,a}=\bigcap_{a\in\Sk^0(\Delta(\tau))}\tau_{l,a}=\tau_{l,\Delta(\tau)}.$$  

Moreover 
 (\ref{item:tau_Delta sub tau_Delta'}) 
(resp. (\ref{item:dim tau_Delta}))  follows from 
Lemma~\ref{lemma:tauDelta=BlDeltavee}~(\ref{item:tau_Delta tau_Delta'}) (resp. (\ref{item:dim_X_Delta})). 
\end{proof}

\begin{cor}\label{cor:Delta=intersection}
Let $\Delta\in\Vor_{l}(\xi^{\natural})$ and $v\in X^{\vee}$. Then 
$\Delta=\bigcap_{\Delta\subset\Sigma_{l}(v)}\Sigma_{l}(v)$. 
\end{cor}
\begin{proof}By Proposition~\ref{prop:Delta to tauDelta corresp}, 
there exists $\tau\in\Fan(\xi^{\natural}_{l})$ 
such that $\Delta=\Delta(\tau)$. 
Let $A=\bigcap_{-v\in\Sk^0(\Cut(\tau_{l,\Delta}))}\Sigma_{l}(v)$ 
and $B=\bigcap_{\Delta\subset\Sigma_{l}(v)}\Sigma_{l}(v)$. 
Then $\Delta\subset B\subset A=\Delta$. Hence $\Delta=B$
\end{proof}

\begin{cor}\label{cor:Sk0 tau_l_Delta}
Let $\Delta\in\Vor_{l}(\xi^{\natural})$, $\Sigma^*_{l}:=
\Conv(v\in X^{\vee}; \Sigma_{l}(v)\cap\Sigma_{l}(0)\neq\emptyset)$, 
$n\in\bN$ and 
$n\Sigma_{l}^*:=\{nx; x\in\Sigma_{l}^*\}$.
Then 
\begin{enumerate}
\item\label{item:Sk0_Cut_Fan=Xvee} 
$\Sk^0(\Cut(\Fan(\xi^{\natural}_{l})))=e(\xi)X^{\vee}$;
\item\label{item Sk0_Cut_tau_Delta}
$\Sk^0(\Cut(\tau_{l,\Delta}))
=\{e(\xi)v;\Delta\subset\Sigma_{l}(-v), v\in X^{\vee}\}$;
\item\label{item:tau_l_a subset Sigma*_l}
$\Cut(\tau_{l,a})-\Cut(\tau_{l,a})=
\Cut(\tau_{l,a})+\Cut(\tau_{l,-a})\subset2\Sigma^*_{l}$\ 
$(\forall a\in X)$;
\item\label{item:n_Sigma*_l}$-\Sigma^*_{l}=\Sigma^*_{l}$ and 
$n\Sigma_{l}^*=\Sigma_{l}^*+\cdots+\Sigma_{l}^*$ $(\text{$n$-times})$.
\end{enumerate}
\end{cor}
\begin{proof}
If $\dim\cut_{l}(\Delta)=0$, then 
$\Delta=\Sigma_{l}(c)$ for some $c\in X^{\vee}$. 
Hence $\Sk^0(\Cut(\Fan(\xi^{\natural}_{l})))=e(\xi)X^{\vee}$ 
by Lemma~\ref{lemma:Cut and Vor_l}~(\ref{item:Cut tau_l_0_u}). 
This proves (\ref{item:Sk0_Cut_Fan=Xvee}).  Let $v\in X^{\vee}$. 
By Proposition~\ref{prop:Delta to tauDelta corresp}~(2a), $e(\xi)v\in\Sk^0(\Cut(\tau_{l,\Delta}))$ iff $\Cut(\tau_{l,\Sigma_{l}(-v)})
\subset\Cut(\tau_{l,\Delta})$ iff $\Sigma_{l}(-v)\supset\Delta$, which proves  (\ref{item Sk0_Cut_tau_Delta}). 
Let $a\in X$. By Corollary~\ref{cor:Sigma_l}, 
we can write $a=2l\mu(u)+a'$ for some $a'\in\Sigma_{l}$ and $u\in X^{\vee}$. Then $\Cut(\tau_{l,-a})=-\Cut(\tau_{l,a})$ and 
$\Cut(\tau_{l,\pm a})\pm e(\xi)u=\Cut(\tau_{l,\pm a'})\subset\Sigma^*_{l}$ by Lemma~\ref{lemma:Cut and Vor_l}~(\ref{item:Cut tau_l_a_u}).  Hence 
$\Cut(\tau_{l,a})-\Cut(\tau_{l,a})=\Cut(\tau_{l,a})+\Cut(\tau_{l,-a})
=\Cut(\tau_{l,a'})+\Cut(\tau_{l,-a'})\subset 2\Sigma^*_{l}$ because $\Sigma^*_{l}$ is convex.  
This proves (\ref{item:tau_l_a subset Sigma*_l}).   
Since $-\Sigma_{l}(v)=\Sigma_{l}(-v)$ $(\forall v\in X^{\vee})$, we have $-\Sigma_{l}^*=\Sigma_{l}^*$. 
Let $n\in\bN$. 
Since $(1/n)\sum_{i=1}^nx_i\Sigma_{l}^*$ 
for $x_i\in \Sigma_{l}^*$ $(\forall i\in [1,n])$, we have $n\Sigma_{l}^*=\Sigma_{l}^*+\cdots+\Sigma_{l}^*$ $(\text{$n$-times})$. This proves (\ref{item:n_Sigma*_l}). 
\end{proof}

\begin{lemma}\label{lemma:torus embedding Zl_Delta}
Let $\alpha\in\Sk^0(\Sigma_{l}(0))$. 
Then
\begin{enumerate}
\item\label{item:tau_alpha}
 $\Cone(\Cut(\tau_{l,\alpha,0}))
=\{u\in X^{\vee}_{\bR};
u(\beta)\geq u(\alpha)\ (\forall\beta\in\Sigma_{l})\}
=\Cone(C_{l}^{\alpha})^{\vee}$;
\item\label{item Xvee ConeCuttau}
$X^{\vee}_{\bR}=\bigcup_{\alpha\in\Sk^0(\Sigma_{l}(0))}
\Cone(\Cut(\tau_{l,\alpha,0}))$;
\item\label{item:Sigma_0} 
$\Sigma_{l}(0)
=\left\{x\in X_{\bR}; u(x)\geq u(\beta)\ 
\left(\begin{matrix}\forall \beta\in\Sk^0(\Sigma_{l}(0))\\ 
\forall u\in X^{\vee}\cap\Cone(\Cut(\tau_{l,\beta,0}))
\end{matrix}\right)
\right\}$.\end{enumerate}
\end{lemma}
\begin{proof}Let $\alpha\in\Sk^0(\Sigma_{l}(0))$ and 
$f_{x,\alpha}(u):=e(\xi)D_{l}(x)+u(x-\alpha)$.  
By Definition~\ref{defn:torus embedding of tildecX} and 
Corollary~\ref{cor:Dlalpha=0 iff alpha in Sigma},
we have $\Cut(\tau_{l,\alpha,0})=
\{u\in X^{\vee}_{\bR}; f_{x,\alpha}(u)\geq 0\ 
(\forall x\in X)\}$.  
Since the cone 
$\Cone(\Cut(\tau_{l,\alpha,0}))$ is bounded only by 
finitely many faces passing 
through the origin ($u=0$) defined by $f_{x,\alpha}=0$, 
it is bounded by $f_{x,\alpha}\geq 0$ 
with $D_{l}(x)=0$. Hence  
$\Cone(\Cut(\tau_{l,\alpha,0}))=\{u\in X^{\vee}_{\bR};u(x-\alpha)\geq 0\ 
(\forall x\in\Sigma_{l})\}=\Cone(C_{l}^{\alpha})^{\vee}$,
which proves (\ref{item:tau_alpha}).  

By Lemma~\ref{lemma:Cut and Vor_l}~(\ref{item:Xvee R}), 
$\bigcup_{\beta\in\Sk^0(\Sigma_{l}(0))}\Cut(\tau_{l,\beta,0})$ contains an open neighborhood of the origin in $X^{\vee}_{\bR}$, which proves (\ref{item Xvee ConeCuttau}). We see (\ref{item:Sigma_0}) by (\ref{item:tau_alpha}):
\begin{gather*}
\Sigma_{l}(0)=
\bigcap_{\beta\in\Sk^0(\Sigma_{l}(0))}
(\beta+\Cone(C_{l}^{\beta}))\\
=\left\{x\in X_{\bR}; u(x-\beta)\geq 0\ 
\left(\begin{matrix}\forall \beta\in\Sk^0(\Sigma_{l}(0))\\ 
\forall u\in X^{\vee}\cap\Cone(\Cut(\tau_{l,\beta,0}))
\end{matrix}\right)
\right\}. 
\end{gather*} 
\end{proof}
\begin{defn}
For $\Delta\in\Vor_{l}$, 
let $Z_{l}(\Delta)$ be the unique reduced closed $\tG_0$-orbit of 
$(W_{l,\Delta})_0$ and $\barZ_{l}(\Delta)$ 
the closure of $Z_{l}(\Delta)$ in $\tP_{l,0}$ with reduced structure. 
See Corollary~\ref{cor:torus emb ZlSigma(0)}.
\end{defn}
\begin{cor}\label{cor:torus emb ZlSigma(0)}
We define a fan in $X^{\vee}$ by
\begin{equation*}
\Fan(\Sigma_{l}(0))=\{
\text{$\Cone(\Cut(\tau_{l,\Delta}))$;  
$\Delta\in\Sk(\Sigma_{l}(0))$}\}.
\end{equation*}  Then  $\barZ_{l}(\Sigma_{l}(0))$ is a torus embedding 
associated with $\Fan(\Sigma_{l}(0))$.
\end{cor}
\begin{proof}This follows from Lemma~\ref{lemma:torus embedding Zl_Delta}~(\ref{item Xvee ConeCuttau}) and Theorem~\ref{thm:W_l,alpha,0}~(\ref{item:boundary_Sigma(0)0 ring str}).
\end{proof}

\subsection{The structure of $P_{l,0}$}
\label{subsec:str of Pl0}
\begin{thm}
\label{thm:str of P_l}  
Let $\pi_{l}:\tP_{l,0}\to P_{l,0}$ be the natural morphism. 
Let $\overline{\Orb}(P_{l,0})$ be the set 
of all $G_0$-stable closed irreducible reduced subschemes of $P_{l,0}$. 
We define a map $\omega_{l}:\overline{\Vor}_{l}\to \overline{\Orb}(P_{l,0})$ by $\omega_{l}(\Delta)=\pi_{l}(\barZ_{l}(\Delta))$. 
Then 
\begin{enumerate}
\item\label{item:omega_l} 
$\omega_{l}$ 
is a bijective correspondence between $\overline{\Vor}_{l}$ and 
$\overline{\Orb}(P_{l,0})$;
\item\label{item:indep l} $\omega_{l}$ is independent of $l$ in the sense that $\psi_{l,ll'}\circ\omega_{l}=\omega_{ll'}\circ\overline{\vor}_{l,ll'}$ $(\forall l'\geq 1)$;
\item\label{item:Irred P_l_0} the set of irreducible components of $P_{l,0}$ consists of $\pi_{l}(\barZ_{l}(\Sigma_{l}(v)))$\ 
$(v\in X^{\vee}/\beta(Y))$ where 
$Z_{l}(\Sigma_{l}(v))=(W_{l,0,v})_0$.
\end{enumerate}
\end{thm}
\begin{proof}
Every irreducible component of $\tP_{l,0}$ is a normal torus embedding, so that every $G_0$-stable irreducible reduced subscheme is the closure of a $G_0$-orbit, and every $G_0$-orbit is of the form $\orb(\tau)$ for some $\tau\in\Fan(\xi^{\natural}_{l})$  with the notation of \cite[Prop.~1.6, p.~10]{Oda85}. 
In our case, $\orb(\tau_{l,\Delta})=Z_{l}(\Delta)$ and any $\tau\in\Fan(\xi^{\natural}_{l})\setminus\{(0)\}$ is of the form $\tau_{l,\Delta}$ for some $\Delta\in\Vor_{l}$ by Proposition~\ref{prop:Delta to tauDelta corresp}~(\ref{item:tau=tau_Delta}). 
Thus $\omega_{l}$ is surjective. Since $\cut_{l}$ 
is bijective, so is 
the quotient map 
$\overline{\Vor}_{l}\to\Cut(\Fan(\xi^{\natural}_{l}))/(S_y;y\in Y)$ 
of $\cut_{l}$. Hence 
$\omega_{l}$ is bijective. This proves (\ref{item:omega_l}). 
(\ref{item:indep l}) follows from Remark~\ref{rem:vo_lll'} and Theorem~\ref{thm:scheme normalization isom}. 
 
 By Theorem~\ref{thm:W_l,alpha,0}~(\ref{item:Sigma(0)0})-(\ref{item:boundary_Sigma(0)0}), 
$P_{l,0}$ consists of $|B|$ irreducible components $(\pi_{l}(\barZ_{l}(\Sigma_{l}(v)));v\in X^{\vee}/\beta(Y))$ where  
$\barZ_{l}(\Sigma_{l}(v))=S_y\barZ_{l}(\Sigma_{l}(v+\beta(y)))$\ $(\forall y\in Y)$. This proves (\ref{item:Irred P_l_0}).  This completes the proof. 
\end{proof}

\begin{cor}\label{cor:Pl0 is reduced} 
$P_{l,0}$ is reduced and Cohen-Macaulay.
\end{cor}
\begin{proof}Since $\tP_{l}$ is an $S$-torus embedding, $\tP_{l}$ is 
Cohen-Macaulay, hence so is $P_{l}$ 
by Lemma~\ref{lemma:toric is excellent}~(\ref{item:C normal CM}), 
so that $P_{l,0}$ is Cohen-Macaulay. 
By Theorem~\ref{thm:W_l,alpha,0}~(\ref{item:Sigma(0)0}), $P_{l,0}$ is regular at the generic point of any irreducible component of $P_{l,0}$. Hence $P_{l,0}$ is reduced everywhere by $(\op{R}_0)$ and $(\op{S}_1)$.  
\end{proof}

\section{Relative compactifications}
\label{sec:relative compactifications}
We use 
Notation~\ref{notation:Notation Rinit Sinit keta=Kmin(xi)}. 
Let $(G,\cL)$ a semiabelian $R_{\init}$-scheme 
 with $G_{0_{\init}}$ a {\it $k(0_{\init})$-split torus} 
and $\cL$ ample symmetric and rigidified, and 
$\cG$ the N\'eron model (over $R_{\init}$) of $G_{\eta_{\init}}$. 
Let $\xi:=FC(G,\cL)=(X,Y,a,b,A,B)$, and let 
$$\xi^{\natural}=(X,Y,\epsilon,b^e,E,\Sigma)$$  
be a N\'eFC kit over $R$ of $\xi^e_{2N}$ 
in \S~\ref{subsec:graded alg R(xinatural)}.
By \cite[III, 8.1, p.~78]{FC90}, 
\begin{equation}\label{eq:Neron subgroup}
(\cG_{0_{\init}}/G_{0_{\init}})(\overline{k(0_{\init})})
\simeq\cG_{0_{\init}}(k(0_{\init})^{\sep})/G_{0_{\init}}(k(0_{\init})^{\sep})
\simeq X^{\vee}/\beta(Y).
\end{equation}

Let $\Sigma_{l}(0)$ be the Voronoi polytope for $c=0$ 
in \S~\ref{subsec:Voronoi cells}.  
Throughout this section,  
we assume (\ref{assump:Sigmal(0) integral}). 
Let 
$\tQ_{l}:=\Proj R_{l}(\xi^{\natural})$, 
and let $(\tP_{l},\cL_{l})$ be the normalization of 
$\tQ_{l}$ with $\tcL_{l}$ the pullback of $\cO_{\tQ_{l}}(1)$ to $\tP_{l}$. 
By Lemma~\ref{lemma:Voronoi_cap_X=Sigmal},  
$0\in\Sigma_{l}\cap\Sigma_{l}(0)^0$.  
By Theorem~\ref{thm:W_l,alpha,0}~(\ref{item:Sigma(0)0}), we define 
$$\tG_{l}=W_{l,0,0}\simeq T_{X,R}:=\Spec R[w^x;x\in X],$$ which 
is independent of $l$. We denote the restriction 
of $\tcL_{l}$ on $\tP_{l}$ to $\tG_{l}$ by the same letter $\tcL_{l}$. 
Then we define an octuple 
\begin{gather*}
\zeta_{l}=(\tG_{l}, X, Y, \iota_{l}, 
\phi_{l}, b, \tcL_{l}, \epsilon_{l}), 
\end{gather*}where $\phi:Y\to X$ is the inclusion, 
$\phi_{l}=2Nl\phi$, and $\iota_{l}: Y\to\tG_{l}(k(\eta))$ 
(resp. $\iota_{l}^e: X^{\vee}\to\tG_{l}(k(\eta))$) 
is a homomorphism (independent of $l$) defined by 
\begin{gather*}
w^x(\iota_{l}(y))=b(y,x)\quad (x\in X, y\in Y),\\
w^x(\iota^e_{l}(u))=b^e(u,x)\quad (x\in X, u\in X^{\vee}).
\end{gather*}

This octuple $\zeta_{l}$ is a split object 
in $\DD_{\ample}$  
in the sense of \cite[p.~57]{FC90}, where  
the abelian part of $\tG_{l}$ is trivial. 
Moreover 
$(\tP_{l}, \tcL_{l}, \tG_{l}, Y, \iota_{l})$ 
is a relatively complete 
model of $\zeta_{l}$ with data (a)-(e) 
satisfying the conditions (1)-(3) of \cite[pp.~60-61]{FC90}. 
Under our notation, (a)-(e) are given as follows:  
\begin{enumerate}
\item[(a)] the normal $S$-torus embedding $\tP_{l}$   
associated with $\Fan(\xi^{\natural}_{l})$, which contains  
$\tG_{l}$ as a Zariski open subset;
\item[(b)] the invertible sheaf $\tcL_{l}$ on $\tP_{l}$ 
extending $\cO_{\tG_{l}}$ on $\tG_{l}$; \footnote{See 
\cite[6.1.2]{Nakamura16}. }
\item[(c)] the action of $\tG_{l}$ on $\tP_{l}$ extending the 
translation action of $\tG_{l}$ on itself;
\item[(d)] the action of the lattice 
$Y$ on $(\tP_{l},\tcL_{l})$ 
extending that on $(\tG_{l,\eta},\cO_{\tG_{l,\eta}})$ 
via $\iota_{l}$, which we denote by $S_y:\tP_{l}\to\tP_{l}$ 
and $\tS_y:S_y^*\tcL_{l}\to\tcL_{l}$;
\item[(e)] the action of $\tG_{l}$ on $\tcL_{l}$ extending 
that on $\cO_{\tG_{l}}$. \footnote{See \S~\ref{subsec:Rl(xinatural)}. 
See also \S~\ref{sec:action of Gl}.}
\end{enumerate}

Let $\delta_u:=\tT_{\iota_{l}^e(u)}$ 
be translation of $\tG_{l}$ by $\iota^e(u)$\ $(u\in X^{\vee})$ and 
$$\tC_{l}:=(\tP_{l}\setminus \bigcup_{u\in X^{\vee}}\delta_u(\tG_{l}))_{\red},\ \tD_{l}:=(\tP_{l}\setminus \bigcup_{y\in Y}S_y(\tG_{l}))_{\red}.
$$ 
Let $(P_{l},\cL_{l})$ (resp. $C_{l}$, $D_{l}$) be the algebraization 
of the formal quotient $(\tP_{l}^{\wedge},\tcL_{l}^{\wedge})/Y$ 
(resp. $\tC_{l}^{\wedge}/Y$, $\tD_{l}^{\wedge}/Y$). We also define  
\begin{equation} 
\label{defn:cGl=union of deltau Gl}
G_{l}=P_{l}\setminus D_{l},\quad 
\cG_{l}=P_{l}\setminus C_{l}
=\bigcup_{u\in X^{\vee}/\beta(Y)}\delta_u(G_{l}),
\end{equation}where  
$G_{l}$ is a semiabelian $S$-scheme 
with $G_{l}^{\wedge}\simeq\tG_{l}^{\wedge}$ 
by \cite[5.7/5.8]{FC90}.

\begin{thm}
\label{thm:summary of Pl Gl}
Assume (\ref{assump:Sigmal(0) integral}). Then 
\begin{enumerate}
\item\label{item:Pl}$P_{l}$ is an irreducible flat 
projective normal Cohen-Macaulay 
$S$-scheme;
\item\label{item:open} 
$G_{l}$ and $\cG_{l}$ 
are open subschemes of $P_{l}$ and semiabelian group $S$-schemes
such that $P_{l,\eta}=G_{l,\eta}=\cG_{l,\eta}$ is an abelian variety  
and $G_{l}=\cG_{l}^0$;
\item\label{item:conn neron model}
 $G_{l}$  is the connected N\'eron model of $G_{l,\eta}$;
\item\label{item:Phi acts}
$(\cG_{l,0}/G_{l,0})(\overline{k(0)})\ 
(\simeq X^{\vee}/\beta(Y))$ acts on $P_{l}$;
\item\label{item:sub of neron model}  
the group $(\cG_{l,0}/G_{l,0})(\overline{k(0)})$ 
is (identified with) a subgroup 
$e(\xi)(X^{\vee}/\beta(Y))$ of the component group $X^{\vee}/e(\xi)\beta(Y)$ of the N\'eron model of $G_{l,\eta}$; 
\item\label{item:codim two}
 $\codim_{P_{l,0}} (P_{l,0}\setminus\cG_{l,0})=1$,  
and any irreducible component of $P_{l,0}$ 
is the closure in $P_{l,0}$ of an irreducible component of $\cG_{l,0}$;
\item\label{item:uniqueness}   
$P_{l}$, $G_{l}$ and $\cG_{l}$ are independent of the choice of $l$;
\item\label{item:pullback conn neron} 
$(G_{l},\cL_{l})\simeq 
(G_S,\cL^{\otimes 2Nl}_S)$ as polarized group $S$-schemes;
\item\label{item:pullback of neron over Rinit}
$\cG_{l}\simeq\cG_S$: the pullback to $S$ of the N\'eron model $\cG$ of 
$G_{\eta_{\init}}$; the quotient of $\cG_{l}$ 
by $\Aut(k(\eta)/k(\eta_{\init}))$ is the N\'eron model of 
$G_{\eta_{\init}}$;
\item\label{item:ample cubical cLl}the restriction of $\cL_{l}$ 
to $\cG_{l}$ is a unique ample cubical invertible sheaf 
on $\cG_{l}$ such that $\cL_{l,\eta}\simeq\cL^{\otimes 2Nl}_{\eta}$.  
\end{enumerate}
\end{thm}
\begin{proof} 
 (\ref{item:Pl})-(\ref{item:conn neron model}) 
are due to   
\cite[p.~66]{FC90}, Proposition~\ref{prop:semi abelian into neron} and 
Eq.~(\ref{defn:cGl=union of deltau Gl}). Note that $P_{l}$ is normal 
 and Cohen-Macaulay by 
Lemma~\ref{lemma:toric is excellent}. 
The (component) group 
$(\cG_{l,0}/G_{l,0})(\overline{k(0)})$ acts on $P_{l}$ via $\delta_u$ 
by Eq.~(\ref{defn:cGl=union of deltau Gl}), 
which proves (\ref{item:Phi acts}). 

The semiabelian $S$-scheme $G_{l}$ is 
the connected N\'eron model of $G_{l,\eta}$ 
by Proposition~\ref{prop:semi abelian into neron}. 
 Let $\cG'_{l}$ be the N\'eron model of $G_{l,\eta}$. 
Since $\tG_{l}=W_{l,0,0}$ by (a), 
in view of \cite[III, 8.1, p.~78]{FC90}   
and Lemma~\ref{lemma:completeness},  
\begin{align*}
\cG'_{l}(R)&
=G_{l}(k(\eta))=P_{l}(k(\eta))=P_{l}(R)=P^{\wedge}_{l}(R)
=\tP_{l}^{\wedge}(R)/Y\\
&=\tP_{l}(R)/Y=\tP_{l}(k(\eta))/Y=\tG_{l}(k(\eta))/\iota_{l}(Y),\\
\cG'_{l}(R)/G_{l}(R)&\simeq(\tG_{l}(k(\eta))/\iota_{l}(Y))/G_{l}(R)
\simeq \tG_{l}(k(\eta))/\tG_{l}(R)\cdot \iota_{l}(Y)\\
&\simeq 
\Hom(X,k(\eta)^{\times}/R^{\times})/e(\xi)\beta(Y)
\simeq X^{\vee}/e(\xi)\beta(Y).
\end{align*}
The same is true for  
any \'etale $R$-algebra $R'$: 
{\small$\cG'_{l}(R')/G_{l}(R')\simeq X^{\vee}/e(\xi)\beta(Y)$.}
Hence by \cite[2.2/13]{BLR90}, 
\begin{equation}\label{eq:comp of cG'l} 
\begin{aligned}
(\cG'_{l,0}/G_{l,0})(\overline{k(0)})
&\simeq \cG'_{l,0}(\overline{k(0)})/G_{l,0}(\overline{k(0)})
\simeq\cG'_{l,0}(k(0)^{\sep})/G_{l,0}(k(0)^{\sep})\\
&\simeq\cG'_{l}(R^{\unr})/G_{l}(R^{\unr})\simeq X^{\vee}/e(\xi)\beta(Y)
\end{aligned}
\end{equation}where 
$R^{\unr}$ is the integral closure of $R$ 
in the maximal unramified extension of $k(\eta)$. 
Since $\delta_u(G_{l})\cdot\delta_v(G_{l})
=\delta_{u+v}(G_{l})$\ 
$(\forall u,v\in X^{\vee}/\beta(Y))$, 
$\cG_{l}$ is a smooth group $S$-scheme. 
Since $\cG_{l,\eta}\simeq\cG'_{l,\eta}$,
there is a canonical $S$-homomorphism $h:\cG_{l}\to \cG'_{l}$ 
by N\'eron mapping property and 
Corollary~\ref{cor:morphism of neron model of gr scheme}. 
By Proposition~\ref{prop:semi abelian into neron},    
$h$ is an isomorphism between the identity components 
of $\cG_{l}$ and $\cG'_{l}$. Let $e_{l}$ be 
the unit section of $G_{l}$. Since 
$\cG_{l,0}/G_{l,0}\simeq X^{\vee}/\beta(Y)$ 
by Eq.~(\ref{defn:cGl=union of deltau Gl}),  there exists, 
for every connected component $Z$ of $\cG_{l,0}$, a section 
$\delta_u(e_{l})\in\cG_{l}(R)$ $(u\in X^{\vee}/\beta(Y))$ 
with $\delta_u(e_{l})\cap Z\neq\emptyset$.  
Hence $h$ is quasi-finite. 
Since $\cG'_{l}$ is smooth over $S$, 
$h$ is an open immersion by Zariski's main theorem.  Hence 
$h$ induces an injection 
$\cG_{l,0}/G_{l,0}\hookrightarrow\cG'_{l,0}/G_{l,0}$,  
whose image is $h(\cG_{l,0}/G_{l,0})=e(\xi)(X^{\vee}/\beta(Y))$ in 
$X^{\vee}/e(\xi)\beta(Y)$  because $v_sb^e(u,x)=e(\xi)u(x)$ by 
Lemma~\ref{lemma:be 2nd}. This proves 
(\ref{item:sub of neron model}). 

By Theorem~\ref{thm:str of P_l}, 
the number of irreducible components of 
$P_{l,0}$ is equal to that of $\cG_{l,0}$, hence 
any irreducible component of $P_{l,0}$ is the closure of an irreducible component of $\cG_{l,0}$. This proves (\ref{item:codim two}).\par

Let $l'\in\bN$. Then $P_{l}\simeq P_{ll'}$ 
by Theorem~\ref{thm:scheme normalization isom}. 
Since $\cG_{l}=P_{l}\setminus C_{l}$ and $C_{l}$ is the singular locus of 
$P_{l,0}$ for any $l$, we have $\cG_{l}\simeq \cG_{ll'}$.

Since  $G_{l,\eta}\simeq P_{l,\eta}
\simeq P_{ll',\eta}\simeq G_{ll',\eta}$, and 
$G_{ll'}$ is the connected N\'eron model 
of $G_{ll',\eta}$ by (\ref{item:conn neron model}), 
we have $G_{l}\simeq G_{ll'}$. This proves (\ref{item:uniqueness}). \par
Next, recall $\FC(G,\cL)=\xi=(X,Y,a,b,A,B)$ over $R_{\init}$ with $\cL$ symmetric. Then $a(y)=a(-y)$ and $a(y)^2=b(y,y)$ $(\forall y\in Y)$ by 
Lemma~\ref{lemma:FC(G_L) L symmetric}.  
Let $(G_{\eta},\cL_{\eta}):=(G_S,\cL_S)\times_S\eta$.  
By Lemma~\ref{lemma:H0(G,Lm) over keta as Fourier},  
 for any $m\geq 1$, 
\begin{equation*}
\Gamma(G_{\eta},\cL_{\eta}^{\otimes 2m})
=\left\{
\theta=\sum_{x\in X}\sigma_x(\theta)w^x;
\begin{matrix}
\sigma_{x+2my}(\theta)
=a(y)^{2m}b(y,x)\sigma_x(\theta)\\
\sigma_x(\theta)\in k(\eta)\ (\forall x\in X, \forall y\in Y)
\end{matrix}
\right\}.
\end{equation*}

Let $N:=|B_{\xi}|$.  
By Lemma~\ref{lemma:theta on Pl}, there exists $m_0$ such that 
for any $m\geq m_0$, 
\begin{equation}
\label{eq:Gamma(Gleta Leta^m) mgeqm0}
\begin{aligned}
\Gamma&(G_{l,\eta},\cL_{l,\eta}^{\otimes m})
=\Gamma(P_{l,\eta},\cL_{l,\eta}^{\otimes m})\\
&=\left\{
\theta=\sum_{x\in X}\sigma_x(\theta)w^x;
\begin{matrix}
\sigma_{x+2Nl my}(\theta)
=a(y)^{2Nl m}b(y,x)\sigma_x(\theta)\\
\sigma_x(\theta)\in k(\eta)\ (\forall x\in X, \forall y\in Y)
\end{matrix}
\right\}
\end{aligned}
\end{equation} 
Hence 
$\Gamma(G_{l,\eta},\cL_{l,\eta}^{\otimes m})=
\Gamma(G_{\eta},\cL_{\eta}^{\otimes 2Nl m}).$ 
Since $\xi=\FC(G,\cL)$, 
$G_{l,0}\simeq T_{X,k(0)}\simeq G_0$ by Theorem~\ref{thm:W_l,alpha,0}. It follows from  Theorem~\ref{thm:G isom G' iff FC same} that $(G_{l},\cL_{l}^{\otimes m})\simeq (G_S,\cL^{\otimes 2Nl m}_S)$\  
as polarized group $S$-schemes $(\forall m\geq m_0)$, so that 
$(G_{l},\cL_{l})\simeq (G_S,\cL^{\otimes 2Nl}_S)$. 
This proves (\ref{item:pullback conn neron}).\par
Next we prove (\ref{item:pullback of neron over Rinit}). 
By Eqs.~(\ref{eq:Neron subgroup})-(\ref{eq:comp of cG'l}) 
\begin{align*}
(\cG_{0}/G_{0})(\overline{k(0)})
&\simeq 
\Hom(X,k(\eta_{\init})^{\times}/R^{\times}_{\init})/\beta(Y)\\
&\overset{\Hom(X,i)\,\text{mod}\,e(\xi)\beta(Y)}{\hookrightarrow}\Hom(X,k(\eta)^{\times}/R^{\times})
/e(\xi)\beta(Y)\\
&=X^{\vee}/e(\xi)\beta(Y)\simeq 
(\cG'_{l,0}/G_{l,0})(\overline{k(0)})
\end{align*}where $i:k(\eta_{\init})^{\times}/R^{\times}_{\init}\to 
k(\eta)^{\times}/R^{\times}$ is the natural injection.   
Hence as subgroups of $(\cG'_{l,0}/G_{l,0})(\overline{k(0)})
\simeq X^{\vee}/e(\xi)\beta(Y)$, we have, by (\ref{item:Phi acts}), 
\begin{equation}\label{eq:Phi equal}
(\cG_{0}/G_{0})(\overline{k(0)})=e(\xi)(X^{\vee}/e(\xi)\beta(Y))=(\cG_{l,0}/G_{l,0})(\overline{k(0)}).
\end{equation}
 
It follows that $\cG_{0}=\cG_{l,0}$ and 
$\cG_S\simeq\cG_{l}$ as subgroup $S$-schemes 
of $\cG'_{l}$. 
Let $\Gamma:=\Gal(k(\eta)/k(\eta_{\init}))$. The group $\Gamma$ acts on both 
$\cG_{l}$ and $G_{l}$. Let $\barG_{l}:=G_{l}/\Gamma$ and 
$\bar\cG_{l}:=\cG_{l}/\Gamma$. 
By Theorem~\ref{thm:summary of Pl Gl}~(\ref{item:pullback conn neron}), 
$G_{l}\simeq G_S$, so that $\barG_{l}\simeq G$. 
Moreover we have a $\Gamma$-invariant morphism 
$h:\cG_{l}\to\cG$ via $\cG_S/\Gamma\simeq\cG$, 
which induces an isomorphism $\cG_{l}/\Gamma\simeq\cG$ 
by Eq.~(\ref{eq:Phi equal}). 
This proves (\ref{item:pullback of neron over Rinit}).\par
Next we prove (\ref{item:ample cubical cLl}). 
Let $\Phi:=X^{\vee}/\beta(Y)$, 
$u:=(u_i)_{i\in[1,n]}\in\Phi^n$ $(n=2,3)$ and 
$\delta_u:=(\delta_{u_i})_{i\in[1,n]}\in\Aut_S(\cG_{l}^n)$. Then 
\begin{equation}\label{eq:Phi invariance of Theta Lambda} 
\delta_u^*\Theta(\cL_{l})=\Theta(\cL_{l}),\ 
\delta_v^*\Lambda(\cL_{l})=\lambda(\cL_{l})\  
(\forall u\in\Phi^3,\forall v\in\Phi^2).
\end{equation}

By (\ref{item:pullback conn neron}), 
$(\cL_{l})_{|G_{l}}$ has a cubical structure $\tau_0$, 
which gives rise to a nonvanishing section 
$\tau_0\in\Gamma(G^3_{l}, \Theta(\cL_{l})^{\otimes(-1)})$. 
Since $G_{l}^3$ is dense open in $\cG^3_{l}$, 
$\tau_0$ induces a meromorphic section $\tau$ of 
$\Theta(\cL_{l})^{\otimes(-1)}$ over $\cG_{l}^3$ with 
$\tau_{|G^3_{l}}=\tau_0$.     
By Eq.~(\ref{eq:Phi invariance of Theta Lambda}), 
$\tau$ is regular and nonvanishing on $\cG^3_{l}$, 
so that $\tau$ gives an isomorphism $\Theta(\cL_{l})\simeq\cO_{\cG^3_{l}}$. 
Let $\xi:\cG^2_{l}\to\cG^2_{l}$ be the permutation $(x,y)\mapsto (y,x)$. 
Since $\Theta(\cL_{l|G_{l}})$ is cubical, 
$\Lambda(\cL_{l|G_{l}})$ is symmetric, {\it i.e.}, $\xi^*\Lambda(\cL_{l})=\Lambda(\cL_{l})$ on $G^2_{l}$, which extends 
to an isomorphism on $\cG^2_{l}$ by 
Eq.~(\ref{eq:Phi invariance of Theta Lambda}). 
Hence $\Lambda(\cL_{l})$ is symmetric, so that 
$\tau$ is a cubical structure of $\cL_{l}$. 
This proves 
(\ref{item:ample cubical cLl}). 
\end{proof}

\begin{cor}\label{cor:neron if exi=1}
If $e(\xi)=1$, then $\cG_{l}$ is the N\'eron model of 
$G_{\eta}=G_{\eta_{\init}}$. 
\end{cor}
\begin{proof}
Clear from Theorem~\ref{thm:summary of Pl Gl}~
(\ref{item:pullback of neron over Rinit}).  
\end{proof}

\begin{cor}\label{cor:rel compact of G}
Let $\Gamma:=\Gal(k(\eta)/k(\eta_{\init}))$. 
The N\'eron model $\cG$ of $G_{\eta_{\init}}$ 
is the quotient $\cG_{l}/\Gamma$,
 which is relatively compactified by the quotient $P_{l}/\Gamma$.
\end{cor}
\begin{proof}
Clear from Theorem~\ref{thm:summary of Pl Gl}~
(\ref{item:pullback of neron over Rinit}).  
\end{proof}

\begin{cor}\label{cor:H0(Pl Llm)}  
Eq.~(\ref{eq:Gamma(Gleta Leta^m) mgeqm0}) is true for any $m\geq 1$. 
In other words,  Lemma~\ref{lemma:theta on Pl} is true for any $m\geq 1$. 
\end{cor}
\begin{proof}By Theorem~\ref{thm:summary of Pl Gl}~(\ref{item:pullback conn neron}), $(P_{l,\eta},\cL_{l,\eta})\simeq (G_{l,\eta},\cL_{l,\eta})
\simeq (G_{\eta},\cL^{\otimes 2Nl}_{\eta})$. Hence Corollary follows from Lemma~\ref{lemma:H0(G,Lm) over keta as Fourier}.
\end{proof}

\section{The action of $G_{l}$ on $P_{l}$}
\label{sec:action of Gl}
The aim of this section is to prove Theorem~\ref{thm:action of Gl on Pl}.

\subsection{Rigid analytic spaces} 
\label{subsec:rigid analytic spaces}
We refer to \cite{Bosch14} for rigid analytic geometry.
Let $K$ be a field with non-trivial complete non-Archimedean valuation.
We briefly recall the definition of the category $\Rig_K$ 
of rigid $K$-spaces \cite[5.3/4]{Bosch14}, 
which is defined as a full subcategory 
of the category $\LGRS_K$ of locally G-ringed 
$K$-spaces \cite[5.3/1]{Bosch14}:
\begin{enumerate}
\item[(a)]
an object of $\LGRS_K$ is a pair $(X,\cO_X)$ of a set $X$ with Grothendieck topology \cite[5.1/1 and the paragraph after 5.1/2]{Bosch14} and a sheaf $\cO_X$ of $K$-algebras on $X$ whose stalk $\cO_{X,x}$ $(\forall x\in X)$ is a local ring;
\item[(b)]
a morphism of $\LGRS_K$ is a pair $(f,f^\sharp):(X,\cO_X)\to(Y,\cO_Y)$ of a map $f:X\to Y$ that is continuous with respect to the Grothendieck topologies and a natural transformation $f^\sharp:\cO_Y\to f_*\cO_X$ whose induced homomorphism $\cO_{Y,f(x)}\to \cO_{X,x}$ $(\forall x\in X)$ is a local homomorphism;
\item[(c)]
an object of $\Rig_K$ is an object of $\LGRS_K$ whose Grothendieck topology satisfies the completeness conditions \cite[5.1/5(G$_0$)-(G$_2$)]{Bosch14} and that is locally isomorphic to the affinoid $K$-space $\Sprig A$ associated with an affinoid $K$-algebra $A$ \cite[the paragraph after 5.3/1]{Bosch14}, whose underlying set is the set $\Max(A)$ of all maximal ideals of $A$. 
\end{enumerate}

A $K$-scheme $X$ is regarded as an object of $\LGRS_K$ whose Grothendieck topology is just the Zariski topology 
of $X$ \cite[the paragraph after 5.1/1]{Bosch14}.

\begin{defn}[{\cite[5.4/3-5]{Bosch14}}]\label{defn:an}
Let $\Sch_K^{\lft}$ be the category of $K$-schemes 
locally of finite type (abbr. lft). 
Let $X$ be an lft $K$-scheme. Then there exists 
 a rigid $K$-space $Y$ together with a morphism $\iota:Y\to X$ of $\LGRS_K$ satisfying the following universal property: any morphism $Z\to X$ of $\LGRS_K$ from a rigid $K$-space $Z$ uniquely factors as $Z\to Y\stackrel{\iota}{\to}X$. 
We denote $Y$ (resp. $\iota$) by $X^{\an}$ (resp. $\iota_X$) and call 
$(X^{\an},\iota_X)$ or simply $X^{\an}$ the 
{\it rigid analytification of $X$.} 
 For a morphism $f:X_1\to X_2$ of lft $K$-schemes, there exists a unique morphism $f^\an:X_1^\an\to X_2^\an$ of rigid $K$-spaces such that 
$f\circ\iota_{X_1}=\iota_{X_2}\circ f^{\an}$. Thus we obtain a functor 
$\an:\Sch_K^{\lft}\to\Rig_K$ sending $X\mapsto X^{\an}$. 
\end{defn}

The rigid analytification of an lft $K$-scheme $X$ defines a functor $F_X:\cF\mapsto \cF^{\an}$ from the category of coherent $\cO_X$-modules to the category of coherent $\cO_{X^\an}$-modules \cite[6.3/11]{Bosch14}.
\begin{thm}[{\cite[6.3/12-13]{Bosch14}, \cite[S\"atze 4.11 und 5.1]{Kopf74}}]
Let $X$ be a proper $K$-scheme.
Then the functor $F_X$ gives a categorical equivalence.
\end{thm}
The kernel of the homomorphism $\cO_X\to f_*\cO_Y$ associated with a closed immersion $f:Y\to X$ of rigid $K$-spaces \cite[6.3/1]{Bosch14} is a coherent ideal of $\cO_X$ \cite[6.1/4 and 3.1/3(i)]{Bosch14}, which gives rise to 
a one-to-one correspondence between the set of closed subspaces of $X$ and the set of coherent ideals of $\cO_X$.
As a consequence of the above theorem, we obtain

\begin{thm}\label{thm:Chow}
Let $X$ be a proper $K$-scheme. 
Then 
\begin{enumerate}
\item\label{item:closed subspace} if $Y$ a closed subspace of $X^\an$, then there exists a unique closed subscheme $Z$ of $X$ such that $Z^\an\simeq Y$ over $X^\an$;
\item\label{item:Hom} if $Y$ and $Z$ are closed subschemes of $X$, then 
 there is a bijection  
$$\Hom_{\Sch_K^{\lft}}(Y,Z)\ni f\mapsto f^{\an}\in \Hom_{\Rig_K}(Y^{\an},Z^{\an}).$$ 
\end{enumerate}
\end{thm}

We denote the valuation ring of $K$ by $R$ 
(not necessarily discrete in this subsection), 
the maximal ideal of $R$ by $I$ 
and the category of formal $R$-schemes locally of topologically finite type by $\op{FSch}_R^{\op{ltft}}$ \cite[7.4/1]{Bosch14}.  
\begin{prop}[{\cite[7.4]{Bosch14}}]\label{prop:generic fiber}
The functor $A\mapsto A\otimes_RK$ from the category of $R$-algebras of topologically finite type \cite[7.3/3]{Bosch14} to the category of affinoid $K$-algebras \cite[3.1/1]{Bosch14} gives rise to a functor 
$$\rig:\op{FSch}_R^{\op{ltft}}\to\Rig_K$$ 
sending $\Spf A$ to $(\Spf A)_{\rig}:=\Sprig A\otimes_RK$.
\end{prop}
\begin{defn}[{\cite[7.4/3]{Bosch14}}]\label{defn:gen}
Let $X$ be a formal $R$-scheme locally of topologically finite type
 (abbr. ltft) and $X_{\rig}:=\rig(X)$. We call $X_{\rig}$ 
the \textit{generic fiber of $X$}.   
\end{defn}
\begin{defn}[{\cite[7.4/1, 7.3/3, 7.3/5]{Bosch14}}]\label{defn:adm}
A formal $R$-scheme $X$ is said to be \textit{admissible} 
if it is ltft and $\cO_X$ 
has no $I$-torsion where $I$ is the maximal ideal of $R$. 
If $R$ is a CDVR, then $X$ is admissible iff it is ltft and $R$-flat.
\end{defn}

\begin{defn}\label{defn:closed pt}
For a locally G-ringed $K$-space $X$ (resp. an lft $K$-scheme $Y$),
we denote the underlying set of $X$ by $|X|$ 
(resp. the set of closed points of $Y$ by $|Y|_{\closed}$). 
Note that $|Y|_{\closed}=Y(\Omega)/G_{\abs}$ where 
$G_{\abs}=\Aut(\Omega/k(\eta))$. For an lft $R$-scheme $Z$, 
we denote $Z_K=Z\otimes_RK$. 
\end{defn}

\begin{defn}[{\cite[8.3/1, 8.3/3, 8.3/8]{Bosch14}}]\label{defn:rigpt}
Let $X$ be an admissible formal $R$-scheme.
A \textit{rig-point $Q$ of $X$} is an admissible formal closed $R$-subscheme of $X$ that is isomorphic to the formal spectrum of  
a local integral domain $R'$ of dimension $1$.  
We denote the set of rig-points of $X$ by $|X|_{\rig}$. 
\end{defn}
\begin{prop}[{\small {\cite[5.4/4, 7.3/3, 7.2/4, 8.3/7]{Bosch14}}}]
\label{prop:rigpt}
The following are true:
\begin{enumerate}
\item any $R$-flat ltft formal $R$-scheme   
is admissible;
\item\label{item:|X|_closed = |X_an|}
if $X$ is an lft $K$-scheme, then 
there is  a canonical bijection $|X|_{\closed}\simeq|X^{\an}|$ 
sending 
a closed point $P$ of $X$ to its rigid analytification $P^{\an}$;
\item\label{item:Xwedge_rig subset Xeta_an}
 if $X$ is an $R$-flat lft $R$-scheme,  
then $X^{\wedge}$  
is an admissible formal $R$-scheme;  $X^{\wedge}_{\rig}$ admits a canonical open immersion $\iota$ into $X_K^{\an}$, 
and if moreover 
$|X^{\wedge}_{\rig}|\simeq|X_K^{\an}|$ 
\footnote{For example, if $X$ is proper over $R$, then $|X^{\wedge}_{\rig}|\simeq|X_K^{\an}|$ by $\iota$. }  by $\iota$, 
then $X^{\wedge}_{\rig}\simeq X_K^{\an}$ 
by $\iota$; 
\item\label{item:|Xrig|=|X|rig} 
if $X$ is an admissible formal $R$-scheme, then there is a canonical 
bijection $|X|_{\rig}\simeq |X_{\rig}|$ 
sending a rig-point of $X$ to its generic fiber. 
\end{enumerate}
\end{prop}

\begin{example}
\label{exam:Tate algebra}
Assume that $R$ is a CDVR.  Let $R[\zeta]$ be 
a polynomial $R$-algebra and $\fraka$ an ideal of $R[\zeta]$. 
We define  
$$R\langle\zeta\rangle=\left\{\sum_{\nu=0}^{\infty}c_{\nu}\zeta^{\nu}\in R[[\zeta]];\lim_{\nu\to\infty} c_{\nu}=0\right\},\ K\langle\zeta\rangle=R\langle\zeta\rangle\otimes_RK.$$  

 Let $X:=\Spec R[\zeta]/\fraka$. Then 
$R\langle\zeta\rangle\simeq R[\zeta]^{\wedge}$  
and $X^{\wedge}\simeq\Spf R\langle\zeta\rangle/\fraka R\langle\zeta\rangle$. 
Moreover $X_K^{\an}$ and $X^{\wedge}_{\rig}$ are given by  
$$X_K^{\an}=\bigcup_{i=0}^{\infty}\Sprig K\langle \zeta/c^i\rangle/\fraka K\langle \zeta/c^i\rangle,\  X^{\wedge}_{\rig}=\Sprig K\langle\zeta\rangle/\fraka K\langle\zeta\rangle$$ 
for any $c\in K^{\times}$ with $|c|>1$. See \cite[2.2]{Bosch14}. 
\end{example}

\subsection{The rigid analytic space $\tP_{l,\eta}^{\an}$}
\label{subsec:adm cover of tPl}
In the remainder of this section, we assume (\ref{assump:Sigmal(0) integral}).
Now we return to \S~\ref{sec:relative compactifications}.  
Let $\tP_{l,\eta}^{\an}$ be 
the rigid analytification of the $k(\eta)$-scheme 
$\tP_{l,\eta}$. 
Let $a\in\Sk^0(\Vor_{l})$, $B_{l,a}:=R[\tau_{l,a}^{\vee}\cap\tX]$, 
$W_{l,a}:=\Spec B_{l,a}$  and 
$R\langle\tau_{l,a}^{\vee}\cap\tX\rangle:=B_{l,a}^{\wedge}$. 
Then $W^{\wedge}_{l,a}$ (resp. $W^{\wedge}_{l,a,\rig}$) is an  
admissible formal $R$-scheme (resp. an affinoid open subset of 
$\tP^{\an}_{l,\eta}$) 
given by 
$W^{\wedge}_{l,a}=\Spf R\langle\tau_{l,a}^{\vee}\cap\tX\rangle$ 
(resp. $W^{\wedge}_{l,a,\rig}=\Sprig k(\eta)\langle \tau_{l,a}^{\vee}\cap\tX\rangle$). By Lemma~\ref{lemma:cutlog} below,  
$(W^{\wedge}_{l,a,\rig};a\in\Sk^0(\Vor_{l}))$ 
is a locally finite affinoid covering of $\tP_{l,\eta}^{\an}$.
 By Lemma~\ref{lemma:completeness} and Proposition~\ref{prop:rigpt}~(\ref{item:|X|_closed = |X_an|})-(\ref{item:|Xrig|=|X|rig}), 
$|\tP_{l,\eta}^{\an}|=|\tP_{l,\eta}|_{\closed}=\tP_{l}(\Omega)/G_{\abs}=\tP_{l}(R_{\Omega})/G_{\abs}
\subset|\tP^{\wedge}_{l}|_{\rig}
=|\tP^{\wedge}_{l,\rig}|\subset |\tP_{l,\eta}^{\an}|$ 
where $\Omega:=\overline{k(\eta)}$. 
Hence we obtain  
\begin{equation}\label{eq:tP_wedeg_rig=P_eta_an}
|\tP^{\wedge}_{l,\rig}|=|\tP_{l,\eta}^{\an}|
=\tP_{l}(R_{\Omega})/G_{\abs},\ 
\tP^{\wedge}_{l,\rig}=\tP_{l,\eta}^{\an}.
\end{equation}

 Let $Q\in |\tP^{\wedge}_{l,\rig}|=\tP_{l}(\Omega)/G_{\abs}$ 
and $\tQ\in\tP_{l}(\Omega)$ 
 a representative of $Q$. Then  $\log|w^{x}(Q)|:=\log|\tQ^*(w^x)|\in\bQ$ $(x\in \tX)$ is well-defined in the sense that 
it is independent of the choice of $\tQ$. 
See Definition~\ref{defn:log and cutlog}. \par

\begin{defn}\label{defn:specialization}
Let $\Omega$ be an algebraic closure of $k(\eta)$ 
and $R_{\Omega}$ the integral closure of $R$ in $\Omega$. 
Let $Q\in\tP_{l}(\Omega)$. 
By Lemma~\ref{lemma:completeness}, 
$\tP_{l}(\Omega)=\tP_{l}(R_{\Omega})$, so that 
$Q\in\tP_{l}(R_{\Omega})$, which is an $R$-flat finite 
formal $R$-scheme. \footnote{As is equivalent, it is 
an admissible finite formal $R$-scheme. 
See \S~\ref{subsec:rigid analytic spaces}.}
The closed fiber $Q_0\in\tP_{l,0}(\overline{k(0)})$ of $Q$ is 
the specialization $\specialization(Q)$ of $Q$ 
in the sense 
\cite[p.~200]{Bosch14}. 
The map $\specialization:\tP_{l}(R_{\Omega})
\to \tP_{l}(\overline{k(0)})$ is surjective by \cite[8.3/8]{Bosch14}.
\end{defn}

\begin{defn}\label{defn:vs}
Let $s$ be a uniformizer of $R$ and define an 
{\it absolute value $|\cdot|$ of $k(\eta)$} 
by $|a|=e^{-v_s(a)}$ for $a\in k(\eta)$. Then  
$v_s$ and $|\cdot|$ can be uniquely extended to $\Omega$ 
\cite[A, Th.~3]{Bosch14}.   Note that $|a|\geq 0$, and $|a|=0$ iff $a=0$.
\end{defn}

\begin{defn}\label{defn:log and cutlog}
Let $Q\in\tP_{l}(\Omega)$ and 
$\log|w^{x}(Q)|:=\log|Q^*(w^x)|\in\bQ$ $(x\in \tX)$.
Since $w^{m_0}(Q)=s$, we define 
\begin{align*}
\log(Q)&=-\sum_{i=0}^g(\log|w^{m_i}(Q)|)f_i\in\tX_{\bQ}^{\vee},
\\
\cutlog(Q)&=-f_0+\log(Q)=-\sum_{i=1}^g(\log|w^{m_i}(Q)|)f_i
\in X_{\bQ}^{\vee}. 
\end{align*} 
\end{defn}

\begin{defn}\label{defn:Sup_vee}
Let $Z$ be a subset of $\tP_{l,0}$. We define 
\begin{align*}
\Sup(Z)&=\{a\in X; Z\subset (W_{l,a})_0\}. 
\end{align*}
\end{defn}

\begin{lemma}\label{lemma:cutlog}Let 
$a\in X$, 
$u\in X^{\vee}$ and $Q\in\tP_{l}(R_{\Omega})$. 
Then 
\begin{enumerate}
\item\label{item:open nbd of q} 
if $Q_0\in(W_{l,a})_0$, 
then $Q\in W_{l,a}(R_{\Omega})$;
\item\label{item:cutlog}
$Q\in W_{l,a}(R_{\Omega})\Leftrightarrow 
\log(Q)\in\tau_{l,a}\Leftrightarrow \cutlog(Q)\in\Cut(\tau_{l,a})$;
\item\label{item:action of delta_u on cutlog} 
$\cutlog(\delta_u(Q))=e(\xi)u+\cutlog(Q)$;
\item\label{item:FPd Pl}
$\delta_u((W_{l,a})_0)\cap (W_{l,a})_0
=\emptyset$\ $(\forall a\in X)$ if 
$e(\xi)u\in X^{\vee}\setminus 2\Sigma^*_{l}$.  
\end{enumerate}
\end{lemma}
\begin{proof}
Suppose $Q_0\in(W_{l,a})_0$. 
Since $W_{l,a}$ is open in $\tP_{l}$, 
the morphism $Q:\Spec R_{\Omega}\to \tP_{l}$ 
factors through $W_{l,a}$, that is, 
$Q\in W_{l,a}(R_{\Omega})$. This proves (\ref{item:open nbd of q}). 
Recall $W_{l,a}=\Spec R[w^x; x\in\tau_{l,a}^{\vee}\cap\tX]$. 
Hence 
\begin{gather*}W_{l,a}(R_{\Omega})=
\{Q\in \tP(R_{\Omega}); 0\leq |w^{x}(Q)|\leq 1\ 
(\forall x\in\tau_{l,a}^{\vee}\cap\tX)\}. 
\end{gather*}
Hence $Q\in W_{l,a}(R_{\Omega})$ iff 
$-\log|w^{x}(Q)|=\langle x,\log(Q)\rangle\geq 0
\ (\forall x\in\tau_{l,a}^{\vee}\cap\tX)$. 
It follows that $Q\in W_{l,a}(R_{\Omega})$ iff $\log(Q)\in\tau_{l,a}$. 
This proves (\ref{item:cutlog}). 
Next (\ref{item:action of delta_u on cutlog}) follows from 
$\langle x,\log(\delta_u(Q))\rangle=-\log|w^x(\delta_u(Q)))|=\langle x,\log(Q)+e(\xi)u\rangle$.  
 \par
Finally let $q\in\delta_u((W_{l,a})_0)\cap (W_{l,a})_0$.
By \cite[8.3/8]{Bosch14}, \footnote{See Definition~\ref{defn:specialization}.} 
there exists $Q\in W_{l,a}(R_{\Omega})$ such that $Q_0=q$ and 
$\delta_u(Q)\in W_{l,a}(R_{\Omega})$. 
Hence $t^*:=\cutlog(Q)\in\Cut(\tau_{l,a})$ and $e(\xi)u+t^*\in\Cut(\tau_{l,a})$. It follows from Corollary~\ref{cor:Sk0 tau_l_Delta}
~(\ref{item:tau_l_a subset Sigma*_l}) that 
$e(\xi)u\in\Cut(\tau_{l,a})-\Cut(\tau_{l,a}) 
\subset 2\Sigma^*_{l}$. This proves (\ref{item:FPd Pl}). 
\end{proof}

\begin{lemma}\label{prop:P_eta_an}
The action of $Y$ on $\tP_{l,\eta}^{\an}$ through $(S_y;y\in Y)$ 
is free and proper, and  
$\tP_{l,\eta}^{\an}/Y\simeq P_{l,\eta}^{\an}=G_{l,\eta}^{\an}$. 
\end{lemma}
\begin{proof}See \cite[3.6]{Mumford72}. 
 By Definition~\ref{defn:twisted Mumford family}, 
$P_{l}^{\wedge}$ is the quotient of  
$\tP_{l}^{\wedge}$ by the free and proper action 
$(S_y^{\wedge}; y\in Y)$ of $Y$. 
Let $\varpi_{l}:\tP_{l}^{\wedge}\to P_{l}^{\wedge}$ 
be the quotient morphism 
and $\varpi_{l,\rig}:\tP_{l,\rig}^{\wedge}\to P_{l,\rig}^{\wedge}$ the natural morphism induced by $\varpi_{l}$.  
Since $S_y=\delta_{\beta(y)}$, 
the action $(S_{y,\rig}^{\wedge};y\in Y)$  of $Y$ 
on $|\tP_{l,\rig}^{\wedge}|\ (=|\tP_{l,\eta}^{\an}|)$ is free 
and proper by Lemma~\ref{lemma:cutlog}~(\ref{item:action of delta_u on cutlog})-(\ref{item:FPd Pl}) and 
$|\tP_{l,\rig}^{\wedge}|/Y=|P_{l,\rig}^{\wedge}|$.  
Hence $\tP_{l,\rig}^{\wedge}/Y\simeq P_{l,\rig}^{\wedge}$.
\footnote{This is proved by taking a suitable refinement of 
$W_{l,a,\rig}^{\wedge}$. We omit the details.  }  
By Eq.~(\ref{eq:tP_wedeg_rig=P_eta_an}), 
$\tP_{l,\rig}^{\wedge}=\tP_{l,\eta}^{\an}$, while 
$P_{l,\rig}^{\wedge}=P_{l,\eta}^{\an}$ by Proposition~\ref{prop:rigpt}~(\ref{item:Xwedge_rig subset Xeta_an}), 
so that $P_{l,\eta}^{\an}$ is the quotient of  
$\tP_{l,\eta}^{\an}$ by $Y$. 
\end{proof}

\subsection{Theorem, Notation and Definitions}
\label{subsec:theorems notation}

In the rest of this section, we prove the following:
\begin{thm}
\label{thm:action of Gl on Pl}
Assume (\ref{assump:Sigmal(0) integral}). Then  
there exists an $S$-morphism $G_{l}\times_SP_{l}\to P_{l}$ that 
extends the group law $G_{l,\eta}\times_\eta G_{l,\eta}\to G_{l,\eta}$ of $G_{l,\eta}$.
\end{thm}

Since $l$ is constant throughout this section, we set as follows:
\begin{gather*}\tG:=\tG_{l},\ \tP:=\tP_{l},\ G:=G_{l},\ P:=P_{l},\\  
W_{\Delta}:=W_{l,\Delta},\ B_{\Delta}:=B_{l,\Delta}\ (\Delta\in\Vor_{l}(\xi^{\natural})).
\end{gather*}

\begin{defn}\label{defn:Z_0 and Z_eta}
Let $Z\in\{\tG,G,\tP,P\}$. Then we define 
\begin{gather*}
Z_{12}=Z\times_SZ,\ 
Z_{123}=Z\times_SZ\times_SZ,\\
\hskip 2.5pt \tL_{12}=\tG\times_S\tP,\, 
\tL_{123}=\tG\times_S\tP\times_S\tP,\\
L_{12}=G\times_SP,\ 
L_{123}=G\times_SP\times_SP.
\end{gather*} For example, $\tP_{123}=\tP\times_S\tP\times_S\tP$.
\end{defn}\label{defn:tilde_rho}

\begin{defn}\label{defn:action of tG and Y}
We denote the action of $\tG$ on $\tP$ 
by $\tilde\rho:\tL_{12}\to\tP$.   The action 
$\tilde\rho$ satisfies the following condition:
\begin{equation}\label{eq:condition on trho}
\trho\circ (\id_{\tG}\times_SS_y)=S_y\circ\trho\  
(\forall y\in Y). \end{equation}   

We also denote 
the action of $X^{\vee}$ (resp. $Y$) on $\tP$ by $(\Phi(u)=\delta_u; u\in X^{\vee})$\ (resp. $(\Psi(y)=S_y; y\in Y)$) where $\Psi(y)=\Phi(\beta(y))$ 
$(\forall y\in Y)$. 
\end{defn}

\begin{defn}\label{defn:quasi prod}
Let $Q,Q'\in\tP_{l}(\Omega)=\tG_{l}(\Omega)$. 
Then $\trho(Q,Q')\in\tP_{l}(\Omega)=\tP_{l}(R_{\Omega})$, 
which itself is also an $R$-flat (hence admissible) formal $R$-scheme. 
Hence we have the {\it specialization} 
$\trho(Q,Q')_0\in\tP_{l,0}(\overline{k(0)})$. 
Note $\trho(S_yQ,Q')=S_y\trho(Q,Q')$ $(\forall y\in Y)$. 
\end{defn}

\begin{defn}\label{defn:Psi_123^I}
Let $i\in [1,2]$, $j\in [1,3]$ and $I$ any subset of $[1,3]$. 
Let $p_i:\tP_{12}\to \tP$ (resp. $q_j:\tP_{123}\to \tP$) be   
the $i$-th (resp. the $j$-th) projection.    
We also define $\Phi^{(i)}_{12}:X^{\vee}\to \Aut_S(\tP_{12})$ and 
$\Phi^{(I)}_{123}:X^{\vee}\to \Aut_S(\tP_{123})$ by 
\begin{gather*}
p_k\circ\Phi^{(i)}_{12}(u)=\begin{cases}\delta_u\circ p_k& (k=i)\\
\id_{\tP}\circ p_k& (k\neq i),\end{cases}\ 
q_k\circ\Phi^{(I)}_{123}(u)=\begin{cases}\delta_u\circ q_k& (k\in I)\\
\id_{\tP}\circ q_k& (k\not\in I)\end{cases}
\end{gather*}for $u\in X^{\vee}$.  
For $J=(12)$ or $(123)$, and 
any proper subset $I$ of $J$, we define $\Psi^{(I)}_{J}:Y\to \Aut_S(\tP_J)$ 
by $\Psi^{(I)}_{J}=\Phi^{(I)}_{J}\circ\beta$. 
For example, 
$$\Phi^{(23)}_{123}(u)=\id_{\tP}\times_S\, \delta_u\times_S \delta_u, \ \Psi^{(2)}_{12}(y)=\id_{\tP}\times_S S_y\ \ (\forall u\in X^{\vee}, \forall y\in Y).$$

In what follows,  
we use the same notation $\Phi^{(i)}_{12}$, $\Phi^{(I)}_{123}$ 
and $\Psi^{(I)}_{J}$ for the induced automorphisms 
on $I$-adic completions or rigid analytifications. 

We denote by $Y^n$ 
the product of $n$-copies of $Y$. 
\end{defn}

\subsection{Proof of Theorem~\ref{thm:action of Gl on Pl} - start}
Let $\trho$  be the action of $\tG$ on $\tP$ and 
$\tGamma:=(\id_{\tL_{12}},\trho):\tL_{12}\to \tL_{123}$ the graph of $\trho$.
Since $\tP$ is separated over $S$, the graph $\tilde\Gamma$ is a closed immersion, which defines a closed subscheme $\tW\ (\simeq \tL_{12})$ of $\tL_{123}$. 
Hence $\tL_{123,\eta}=\tP_{123,\eta}$ and 
$\tW\cap\tP_{123,\eta}=\tW\cap\tL_{123,\eta}=\tW_{\eta}$.
Since $\tP_{123}$ is locally noetherian, the immersion 
$j:\tW_{\eta}\to\tP_{123}$ is quasi-compact.
Therefore we can take the scheme-theoretic image $\tV$ of $j$ 
\cite[2.5, p.~55]{BLR90}, which is the closed subscheme 
of $\tP_{123}$ defined by the kernel of 
the homomorphism $j_*j^*:\cO_{\tP_{123}}\to j_*\cO_{\tW_{\eta}}$.
By \cite[IV$_2$, 2.8.2]{EGA}, $\tV$ is $S$-flat with 
$\tV_{\eta}=\tW_{\eta}$.  Since $\tW$ is $S$-flat, we have 
by \cite[IV$_2$, 2.8.5]{EGA} 
\begin{equation}\label{eq:tV_cap_tL123=tW}
\tV\cap\tL_{123}=\tW.
\end{equation}

Since $\tilde\rho$ is ($Y$-)equivariant with respect to $\Psi_{12}^{(2)}$ 
and $\Psi$ by Eq.~(\ref{eq:condition on trho}), 
$\tilde\rho$ descends to an action 
$\hat\rho:L_{12}^{\wedge}\to P^{\wedge}$ of $G^\wedge$ on $P^\wedge$. 
Let $q_\Psi:\tP^{\wedge}\to\tP^\wedge/Y=P^{\wedge}$ be 
the quotient by $Y$ via $\Psi$.
Let $\tilde\Lambda$ be the composite
\begin{equation}\label{eq:tilde_Lambda}
 \tilde\Lambda:\tL_{12}^\wedge\xrightarrow{\tilde\Gamma^\wedge}\tL_{123}^\wedge=\tL_{12}^\wedge\times_{S^{\wedge}}\tP^\wedge\xrightarrow{\id_{\tL_{12}^\wedge}\times_{S^{\wedge}} q_\Psi}
\tL_{12}^{\wedge}\times_{S^{\wedge}}P^{\wedge}
\end{equation} Since $\tGamma$ is 
equivariant with respect to $\Psi_{12}^{(2)}$ and $\Psi_{123}^{(23)}$, so is 
$\tilde\Lambda$ with respect to $\Psi_{12}^{(2)}$ and $\Psi_{12}^{(2)}\times_{S^\wedge}\id_{P^\wedge}$.
Since $\tG^{\wedge}\simeq G^{\wedge}$, we have    
$L_{12}^{\wedge}\simeq\tL_{12}^{\wedge}/Y$, so that 
$\tilde\Lambda$ induces a morphism $\hat\Gamma$: the graph of $\hat\rho$. 
Since $P^{\wedge}$ is separated over $S^\wedge$,  $\hat\Gamma$ is a closed immersion, which defines a formal closed subscheme $\mathcal{W}$ 
of $L_{123}^{\wedge}$:
\begin{equation}\label{eq:hatGamma}
\hat\Gamma:L_{12}^\wedge\simeq\tL_{12}^\wedge/Y\overset{\simeq}{\longrightarrow}\cW\hookrightarrow(\tL_{12}^\wedge\times_{S^{\wedge}} P^\wedge)/Y\simeq L_{123}^\wedge. 
\end{equation}  In the following, we algebraize $\hat\rho$ by ``compactifying'' $\mathcal{W}$ in $P_{123}^{\wedge}$, to be more precise, 
by constructing a closed 
subscheme  $W$ (resp. $V$) of $L_{123}$ (resp. $P_{123}$) such that 
$W=V\cap L_{123}$ and $W^{\wedge}=\cW$.

\begin{lemma}\label{lemma:intersection tV+y_cap_tL123}
Let $Z\in\{\tW_0,\tW^{\wedge}\}$ and 
$Z'\in\{\tV_0,\tV^{\wedge},\tV^{\wedge}_{\rig}\}$ . Then 
\begin{enumerate}
\item\label{item:Psi1ytV_cap_tL=Psi3-ytV} 
$\Psi_{123}^{(1)}(y)(\tV)\cap\tL_{123}=
\Psi_{123}^{(3)}(-y)(\tW)$\ $(\forall y\in Y)$;
\item\label{item:Phi3uZ_cap_Z=empty iff uneq0} $\Phi_{123}^{(3)}(u)(Z)\cap Z=\emptyset$ iff $u\in X^{\vee}\setminus(0)$;
\item\label{item:Phi3u_tV0_cap_tV0=nonempty finite u} 
$\Phi_{123}^{(i)}(u)(Z')\cap Z'=\emptyset$\ $(i=1,3)$\ 
if $e(\xi)u\in X^{\vee}\setminus 6\Sigma^*_{l}$.
\end{enumerate}
\end{lemma}
\begin{proof}
By Eq.~(\ref{eq:tV_cap_tL123=tW}), we obtain
$\tW=\tV\cap\tL_{123}=\Psi_{123}^{(13)}(y)(\tV)
\cap\Psi_{123}^{(3)}(y)\tL_{123}=\Psi_{123}^{(3)}(y)\left(
\Phi_{123}^{(1)}(y)(\tV)\cap\tL_{123}\right)$,
which proves (\ref{item:Psi1ytV_cap_tL=Psi3-ytV}).\par 
Next we prove (\ref{item:Phi3uZ_cap_Z=empty iff uneq0}). 
Suppose $(\Phi_{123}^{(3)}(u)(\tW_0)\cap\tW_0)(\overline{k(0)})\neq\emptyset$. 
Then there exist two closed points $(h_0,w_0,v_0)$ and $(g_0,z_0,x_0)$ of $\tW_0$ such that $(h_0,w_0,v_0)=(g_0,z_0,\delta_ux_0)$. Hence 
$x_0=\trho_0(g_0,z_0)=\trho_0(h_0,w_0)=v_0=\delta_ux_0$ 
where $\trho_0=\trho\times_S0$. 
Since $\delta_u$ has a fixed point $x_0$, we have $u=0$. 
This proves (\ref{item:Phi3uZ_cap_Z=empty iff uneq0}) 
for $Z=\tW_0$.  
The case $Z=\tW^{\wedge}$ follows from the case $Z=\tW_0$.\par
Finally we prove (\ref{item:Phi3u_tV0_cap_tV0=nonempty finite u}). 
Since the case $i=1$ is proved in the same manner as $i=3$, 
we shall first prove the case $i=3$ for $Z'=\tV_0$. 
Suppose $\Phi_{123}^{(3)}(u)(\tV_0)\cap\tV_0\neq\emptyset$ 
 for some $u\in X^{\vee}$. 
By \cite[8.3/8]{Bosch14}, we can find 
$Q:=(h,w,v),\ Q':=(g,z,x)\in\tV(R_{\Omega})$  
with $Q_0=\Phi_{123}^{(3)}(u)Q'_0$. 
Hence
$(h_0,w_0,v_0)=(g_0,z_0,\delta_ux_0)$, so that $(h_0,w_0)=(g_0,z_0)$, 
$v_0=\trho(h,w)_0$, $x_0=\trho(g,z)_0$ and $v_0=\delta_ux_0$. 
Let $a\in\Sup(h_0)=\Sup(g_0)$, $b\in\Sup(w_0)=\Sup(z_0)$, $d\in\Sup(v_0)$,  
$\Gamma^*:=\Cut(\tau_{l,a})+\Cut(\tau_{l,b})$ and 
$\Gamma^{**}:=\Cut(\tau_{l,d})$. 
We have
\begin{align*}
t_1^*&:=\cutlog(\trho(h,w))=\cutlog(h)+\cutlog(w)\in\Gamma^*,\\
t_2^*&:=\cutlog(\trho(g,z))=\cutlog(g)+\cutlog(z)\in\Gamma^*. 
\end{align*}  
By $\Sup(v_0)=\Sup(\delta_ux_0)$, we have    
$t_1^*\in \Gamma^{**}$ and 
$t_2^*\in -e(\xi)u+\Gamma^{**}$ 
in view of Lemma~\ref{lemma:cutlog}~(\ref{item:action of delta_u on cutlog}). 
 By Corollary~\ref{cor:Sk0 tau_l_Delta}
~(\ref{item:tau_l_a subset Sigma*_l})-(\ref{item:n_Sigma*_l}), 
$t_1^*-t_2^*\in e(\xi)u+\Gamma^{**}-\Gamma^{**}\subset e(\xi)u+2\Sigma^*_{l}$, so that $e(\xi)u\in 2\Sigma^*_{l}+t_1^*-t_2^*\subset 2\Sigma^*_{l}+\Gamma^*-\Gamma^*\subset 6\Sigma^*_{l}$. This proves (\ref{item:Phi3u_tV0_cap_tV0=nonempty finite u}) 
for $Z'=\tV_0$. The other cases of 
(\ref{item:Phi3u_tV0_cap_tV0=nonempty finite u}) follow from it. 
\end{proof}

\subsection{The condition $(\FP)_Y$}
We consider the condition:
\begin{equation*}
(\FP)_Y\qquad \Psi_{123}^{(i)}(y)(\tV_0)\cap\tV_0=\emptyset\ \ 
(i=1,3, \forall y\in Y\setminus(0)).
\end{equation*} 

Suppose that we have proved Theorem~\ref{thm:action of Gl on Pl} 
under $(\FP)_Y$. Then  
Theorem~\ref{thm:action of Gl on Pl} 
in the general case follows from it.  
Indeed,  
there exists $n\in\bN$ 
such that $n\beta(Y)\cap 6\Sigma^*_{l}=\emptyset$. 
Hence $(\FP)_{nY}$ is true 
by Lemma~\ref{lemma:intersection tV+y_cap_tL123}~(\ref{item:Phi3u_tV0_cap_tV0=nonempty finite u}), while the action of $nY$ on $P_{l}$ is free and proper 
by Lemmas~\ref{prop:P_eta_an}/\ref{lemma:cutlog}~
(\ref{item:action of delta_u on cutlog})\,-(\ref{item:FPd Pl}). 

 Then we construct $P^{(n)}_{l}$ and $G^{(n)}_{l}$ 
in a manner similar to Theorem~\ref{thm:summary of Pl Gl} 
by algebraizing the formal quotient $\tP^{\wedge}_{l}/nY$. 
Let $h:P^{(n)}_{l}\to P_{l}$ be the natural 
finite surjective morphism. Then we have $G_{l}=h(G^{(n)}_{l})$. 
By assumption, there exists an action $\rho^{(n)}$ 
of $G^{(n)}_{l}$ on $P^{(n)}_{l}$ extending the group law of 
$G^{(n)}_{l}$.   
This induces in the natural manner 
an action $\rho$ of $G_{l}$ on $P_{l}$ 
extending the group law of $G_{l}$. 
Thus it suffices to prove Theorem~\ref{thm:action of Gl on Pl} under 
$(\FP)_Y$. \par
{\it From now throughout}\, \S~\ref{sec:action of Gl}
{\it we assume $(\FP)_Y$.} 
Let $\tcV:=\bigcup_{y\in Y}\Psi_{123}^{(1)}(y)(\tV^\wedge)$.
By $(\FP)_Y$, $\tcV$ is an {\it admissible} 
formal closed subscheme of $\tP_{123}^{\wedge}$.

\begin{defn}\label{defn:cV and V}
Let $\tilde\iota:\tcV\to\tP_{123}^\wedge$ be the natural closed immersion.
Since $\tilde\iota$ is equivariant 
with respect to $\prod_{i\in[1,3]}\Psi_{123}^{(i)}$, 
$\tilde\iota$ induces a closed immersion 
$\tilde{\iota}_{\quot}: \tcV/Y^3\to
\tP_{123}^{\wedge}/Y^3=P_{123}^{\wedge}$, which defines a {\it formal closed subscheme $\cV$ of} $P_{123}^{\wedge}$.
We define $V$ to be a {\it unique closed subscheme of $P_{123}$} 
which algebraizes $\cV$ by \cite[III$_1$, 5.4.1/5.4.5]{EGA},  
 that is, $V^{\wedge}\simeq\cV$ over $P^{\wedge}_{123}$. \par
Let $F:=V\cap(D_{l}\times_SP_{l}\times_S P_{l})$ 
by Eq.~(\ref{defn:cGl=union of deltau Gl}). 
It is a closed $S$-subscheme of $V$ such that 
$V\setminus F=V\cap L_{123}$ because 
$L_{123}=(P_{l}\setminus D_{l})\times_{S}
P_{l}\times_{S}P_{l}$. 
\end{defn}

\subsection{Proof  continued}
\label{subsec:proof continued}
\begin{lemma}\label{lem:P} Let $W:=V\setminus F=V\cap L_{123}$. 
 Then the following are true:
\begin{enumerate}
\item\label{item:graph of action 0} 
$W_0$ is the graph of $\rho_0:=\rho\times_S0$ and  
$W_0=\cW_0\simeq G_0\times_0P_0$;
\item\label{item:graph of action eta} 
$W_{\eta}$ is the graph of the group law of $G_{\eta}$;
\item\label{item:integral graph} $W$ is an integral scheme.
\end{enumerate}
\end{lemma}
\begin{proof}By $(\FP)_Y$ and Lemma~\ref{lemma:intersection tV+y_cap_tL123}~(\ref{item:Psi1ytV_cap_tL=Psi3-ytV})-(\ref{item:Phi3uZ_cap_Z=empty iff uneq0}),   
$\cV\cap L^{\wedge}_{123}
\simeq(\tcV\cap\tL^{\wedge}_{123})/Y^3
=\coprod_{y\in Y}\Psi_{123}^{(3)}(-y)(\tW^{\wedge})/Y^3\simeq\tW^{\wedge}/Y^2\simeq\cW$, whence $W_0=V_0\setminus F_0= 
\cV_0\cap L^{\wedge}_{123,0}=\cW_0$. Since $\cW_0=\hat\Gamma_0(L_{12,0}^{\wedge})$ by Eq.~(\ref{eq:hatGamma}), this proves (\ref{item:graph of action 0}). \par
Next we shall  prove  (\ref{item:graph of action eta}). 
The formal scheme  
$\tV^{\wedge}$ is an admissible formal closed subscheme of 
$\tP_{123}^{\wedge}$ because $\tV^{\wedge}$ is flat over $S^{\wedge}$. 
Hence by Proposition~\ref{prop:rigpt}~(\ref{item:Xwedge_rig subset Xeta_an})-(\ref{item:|Xrig|=|X|rig}), 
$|\tV^{\wedge}_{\rig}|=|\tV^{\wedge}|_{\rig}\subset|\tV^{\an}_{\eta}|$.    
By Eq.~(\ref{eq:tP_wedeg_rig=P_eta_an}), 
we have 
$|\tV^{\wedge}_{\rig}|=|\tV^{\an}_{\eta}|\cap|\tP_{123,\rig}^{\wedge}|
=|\tV^{\an}_{\eta}|\cap|\tP_{123,\eta}^{\an}|
=|\tV^{\an}_{\eta}|$, so that $\tV^{\wedge}_{\rig}\simeq 
\tV^{\an}_{\eta}=\tW^{\an}_{\eta}$. Meanwhile 
since $\tV^{\wedge}_{\rig}$ is $\Psi^{(13)}_{123}$-stable, 
by $(\FP)_Y$, we obtain   
$\cV_{\rig}\simeq\tcV_{\rig}/Y^3
=\coprod_{y\in Y}\Psi_{123}^{(3)}(y)(\tV^{\wedge}_{\rig})/Y^3
\simeq\tV^{\wedge}_{\rig}/Y^2$. 
Since $V_{\eta}=(V\setminus F)_{\eta}=W_{\eta}$, 
we obtain 
\begin{align}\label{eq:Weta_ansimeq tWeta_an/YY}
W^{\an}_{\eta}&=V_{\eta}^{\an}\simeq   
V^{\wedge}_{\rig}\simeq\cV_{\rig}\simeq\tV^{\wedge}_{\rig}/Y^2\simeq\tV^{\an}_{\eta}/Y^2=\tW^{\an}_{\eta}/Y^2.  
\end{align} 

Let $\sigma:\tW^{\an}_{\eta}/Y^2\to W^{\an}_{\eta}$ 
be the isomorphism (\ref{eq:Weta_ansimeq tWeta_an/YY}). 
Since 
$V$ is a closed $S$-subscheme of $P_{123}$, 
$W_{\eta}^{\an}$ is a closed (rigid analytic) subspace 
of $L_{123,\eta}^{\an}$ by Eq.~(\ref{eq:Weta_ansimeq tWeta_an/YY}). 
Meanwhile, 
by Lemma~\ref{prop:P_eta_an}, 
we have an isomorphism $\tau:L_{12,\eta}^{\an}
\simeq P_{12,\eta}^{\an}
\simeq\tP_{12,\eta}^{\an}/Y^2\simeq\tL_{12,\eta}^{\an}/Y^2$.   

Let $\tp_{12}:\tL_{123}\to \tL_{12}$ (resp. $p_{12}:L_{123}\to L_{12}$)  be  
the projection to the first two factors of $\tL_{123}$  (resp. $L_{123}$).  
Since $\tP_{\eta}=\tG_{\eta}$ is a group $k(\eta)$-variety, we have 
$\tL_{12,\eta}=\tG_{\eta}\times_{\eta}\tG_{\eta}$ and   
$\tGamma_{\eta}:\tL_{12,\eta}\simeq\tW_{\eta}$, 
which induces an isomorphism 
$\tGamma_{\eta}^{\an}:\tL_{12,\eta}^{\an}\simeq\tW_{\eta}^{\an}$ with 
$\tp_{12,\eta}^{\an}\circ\tGamma_{\eta}^{\an}=\id_{\tL_{12,\eta}^{\an}}$.  
Since 
$\tGamma_{\eta}^{\an}$ is $Y^2$-equivariant with respect 
to $\Psi_{12}^{(1)}\times \Psi_{12}^{(2)}$ and 
$\Psi_{123}^{(13)}\times \Psi_{123}^{(23)}$ 
by Eq.~(\ref{eq:condition on trho}), $\tGamma_{\eta}^{\an}$ induces a morphism 
$\Gamma_{\an}:\tL_{12,\eta}^{\an}/Y^2
\xrightarrow{}\tL_{123,\eta}^{\an}/Y^3$.  
Meanwhile as $\trho_{\eta}^{\an}$  
is $Y$-equivariant with respect to $\Psi_{12}^{(2)}$ 
and $\Psi$,   
$\trho_{\eta}^{\an}$ induces a morphism 
$\rho_{\an}:\tL_{12,\eta}^{\an}/Y^2\to \tP_{\eta}^{\an}/Y$. 
Then $\Gamma_{\an}$ is the graph of $\rho_{\an}$. 
Hence  we have an isomorphism 
\begin{equation}
\label{eq:final isom}
\sigma\circ\Gamma_{\an}\circ\tau:L_{12,\eta}^{\an}\overset{\tau}{\to}\tL_{12,\eta}^{\an}/Y^2\overset{\Gamma_{\an}}{\to}\tW_{\eta}^{\an}/Y^2\overset{\sigma}{\to} W^{\an}_{\eta}\ (\subset L_{123,\eta}^{\an})
\end{equation} 
with $p_{12,\eta}^{\an}\circ\sigma\circ\Gamma_{\an}\circ\tau
=\id_{L_{12,\eta}^{\an}}$.  
It follows that $\sigma\circ\Gamma_{\an}\circ\tau$ is (regarded as)
a closed immersion of $L_{12,\eta}^{\an}$ into $L_{123,\eta}^{\an}$.  

In what follows, we regard $\sigma=\id_{W_{\eta}^{\an}}$ 
and $\tau=\id_{L_{12,\eta}^{\an}}$ for simplicity. 
 Let $q_3:L_{123}\to P$  
be the projection to the third factor, $p:=(p_{12})_{|W}$ and 
$q:=(q_3)_{|W}$. 
Since $L_{12,\eta}$ and $W_{\eta}$ are (regarded as)
closed subschemes of the proper $k(\eta)$-scheme $L_{123,\eta}$, 
we can apply Theorem~\ref{thm:Chow}~(\ref{item:Hom}). 
Since $p_{\eta}^{\an}=(p_{12,\eta}^{\an})_{|W^{\an}_{\eta}}$ 
 is an isomorphism by Eq.~(\ref{eq:final isom}), so is 
$p_{\eta}$. Since $\Gamma_{\an}:L_{12,\eta}^{\an}\to W_{\eta}^{\an}$ 
is an isomorphism with 
$p_{\eta}^{\an}\circ\Gamma_{\an}=\id_{L_{12,\eta}^{\an}}$, 
there exists 
an isomorphism $\Gamma:L_{12,\eta}\to W_{\eta}$ of $k(\eta)$-schemes  
such that $\Gamma^{\an}=\Gamma_{\an}$     
and $p_{\eta}\circ\Gamma=\id_{L_{12,\eta}}$.

Now we define a $k(\eta)$-morphism  
\begin{equation}\label{eq:rho_eta}\rho^{\eta}=q_{\eta}\circ \Gamma:L_{12,\eta}\ \to W_{\eta}\to P_{\eta} 
\end{equation}where 
$(\rho^{\eta})^{\an}=q_{\eta}^{\an}
\circ \Gamma^{\an}=q_{\eta}^{\an}\circ
\Gamma_{\an}=\trho_{\eta}^{\an}$ mod $Y^2$. Since $\trho_{\eta}$ is an action of $\tG_{\eta}$ on $\tP_{\eta}$, $(\rho^{\eta})^{\an}$ is an action of $G^{\an}_{\eta}$ on $P^{\an}_{\eta}=G^{\an}_{\eta}$, which is just the group law of $G^{\an}_{\eta}$. The group law of $G_{\eta}^{\an}$ 
yields that of $G_{\eta}$ by Theorem~\ref{thm:Chow}~(\ref{item:Hom}), 
whence $\rho^{\eta}:G_{12,\eta}=P_{12,\eta}\to P_{\eta}=G_{\eta}$ 
defines a group law of $G_{\eta}$. 
Since a group $k(\eta)$-scheme structure on $G_{\eta}$ is uniquely determined by the choice of a $k(\eta)$-rational point on $G_{\eta}$ as {\it unit}, the group law $\rho^{\eta}$ of $G_{\eta}$ is equal to that of the abelian variety $G_{\eta}$ over $k(\eta)$. Hence $W_{\eta}$ is the graph of 
the group law of $G_{\eta}$. This proves (\ref{item:graph of action eta}).

Let us prove (\ref{item:integral graph}).  
Let $w\in W_0$. Since $W_0=\cW_0$ by Lemma~\ref{lem:P},  
there exists a point $\tilde{w}$ of $\tW_0$ above $w$ 
by $\tW_0/Y=\tW_0^{\wedge}/Y\simeq \cW_0$.
Let $B:=\cO_{\tW,\tw}$ and $C:=\cO_{W,w}$. 
Since $Y$ acts freely on $\tW^{\wedge}$, we have 
$B^{\wedge}\simeq C^{\wedge}$.  
Since $\tW\simeq\tL_{12}$, $B$ is reduced, so that $C$ is reduced by Corollary~\ref{cor:if Chat=Bhat Rk Sk}. 
Hence $W$ is reduced at any point on $W_0$.  
Since $B$ is $R$-flat, so are $B^{\wedge}$ and $C^{\wedge}$.
Hence $C$ is $R$-flat, so that $W$ is $R$-flat at any point on $W_0$. Hence $W$ is the closure of $W_{\eta}$ in $L_{123}$ by \cite[IV$_2$, 2.8.5]{EGA}. Since $W_{\eta}$ is irreducible and reduced by $W_{\eta}\simeq L_{12,\eta}$,  
so is $W$. This proves (\ref{item:integral graph}).
\end{proof}

\subsection{Completion of the proof} 
Let $p=(p_{12})_{|W}:W\to L_{12}$ and $q=(q_3)_{|W}:W\to P$ as 
in \S~\ref{subsec:proof continued}. 
First we prove that $p$ is an $S$-isomorphism. 
By Lemma \ref{lem:P}, 
$p$ is a bijective birational morphism between integral schemes.
Since $L_{12}$ is normal, Zariski's main theorem shows 
that $p$ is an isomorphism \cite[III$_1$, 4.4.9]{EGA}. 
Since $W$ is a closed subscheme of $L_{123}$, 
this implies by (\ref{item:graph of action eta}) 
that $W$ is the graph of an $S$-morphism 
$\rho:=q\circ p^{-1}:L_{12}\to P$ extending 
the group law $\rho^{\eta}$ in Eq.~(\ref{eq:rho_eta}) of $G_{\eta}$ over $k(\eta)$, that is, $\rho_{\eta}=\rho^{\eta}$. Hence $\rho$ is an action of $G$ on $P$ which extends the group law $G_{\eta}$. This completes the proof of 
Theorem~\ref{thm:action of Gl on Pl}. 

\begin{rem}\label{rem:group scheme G by same proof}
If we take $\tG_J$ (resp. $G_J$) instead of $\tL_J$ (resp. $L_J$) 
in the above proof, we can define a group $S$-scheme structure of $G$. 
This implies that 
the action $\rho$ of $G$ on $P$ also extends the group law of $G$. 
There seems to be no available literature 
for Theorem~\ref{thm:action of Gl on Pl}. 
\end{rem}

\section{Proof of Theorem~\ref{thm:main thm} -- the split case}
\label{sec:proof of main thm}

\subsection{Our plan}
\label{subsec:our plan}
 
We prove Theorem~\ref{thm:main thm}   
in \S\S~\ref{sec:proof of main thm}\,-\ref{sec:nonsplit case}. 
Let $R$ be a CDVR, $k(\eta)$ its fraction field, $S:=\Spec R$ and 
$\eta$ (resp. $0$) the generic (resp. closed) point of $S$. 
We use 
Notation~\ref{notation:Notation Rinit Sinit keta=Kmin(xi)} freely. 
Therefore $k(\eta)$ (resp. $R$, $S$, $0$) 
in Theorem~\ref{thm:main thm} is denoted in 
\S\S~\ref{sec:proof of main thm}\,-\ref{sec:nonsplit case} 
by $k(\eta_{\init})$ (resp. $R_{\init}$, $S_{\init}$, $0_{\init}$). 
Let $(G,\cL)$ be a semiabelian $S_{\init}$-scheme 
with $\cL$ symmetric ample cubical invertible 
(by Remark~\ref{rem:symmetric L}). 
After a finite \'etale Galois base change $S'\to S_{\init}$ 
by \cite[X, \S~1]{SGA3}, we have a 
semiabelian $S'$-scheme $(G_{S'},\cL_{S'})$ with closed fiber 
$G_{0'}$ a split $k(0')$-torus where 
$R':=\Gamma(\cO_{S'})$ is a CDVR and 
$0'$ is the closed point of $S'$. 
Let $k(\eta')$ be the fraction field of $R'$. 
By Theorem~\ref{thm:FC datum}, $\FC(G_{S'},\cL_{S'})$ is a 
{\it totally degenerate} FC datum. 

Then we can find an eFC datum $(X,Y,a,b^e,A,B)$ 
by Lemma~\ref{lemma:be 2nd} if we enlarge $k(\eta')$ 
to a suitable normal extension $K'$ of $k(\eta_{\init})$, which we denote 
by $k(\eta_{\min})$ by abuse of notation. 
Let $R_{\min}$ be the integral closure of $R$ in $k(\eta_{\min})$ and 
$\eta_{\min}$ the generic point of $\Spec R_{\min}$. 

Let $k(\eta_{\spl})$ be the minimal {\it unramified  Galois} extension 
of $k(\eta_{\init})$ in $k(\eta_{\min})$  such that 
$G_{0_{\spl}}$ is a $k(0_{\spl})$-split torus,     
$R_{\spl}$ the integral closure of $R_{\init}$ in 
$k(\eta_{\spl})$ and $S_{\spl}:=\Spec R_{\spl}$ where $0_{\spl}$ is 
the closed point of $S_{\spl}$. Hence we have a sequence of extensions:
\begin{equation}\label{eq:field extensions}
k(\eta_{\min})\supset k(\eta_{\spl})
\supset k(\eta_{\init}).
\end{equation}

Let $\xi:=\FC(G_{S_{\spl}},\cL_{S_{\spl}})=(X,Y,a,b,A,B)$ and  
$\xi^e:=(X,Y,a,b^e,A,B)$. We denote $S_{\min}(\xi)$ 
(resp. $R_{\min}(\xi)$, $\eta_{\min}$) by $S$ 
(resp. $R$, $\eta$) in this section. 

Theorem~\ref{thm:summary of Pl Gl} constructs a quadruple 
$(P_{l},\cL_{l},G_{l},\cG_{l})$ over $R$. 
To prove Theorem~\ref{thm:main thm}, 
we descend the quadruple to $R_{\init}$ following 
Eq.~(\ref{eq:field extensions}). 
Therefore we suffice to consider the following two cases: 
\begin{enumerate}
\item[]Case 1.\ descent from $S_{\min}$ to $S_{\spl}$\ 
(\S~\ref{subsec:Proof in Case 1});
\item[]Case 2.\ descent from $S_{\spl}$ to $S_{\init}$ 
\ (\S~\ref{sec:nonsplit case}).
\end{enumerate}

\subsection{Galois descent}
\label{subsec:Galois descent}
Let $R$ be a CDVR, 
$K:=k(\eta)$ its fraction field, $L$ a finite 
Galois extension of $K$, $R_L$ the integral closure of $R$ 
in $L$ and $\Gamma:=\Gal(L/K)=\Gal(R_L/R)$. 
Hence $R_L$ is an $R$-free module of finite rank. 
Note that we do not assume  that $L$ is unramified over $K$. 
Our basic references for this subsection are 
\cite{Conrad}, \cite[Chap.~17]{W79} and \cite[Descent]{Stacks}. 
\begin{defn}
\label{defn:Gamma str}(\cite[\S~2]{Conrad})
Let $V$ be an $R_L$-module, $\sigma\in\Gamma$ and $\psi\in\Hom_R(V, V)$. 
We call $\psi$ {\it $\sigma$-semilinear} if $\psi(ax)=\sigma(a)\psi(x)$ $(\forall a\in R_L, \forall x\in V)$. 
A {\it $\Gamma$-structure of $V$} is defined to be the set 
$(\psi_{\sigma};\sigma\in\Gamma)$ such that 
\begin{enumerate}
\item[(i)] $\psi_{\sigma}\in\Aut_R(V)$ and $\psi_{\sigma}$ is $\sigma$-semilinear\ $(\forall\sigma\in\Gamma)$;
\item[(ii)]  $\psi_{\id_L}=\id_V$, 
$\psi_{\sigma}\psi_{\tau}=\psi_{\sigma\tau}$ 
$(\forall \sigma,\tau\in\Gamma)$.
\end{enumerate} 

The pair $(V,(\psi_{\sigma};\sigma\in\Gamma))$
 is called a {\it descent datum}, or we say that 
$(\psi_{\sigma};\sigma\in\Gamma)$ is a descent datum for $V$. 
\end{defn}

\begin{defn}\label{defn:Gamma str can}
Let $V$ be an $R_L$-module. 
An \textit{$R$-form of $V$} is defined to be an $R$-submodule $U$ of $V$ such that the $R$-homomorphism $U\otimes_RR_L\to V$ sending $U\otimes_RR_L\ni u\otimes a\mapsto au\in V$ 
is bijective.  
Then we define a $\Gamma$-structure 
$\underline{\can}=(\psi^{\can}_{\sigma}\in\Aut_R(V);\sigma\in\Gamma)$ 
of $V$ 
by 
$$\psi^{\can}_{\sigma}(au)=\sigma(a)u\ (\forall a\in R_L,\forall u\in U),$$
which we call the canonical descent datum of $U\otimes_RR_L$. 
\end{defn}

\begin{thm}\label{thm:Galois descent/RL}With the same notation as above, 
the following two categories are equivalent:
\begin{enumerate}
\item[(i)] the category of $R$-modules $U$;
\item[(ii)] the category of descent data $(V,(\psi_{\sigma};\sigma\in\Gamma))$ for $R_L$-modules $V$.
\end{enumerate}

The equivalence of the categories is given as follows:
\begin{align*}
(V,(\psi_{\sigma};\sigma\in\Gamma))&\mapsto U:=V^{\Gamma}:=
\{x\in V: \psi_{\sigma}(x)=x\ (\forall \sigma\in\Gamma)\},\\
U&\mapsto (U\otimes_RR_L,\underline{\can}).
\end{align*}
\end{thm}

In the above correspondence, the $R$-form 
$U=V^{\Gamma}$ of $V$ is called a {\it descent to $R$ of 
$(V,(\psi_{\sigma};\sigma\in\Gamma))$} 
or a {\it descent of $V$ to $R$} if no confusion is possible.  

A direct proof of Theorem~\ref{thm:Galois descent/RL} is given 
in \cite{Conrad} when $(R,R_L)=(K,L)$. 
Theorem~\ref{thm:Galois descent/RL} 
is a corollary of faithfully flat descent 
\cite[17.2, 17.5-17.7]{W79} in the general case when 
$L$ is a finite Galois extension of $K$.
See also \cite[Descent, tags~039W, ~023N, ~0D1V]{Stacks}.  
Since there is somewhat notational difference between them 
 and here, we explain about it.\par 

Let $A:=R_L$, $Y:=\Spec R$, $X:=\Spec A$, and let $\pi:X\to Y$ be 
the morphism induced from the natural inclusion $R\subset A$. 
Since $\pi$ is finite Galois surjective,  
we can apply \cite[17.2, 17.7]{W79} to $\pi:X\to Y$.
See also \cite[Descent, tag~0D1V]{Stacks}. 

\begin{defn}\label{defn:ssigma tsigma}
Let $\sigma,\tau\in\Gamma$. 
We define   
an $R$-homomorphism $s_{\sigma}^*:A\to A$ by  
$s_{\sigma}^*(a)=\sigma(a)$ $(\forall a\in A)$, 
and let 
$s_{\sigma}:X\to X$ be the $Y$-morphism 
induced from $s_{\sigma}^*$ \footnote{One also writes 
$s_{\sigma}=\Spec(\sigma)$.} 
and $t_{\sigma}:=s_{\sigma}^{-1}$.  
Then  $s_{\tau}\circ s_{\sigma}=s_{\sigma\tau}$ 
and $t_{\sigma}\circ t_{\tau}=t_{\sigma\tau}$. 
\end{defn}

Let $V$ be an $A$-module and  
$t_{\sigma}^*V:=V\otimes_{A_2;A_2\overset{t_{\sigma}^*}{\to}A_1}A_1$  
where $A_i:=A$ $(i\in [1,2])$. 
Let $\cF:=\widetilde{V}$ be 
the quasi-coherent $\cO_X$-module associated with $V$ 
and $t_{\sigma}^*\cF:=\cF\otimes_{\cO_X}t_{\sigma}^*(\cO_X)$. Then 
$V=H^0(X,\cF)$ and $t_{\sigma}^*V=H^0(X,t_{\sigma}^*\cF)$. 
The $A\ (=1\otimes A_1)$-action 
on $t_{\sigma}^*V$ is given by  
$a(v\otimes 1):=(1\otimes a)(v\otimes 1)=v\otimes a=\sigma(a)v\otimes 1$\ $(\forall a\in A, \forall v\in V).$  
Meanwhile, we define an $A$-module $(V_{\sigma},\cdot_{\sigma})$ by 
setting $V_{\sigma}=V$ (as an $R$-module) 
and an action $\cdot_{\sigma}$ of $A$ on $V$ by 
$a\cdot_{\sigma} v=\sigma(a)v$ 
\ $(\forall a\in A, \forall v\in V)$. 
 We denote $(V,\cdot_{\sigma})$ by  
$V_{\sigma}$ for brevity. When we write $av$, we mean $av=a\cdot_{\id_A} v$ 
in $V$ via $V_{\sigma}=V=V_{\id_A}$. 

Then we define an isomorphism of $A$-modules: 
\begin{equation}\label{eq:kappa}
\kappa_{\sigma}:t_{\sigma}^*V\simeq V_{\sigma}
\end{equation}
by $\kappa_{\sigma}(\sigma(a)v\otimes 1)=\kappa_{\sigma}(v\otimes a):=a\cdot_{\sigma}v=\sigma(a)v$ $(\forall a\in A_1, \forall v\in V)$. 
 
\begin{lemma}\label{lemma:phi_sigma_tau iff varphi_sigma_tau}
Let $\sigma\in\Gamma$, $\phi_{\sigma}\in\Hom_R(V,t_{\sigma}^*V)$, 
$\psi_{\sigma}:=\kappa_{\sigma}\circ\phi_{\sigma}$ and let 
$\varphi_{\sigma}:\cF\to t_{\sigma}^*\cF$ be 
the {\it $\pi^{-1}(\cO_Y)$-homomorphism} associated with $\phi_{\sigma}$. Then 
\begin{enumerate}
\item\label{item:semilinear} $\psi_{\sigma}$ is $\sigma$-semilinear 
$\Leftrightarrow\phi_{\sigma}\in\Hom_A(V,t_{\sigma}^*V)
\Leftrightarrow \varphi_{\sigma}\in
\Hom_{\cO_X}(\cF,t_{\sigma}^*\cF)$;
\item\label{item:cocycle}  
$\psi_{\sigma\tau}=\psi_{\sigma}\circ\psi_{\tau}
\Leftrightarrow\phi_{\sigma\tau}=t_{\tau}^*\phi_{\sigma}\circ\phi_{\tau}
\Leftrightarrow
\varphi_{\sigma\tau}
=t_{\tau}^*\varphi_{\sigma}\circ\varphi_{\tau}$ 
$(\forall\sigma,\tau\in\Gamma)$.
\end{enumerate}
\end{lemma}
\begin{proof}By (\ref{eq:kappa}),  
$\phi_{\sigma}(v)=\kappa_{\sigma}^{-1}\circ\psi_{\sigma}(v)
=\psi_{\sigma}(v)\otimes 1\ (\forall v\in V)$. Hence  
$\phi_{\sigma}\in\Hom_A(V,t_{\sigma}^*V)$ iff  
$\phi_{\sigma}(av)=\sigma(a)\phi_{\sigma}(v)$ $(\forall a\in A)$ 
iff $\psi_{\sigma}(av)=\sigma(a)\psi_{\sigma}(v)$ $(\forall a\in A)$,
which proves (\ref{item:semilinear}).   
Since $V_{\sigma}=V$,  we have 
$\psi_{\sigma}\in\Hom_A(V,V)$. 
If $\psi_{\sigma\tau}=\psi_{\sigma}\circ\psi_{\tau}$, then 
\begin{align*}
\phi_{\sigma\tau}&=\psi_{\sigma\tau}\otimes 1=(\psi_{\sigma}\circ\psi_{\tau})\otimes 1:V\overset{\psi_{\sigma\tau}}{\to} V\overset{\kappa_{\sigma\tau}^{-1}}{\to} t_{\sigma\tau}^*V\\
&=t_{\tau}^*(\psi_{\sigma}\otimes 1)\circ(\psi_{\tau}\otimes 1)=t_{\tau}^*\phi_{\sigma}\circ\phi_{\tau}:V\to t_{\tau}^*V\to 
t_{\tau}^*t_{\sigma}^*V.
\end{align*}
Thus we see 
$\psi_{\sigma\tau}=\psi_{\sigma}\circ\psi_{\tau}
 \Leftrightarrow 
\phi_{\sigma\tau}=t_{\tau}^*\phi_{\sigma}\circ\phi_{\tau}     
 \Leftrightarrow 
H^0(\varphi_{\sigma\tau})
=H^0(t_{\tau}^*\varphi_{\sigma}\circ\varphi_{\tau})
\Leftrightarrow \varphi_{\sigma\tau}
=t_{\tau}^*\varphi_{\sigma}\circ\varphi_{\tau}$
because $\phi_{\rho}=H^0(\varphi_{\rho})$ 
$(\forall\rho\in\Gamma)$. This proves (\ref{item:cocycle}).  
\end{proof}

\begin{rem}Let $\varphi'_{\sigma}:=s_{\sigma}^*(\varphi_{\sigma}^{-1})$.
Then $\varphi_{\sigma\tau}
=t_{\tau}^*\varphi_{\sigma}\circ\varphi_{\tau}\ 
(\forall \sigma,\tau\in\Gamma)\Leftrightarrow \varphi'_{\sigma\tau}=s_{\sigma}^*\varphi'_{\tau}\circ \varphi'_{\sigma}\ (\forall \sigma,\tau\in\Gamma)$. 
See \cite[0D1V]{Stacks} for the second equality. 
\end{rem}

\subsection{Quotients by $\bf\mu$}
\begin{defn}Let $A$ be a commutative ring and $S=\Spec A$. 
A finite {\it split multiplicative} group $S$-scheme 
$\bf\mu$ is a finite commutative group $S$-scheme 
which is the Cartier dual $D(Z)$ 
of some constant finite abelian group $S$-scheme $Z$,  
say $Z=\bigoplus_{i=1}^n(\bZ/e_i\bZ)_S$ for some $e_i\geq 1$.  
Let $B:=\Gamma(\cO_{\bf\mu})$. 
Then $B$ is a Hopf algebra over $A$ with 
$A$-algebra homomorphisms
$\Delta:B\to B\otimes_AB$ comultiplication, 
$\epsilon:B\to A$ counit and 
$\iota:B\to B$ coinverse such that 
\begin{gather*}
(\id_B\otimes\Delta)\circ\Delta=(\Delta\otimes\id_B)\circ\Delta,\\
(\epsilon\otimes\id_B)\circ\Delta=\id_B,\ \ 
(\iota,\id_B)\otimes\Delta=j\circ\epsilon 
\end{gather*}where $j:A\to B$ is the inclusion. We call $B$ 
the {\it Hopf algebra of} $\bf\mu$. 
\end{defn}
\begin{defn}\cite[1.6.2]{Mont93}\ \ 
An $A$-module $M$ is a left $B$-{\it comodule} 
if $M$ is $A$-flat and 
there is an $A$-homomorphism 
$\rho_M:M\to B\otimes_AM$ such that 
\begin{equation}\label{eq:M is B-comodule}
\begin{aligned}(\id_{B}\otimes\rho_M)\circ\rho_M
&=(\Delta\otimes\id_M)\circ\rho_M,\\
(\epsilon\otimes\id_M)\circ\rho_M&=1_A\otimes\id_M.
\end{aligned}
\end{equation}

We call  $\rho_M$ the {\it coaction of $B$ on $M$}, which we   
denote by $\rho^B_M$ if necessary.
We say that {\it $\mu$ acts on $M$} if $M$ is a left $B$-comodule. 
Note that $B$ is a left $B$-comodule with $\rho_B=\Delta$.
\end{defn}

\begin{lemma}\label{lemma:decomposition of M}Let 
$G$ be either a finite split multiplicative group $S$-scheme 
or a split $S$-torus. Let $B:=\Gamma(\cO_{G})$ and 
$D(G):=\Hom_S(G,\bG_{m,S})$. Then  
\begin{enumerate}
\item\label{item:B}
$B=\bigoplus_{\chi\in D(G)}Aw^{\chi}$ with 
$w^{\chi}\cdot w^{\chi'}=w^{\chi+\chi'}$\ 
$(\forall\chi,\chi'\in D(G))$;
\item\label{item:decompose M}for 
a left $B$-comodule $M$, 
there is a unique direct sum decomposition $M=\bigoplus_{\chi\in D(G)}M^{\chi}
$ where $M^{\chi}:=\{m\in M; \rho_M(m)=w^{\chi}\otimes m\}$. 
\end{enumerate}
\end{lemma}

We say that $m\in M$ is $\bf\mu$-{\it invariant} if $m\in M^{0}$. We denote  
$M^{\chi}$ (resp. $M^{0}$) by $M_B^{\chi}$ (resp. $M_B^{0}$ or 
$M^{\bf\mu\text{-}\inv}$) 
if necessary.

\begin{proof}(\ref{item:B}) is clear. We prove (\ref{item:decompose M}).
Let $m\in M$ and  $\rho_M(m)=\sum_{\chi\in D(G)}w^{\chi}\otimes m_{\chi}$ for some $m_{\chi}\in M$. 
By Eq.~(\ref{eq:M is B-comodule}), 
$m=(\epsilon\otimes\id_M)\circ\rho_M(m)=\sum_{\chi\in D(G)}m_{\chi}$ and 
$\rho_M(m_{\chi})=w^{\chi}\otimes m_{\chi}$ 
 $(\forall \chi\in D(G))$. Since $(w^{\chi};\chi\in D(G))$ 
is a free $A$-basis of $B$, it is obvious that  
$M^{\chi}\cap M^{\chi'}=0$ if $\chi\neq\chi'$, and hence 
$m_{\chi}$ is uniquely determined by $m$ and $\chi$.  
This proves Lemma. \end{proof}

\begin{defn}\label{defn:B-comodule C}
Let $q\in\bN$, $B:=\Gamma(\cO_{\mu_{q,S}})=A[w]/(w^q-1)$,  
$\beta\in A^{\times}$ and $C:=A[x]/(x^q-\beta)$. 
We define an $A$-homomorphism $\rho_C:C\to B\otimes_AC$ by 
$\rho_C(x)=w\otimes x$. Then  
$\rho_C$ is a coaction of $B$ on $C$. 
\end{defn}

\begin{lemma}\label{lemma:mu acts freely Smin}
Let $V:=\Spec\ C$. Then 
$\mu_{q,S}$ acts on $V$ freely in the sense that 
the natural homomorphism $f^*:C\otimes_AC\to B\otimes_AC$ 
sending $c\otimes d\mapsto \rho_C(c)\cdot(1\otimes d)$ $(c,d\in C)$ 
is surjective.
\end{lemma}
\begin{proof}Let $B=A[w^i; i\in[0,q-1]]/(w^q-1)$. 
We check 
$f^*(x\otimes x^{q-1}/\beta)=(w\otimes x)\cdot(1\otimes x^{q-1}/\beta)
=w\otimes 1$. Hence $f^*$ is surjective.
\end{proof}

\begin{cor}\label{cor:free action of mu on Smin schemes}
The finite split multiplicative group $S_{\min}$-scheme 
$\mu_{S_{\min}}$ acts freely 
(in the sense of Lemma~\ref{lemma:mu acts freely Smin}) on 
$P_{l_0l}$, $W_{l_0l,a}$, $\cG_{l_0l}$ and $G_{l_0l}$.  
\end{cor}

Hence \cite[\S~12, Thm.~1, p.~104]{Mumford12} is true for these 
$S_{\min}$-schemes. This is an analogue to 
Theorem~\ref{thm:Galois descent/RL}  
when $\Gamma$ is replaced by $\mu$.

\subsection{Proof of Theorem~\ref{thm:main thm} - Case 1.}
\label{subsec:Proof in Case 1}\qquad\par
\noindent{\it Step 1.}\ \ 
Let $(G,\cL)$ be a semiabelian $S_{\spl}$-scheme with $\cL$ 
symmetric ample cubical invertible.   
By definition, 
$G_{0_{\spl}}$ is a split $k(0_{\spl})$-torus:
\begin{equation}\label{eq:torus G0}G_{0_{\spl}}
=\Spec k(0_{\spl})[w^x; x\in X]
\end{equation} for some lattice $X$ of rank $g$. 
 Let $\cG$ be the N\'eron model of $G_{\eta_{\spl}}$. 
By Proposition~\ref{prop:semi abelian into neron}, we (can) 
identify $G$ with the identity component of $\cG$. \par
 Let $\xi:=\FC(G,\cL)=(X,Y,a,b,A,B)$. Recall  
$k(\eta)=k_{\min}(\xi)$, $R=R_{\min}(\xi)$ and $S=S_{\min}(\xi)$.  
By the proof of Lemma~\ref{lemma:be 2nd}, 
there exist a N\'eFC kit $\xi^{\natural}$ over $S$ 
of $\xi^e_{2N}$. Let 
$\tilde\Gamma:=\Aut(k(\eta)/k(\eta_{\spl}))
=\Aut(S/S_{\spl})$. If $\chara k(\eta)=0$, then 
$\tilde\Gamma=\Gal(S/S_{\spl})$.   
If $\chara k(\eta)=p>0$, then $\tilde\Gamma$ 
is a finite group $S_{\spl}$-scheme, and  
the identity 
component $\mu$ of $\tilde\Gamma$ is a normal subgroup $S_{\spl}$-scheme 
of $\tilde\Gamma$ such that 
$\tilde\Gamma/\mu\simeq\tilde\Gamma_{\red}$ : 
 the reduced part of $\tilde\Gamma$, which is 
a subgroup scheme of $\tilde\Gamma$. 
By the proof of 
Lemma~\ref{lemma:be 2nd} and \cite[4.24/4.25]{Nakamura24}, 
$\mu\simeq\prod_{i=1}^g\mu_{q_i,S}$ where 
$q_i$ is some power of $p$ dividing $e_i$ $(i\in[1,g])$, 
and  there exists a finite group $\Gamma$ 
such that $\tilde\Gamma_{\red}\simeq\Gamma_S$. 
See \cite[\S~4.6]{Nakamura24} for more detail. Regardless of $\chara k(\eta)$, 
$S/\tilde\Gamma\simeq S_{\spl}$ in any case. 

Let $\xi^{\natural}_{l}$\ $(l\geq 1)$ be the $l$-th N\'eFC kit over $S$
induced from $\xi^{\natural}$ in Definition~\ref{defn:xinatural_l} 
with $\xi^{\natural}_1=\xi^{\natural}$:
$$\xi^{\natural}_{l}=(X,Y,\epsilon_{l},b^e_{l},E_{l},\Sigma_{l}).
$$

By Lemma~\ref{lemma:constructing integral Sigma_n}, 
there exists $l_0\in\bN$ such that 
\begin{equation}
\label{eq:integral Sigmal(0)}
\text{$\Sigma_{l_0l}(0)$ is integral for any integer $l\geq 1$.}
\end{equation}

By Theorems~\ref{thm:summary of Pl Gl}/\ref{thm:action of Gl on Pl}, we obtain
\begin{claim}\label{claim:Pl0l Gl0l cGl0l GS cGS}
Let $S=S_{\min}(\xi)$ and 
$(P_{l_0l},\cL_{l_0l}):
=(P_{l_0l}(\xi^{\natural}),\cL_{l_0l}(\xi^{\natural}))$.
Then 
\begin{enumerate}
\item[(i')]\label{item:Cohen-Macaulay} 
$(P_{l_0l},\cL_{l_0l})$ is an irreducible flat 
projective Cohen-Macaulay $S$-scheme;
\item[(ii')]\label{item:Xvee/betaY}there exist 
open $S$-subschemes $G_{l_0l}$ and $\cG_{l_0l}$ of $P_{l_0l}$ 
such that $(G_{l_0l},\cL_{l_0l})\simeq (G_S,\cL^{2Nl_0l}_S)$, 
$\cG_{l_0l}\simeq\cG_S$ and  
 $(\cG_{l_0l,0}/G_{l_0l,0})(\overline{k(0)})\simeq X^{\vee}/\beta(Y)$;
\item[(iii')]\label{item:codim 2} $P_{l_0l}\setminus \cG_{l_0l}$ is 
of codimension two in $P_{l_0l}$;
\item[(iv')]\label{item:action of cGl on Pl}
 $\cG_{l_0l}$ acts on $P_{l_0l}$ extending 
multiplication of  $\cG_{l_0l}$;
\item[(v')]\label{item:cubical sheaves cLl}
 $(\cL_{l_0l})_{|\cG_{l_0l}}$ is 
ample cubical with $\cL_{l_0l,\eta}=\cL^{2Nl_0l}_{\eta}$.
\end{enumerate}
\end{claim}

{\it Step 2.}\ \  Let $s$ (resp. $t$) 
be a uniformizer of $R$ (resp. $R_{\spl}$) 
and $I:=tR_{\spl}$. 
Since $G_{0_{\spl}}$ is a split $k(0_{\spl})$-torus, 
$G^{\wedge}$ is a formal $S_{\spl}^{\wedge}$-torus: 
$G^{\wedge}=\Spf R_{\spl}[w^x; x\in X]^{I\op{-adic}}$.  
Since $R$ is finite flat over $R_{\spl}$ (hence $R_{\spl}$-free), 
$\Gamma(\cO_{G_{S}^{\wedge}})=\Gamma(\cO_{G^{\wedge}})\otimes_{R_{\spl}}R  
=R[w^x; x\in X]^{IR\op{-adic}}$. Let $M:=\Gamma(\cO_{G_{S}^{\wedge}})$ and 
$\rho_M$ the coaction of $\mu$ on $M$. 
The group $S$-scheme $\tilde\Gamma$ acts on $M$ as follows: 
\begin{gather*}
\rho_M(cw^x)=\rho_M(c)\cdot (1\otimes w^x)
\in\Gamma(\cO_{\bf\mu})\otimes_{R} M,\\
\sigma(cw^x)=\sigma(c)w^x\ 
(\forall \sigma\in\Gamma, \forall c\in R, \forall x\in X).
\end{gather*}

Let $V_{m,\spl}:=\Gamma(G_{\eta_{\spl}},\cL_{\eta_{\spl}}^{\otimes m})$ 
and $V_{m}:=\Gamma(G_{\eta},\cL_{\eta}^{\otimes m})=V_{m,\spl}\otimes_{R_{\spl}}R$ $(m\geq 1)$. Then $\tilde\Gamma$ acts 
on the $k(\eta)$-vector subspace $V_{m}$ 
of $\Gamma(\cO_{G^{\wedge}})\otimes_{R_{\spl}}k(\eta)$ 
via  Eq.~(\ref{eq:Gamma(Geta,Leta) as subspace}). Hence 
by Lemma~\ref{lemma:H0(G,Lm) over keta as Fourier},
\begin{equation}\label{eq:subspce of invariants of Vm}  
V_{m}^{\tilde\Gamma}=V_{m,\spl}
\end{equation}

Let $M_{l_0l}:=R[X][\theta_{l_0l}]$. 
We define the action of $\tilde\Gamma$ on $M_{l_0l}$  by 
$\rho_{M_{l_0l}}(cw^x\theta_{l_0l})
=\rho_M(cw^x)(1\otimes\theta_{l_0l})$ and 
$\psi_{\sigma}(cw^x\theta_{l_0l})=\sigma(c)w^x\theta_{l_0l}$ 
$(\forall c\in R, \forall x\in X)$. 
Recall 
\begin{gather*}
R_{l_0l}(\xi^{\natural})=R(\xi^{\natural}_{l_0l})
=R[\xi_{l_0l,\alpha,u}\theta_{l_0l};
\alpha\in\Sigma_{l_0l},u\in X^{\vee}],\\
\ \ \xi_{l_0l,\alpha,u}=\ep_{l}(u)b^e(u,\alpha)
w^{\alpha+2l_0l\mu_{l_0l}(u)}.
\end{gather*}
Then 
$R_{l_0l}(\xi^{\natural})$ is an $R$-subalgebra of 
$k(\eta)[X][\theta_{l_0l}]$ 
stable under $\tilde\Gamma$ because   
$\psi_{\sigma}(\xi_{l_0l,\alpha,u})$ is an $R^{\times}$-multiple 
of $\xi_{l_0l,\alpha,u}$ and $\rho_{M_{l_0l}}(R_{l_0l}(\xi^{\natural}))
\subset\Gamma(\cO_{\bf\mu})\otimes_{R_{\spl}} R_{l_0l}(\xi^{\natural})$. 
Now we define 
an $R_{\spl}$-subalgebra of $k(\eta_{\spl})[X][\theta_{l_0l}]$ by 
\begin{equation}\label{eq:R form of R'l_xinatural}
R_{l_0l,\spl}(\xi^{\natural})
=R_{\spl}[\zeta_{l_0l,x}\theta_{l_0l};x\in X],\ 
\zeta_{l_0l,x}:=t^{D_{l_0l}(x)}w^{x}.
\end{equation}  

By Definition~\ref{defn:Dl(x)},     
$\xi_{l_0l,\alpha,u}$ 
is an $R^{\times}$-multiple of 
$\zeta_{l_0l,x}$ 
if $x=\alpha+2l_0l\mu(u)$ for some 
$(\alpha,u)\in\Sigma_{l_0l}\times X^{\vee}$.  
Since $\zeta_{l_0l,x}$ is $\tilde\Gamma$-invariant,  
$R_{l_0l,\spl}(\xi^{\natural})$ is 
an $R_{\spl}$-subalgebra of $R_{l_0l}(\xi^{\natural})$ 
such that 
\begin{equation*}\label{eq:R form of R'l_xinatural (2)}
R_{l_0l}(\xi^{\natural})
=R_{l_0l,\spl}(\xi^{\natural})\otimes_{R_{\spl}}R,\  
R_{l_0l}(\xi^{\natural})^{\tilde\Gamma}=R_{l_0l,\spl}(\xi^{\natural}).
\end{equation*}
 
We see that $S_y^*$ acts on $R_{l_0l,\spl}(\xi^{\natural})$, while 
$\delta_u^*$ does not keep $R_{l_0l,\spl}(\xi^{\natural})$ stable in general. 
Let $\tQ^*_{l_0l}:=\Proj R_{l_0l,\spl}(\xi^{\natural})$  and 
$(\tP^*_{l_0l},\tcL^*_{l_0l})$ be 
the normalization of $\tQ^*_{l_0l}$ with  
$\tcL^*_{l_0l}$ the pullback of $\cO_{\tQ^*_{l_0l}}(1)$ 
to  $\tP^*_{l_0l}$. Let 
\begin{equation*}\label{eq:tC*l tD*l}
\begin{aligned}
W_{l_0l,a}&:=\Spec R[\tau_{l_0l,a}^{\vee}\cap \tX],\\  
W^*_{l_0l,a}&:=W_{l_0l,a}/\tilde\Gamma
=\Spec R_{\spl}[\tau_{l_0l,a}^{\vee}\cap \tX],\\
\tC^*_{l_0l}&:=(\tP^*_{l_0l}\setminus \bigcup_{u\in X^{\vee}}
W^*_{l_0l,2l_0l\mu(u)})_{\red},\\
\tD^*_{l_0l}&:=(\tP^*_{l_0l}\setminus 
\bigcup_{y\in Y}W^*_{l_0l,2Nl_0l y})_{\red}. 
\end{aligned}
\end{equation*}

By Corollary~\ref{cor:Wlalphau as covering}, 
$(W_{l_0l,a}; a\in X)$ (resp. $(W^*_{l_0l,a}; a\in X)$) 
is an open covering of $\tP_{l_0l}$ (resp. $\tP^*_{l_0l}$). 
Let $(P^*_{l_0l},\cL^*_{l_0l})$ be the algebraization 
of the formal quotient 
$(\tP_{l_0l}^{*\wedge},\tcL_{l_0l}^{*\wedge})/Y$ 
and $C^*_{l_0l}$ (resp. $D^*_{l_0l}$)  
the algebraization 
of the formal quotients $\tC_{l_0l}^{*\wedge}/Y$ (resp. 
$\tD_{l_0l}^{*\wedge}/Y$). We also define  
\begin{equation}
\label{defn:Gl* cG*l}
\cG^*_{l_0l}=P^*_{l_0l}\setminus C^*_{l_0l},\quad 
G^*_{l_0l}=P^*_{l_0l}\setminus D^*_{l_0l}.
\end{equation}

{\it Step 3.}\ \ Finally we prove that  $(P^*_{l_0l},\cL^*_{l_0l},G^*_{l_0l},\cG^*_{l_0l})$ is the desired descent of 
$(P_{l_0l},\cL_{l_0l},G_{l_0l},\cG_{l_0l})$ 
to $S_{\spl}$. 
Note first that $P^*_{l_0l}$ (resp. $G^*_{l_0l}$, $\cG^*_{l_0l}$) 
is uniquely determined by $\xi$. Moreover $G^*_{l_0l}$ and 
$\cG^*_{l_0l}$ are group $S_{\spl}$-schemes, as is proved in parallel to 
Theorem~\ref{thm:action of Gl on Pl} 
by Remark~\ref{rem:group scheme G by same proof}.

\begin{claim}\label{claim:P*l}The following are true over $S_{\spl}$:
\begin{enumerate}
\item[(i$^*$)] $(P^*_{l_0l},\cL^*_{l_0l})$ is an irreducible flat 
projective Cohen-Macaulay $S_{\spl}$-scheme;
\item[(ii$^*$)]  
$(G^*_{l_0l},\cL^*_{l_0l})\simeq 
(G,\cL^{\otimes 2Nl_0l})$ and $\cG^*_{l_0l}\simeq\cG$; 
\item[(iii$^*$)] $P^*_{l_0l}\setminus\cG^*_{l_0l}$ is of codimension two;
\item[(iv$^*$)] $\cG^*_{l_0l}$ acts on $P^*_{l_0l}$ 
extending multiplication of $\cG^*_{l_0l}$;
\item[(v$^*$)]
 $(\cL^*_{l_0l})_{|\cG^*_{l_0l}}$ is 
ample cubical with $\cL^*_{l_0l,\eta_{\spl}}=\cL^{2Nl_0l}_{\eta_{\spl}}$.
\end{enumerate}
\end{claim}
\begin{proof}There is a finite morphism 
$\pi:(P_{l_0l},\cL_{l_0l})\to (P^*_{l_0l},\cL^*_{l_0l})$ 
by the above construction such that 
$W_{l_0l,a}^{\wedge}=(\pi^{\wedge})^{-1}(W^{*\wedge}_{l_0l,a})$\ 
 $(\forall a\in X)$,  
so that $G_{l_0l}=\pi^{-1}(G^*_{l_0l})$, 
$\cG_{l_0l}=\pi^{-1}(\cG^*_{l_0l})$ and 
$P_{l_0l}\setminus\cG_{l_0l}
=\pi^{-1}(P^*_{l_0l}\setminus\cG^*_{l_0l})$. 
Hence (i$^*$) (resp. (iii$^*$)) follows 
from Claim~\ref{claim:Pl0l Gl0l cGl0l GS cGS}~(i') (resp. (iii')). 

 Next we prove (ii$^*$). By Eq.~(\ref{eq:subspce of invariants of Vm}), 
\begin{equation}\label{eq:Gamma(G* cL*)=Gamma(G,cL)_Gamma}
\Gamma(G^*_{l_0l,\eta_{\spl}},\cL^{*}_{l_0l,\eta_{\spl}})
=\Gamma(G_{l_0l,\eta},\cL_{l_0l,\eta})^{\tilde\Gamma}
=V_{2Nl_0l}^{\tilde\Gamma}=V_{2Nl_0l,\spl}.
\end{equation}  Hence  
$\FC(G^*_{l_0l},\cL^*_{l_0l})=\FC(G,\cL^{\otimes 2Nl_0l})$ 
over $S_{\spl}$ 
by Lemma~\ref{lemma:H0(G,Lm) over keta as Fourier}.  
By Theorem~\ref{thm:G isom G' iff FC same},   
$(G^*_{l_0l},\cL^*_{l_0l})\simeq 
(G,\cL^{\otimes 2Nl_0l})$. 
By Eq.~(\ref{defn:Gl* cG*l}), we also see 
 {\small\begin{align*}
&\cG^*_{l_0l,0_{\spl}}/G^*_{l_0l,0_{\spl}}
\simeq\coprod_{u\in X^{\vee}}(W^*_{l_0l,2l_0l\mu(u)})_{0_{\spl}}/
\coprod_{y\in Y}(W^*_{l_0l,2Nl_0l y})_{0_{\spl}}\simeq X^{\vee}/\beta(Y)
\end{align*}}where 
$W^*_{l_0l,2l_0l\mu(u)}=\Spec R_{\spl}[t^{u(x)}w^x; x\in X]
\simeq T_{X,R_{\spl}}$. It follows from Theorem~\ref{thm:summary of Pl Gl}~(\ref{item:conn neron model})/(\ref{item:pullback conn neron})/(\ref{item:pullback of neron over Rinit}) that  
\begin{align*}
(\cG^*_{l_0l,0_{\spl}}/G^*_{l_0l,0_{\spl}})(\overline{k(0_{\spl})})
&\simeq (\cG_{l_0l,0}/G_{l_0l,0})(\overline{k(0)})
\simeq (\cG_{0_{\spl}}/G_{0_{\spl}})(\overline{k(0_{\spl})}).
\end{align*}  

Hence $\cG^*_{l_0l}$ is the N\'eron model 
of $G^*_{l_0l,\eta_{\spl}}\simeq G_{\eta_{\spl}}$, 
so that $\cG^*_{l_0l}\simeq\cG$. This proves (ii$^*$). 
(iv$^*$) is proved in parallel 
to Theorem~\ref{thm:action of Gl on Pl}. 
Since $\cL_{l_0l}$ is cubical, 
we have a unique (rigidified) global nonvanishing section 
$\tau\in\Gamma(\cG_{l_0l}^3,\Theta(\cL_{l_0l})^{\otimes(-1)})$, 
which is $\tilde\Gamma$-invariant, so that 
it descends to a global nonvanishing section 
$\tau_{\spl}\in\Gamma((\cG^*_{l_0l})^3,
\Theta(\cL^*_{l_0l})^{\otimes(-1)})$.
Since $\lambda(\cL^*_{l_0l})$ is symmetric, 
$\cL^*_{l_0l}$ is cubical on $\cG^*_{l_0l}$. The rest is clear. 
This proves (v$^*$).
\end{proof}

This completes the proof of Theorem~\ref{thm:main thm} 
in Case~1.

\section{Proof of Theorem~\ref{thm:main thm} -- the non-split case}
\label{sec:nonsplit case}
In this section, we prove Theorem~\ref{thm:main thm} 
in Case~2. 
when $G_{0_{\init}}$ is a non-split $k(0_{\init})$-torus. 
We denote $Z_{\spl}$ and $Z_{\init}$ 
by $Z^*$ and $Z$ respectively for $Z\in\{S,R,\eta,0\}$, 
{\it e.g.}, $S^*=S_{\spl}$ and $S=S_{\init}$.

\subsection{The Hopf algebra $\Gamma(\cO_{G^{\wedge}})$}
\label{subsec:Hopf algebra}
Let $R$ be a CDVR with $s$ uniformizer, $I:=sR$, $S:=\Spec R$ and  
$G$ an $S$-flat group $S$-scheme of finite type 
such that $G_0$ is an affine $k(0)$-scheme. 
Let $S_n:=\Spec R/I^{n+1}$ and $G_n:=G\times_SS_n$ $(n\geq 0)$.   
Since $G_0$ is affine, so is $G_n$  
by \cite[Rem.~1.1]{Artin76}. 
\footnote{The proof of \cite[Rem.~1.1]{Artin76} 
is applied to our case as well.}  
Hence $G_n$ is an affine group $S_n$-scheme.  
Let $\Delta_n$ (resp. $\epsilon_n$ and $\iota_n$)  
be comultiplication (resp. counit and coinverse) of $G_n$.  
 Since $G_n\simeq G_{n+1}\times_{S_{n+1}}S_n$ $(\forall n\geq 0)$, 
we have a projective system 
$(\Gamma(\cO_{G_n}),\Delta_n, \epsilon_n, \iota_n;n\geq 0)$ of 
Hopf algebras. 
Then 
$\Gamma(\cO_{G^{\wedge}})=\projlim\Gamma(\cO_{G_n})$, which  
is a Hopf algebra over $R$ whose comultiplication $\Delta$  
(resp. counit $\epsilon$, coinverse $\iota$) is the projective limit of 
$\Delta_n$ (resp. $\epsilon_n$, $\iota_n$). 
 
Let $k(\eta^*)$ be a finite Galois extension 
of $k(\eta)$, $R^*$ the integral closure of $R$ in $k(\eta^*)$, 
$S^*:=\Spec R^*$,  
 $\varpi:S^*\to S$ the associated finite 
flat morphism, 
$I^*:=IR^*$, 
$0^*$ the closed point of $S^*$, 
$\eta^*$ the generic point of $S^*$ and  
$k(0^*)$ the residue field of $R^*$. 
Let $G_{S^*}:=G\times_S{S^*}$, and  
let   
$G_{S^*}^{\wedge}$ be the $I^*$-adic completion of $G_{S^*}$. 
Then $\Gamma(\cO_{G_{S^*}^{\wedge}})$ is 
the $I^*$-adic completion of 
$\Gamma(\cO_{G^{\wedge}})\otimes_RR^*$. 
Since $R^*$ is finite and flat over $R$, 
$R^*$ is a finite $R$-free module. Hence  
$\Gamma(\cO_{G_{S^*}^{\wedge}})
=\Gamma(\cO_{G^{\wedge}})\otimes_RR^*$, 
which  is also a Hopf algebra over $R^*$, 
whose comultiplication $\Delta^{*}$, counit $\epsilon^{*}$ 
and coinverse $\iota^{*}$ are given by the pullbacks: 
$\Delta^{*}=\Delta\otimes_R\id_{R^*}$,  
$\epsilon^{*}=\epsilon\otimes_R\id_{R^*}$ and   
$\iota^{*}=\iota\otimes_R\id_{R^*}$.

\subsection{The Galois action}
\label{subsec:Galois action}
We restart from \S~\ref{subsec:our plan} Case~2. Now we set 
$S^*=S_{\spl}$ and $S=S_{\init}$ to simplify the notation. 
  Let $(G,\cL)$ be a semiabelian $S$-scheme 
with $\cL$ symmetric and rigidified such that $G_0$ is a $k(0)$-torus. 
Let $\lambda(\cL_{\eta}):G_{\eta}\to G_{\eta}^t$ be 
the polarization morphism, 
$K(\cL_{\eta}):=\ker(\lambda(\cL_{\eta}))$ and $K(\cL)$ 
the scheme-theoretic closure of $K(\cL_{\eta})$ in $G$.   
 By \cite[IV, 2.4]{MB85},  
$K(\cL)$ is a closed subgroup scheme of $G$ flat quasi-finite over $S$. 
Let $G^t:=G/K(\cL)$ \cite[p.~35]{FC90},  
which we call the dual of $G$. 
\par
Let $\varpi:S^*\to S$ be the natural morphism. It is  
a finite \'etale Galois cover 
such that $G_{0^*}:=G_{S^*}\times_{S^*}0^*$ is a split $k(0^*)$-torus:
\begin{equation}\label{eq:torus G'0}G_{0^*}=\Spec k(0^*)[w^x; x\in X]
\end{equation} for some lattice $X$ of rank $g$. 
By our choice of $k(\eta^*)=k(\eta_{\spl})$ in \S~\ref{subsec:our plan}, 
\begin{equation}
\label{eq:minimality assumption}
\text{$S^*$ is the minimal \'etale Galois cover of $S$ satisfying 
(\ref{eq:torus G'0}).} 
\end{equation}

Since $G_{0^*}$ is a split $k(0^*)$-torus, so is 
the closed fiber $G^t_{0^*}$:
\begin{equation*}\label{eq:torus G't0}
G^t_{0^*}=\Spec k(0^*)[w^y; y\in Y]
\end{equation*}for a sublattice $Y$ of $X$ of finite index  
by \cite[p.~35]{FC90}. \par

 Let $\Gamma:=\Gal(S^*/S)$. Since $\varpi$ is \'etale, 
$s$ is a uniformizer of $R^*$, so that 
$k(0^*)=R^*/sR^*=(R/sR)\otimes_RR^*=k(0)\otimes_RR^*$. 
Hence $\dim_{k(0)}k(0^*)=\rank_RR^*=\dim_{k(\eta)}k(\eta^*)$.  
By Hensel's lemma, 
\begin{equation}\label{eq:Hensel}
\Gal(k(0^*)/k(0))=\Gal(S^*/S)=\Gal(k(\eta^*)/k(\eta)). 
\end{equation}

\begin{defn}\label{defn:phi_sigma f_sigma} 
Let $\sigma\in\Gamma$. 
We define an $R$-homomorphism $s_{\sigma}^*:R^*\to R^*$ by  
$s_{\sigma}^*(a)=\sigma(a)$ $(\forall a\in R^*)$, 
and let $s_{\sigma}:S^*\to S^*$  the $S$-morphism 
induced from $s_{\sigma}^*$. 
Since $\sigma$ acts on $S^*$, so does $\sigma$ on $G_{S^*}$ and the $I^*$-adic completion $G_{S^*}^{\wedge}$ of $G_{S^*}$ where
$I^*=IR^*$.  We define 
$\psi_{\sigma}\in\Aut_R(\Gamma(\cO_{G_{S^*}^{\wedge}}))$ 
by $\psi_{\sigma}(z\otimes a)=z\otimes\sigma(a)$\ 
$(\forall a\in R^*,\forall z\in\Gamma(\cO_{G^{\wedge}}))$ via the 
identification $\Gamma(\cO_{G_{S^*}^{\wedge}})
=\Gamma(\cO_{G^{\wedge}})\otimes_RR^*$.  Then $\psi_{\sigma}$ is 
$\sigma$-{\it semilinear}   
and $\psi_{\sigma\tau}=\psi_{\sigma}\circ\psi_{\tau}$ $(\forall \sigma,\tau\in\Gamma)$. Thus  $(\psi_{\sigma};\sigma\in\Gamma)$ is a $\Gamma$-structure of 
$\Gamma(\cO_{G_{S^*}^{\wedge}})$ and $\Gamma(\cO_{G_{S^*}^{\wedge}})^{\Gamma}
=\Gamma(\cO_{G^{\wedge}})$. 
Let $\cL_{S^*}$ be the pullback of $\cL$ to $G_{S^*}$. Then  
 $\Gamma$ 
acts on $\Gamma(G_{S^*},\cL_{S^*}^{\otimes m})$ $(\forall m\geq 1)$,  
so that $\Gamma(G_{S^*},\cL_{S^*}^{\otimes m})\simeq \Gamma(G,\cL^{\otimes m})\otimes_RR^*$ as $\Gamma$-modules.  
\end{defn}

\begin{lemma}
\label{lemma:nomal form of Gamma action}
The following are true: 
\begin{enumerate}
\item\label{item:Hom(Gamma,Aut)} 
there exists $\lambda\in\Hom(\Gamma,\Aut_{\bZ}(X))$ 
such that $\lambda_{\sigma}(Y)=Y$ and 
$\psi_{\sigma}(w^x)=w^{\lambda_{\sigma}(x)}$\ 
$(\forall x\in X, \forall \sigma\in\Gamma)$  
where $\lambda_{\sigma}:=\lambda(\sigma)$;
\item\label{item:Gamma effective}
$\Gamma$ acts effectively on $X$ via $(\lambda_{\sigma};\sigma\in\Gamma)$.   
\end{enumerate}
\end{lemma}
\begin{proof}By Eq.~(\ref{eq:torus G'0}), 
$G_{S^*}^{\wedge}$ is a formal $R^*$-torus, 
so that $\Gamma(\cO_{G_{S^*}^{\wedge}})=R^*[w^x; x\in X]^{I^*\op{-adic}}$. 
Let  $\cB:=\Gamma(\cO_{G^{\wedge}})$ and 
$\cB^*:=\Gamma(\cO_{G_{S^*}^{\wedge}})=\cB\otimes_RR^*$.   
Recall that $\Delta^*$ (resp. $\epsilon^*$ 
and $\iota^*$) is comultiplication  
(resp. counit and coinverse) of $\cB^*$. 
Let $\sigma\in\Gamma$. Since $\psi_{\sigma}$ commutes 
with $\Delta^*$, $\epsilon^*$  
and $\iota^*$, we obtain
$\Delta^*(\psi_{\sigma}(w^x))
=(\psi_{\sigma}\otimes\psi_{\sigma})(w^x\otimes w^x)$, 
$\epsilon^*(\psi_{\sigma}(w^x))=1$  and 
$\iota^*(\psi_{\sigma}(w^x))=\psi_{\sigma}(w^{-x})$ $(\forall x\in X)$.

Since $\cB^*$ is an $R^*$-free module with basis $(w^y;y\in X)$, 
there is a unique expression $\psi_{\sigma}(w^x)=\sum_{y\in X}c_yw^y$ for
some $c_y\in R^*$, so that $\Delta^*(\psi_{\sigma}(w^x))=\sum_{y\in X}c_yw^y\otimes w^y$. Hence 
\begin{align*}
\sum_{y\in X}c_yw^y\otimes w^y&=
(\psi_{\sigma}\otimes\psi_{\sigma})(w^x\otimes w^x)=\sum_{y,z\in X}c_yc_zw^y\otimes w^z. 
\end{align*}  
Hence $c_y^2=c_y$ and $c_yc_z=0$\ $(\forall y,z\in X, y\neq z)$
because $(w^y\otimes w^z;y,z\in X)$ 
is a free $R^*$-basis of $\cB^*\otimes_{R^*} \cB^*$. Hence 
$\psi_{\sigma}(w^x)=w^y$ for a unique $y\in X$, which we denote 
by $\lambda_{\sigma}(x)$. Hence 
$\psi_{\sigma}(w^x)=w^{\lambda_{\sigma}(x)}$. Hence  
$\lambda_{\sigma}\in \Aut_{\bZ}(X)$ and $\lambda_{\sigma\tau}=\lambda_{\sigma}\lambda_{\tau}$\ 
$(\forall\sigma,\tau\in\Gamma)$, so that 
$\lambda\in\Hom(\Gamma,\Aut_{\bZ}(X))$. 
Since the quotient morphism $G\times_00^*\to G^t\times_00^*$ is $\Gamma$-equivariant,  $\Gamma(\cO_{G^t_{0^*}})$ is a $\Gamma$-submodule of 
$\Gamma(\cO_{G_{0^*}})$ via \cite[p.~35]{FC90}, so that 
$\lambda_{\sigma}(Y)=Y$.   
This proves (\ref{item:Hom(Gamma,Aut)}). 
\par

Let $\Gamma_0:=\{\sigma\in\Gamma;\lambda_{\sigma}(x)=x\ (\forall x\in X)\}$, $k(\eta'):=k(\eta^*)^{\Gamma_0}$, $R'$ 
the integral closure of $R$ in $k(\eta')$, 
$S':=\Spec R'$, $0'$ the closed point of $S'$ and $k(0')$ 
the residue field of $R'$.  
Since $k(0')=k(0^*)^{\Gamma_0}$ by  Eq.~(\ref{eq:Hensel}),  
$$G_{0'}\simeq G_{0^*}/\Gamma_0\simeq \Spec (k(0^*)[w^x;x\in X])^{\Gamma_0}= 
\Spec k(0')[w^x;x\in X].$$
By (\ref{eq:minimality assumption}), $S'=S^*$ and 
$\Gamma_0=\{1\}$.  
This proves (\ref{item:Gamma effective}). 
\end{proof}

\subsection{The $S^*$-scheme $P^*_{l_0l}$}
\label{subsec:P*l0l}
Now we return to the proof of Theorem~\ref{thm:main thm}. 
With the notation in \S~\ref{subsec:Galois action},
there exists a finite {\it \'etale} Galois cover $\varpi:S^*\to S$ such that 
$G_{0^*}$ is a split $k(0^*)$-torus. Then we obtain 
an FC datum 
$$\xi:=\FC(G_{S^*},\cL_{S^*})=(X,Y,a,b,A,B)$$ with $\cL_{S^*}$ symmetric. 
Let $N:=|B_{\xi}|=|B|$. 
There exist a  N\'eFC kit $\xi^{\natural}$ over $S_{\nor}(\xi)$ 
of $\xi^e_{2N}$ and the induced $l$-th N\'eFC kit from $\xi^{\natural}$
$$\xi^{\natural}_{l}:=(X, Y, \epsilon_{l},b^e_{l},E_{l}, \Sigma_{l})\  
(l\geq 1).$$
Let $l_0\in\bN$ such that 
$\Sigma_{l_0}(0)$ is integral and 
\begin{align*}
R^*_{l_0l}(\xi^{\natural})&:=
R^*[\zeta_{l_0l,x}\theta_{l_0l}; x\in X],\ 
\zeta_{l_0l,x}:=s^{D_{l_0l}(x)}w^{x}\ (l\geq 1). 
\end{align*} 

We recall the following $S^*$-schemes 
from Case~1 Step~2 and 
Eq.~(\ref{eq:R form of R'l_xinatural}):
\begin{gather*}
\tQ^*_{l_0l}:=\Proj R^*_{l_0l}(\xi^{\natural}),\   
(\tP^{*}_{l_0l},\tcL^{^*}_{l_0l}),\  
W^*_{l_0l,a},\ (P^*_{l_0l},\cL^*_{l_0l}).
\end{gather*}

Note that Claim~\ref{claim:P*l} with $S^*=S_{\spl}$ 
is true. Hence $P^*_{l_0l}$ is an $S^*$-scheme with open subschemes  
$G^*_{l_0l}$ and $\cG^*_{l_0l}$ 
satisfying  
 (i$^*$)/(iii$^*$)-(v$^*$) and 
(ii$^{**}$)\  $(G^*_{l_0l},\cL^*_{l_0l} )\simeq (G_{S^*}, \cL_{S^*}^{\otimes 2Nl_0l})$,\  
$\cG^*_{l_0l}\simeq \cG_{S^*}$.

\begin{lemma}
\label{lemma:a b lambda_sigma}Under the notation 
in Lemma~\ref{lemma:nomal form of Gamma action}, we have  
\begin{gather*}
\sigma(a(y))=a(\lambda_{\sigma}(y)),\ \sigma(b(y,x))=b(\lambda_{\sigma}(y),\lambda_{\sigma}(x))\\ 
(\forall y\in Y, \forall x\in X, \forall\sigma\in\Gamma).
\end{gather*}
\end{lemma}
\begin{proof}
Let $G^*:=G_{S^*}$, $\cL^*:=\cL_{S^*}$, $G^*_{\sigma}:=G^*$  
 and $g_{\sigma}:=\id_G\times_Ss_{\sigma}:G^*_{\sigma}\to G^*$.
\footnote{See Definition~\ref{defn:ssigma tsigma} for $s_{\sigma}$.} 
Then $g_{\sigma}$ is an $S$-isomorphism, which induces an $R$-isomorphism 
$(g^{\wedge}_{\sigma})^*:\Gamma(\cO_{G^{*\wedge}})
\to\Gamma(\cO_{G_{\sigma}^{*\wedge}}).$  
 Then $(g^{\wedge}_{\sigma})^*=\psi_{\sigma}$ and   
$w_{\sigma}^x:=(g^{\wedge}_{\sigma})^*(w^x)
=w^{\lambda_{\sigma}(x)}$ $(\forall x\in X)$ by Lemma~12.2. 
Let $\cL_{\sigma}^*:=g_{\sigma}^*\cL^*$ be 
the pullback of $\cL^*$ to $G^*_{\sigma}$. Then 
by Theorem~3.1, we obtain via the isomorphism $(g^{\wedge}_{\sigma})^*$
\begin{align*}
\Gamma(G^*_{\eta^*},\cL^*_{\eta^*})
&=\left\{
\begin{matrix}
\theta=\sum_{x\in X}c_x(\theta)w^x\in \Gamma(\cO_{G^{*\wedge}});\\
 c_{x+y}(\theta)=a(y)b(y,x)c_x(\theta)\\ 
 c_x(\theta)\in k(\eta^*)\ (\forall x\in X, \forall y\in Y)
\end{matrix}\right\},
\end{align*}
\begin{align*}
\Gamma(G^*_{\sigma,\eta^*},\cL^*_{\sigma,\eta^*})
&=\left\{
\begin{matrix}
\theta_{\sigma}=\sum_{x\in X}c_{\sigma,x}(\theta_{\sigma})w_{\sigma}^x
\in\Gamma(\cO_{G_{\sigma}^{*\wedge}});\\
 c_{\sigma,x+y}(\theta_{\sigma})=\sigma(a(y)b(y,x))c_{\sigma,x}(\theta_{\sigma}
\\ 
 c_{\sigma,x}(\theta_{\sigma})\in k(\eta^*)\ (\forall x\in X, \forall y\in Y)
\end{matrix}\right\}.
\end{align*}

Since $\Gamma(G^*_{\sigma,\eta^*},\cL^*_{\sigma,\eta^*})=\Gamma(G^*_{\eta^*},\cL^*_{\eta^*})$, which is $R^*$-free of rank $|X/Y|$,  
we have the same relations between the coefficients of 
$w^{\lambda_{\sigma}(x+y)}$ and $w^{\lambda_{\sigma}(x)}$: 
$$a(\lambda_{\sigma}(y))b(\lambda_{\sigma}(y),\lambda_{\sigma}(x))=\sigma(a(y)b(y,x))\ (\forall x\in X, \forall y\in Y), 
$$whence  
we obtain $\sigma(a(y))=a(\lambda_{\sigma}(y))$ for $x=0$, 
so that $\sigma(b(y,x))
=b(\lambda_{\sigma}(y),\lambda_{\sigma}(x))$. This completes the proof.
\end{proof}

\begin{defn}\label{defn:psi_sigma,tf_sigma}
Let $\sigma\in\Gamma$. 
We define 
$\tpsi^{\sharp}_{\sigma}\in\Aut(k(\eta^*)[X][\theta_{l_0l}])$ by 
\begin{equation*}\label{eq:tilde_psi_sigma}
\tpsi^{\sharp}_{\sigma}(\theta_{l_0l})
=\theta_{l_0l},\ \tpsi^{\sharp}_{\sigma}(cw^x)=
\sigma(c)w^{\lambda_{\sigma}(x)}\ 
(\forall c\in k(\eta^*), \forall x\in X).
\end{equation*} 
\end{defn}

Regardless of $\deg(\theta_{l_0l})=l_0l$, 
we define here the $*$-degree by $\deg^*(\theta_{l_0l})=1$   
in order to use the notation of \cite[II, 2.1-2.5]{EGA}. Let $\cA:=R^*_{l_0l}$, $\cA_k$ the subset of $\cA$ consisting of all 
homogeneous elements of $*$-degree $k$, 
$V_{l_0l}:=\cA(1)$ the graded $\cA$-module whose subset $(V_{l_0l})_k$ consisting of all homogenous elements of $*$-degree $k$ 
is $\cA_{k+1}$ $(\forall k\geq 0)$ and 
$\widetilde{V}_{l_0l}$ the $\cO_{\tP_{l_0l}^*}$-module 
associated with $V_{l_0l}$. Note that $\tcL^{*}_{l_0l}=\cO_{\tP^*_{l_0l}}(1)=\widetilde{V}_{l_0l}$ and $\tpsi^{\sharp}_{\sigma}(\cA)=\cA$ $(\forall \sigma\in\Gamma)$. 

 Let $\cA^i:=\cA$, $V_i:=V_{l_0l}$  
$(i\in[1,2])$ and  $\tpsi_{\sigma}:\cA^2\to \cA^1$ the restriction to $\cA^2$ 
of $\tpsi^{\sharp}_{\sigma}$. Then $\tpsi_{\sigma}$ is 
a $\sigma$-semilinear $R$-homomorphism, which induces 
$$\tf_{\sigma}\in\Isom_S(\Proj \cA^1, \Proj \cA^2),\ \tLambda_{\sigma}\in\Isom_{R}(\tf_{\sigma}^*\widetilde{V}_2,\widetilde{V}_1)$$ such that 
$\tf_{\sigma}^*=\tpsi_{\sigma}\in\Isom_R(\cA^2, \cA^1)$ and 
 $\tLambda_{\sigma}$ is the 
$R$-isomorphism of $\cO_{\Proj \cA^1}$-modules  
induced from $\tpsi_{\sigma|V_2}$. 
Hence we have $\tf_{\sigma}\in\Aut_S(\tP_{l_0l}^*)$ 
and $\tLambda_{\sigma}\in\Isom_R(\tf_{\sigma}^*\tcL_{l_0l}^*,\tcL_{l_0l}^*)$. 
By $\tpsi_{\sigma\tau}=\tpsi_{\sigma}\circ\tpsi_{\tau}$, we have  
\begin{equation}\label{eq:cocycle tLambda_sigma}
\tf_{\sigma\tau}=\tf_{\tau}\circ\tf_{\sigma},\ \tLambda_{\sigma\tau}=\tLambda_{\sigma}\circ \tf_{\sigma}^*\tLambda_{\tau}\ (\forall\sigma, \tau\in\Gamma).
\end{equation}   

Let $\tf_{\sigma}^{\wedge}$ 
(resp. $\tpsi_{\sigma}^{\wedge}$,   
$\tLambda_{\sigma}^{\wedge}$) 
be the $I^*$-adic completion of 
$\tf_{\sigma}$ (resp. $\tpsi_{\sigma}$,
 $\tLambda_{\sigma}$).

\begin{claim}
\label{claim:descent of f_sigma Lambda_sigma}  The following are true: 
\begin{enumerate}
\item\label{item:descent of f_sigma} $\tf_{\sigma}^{\wedge}$  descends to an 
$S^{\wedge}$-automorphism of $\tP_{l_0l}^{*\wedge}/Y$, 
which is algebraized into 
$f_{\sigma}\in\Aut_S(P_{l_0l}^*)$  such that 
$f_{\sigma\tau}=f_{\tau}\circ f_{\sigma}$ 
$(\forall \sigma,\tau\in\Gamma)$;
\item\label{item:descent of Lambda_sigma} 
$\tLambda_{\sigma}^{\wedge}$ descends to an 
$R$-isomorphism $\tLambda_{\sigma}^{\wedge}\op{mod} Y:
(\tf_{\sigma}^{*}\tcL^{*}_{l_0l})^{\wedge}/Y\to\tcL^{*\wedge}_{l_0l}/Y$, 
which is algebraized  
 into an $R$-isomorphism $\Lambda_{\sigma}:f_{\sigma}^*\cL_{l_0l}^*\to\cL_{l_0l}^*$ such that $\Lambda_{\sigma\tau}=\Lambda_{\sigma}\circ f_{\sigma}^*\Lambda_{\tau}$\ $(\forall\sigma, \tau\in\Gamma)$.
\end{enumerate} 
\end{claim}
\begin{proof}For any $(y,\sigma)\in Y\times\Gamma$, by 
Lemma~\ref{lemma:a b lambda_sigma}, we see
\begin{align*}
\tpsi_{\sigma} S_y^*(\theta_{l_0l})
&=a(\lambda_{\sigma}(y))^{2Nl_0l}
w^{2Nl_0l\lambda_{\sigma}(y)}\theta_{l_0l}
=S_{\lambda_{\sigma}(y)}^*\tpsi_{\sigma}(\theta_{l_0l}),\\
\tpsi_{\sigma} S_y^*(cw^x)
&=\sigma(c)b(\lambda_{\sigma}(y),\lambda_{\sigma}(x))w^{\lambda_{\sigma}(x)}
=S_{\lambda_{\sigma}(y)}^*\tpsi_{\sigma}(cw^x). 
\end{align*}
  
Hence $\tpsi_{\sigma}^{\wedge}\circ S_y^*=S_{\lambda_{\sigma}(y)}^*\circ
\tpsi_{\sigma}^{\wedge}$, so that  
$S_{y}\circ
\tf_{\sigma}^{\wedge}=\tf_{\sigma}^{\wedge}\circ S_{\lambda_{\sigma}(y)}$ and 
$\tLambda_{\sigma}^{\wedge}\circ S_y^*=S_{\lambda_{\sigma}(y)}^*\circ
\tLambda_{\sigma}^{\wedge}$. Hence  $\tf_{\sigma}^{\wedge}$  descends to an 
$S^{\wedge}$-automorphism of $\tP_{l_0l}^{*\wedge}/Y$, which is algebraized 
 into $f_{\sigma}\in\Aut_S(P_{l_0l}^*)$ by \cite[III, 5.4.1]{EGA}. By  Eq.~(\ref{eq:cocycle tLambda_sigma}), 
$f_{\sigma\tau}=f_{\tau}\circ f_{\sigma}$ 
$(\forall \sigma,\tau\in\Gamma)$. 
This proves (\ref{item:descent of f_sigma}).

 Now we identify 
$(\tP^{*\wedge}_{l_0l},\tcL^{*\wedge}_{l_0l})/Y$ with 
$(P^{*\wedge}_{l_0l},\cL^{*\wedge}_{l_0l})$. 
Let $Z:=\bP(\cL_{l_0l}^*)\simeq P_{l_0l}^{*}$, 
$\cF:=\cO_Z\oplus\cL_{l_0l}^*$,    
$W:=\bP(\cF)$ and $W':=W\setminus Z$. 
Then $Z$ is a closed subscheme of $W$ 
and $W'=(\bL_{l_0l}^*)^{\otimes(-1)}$ where $\bL_{l_0l}^*$ is 
the underlying line bundle of $\cL_{l_0l}^*$.    
Let $Z^{\wedge}$ (resp. $W^{\wedge}$, $\cF^{\wedge}$) be 
the $I^*$-adic completion of $Z$ (resp. $W$, $\cF$). 
Let $\Gamma_{\sigma}:=\id_{\cO_{Z^{\wedge}}}
\oplus\tLambda_{\sigma}^{\wedge}
\in\Hom((f_{\sigma}^*\cF)^{\wedge},\cF^{\wedge})$. 
Then $\Gamma_{\sigma}$ induces 
\begin{align*}\hat{\Delta}_{\sigma}&:=\bP(\Gamma_{\sigma}\op{mod}\,Y)\in
\Isom_{S^{\wedge}}(W^{\wedge}, (f_{\sigma}^{*}W)^{\wedge}),\\
\hat{\Delta}_{\sigma|(W')^{\wedge}}&\in
\Isom_{S^{\wedge}}((\bL_{l_0l}^{*\wedge})^{\otimes(-1)},((f_{\sigma}^*\bL_{l_0l}^*)^{\wedge})^{\otimes(-1)})\\
&=\Isom_{S^{\wedge}}((f_{\sigma}^*\bL_{l_0l}^*)^{\wedge},\bL_{l_0l}^{*\wedge}),
\end{align*}
which induce the $R$-isomorphism 
$\tLambda_{\sigma}^{\wedge}\op{mod} Y:(f_{\sigma}^*\cL_{l_0l}^*)^{\wedge}\to(\cL_{l_0l}^*)^{\wedge}$.  
By \cite[III, 5.4.1]{EGA}, there exists 
$\Delta_{\sigma}\in\Isom_S(W,f_{\sigma}^*W)$ such that 
$(\Delta_{\sigma})^{\wedge}=\hat{\Delta}_{\sigma}$. 
 Hence 
$\Delta_{\sigma|Z}=\id_Z$ and   
$(\Delta_{\sigma|W'})^{\wedge}=\tLambda_{\sigma}^{\wedge}\op{mod} Y$, so that 
$$\Delta_{\sigma|W'}\in\Isom_S(W',f_{\sigma}^*(W'))=\Isom_S(f_{\sigma}^*\bL_{l_0l}^*,\bL_{l_0l}^*).$$ 

Let $\Lambda_{\sigma}\in\Isom_R(f_{\sigma}^*\cL_{l_0l}^*,\cL_{l_0l}^*)$ 
be the $R$-isomorphism of $\cO_{P_{l_0l}^*}$-modules induced by $\Delta_{\sigma|W'}$. Then 
$\Lambda_{\sigma}$ algebraizes $\tLambda_{\sigma}^{\wedge}\op{mod} Y$:
$\Lambda_{\sigma}^{\wedge}=(\Delta_{\sigma|W'})^{\wedge}=\tLambda_{\sigma}^{\wedge}\op{mod} Y$.  
This proves (\ref{item:descent of Lambda_sigma}) 
by Eq.~(\ref{eq:cocycle tLambda_sigma}).  
\end{proof}

\begin{claim}\label{claim:f_sigma varphi_sigma}
Let $\sigma\in\Gamma$, $h_{\sigma}:=f_{\sigma^{-1}}$ 
and $\varphi_{\sigma}:=h_{\sigma}^*\Lambda_{\sigma}$. Then 
\begin{enumerate}
\item\label{item:varphi_sigma isom} $\varphi_{\sigma}:\cL^*_{l_0l}\to h_{\sigma}^*(\cL^*_{l_0l})$ 
is an isomorphism of $\cO_{P^*_{l_0l}}$-modules; 
\item\label{item:varphisigma_tau} 
$h_{\sigma\tau}=h_{\sigma}\circ h_{\tau}$ and $\varphi_{\sigma\tau}=h_{\tau}^*\varphi_{\sigma}\circ\varphi_{\tau}$\ $(\forall\tau\in\Gamma)$;
\item\label{item:f_sigma maps} $h_{\sigma}$ maps $G^{*}_{l_0l}$ 
(resp. $\cG^{*}_{l_0l}$) onto itself.
\end{enumerate} 
\end{claim}
\begin{proof}Since $\tpsi_{\sigma}$ is $\sigma$-semilinear,
so are  $\tLambda_{\sigma}$, $\Lambda_{\sigma}$, 
$\tf_{\sigma}^*$ and $f_{\sigma}^*$.  Hence  
$\varphi_{\sigma}$ is an $\cO_{P_{l_0l}^*}$-homomorphism, 
 as is shown in the same manner as 
in Lemma~\ref{lemma:phi_sigma_tau iff varphi_sigma_tau}. 
This proves (\ref{item:varphi_sigma isom}). 
The first half of (2) is clear.
Then we see 
$\varphi_{\sigma\tau}=h_{\sigma\tau}^*\Lambda_{\sigma\tau}
=(h_{\tau}^*\circ h_{\sigma}^*)
(\Lambda_{\sigma}\circ f_{\sigma}^*\Lambda_{\tau})
=(h_{\tau}^*\circ h_{\sigma}^*)\Lambda_{\sigma}\circ  
(h_{\tau}^*\circ h_{\sigma}^*\circ f_{\sigma}^*)\Lambda_{\tau}
=h_{\tau}^*\varphi_{\sigma}\circ \varphi_{\tau}, $ 
which proves the second half of (\ref{item:varphisigma_tau}).  
The rest is clear. 
\end{proof}

\subsection{Descent to $S$}
\label{subsec:descent to S}
By Claim~\ref{claim:f_sigma varphi_sigma}, 
$(h_{\sigma};\sigma\in\Gamma)$ is a finite group of $S$-automorphisms of 
$P_{l_0l}^*$, and $(\varphi_{\sigma};\sigma\in\Gamma)$ 
is a descent datum for $\cL_{l_0l}^*$.  
Hence we have the quotients as $S$-schemes 
of $P^*_{l_0l}$, $G^*_{l_0l}$ and 
$\cG^*_{l_0l}$ by $\Gamma$  (see \cite[Theorem, p.~63]{Mumford12}), 
while $\cL_{l_0l}^*$ descends to an invertible sheaf of the quotient of 
$P^*_{l_0l}$ by $\Gamma$ by \cite[Descent, 0D1V]{Stacks} 
and Claim~\ref{claim:f_sigma varphi_sigma}. Thus   
we define 
\begin{equation}\label{eq:descents}
\begin{aligned}
P^{\flat}_{l_0l}:=P^*_{l_0l}/\Gamma,&\ \ \cL^{\flat}_{l_0l}
:=\text{the descent of $\cL^*_{l_0l}$ to 
$P^{\flat}_{l_0l}$},\\   
G^{\flat}_{l_0l}&:=G^*_{l_0l}/\Gamma,\ 
\cG^{\flat}_{l_0l}:=\cG^*_{l_0l}/\Gamma,
\end{aligned}
\end{equation}where $G^{\flat}_{l_0l}$ and $\cG^{\flat}_{l_0l}$ 
are open in $P^{\flat}_{l_0l}$ by  \cite[V, 4.1/7.1 (ii)]{SGA3}.

Let $P^*:=P^*_{l_0l}$,  
$G^*:=G^*_{l_0l}$, $\cG^*:=\cG^*_{l_0l}$, 
$P^{\flat}:=P^{\flat}_{l_0l}$, 
$G^{\flat}:=G^{\flat}_{l_0l}$ and 
$\cG^{\flat}:=\cG^{\flat}_{l_0l}$. 
Let  $q:P^*\to P^{\flat}$ be the quotient morphism, 
 $U$ any affine open $S$-subscheme of $P^{\flat}$ and $V:=q^{-1}(U)$. 
Since $q$ is affine,  
$V$ is a $\Gamma$-invariant affine open set of $P^*$ 
such that $U=V/\Gamma$ as $S$-scheme.    
$\Gamma(V,\cO_{P^*})$ is an $R^*$-module with some descent datum 
$(\rho_{\sigma};\sigma\in\Gamma)$ 
by Theorem~\ref{thm:Galois descent/RL}, so that   
$\Gamma(V,\cO_{P^*})=\Gamma(U,\cO_{P^{\flat}})\otimes_RR^*$ and  
$\Gamma(U,\cO_{P^{\flat}})=\Gamma(V,\cO_{P^*})^{\Gamma}$.
The same is true for any affine open $S$-subscheme $U$ of $G^{\flat}$ (resp. $\cG^{\flat}$) since $G^{\flat}$ (resp. $\cG^{\flat}$) is an open $S$-subscheme  of $P^{\flat}$. 
 
\begin{defn}\label{defn:T_valued points}
Recall that $\Gamma$ acts on $S^*=\Spec R^*$ 
by $(s_{\sigma}; \sigma\in\Gamma)$. 
Let $T$, $Z$ be $S$-schemes, $Z(T):=\Hom_S(T,Z)$, 
 $Z_{S^*}(T_{S^*}):=\Hom_S(T_{S^*},Z_{S^*})$ and $Q^*\in Z_{S^*}(T_{S^*})$. 
\footnote{Caution : We do not define 
$Z_{S^*}(T_{S^*})=\Hom_{S^*}(T_{S^*},Z_{S^*})$.} 
We define an action of $\Gamma$ on $Z_{S^*}(T_{S^*})$ by 
\begin{equation*}\label{defn:sigma(Q*)}
\sigma(Q^*)=(\id_Z\times_Ss_{\sigma})^{-1}\circ Q^*\circ
(\id_T\times_Ss_{\sigma}). 
\end{equation*} 
Then $\sigma(Q^*)\in Z_{S^*}(T_{S^*})$ and 
$(\sigma\tau)(Q^*)=\sigma(\tau(Q^*))$ $(\forall\sigma,\tau\in\Gamma)$. 
Note that if $Z$ is an $S$-subscheme of $P^{\flat}$, then $Z_{S^*}$ is an $S^*$-subscheme of $P^*$ and $\id_Z\times_Ss_{\sigma}=f_{\sigma|Z_{S^*}}$.  
We also define 
$i_{Z,T}:Z(T)\to Z_{S^*}(T_{S^*})$ 
by $i_{Z,T}(Q)=Q_{S^*}:=Q\times_S\id_{S^*}$. Then $i_{Z,T}$ is injective, 
so that we identify  $Z(T)$ with the subset 
$i_{Z,T}(Z(T))$ of $Z_{S^*}(T_{S^*})$. 
Note that $Z(T)$ is fixed elementwise by $\Gamma$ because 
$\sigma(Q_{S^*})=(\id_Z\times_Ss_{\sigma})^{-1}\circ (Q\times_S\id_{S^*})\circ
(\id_T\times_Ss_{\sigma})=Q\times_S (s_{\sigma}^{-1}\circ\id_{S^*}\circ s_{\sigma})=Q_{S^*}$.  
We denote by $Z_{S^*}(T_{S^*})^{\Gamma}$ the subset of $Z_{S^*}(T_{S^*})$ consisting of all $\Gamma$-invariant elements. 
\end{defn}

\begin{claim}\label{claim:T val points}
Let $T$ be any $S$-scheme, $\tau_{\sigma}:=\id_T\times s_{\sigma}$  and 
$Z\in\{P^{\flat}, G^{\flat}, \cG^{\flat}\}$. Then  
$Z(T)=Z_{S^*}(T_{S^*})^{\Gamma}:=
\{h\in Z_{S^*}(T_{S^*}); h\circ\tau_{\sigma}=f_{\sigma}\circ h\ (\forall\sigma\in\Gamma)\}.$
\end{claim}
\begin{proof}We shall prove the equality when $Z=G^{\flat}$, in which case, 
$Z_{S^*}=G^*$. Let $T^*:=T_{S^*}$. By $i_{G^{\flat},T}$, (we regard)    
$G^{\flat}(T)\subset G^*(T^*)^{\Gamma}$.  
We prove the converse inclusion.  
Since $G^{\flat}$ is open in $P^{\flat}$, 
there exists a finite  cover 
$(U_{\lambda}; \lambda\in\Lambda)$ 
of $G^{\flat}$ by affine open $S$-schemes $U_{\lambda}$  
such that $V_{\lambda}:=q^{-1}
(U_{\lambda})\ (\simeq U_{\lambda}\times_SS^*)$ 
are affine $S^*$-schemes and 
\begin{equation*}
\label{eq:basic descent data}
\begin{aligned}
A^*_{\lambda}&:=\Gamma(\cO_{V_{\lambda}})
=\Gamma(\cO_{U_{\lambda}})\otimes_RR^*,\  
A_{\lambda}:=\Gamma(\cO_{U_{\lambda}})
=\Gamma(\cO_{V_{\lambda}})^{\Gamma}\ (\forall \lambda\in\Lambda).
\end{aligned}
\end{equation*}

Let $h\in G^*(T^*)^{\Gamma}$. Let $(T_i;i\in I)$ be a cover of $T$ by 
affine open $S$-schemes such that $h(T^*_i)\subset V_{\lambda(i)}$
for any $i\in I$ 
and some $\lambda(i)\in\Lambda$ where $T^*_i:=T_i\times_SS^*$. 
Let $B_i:=\Gamma(\cO_{T_i})$ and $B^*_i:=B_i\otimes_RR^*$. 
Let $h_i$ 
be the restriction of $h$ to $T^*_i$, 
which belongs to $V_{\lambda(i)}(T^*_i)$. 
Since $h\circ\tau_{\sigma}=f_{\sigma}\circ h$ 
$(\forall\sigma\in\Gamma)$, we have 
$h_i\circ\tau_{\sigma}=f_{\sigma}\circ h_i$. 
Thus $h_i\in V_{\lambda(i)}(T^*_i)^{\Gamma}$.  

$\sigma\in\Gamma$ acts on $\Hom_R(A_{\lambda(i)},B^*_i)$ 
through the action $\id_{B_i}\otimes\sigma$ on $B_i\otimes_RR^*$. 
Since $R^*$ is $R$-free, there are the following canonical identifications:
\begin{equation*}
\begin{aligned}
V_{\lambda(i)}(T^*_i)&=\Hom_{R^*}(A^*_{\lambda(i)},B^*_i)
=\Hom_R(A_{\lambda(i)},B^*_i),\\
V_{\lambda(i)}(T^*_i)^{\Gamma}
&=\{\phi^*\in\Hom_R(A^*_{\lambda(i)},B^*_i) ; 
\tau_{\sigma}^*\circ\phi^*=\phi^*\circ f_{\sigma}^*\}\\
&=\left\{
\phi^*\in \Hom_R(A_{\lambda(i)},B^*_i)
;\begin{matrix}
(\id_{B_i}\otimes\sigma )(\phi^*(x))=\phi^*(x)\\
(\forall x\in A_{\lambda(i)}, 
\forall\sigma\in\Gamma)
\end{matrix}
\right\}\\
&=\Hom_R(A_{\lambda(i)},B_i)
=\Hom_S(T_i,U_{\lambda(i)}). 
\end{aligned}
\end{equation*}

Hence for each $h_i\in V_{\lambda(i)}(T^*_i)^{\Gamma}$,  
there exists a {\it unique} descent 
$g_i\in U_{\lambda(i)}(T_i)$ of $h_i$. 
Since $(g_i)_{|T_i\cap T_j}=(g_j)_{|T_i\cap T_j}$ 
by the uniqueness of descent, $(g_i;i\in I)$ gives $g\in G^{\flat}(T)$. 
The other cases are proved similarly.  
\end{proof}

\begin{claim}\label{claim:action T val points}
The group $S^*$-scheme $G^*$ (resp. $\cG^*)$ descends to a group $S$-scheme 
$G^{\flat}$ (resp. $\cG^{\flat})$, and the following are true:
\begin{enumerate}
\item[(i$^{\flat}$)] $(P^{\flat},\cL^{\flat}_{l_0l})$ is an irreducible 
flat projective Cohen-Macaulay $S$-scheme;
\item[(ii$^{\flat}$)] $G^{\flat}$ (resp. $\cG^{\flat}$) 
is a group $S$-scheme 
such that $G^{\flat}\simeq G$ (resp. $\cG^{\flat}\simeq\cG$);
\item[(iii$^{\flat}$)] $P^{\flat}\setminus\cG^{\flat}$ is of codimension two;
\item[(iv$^{\flat}$)] there exists a morphism $G^{\flat}\times P^{\flat}\to P^{\flat}$ 
(resp. $\cG^{\flat}\times P^{\flat}\to P^{\flat}$) extending 
multiplication of $G^{\flat}$ (resp. $\cG^{\flat})$;
\item[(v$^{\flat}$)] $(\cL^{\flat}_{l_0l})_{|\cG^{\flat}}$ is ample cubical with $\cL^{\flat}_{l_0l,\eta}=\cL^{2Nl_0l}_{\eta}$.
\end{enumerate}
\end{claim}
\begin{proof}Since $q:P^*\to P^{\flat}$ 
is finite and surjective, 
and $q^*\cL^{\flat}_{l_0l}$ is ample,  
so is  $\cL^{\flat}_{l_0l}$ by \cite[tag~0B5V]{Stacks}.  
Since $P^*$ is irreducible, 
so is $P^{\flat}$. 
Since $P^*$ is 
flat and Cohen-Macaulay, so is 
$P^{\flat}$ because  $q$ is \'etale. \footnote{It is enough if $q$ is 
faithfully flat.} 
This proves (i$^{\flat}$).  

Next we shall prove (ii$^{\flat}$) for $G^{\flat}$. 
First note that 
$G^{\flat}\simeq G$ as $S$-schemes 
since both $G^{\flat}$ and $G$ are descents of $G_{S^*}$ 
with respect to the canonical descent datum.  
It remains to prove  
that $G^{\flat}\simeq G$ as group $S$-schemes.

Let $T$ be any locally noetherian $S$-scheme, $T_0:=T\times_S0$, 
$T^*:=T\times_SS^*$, $T^*_{0^*}:=T\times_00^*$ and  
$T^*_{\eta^*}:=T\times_{\eta}\eta^*$. 
Let $g_i\in G^{\flat}(T)$ $(i\in[1,2])$. 
 We define the product of $g_1$ and $g_2$ 
functorially in $T$. 
Let $h_i:=g_i\times_SS^*$ $(i\in[1,2])$.  
Since $h_i\in G^*(T^*)$, 
we have the product $h:=h_1\cdot h_2\in G^*(T^*)$.  

\smallskip
\noindent{\it Step 1.}  Let $p\in T_0$, 
$\barp:=h(p)\in G^*_{0^*}$, $V=\Spec\cO_{T,p}$, 
$V^*=V\times_SS^*$, $p^*:=p\times_00^*$ and $\cO_{T^*,p^*}^{\wedge}$ 
the $I^*$-adic completion  of $\cO_{T^*,p^*}$. 
We denote by the same $g_i$ the restriction of $g_i$ to $V$, 
$h_i=g_i\times_SS^*$ and $h_i^{\wedge}$ be 
 the $I^*$-adic completion of $h_i$. 
Since $G^*_{0^*}$ is a split $k(0^*)$-torus (\ref{eq:torus G'0}), $G^{*\wedge}$ is a formal $R^*$-torus : $\cB^*:=\Gamma(\cO_{G^{*\wedge}})=R^*[w^x;x\in X]^{I^*\op{-adic}}$. Then the $I^*$-adic completion 
$h^{\wedge}$ of 
the product $h$ is computed via 
$(h_i^{\wedge})^*\in\Hom(\Gamma(\cO_{G^{*\wedge}}),\Gamma(\cO_{V^{*\wedge}}))$
by $(h^{\wedge})^*(w^x)
=(h_1^{\wedge})^*(w^x)\cdot (h_2^{\wedge})^*(w^x)$\ $(\forall x\in X)$. 
 Since $h_i\in G^*(V^*)^{\Gamma}$ 
by Claim~\ref{claim:T val points}, we have 
$h_i\circ\tau_{\sigma}=f_{\sigma}\circ h_i$\ 
$(\forall\sigma\in\Gamma)$, 
we see 
\begin{align*}
((h\circ \tau_{\sigma})^{\wedge})^*(w^x)
&=(\tau_{\sigma}^{\wedge})^*\circ (h^{\wedge})^*(w^x)
=(\tau_{\sigma}^{\wedge})^*\circ
((h_1^{\wedge})^*(w^x)\cdot (h_2^{\wedge})^*(w^x))\\
&=(h_1^{\wedge})^*\circ (f_{\sigma}^{\wedge})^{*}(w^x)
\cdot (h_2^{\wedge})^*\circ (f_{\sigma}^{\wedge})^{*}(w^x)\\
&=(h^{\wedge})^*\circ (f_{\sigma}^{\wedge})^*(w^x)
=((f_{\sigma}\circ h)^{\wedge})^*(w^x).
\end{align*} 

Since $\cO_{T^*,p^*}$ is 
noetherian local and 
$I^*\cO_{T^*,p^*}$ is a proper ideal, we have 
$\bigcap_{n\geq 0}(I^*)^n\cO_{T^*,p^*}=0$, so that 
$\cO_{T^*,p^*}\subset \cO_{T^*,p^*}^{\wedge}$. 
Let $\cB^*_{\barp}$ be the localization 
by the prime ideal of $\cB^*$ corresponding to $\barp$, and 
$\cO_{G^*,\barp}^{\wedge}$
the $I^*$-adic completion of $\cO_{G^*,\barp}$.
Then $\cO_{G^*,\barp}\subset \cO_{G^*,\barp}^{\wedge}=\cB^*_{\barp}$,  
which is the local ring of $G^{*{\wedge}}$ at $\barp$. 

Then the above computation shows $(h\circ\tau_{\sigma})^{\wedge}
=(f_{\sigma}\circ h)^{\wedge}$ as elements of 
{\small $\Hom_{R}(\cO_{G^*,\barp}^{\wedge},\cO_{T^*,p^*}^{\wedge})$}. 
Hence $h\circ\tau_{\sigma}=f_{\sigma}\circ h$ as 
elements of {\small$\Hom_{R}(\cO_{G^*,\barp},\cO_{T^*,p^*})$}.
Since $G$ is of finite type over $S$,    
$h\circ\tau_{\sigma}=f_{\sigma}\circ h$ 
on an open neighborhood $W''$ of $V^*$ in $T^*$. 
Let $W':=\bigcap_{g\in\Gamma}g(W'')$.  
Then $W'$ is a $\Gamma$-invariant open subset of $T^*$
which contains $V^*$. 
Then $h\in G^*(W')^{\Gamma}$. 
Since $W'$ is quasi-projective,  
by \cite[V, 7.1~(i)-(ii)]{SGA3}, 
there exists a quotient $S$-scheme 
$W:=W'/\Gamma$, 
which is an open set of 
$T$  such that 
$V\subset W$ and $W'=W\times_SS^*$. 
Since $p$ can be any point of $T_0$, there exists an open 
subset $U$ of $T$ which contains $T_0$ such that 
$h\in G^*(U^*)^{\Gamma}$ where $U^*:=U\times_SS^*$.

\smallskip
\noindent{\it Step 2.}  Next we consider 
the restriction of $g_i$ to $T_{\eta}$ 
(resp. $h_i$ to $T^*_{\eta^*}$), 
which we denote by $g_i$ (resp. $h_i$). 
Since $G^{\flat}_{\eta}=G^*_{\eta^*}/\Gamma
=(G_{\eta}\times_{\eta}{\eta^*})/\Gamma=G_{\eta}$, 
$G^{\flat}_{\eta}$ is an abelian variety over $k(\eta)$. 
Hence we have $g_1\cdot g_2\in G^{\flat}(T_{\eta})$.  
Since $h_i=g_i\times_{\eta}\eta^*$ $(i\in [1,2])$, 
we have $h=(g_1\cdot g_2)\times_{\eta}\eta^*\in G^*(T^*_{\eta^*})^{\Gamma}$. 

\smallskip 
Now we return to the $S$-scheme $T$. Let $U$ be 
the open set of $T$ chosen in Step~1. 
Let $r_{U^*}:G^*(T^*)\to G^*(U^*)$
(resp. $r_{\eta^*}:G^*(T^*)\to G^*(T^*_{\eta^*})$) 
be the restriction 
to $U^*$ (resp. $T^*_{\eta^*})$. By Step 1, 
 $r_{U^*}(h)\in G^*(U^*)^{\Gamma}$,   
while  $r_{\eta^*}(h)\in G^*(T^*_{\eta^*})^{\Gamma}$ by Step~2. 
Since $T^*=U^*\bigcup T^*_{\eta^*}$, we have $h\in G^*(T^*)^{\Gamma}$, 
so that there exists a unique descent $g\in G(T)$ of $h$ 
by Claim~\ref{claim:T val points}, which we denote by $g_1\cdot g_2$. We call 
$g_1\cdot g_2$ the product of $g_1$ and $g_2$.  
It is easy to prove that the product 
$(g_1,g_2)\mapsto g_1\cdot g_2$ is functorial in $T$, and it is associative. 
Similarly we can also prove by using Claim~\ref{claim:T val points} and \S~\ref{subsec:P*l0l}~(ii$^{**}$) that the product of $G^{\flat}$ satisfies the axioms of group $S$-schemes. It follows that 
$G^{\flat}$ is a group $S$-scheme.  
Since $G^{\flat}\simeq G$ as $S$-schemes with 
the same unit, they are isomorphic as group $S$-scheme 
as in Theorem~\ref{thm:G isom G' iff FC same}. 

Since $\cG^{\flat}_0(k(0)^{\sep})/G^{\flat}_0(k(0)^{\sep})=\cG^*_{0^*}(k(0^*)^{\sep})/G^*_{0^*}(k(0^*)^{\sep})$ and 
$\cG^*$ is the N\'eron model of $G^*_{\eta}$, 
 $\cG^{\flat}$ is the N\'eron model of $G^{\flat}_{\eta}$, 
so that $\cG^{\flat}\simeq\cG$. 
This proves (ii$^{\flat}$). Let $Z^*_{l_0l}:=P^*\setminus\cG^*$ 
and $Z^{\flat}_{l_0l}:=P^{\flat}\setminus\cG^{\flat}$. Since $q$ 
is finite, 
$Z^{\flat}_{l_0l}=q(Z^*_{l_0l})$ is of codimension two in $P^{\flat}$. 
This proves (iii$^{\flat}$).  (iv$^{\flat}$) (resp. (v$^{\flat}$)) 
is proved in a manner similar to (ii$^{\flat}$) (resp. Claim~\ref{claim:P*l}~(v$^*$)) .
 \end{proof} 

This completes the proof of Theorem~\ref{thm:main thm} in Case~2.

\section{Examples}
\label{sec:examples}
Let $R$ be a CDVR with uniformizer $s$, $S:=\Spec R$ and $(P,\cL)$  
a Mumford family over $S$ \cite{AN99} 
of principally polarized abelian surfaces.
The purpose of this section is 
to compare the other three  
families
$$(P^{\flat},\cL^{\flat}),\ 
(P^{\natural},\cL^{\natural}),\ (P_{l},\cL_{l})\ (l\geq 1)
$$with $(P,\cL)$. 
Each of them is a relatively compactification 
of the N\'eron model $\cG$ of $P_{\eta}$ with generic fibers given as follows: 
\begin{gather*}
(P^{\flat}_{\eta},\cL^{\flat}_{\eta})
\simeq (P_{\eta},\cL_{\eta}^{\otimes 2}),\ 
(P^{\natural}_{\eta},\cL^{\natural}_{\eta})
\simeq (P_{\eta},\cL_{\eta}^{\otimes 3}),\  
(P_{l,\eta},\cL_{l,\eta})\simeq (P_{\eta},\cL_{\eta}^{\otimes 6l}),
\end{gather*}where the restriction of 
$\cL_{l}$ and $\cL^{\natural}$ (resp. $\cL^{\flat}$) to $\cG$ is cubical 
(resp. non-cubical). 
We see $P\not\simeq P^{\flat}$,\ $P\not\simeq P^{\natural}$,\  
$P^{\flat}\not\simeq P^{\natural}$, while  
$P^{\natural}\simeq P_{l}$ $(\forall l\geq 1)$. 

\subsection{The $S$-scheme $(P,\cL)$}
\label{subsection:2dim PSQAS}\quad
Let $(m_i;i=1,2)$ be a basis 
of $X:=\bZ^{\oplus 2}=\bZ m_1\oplus\bZ m_2$,  
$m_3:=-(m_1+m_2)$, $\xi=(X,a,b,A,B)$ 
an FC datum with $a(x)=s^{A(x)}$, 
$b(x,y)=s^{B(x,y)}$ and  
\begin{gather*}
A(x)=x_1^2-x_1x_2+x_2^2,\ 
B(x,y)=2x_1y_1-x_1y_2-x_2y_1+2x_2y_2
\end{gather*}where $x=x_1m_1+x_2m_2$ and $y=y_1m_1+y_2m_2$. 
Let $w_i:=w^{m_i}$ $(i\in[1,3])$. 
\par
 
Let $R^{\flat}(\xi)$ and $(S_y;y\in X)$ be the same 
as in Definition~\ref{defn:graded alg Aflat}, 
$(\tP,\tcL):=(\Proj R^{\flat}(\xi),\cO_{\tP}(1))$ and  
$U_{0,0}:=\Spec R[sw_1^{\pm 1}, sw_2^{\pm 1}, sw_3^{\pm 1}]$. 
The $S$-scheme $\tP$ is covered by affine open sets 
$S_y(U_{0,0})$\ $(y\in X)$. 
As in \S~4, 
let $(P,\cL):=(P(\xi),\cL(\xi))$. Then 
the generic fiber $(P_{\eta},\cL_{\eta})$ is 
 a principally polarized abelian surface over $k(\eta)$, 
while the closed fiber $P_0$ of $P$ 
is a union of 
2 copies of $\bP^2_{k(0)}$, whose configuration 
is given by the Delaunay decomposition $\Del_B$ (Fig.~1/left) 
modulo $X$ by \cite[\S6, B, p.~418]{AN99}. 
Note that $P$ is not a relative compactification
of the N\'eron model of $P_{\eta}$.

\subsection{A twist $(P^{\flat},\cL^{\flat})$ of $(P,\cL)$}
\label{subsection:twisted PSQAS}
Let $\xi$ be the same as in \S~\ref{subsection:2dim PSQAS}  
and $\theta^{\flat}:=\vartheta^2\in R^{\sharp}$ 
with the notation in \S~\ref{subsec:a Mumford family}.  
We define an $R$-subalgebra  $R^{\flat}$ of $R^{\sharp}$ by
\begin{equation*}
\label{eq:RD}
R^{\flat}=R[\xi_{\alpha,y}\theta^{\flat}; y\in X, \alpha\in\Sigma]
\end{equation*}where $\xi_{\alpha,y}:=b(\alpha,y)a(y)^2w^{\alpha+2y}$ and 
$\Sigma:=\{0, \pm m_1,\pm m_2, \pm m_3\}$. The $R$-automorphism 
$S_y^*$ in (\ref{eq:Sy on Rsharp}) induces 
$S_y^*\in\Aut_R(R^{\flat})$ $(y\in X)$ by
\begin{equation}
\label{eq:action Sy on RD}
S_y^*(\xi_{\alpha,z}\theta^{\flat})=\xi_{\alpha,y+z}
\theta^{\flat}\quad (z\in X). 
\end{equation}

 Let $(\tP^{\flat},\tcL^{\flat}):=(\Proj R^{\flat},\cO_{\tP^{\flat}}(1))$. 
Let $S_y$ be the $S$-automorphism 
of $\tP^{\flat}$ induced from $S_y^*$. The $S$-scheme $\tP^{\flat}$ is covered by affine open sets $U_{\alpha,y}^{\flat}$ $(\alpha\in\Sigma, y\in Y)$ where 
$U_{\alpha,y}^{\flat}:=\Spec R[\xi_{\beta,z}/\xi_{\alpha,y}; \beta\in\Sigma,z\in X]$. For example, 
\begin{equation*}
\label{eq:charts of hexagon H111 0}
\begin{aligned}
U^{\flat}_{0,0}&=\Spec R[w_1, w_1^{-1}, w_2, w_2^{-1}],\\
U^{\flat}_{m_3,0}&=\Spec R[w_1, sw_1^{-1}, w_2, sw_2^{-1}],\\
U^{\flat}_{-m_2,0}&=\Spec R[w_1^{-1}, sw_1, w_3^{-1}, sw_3],\\
U^{\flat}_{m_1,m_3}&=\Spec R[s^{-1}w_3^{-1}, s^2w_2^{-1}, s^{-1}w_2, s^2w_3].
\end{aligned}
\end{equation*} 

Let $(P^{\flat},\cL^{\flat})$ be the algebraization 
of  the quotient $(\tP^{\flat\wedge},\tcL^{\flat\wedge})/X$. 
The generic fiber $(P^{\flat}_{\eta},\cL^{\flat}_{\eta})$ is 
$k(\eta)$-isomorphic to 
$(P_{\eta},\cL^{\otimes 2}_{\eta})$ 
by Corollary~\ref{cor:FC(GLm)}, 
while 
the closed fiber 
$P^{\flat}_0$ is the union of a 3-times blowing-up of $\bP^2_{k(0)}$ 
and two copies of $\bP^2_{k(0)}$, whose configuration is given 
by Fig.~1 modulo $X$ by \cite[\S~7]{Mumford72}. 
Let $(P,\cL)$ be the simplified Mumford 
family in \S~\ref{subsection:2dim PSQAS}. 
Then there is 
an $S$-morphism $\pi:P^{\flat}\to P$ 
contracting one of $\bP^2_{k(0)}$ 
such that $\pi_{\eta}$ is an isomorphism. Indeed, 
if we let $U^{\flat}:=U^{\flat}_{m_3,0}\cup U^{\flat}_{-m_2,0}\cup U^{\flat}_{m_1,m_3}$,  then the restriction 
$\pi_{|U^{\flat}}:U^{\flat}\to U_{0,0}$ is explicitly given by 
$$\pi^*(sw_1,sw_3^{-1},sw_2,sw_1^{-1},sw_3,sw_2^{-1})
=(sw_1,w_3^{-1},w_2,sw_1^{-1},s^2w_3,s^2w_2^{-1}),
$$ 
whose right hand side is understood, for example, on $U^{\flat}_{m_3,0}$ as 
$$(s\cdot w_1,w_1\cdot w_2,w_2,sw_1^{-1}, sw_1^{-1}\cdot sw_2^{-1},s\cdot sw_2^{-1}). 
$$

The restriction 
$\pi_{|U^{\flat}}$ contracts the surface $\{s=sw_1^{-1}=w_2=0\}$ of $U^{\flat}_{m_3,0}$. 

\begin{figure}[ht]
\label{fig:2dim Voronoi}
\caption{}
\vspace*{0.7cm}
\begin{picture}(80,70)(30,0)
   \put(15,0){\line(-1,0){10}}
   \put(15,0){\line(1,0){70}}
   \put(15,30){\line(-1,0){10}}
   \put(15,30){\line(1,0){70}}
   \put(15,60){\line(-1,0){10}}
   \put(15,60){\line(1,0){70}}
   \multiput(15,60)(30,0){3}{\circle*{3}}
   \multiput(15,30)(30,0){3}{\circle*{3}}
   \multiput(15,0)(30,0){3}{\circle*{3}}
   \put(45,0){\line(0,-1){10}}
   \put(75,0){\line(0,-1){10}}
   \put(15,0){\line(0,-1){10}}
   \put(15,0){\line(0,1){70}}
   \put(45,0){\line(0,1){70}}
   \put(75,0){\line(0,1){70}}
   \put(15,0){\line(-1,-1){10}}
   \put(15,0){\line(1,1){70}}
   \put(45,0){\line(1,1){40}}
   \put(45,0){\line(-1,-1){10}}
   \put(75,0){\line(1,1){10}}
   \put(75,0){\line(-1,-1){10}}
   \put(30,15){\line(1,1){15}}
   \put(60,15){\line(1,1){15}}
   \put(15,30){\line(-1,-1){10}}
   \put(15,30){\line(1,1){40}}
   \put(75,30){\line(1,1){10}}
   \put(15,60){\line(1,1){10}}
   \put(15,60){\line(-1,-1){10}}
   \put(35,-30){$\Del_B$}
   \end{picture}
\vspace*{0.2cm}
\begin{picture}(80,70)(-22,0)
   \put(15,0){\line(-1,0){10}}
   \multiput(15,0)(15,0){5}{\circle*{3}}
   \multiput(15,15)(30,0){3}{\circle*{3}}
   \multiput(15,30)(15,0){5}{\circle*{3}}
   \multiput(15,45)(30,0){3}{\circle*{3}}
   \multiput(15,60)(15,0){5}{\circle*{3}}
   \put(15,0){\line(1,0){70}}
   \put(15,30){\line(1,0){70}}
   \put(15,60){\line(1,0){70}}
   \put(15,30){\line(-1,0){10}}
   \put(15,60){\line(-1,0){10}}
   \put(15,30){\line(1,0){70}}
   \put(15,60){\line(1,0){70}}
   \put(15,0){\line(0,-1){10}}
   \put(45,0){\line(0,-1){10}}
   \put(75,0){\line(0,-1){10}}
   \put(15,0){\line(0,1){70}}
   \put(45,0){\line(0,1){70}}
   \put(75,0){\line(0,1){70}}
   \put(30,0){\line(1,1){55}}
   \put(30,0){\line(-1,-1){10}}
   \put(60,0){\line(1,1){25}}
   \put(60,0){\line(-1,-1){10}}
   \put(15,15){\line(-1,-1){10}}
   \put(15,15){\line(1,1){55}}
   \put(15,45){\line(-1,-1){10}}
   \put(15,45){\line(1,1){25}}
   \put(37,-30){$\tP^{\flat}_0$}
\end{picture}
\vspace*{1.3cm}
\end{figure}

Now we explain how $(P^{\flat},\cL^{\flat})$ is obtained from  
\cite[\S~7]{Mumford72}. 

Let $A:=k[[a,b,c]]$ and $\fm:=(a,b,c)$. 
We define  $\tP^D=\Proj R_{\phi,\Sigma}$ where 
\begin{gather*}
R_{\phi,\Sigma}:=A[\fb(\alpha+y, y)w^{\alpha+2y}\theta^D; 
\alpha\in\Sigma, y\in X],\\
\fb(x, y):=(bc)^{x_1y_1}c^{-(x_1y_2+x_2y_1)}(ac)^{x_2y_2}.
\end{gather*} 

Let 
$\tL^D:=\cO_{\tP^D}(1)$.
Then by taking the $\fm$-adic completion 
$(\tP^D_{\formal},\tL^D_{\formal})$, and then 
by algebraizing the formal quotient of  
$(\tP^D_{\formal},\tL^D_{\formal})$  by the action 
$(S_y;y\in X)$ similar to Eq.~(\ref{eq:action Sy on RD}), we obtain a projective $A$-scheme $(P^D,L^D)$. This $(P^D,L^D)$ 
is the same as the example given in \cite[\S~7]{Mumford72}.
Let $R:=k[[s]]$, and let $\gamma:\Spec R\to \Spec A$ 
be the {\it diagonal} morphism given by $\gamma^*(a,b,c)=(s,s,s)$.
Then 
it is clear that 
$$R^{\flat}=R_{\phi,\Sigma}\otimes_AR,\ 
(P^{\flat},\cL^{\flat})\simeq \gamma^*(P^D,L^D).$$

\subsection{The $S$-scheme $(P^{\natural},\cL^{\natural})$}
\label{subsec:Pdagger}
Let $\xi$ be the same as in \S~\ref{subsection:2dim PSQAS}. 
Since $B=B_{\xi}$ is even, there exists a unique eFC datum 
$\xi^e$ over $R$ with $\xi(\xi^e)=\xi$.  
By Lemma~\ref{lemma:Be and Bephi}, $\xi^e=(X,a,b^e,A,B)$ 
with $b^e(u,x)=s^{u(x)}$ $(\forall x\in X, \forall u\in X^{\vee})$. 
Let $\xi^{\natural}$ be the N\'eFC kit 
of $\xi^e_{2N}=\xi^e_6$ and 
$(P_{l},\cL_{l}):=(P_{l}(\xi^{\natural}),\cL_{l}(\xi^{\natural}))$
\ $(l\geq 1)$.

We shall construct another relative compactification 
$(P^{\natural},\cL^{\natural})$ different from  $(P_{l},\cL_{l})$. 
 Let $m_3=-(m_1+m_2)$, 
$\Sigma=\{0,\pm m_i;i\in[1,3]\}$ and $\Psi:=\{0, \pm m_1\}$. 
We define 
$R^{\natural}$ as follows:
\begin{gather*}
R^{\natural}=R[\epsilon_{+}(u)b^e(u,\alpha)w^{\alpha+\mu(u)}
\theta_{+}; u\in X^{\vee},\alpha\in\Sigma] 
\end{gather*}
where 
\begin{gather*}
\theta_{+}:=\prod_{y\in\Psi}\vartheta_y,\ 
\epsilon_{+}(u):=s^{u_1^2+u_1u_2+u_2^2},\ 
b^e(u,y)=s^{u_1y_1+u_2y_2},\\
\beta(y):=(2y_1-y_2)f_1+(-y_1+2y_2)f_2,\\
\mu(u):=(2u_1+u_2)m_1+(u_1+2u_2)m_2,\\ 
E_{+}(u):=v_s\epsilon_{+}(u)=u(\mu(u))/2=u_1^2+u_1u_2+u_2^2,\\ 
u=u_1f_1+u_2f_2,\ y=y_1m_1+y_2m_2.
\end{gather*}

Note $\epsilon_{+}(\beta(y))=a(y)^3$. 
We define in a manner similar 
to Definition~\ref{defn:Voronoi polytopes} 
$$\Sigma(0)=\{x\in X_{\bR};E_{+}(u)+u(x)\geq 0\ (\forall u\in X^{\vee})\}.$$ 
Then $\Sigma=\Sigma(0)\cap X$ and  $\Sigma(0)$ is integral.  
See Fig.~2. Similarly we define 
$\Vor=\{\Sigma(v);v\in X^{\vee}\ \text{and their faces}\}$ where 
$\Sigma(v)=\Sigma(0)+\mu(v)$.

\begin{figure}
\label{fig:Sigma and Sigma0}
   \begin{picture}(110,90)(56.5,-10)
   \multiput(90,10)(20,0){2}{\circle*{5}}
   \multiput(90,30)(20,0){3}{\circle*{5}}
   \multiput(110,50)(20,0){2}{\circle*{5}}
   \put(77,-2){$m_3$}
   \put(89,52){$m_2$}
   \put(135,35){$m_1$}
   \put(90,10){\line(0,1){20}}
   \put(110,-10){\line(0,1){80}}
   \put(130,30){\line(0,1){20}}
   \put(90,10){\line(1,0){20}}
   \put(70,30){\line(1,0){80}}
   \put(110,50){\line(1,0){20}}
   \put(90,30){\line(1,1){20}}
   \put(110,10){\line(1,1){20}}
   \end{picture} 
\caption{$\Sigma$ and $\Sigma(0)$}
\end{figure}
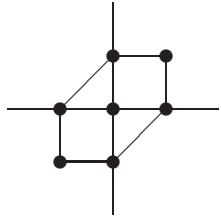

  The algebra $R^{\natural}$ is equipped with endomorphisms $\delta_u^*$ in Eq.~(\ref{eq:defn of delta_y}) and $S_y^*=\delta_{\beta(y)}^*$\ $(u\in X^{\vee}, y\in X)$. Let $(\tP^{\natural},\tcL^{\natural}):=(\Proj R^{\natural}(\xi^e),\cO_{\tP^{\natural}}(1))$. The quotient by $(S_y;y\in X)$ of 
$(\tP^{\natural,\wedge},\tcL^{\natural,\wedge})$  
is algebraized into an $S$-scheme 
$(P^{\natural},\cL^{\natural})$.    

It is shown that 
the $S$-scheme $(P^{\natural},\cL^{\natural})$ has an open subscheme $G$ (resp. $\cG$) which is the connected N\'eron model (resp. the N\'eron model) of $G_{\eta}=\cG_{\eta}=P^{\natural}_{\eta}$.  We also see 
$(P^{\natural}_{\eta},\cL^{\natural}_{\eta})
\simeq (G_{\eta},\cL^{\otimes 3}_{\eta})$. 
The closed fiber $P^{\natural}_0$ consists of 3 irreducible 
components $X_i$ $(i\in [0,2])$ with 
$\delta_{-m_1}(X_0)=X_1,\ \delta_{-m_2}(X_0)=X_2,$ 
each being a 3-times blowing up of $\bP^2_{k(0)}$, 
whose configuration is given by $\op{Vor}\hskip -0.5mm/3X$.  
Moreover $P^{\natural}$ is regular. To see this,  
by Corollary~\ref{cor:if Chat=Bhat Rk Sk},  
it suffices to prove that $\tP^{\natural}$ is regular. 
$\tP^{\natural}$ is covered by affine open sets $W_{\alpha,u}=
\Spec B_{\alpha,u}$\ $(\alpha\in\Sigma, u\in X^{\vee})$ where 
$B_{\alpha,u}=R[\xi_{\beta,v}/\xi_{\alpha,u}; \beta\in\Sigma,v\in X^{\vee}]$.   We see 
$$B_{0,0}=R[w^x;x\in X],\  B_{\epsilon m_i,0}=R[sw^{\epsilon }_i, 
w^{\epsilon }_{i+1}, w^{\epsilon }_{i+2}] $$
where $\epsilon=\pm 1$  and $w_{i}=w_{i+3}$ $(\forall i\in\bZ)$.  
Hence $B_{\alpha,u}=\delta_u^*B_{\alpha,0}$ is regular,  
so that $P^{\natural}$ is regular.  
Note that $(P_{l},\cL_{l})\simeq (P^{\natural},(\cL^{\natural})^{\otimes 2l})$\ $(\forall l\geq 1)$.  
\begin{rem}\label{rem:no cubic str on cLflat on cG}
$\cG:=P^{\flat}\setminus\Sing(P^{\flat}_0)$ is the N\'eron model of $P_{\eta}=G_{\eta}$ where $\Sing(P^{\flat}_0)$ is the singular locus of  $P^{\flat}_0$.  
Moreover $\cL^{\flat}_{|\cG}$ has no cubical structure. Indeed, otherwise, 
$(\cL^{\flat})^{\otimes 3}$ is 
cubical on $\cG$, so that 
$(\cL^{\flat})^{\otimes 3}=\cL^{\natural}$  on $\cG$ 
by Lemma~\ref{lemma:cub str}~(\ref{item:cubic str on G}). Hence 
by the uniqueness of a cubical compactification \cite[1.2]{Nakamura24}, 
$P^{\flat}\simeq P^{\natural}$, which contradicts 
$P^{\flat}_0\not\simeq P^{\natural}_0$.    
\end{rem}
\begin{rem}\label{rem:Jac_phi}Let 
$C$ be a stable curve which is the union 
of two nonsingular rational curves meeting 
at three distinct points over an algebraically closed field. 
The generalized Jacobian variety $\Pic^0_{C}$
of $C$ may be identified with $\bG_{m,k(0)}^{\oplus 2}$. 
\cite{OS79} constructs ``compactifications'' $\Jac_{\phi}(C)$ of 
$\Pic^0_{C}$. By \cite[pp.~84-85]{OS79}, 
$\Jac_{\phi}(C)$ is either $P_0$\ $(\phi=(v_1-v_2)/2)$ 
or $P^{\flat}_0$\ $(\phi=0)$. 
See \S\S~\ref{subsection:2dim PSQAS}-\ref{subsection:twisted PSQAS} 
for $P_0$ and $P^{\flat}_0$. 
No $\Jac_{\phi}(C)$ is isomorphic to $P^{\natural}_0$. 
\end{rem}

\end{document}